\documentclass[
a4paper,    
11pt,
english,
reqno
]{amsbook}

\usepackage{mystyle}

\title{Théories à types dépendants et algèbre de dimension supérieure}
\author[C. Leena Subramaniam]{Chaitanya Leena Subramaniam}

\begin{document}

\begin{titlepage}
\begin{changemargin}{-.8in}{-.8in}
  \includegraphics[width=0.4\textwidth]{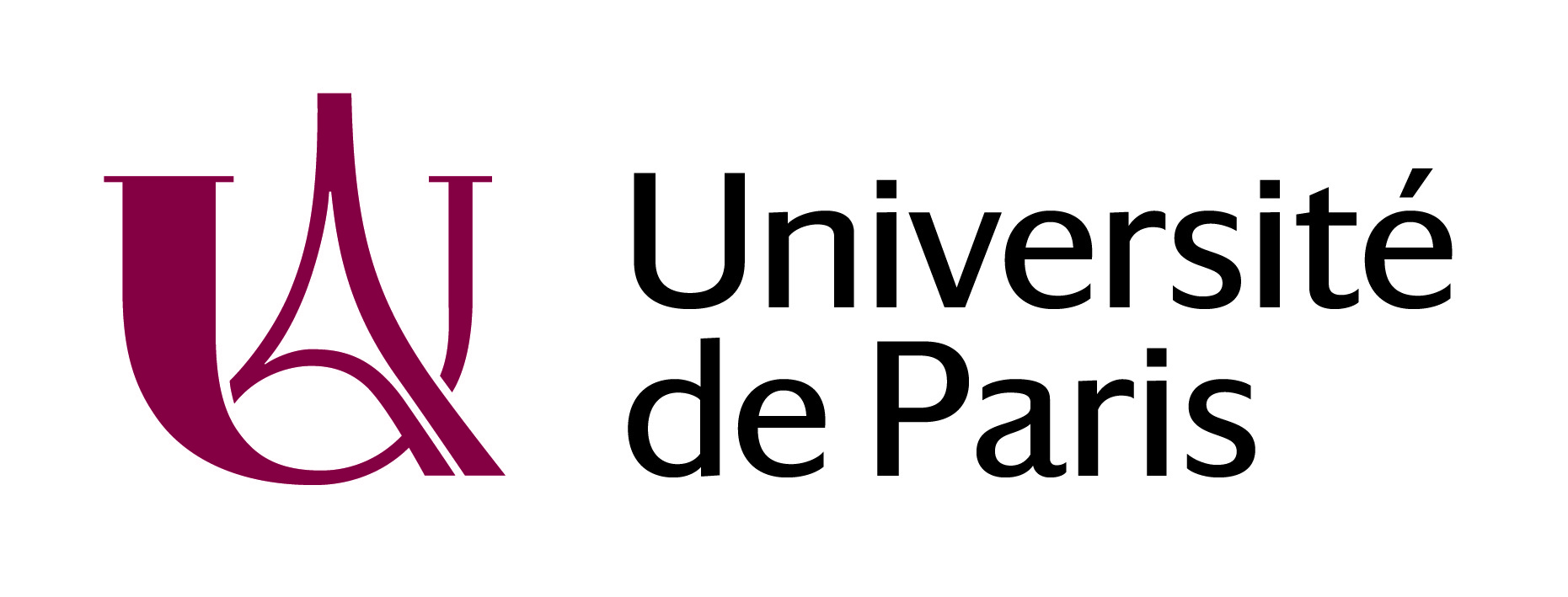}
  
  \vspace*{1cm}
  \raggedleft
  {
    \scshape{\textbf{Université de Paris}} \par}
  \vspace{0.25cm}
  {\large École de Sciences Mathématiques Paris Centre (ED 386)}
  
  {Institut de Recherche Fondamentale en Informatique (IRIF)}

  \vspace*{.25cm}
  {\Large \scshape \textbf{Thèse de doctorat en Informatique}}

  {\large Dirigée par Paul-André Melliès}

  \vspace*{2cm}
  \hrulefill
  \vspace*{0.25cm}
  
  \centering
  {\Huge \scshape \textbf{From dependent type theory to higher algebraic structures}\par}

  \vspace*{0.25cm}
  \hrulefill
  
  \vspace{1cm}

  {\Large \textbf{par Chaitanya Leena Subramaniam}}

  \vspace*{2cm}
  {\slshape Présentée et soutenue publiquement le 28 septembre 2021 \\
    devant un jury composé de : \par}
  \vspace{0.8cm}
  
  \begin{tabular}{ll}
    \textbf{M. Michael Shulman}* (Assoc. Prof., University of San Diego) & Rapporteur \\
    \textbf{M. Carlos Simpson} (DR CNRS, Université de Nice) & Rapporteur \\
    \textbf{M. Tom Hirschowitz}* (DR CNRS, Université Savoie Mont Blanc) & Examinateur \\
    \textbf{Mme Muriel Livernet} (Prof., Université de Paris) & Examinatrice \\
    \textbf{M. Samuel Mimram} (Prof., École Polytechnique) & Examinateur \\
    \textbf{Mme Emily Riehl}* (Assoc. Prof., Johns Hopkins University) & Examinatrice \\
    \textbf{M. Paul-André Melliès} (DR CNRS, Université de Paris) & Directeur de Thèse \\
    \textbf{M. Gilles Dowek} (DR Inria, Université ENS Paris--Saclay) & Membre invité\\ 
    \textbf{M. Thomas Streicher}* (Prof., Technische Universität Darmstadt) & Membre invité\\

    \vspace{0.5cm}
  {\small *en visioconférence}
  \end{tabular}
 
\end{changemargin}
\end{titlepage}

\maketitle

\cleardoublepage
\thispagestyle{empty}
\vspace*{13.5pc}
\begin{center}
  \emph{For Amma and Appa}
\end{center}
\cleardoublepage

\begin{abstract}
  In the first part of this dissertation, we give a definition of ``dependently
  typed/sorted algebraic theory'', generalising the ordinary multisorted algebraic
  theories of Lawvere--Bénabou. Dependently sorted algebraic theories in our
  sense form a strict subclass of the ``generalised algebraic theories'' of
  Cartmell. We prove a classification theorem for dependently sorted algebraic
  theories, and use this theorem to prove the existence of dependently sorted
  algebraic theories for a number of varieties of algebraic structures, such as
  small categories, n-categories, strict and weak omega-categories, planar
  coloured operads and opetopic sets. We also prove a Morita equivalence between
  dependently sorted algebraic theories and essentially algebraic theories,
  showing that every locally finitely presentable category is the category of
  models of some dependently sorted algebraic theory. We give a definition of
  strict and weak homotopical models of a dependently sorted algebraic theory,
  and prove a rigidification theorem in a particular case. We also study the
  opetopic case in detail, and prove that a number of varieties of algebraic
  structures such as small categories and coloured planar operads can be ``typed'' 
  by the category of opetopes.

  The second part of this dissertation concerns accessible reflective
  localisations of locally presentable infinity-categories. We give a definition
  of ``pre-modulator'', and prove a canonical correspondence between
  pre-modulators and accessible orthogonal factorisation systems on a locally
  presentable infinity-category. Moreover, we show that every such factorisation
  system can be generated from a pre-modulator by a transfinite iteration of a
  ``plus-construction''. We give definitions of ``modulator'' and ``left exact
  modulator'', and prove that they correspond to those factorisation systems
  that are modalities and left-exact modalities respectively. Finally we obtain
  a correspondence between left-exact localisations of infinity-topoi and
  left-exact modulators.

  
  \newabstract{french} Dans la première partie de cette thèse, nous proposons une
  définition de \og théorie algébrique à types/sortes dépendants \fg{} qui
  généralise les théories algébriques ordinaires de Lawvere--Bénabou. Les
  théories algébriques à sortes dépendantes, en notre sens, forment une
  sous-classe stricte des \og théories algébriques généralisées\fg{} de
  Cartmell. Nous démontrons un théorème de classification pour les théories
  algébriques à sortes dépendantes, et nous utilisons ce théorème pour montrer
  l'existence de plusieurs de ces théories --- parmi elles, les théories des
  petites catégories, des n-catégories, des omega-catégories strictes et
  faibles, des opérades planaires colorées, et des ensembles opétopiques. Nous
  étudions le cas opétopique en détail. Nous montrons également une équivalence
  de Morita entre les théories algébriques à sortes dépendantes et les théories
  essentiellement algébriques, et nous concluons que chaque catégorie localement
  finiment présentable admet une description comme catégorie des modèles d'une
  théorie algébrique à sortes dépendantes. Nous proposons également les
  définitions des modèles homotopiques strictes et faibles d'une théorie
  algébrique à sortes dépendantes, et nous montrons un théorème de
  rigidification dans un cas particulier.

  La deuxième partie de cette dissertation concerne les localisations
  refléctives accessibles des infini-catégories localement présentables. Nous
  donnons une définition de \og pré-modulateur\fg{} et montrons une
  correspondance entre les pré-modulateurs et les systèmes de factorisations
  accessibles sur une infini-catégorie localement présentable. Nous montrons
  également que chaque tel système de factorisation est engendré à partir d'un
  pré-modulateur par itération transfinie d'une \og construction plus \fg{}.
  Nous proposons les définitions de \og modulateur\fg{} et de \og modulateur
  exact à gauche\fg{} et nous montrons des correspondances avec les modalités et
  les modalités exactes à gauche respectivement.

  \smallskip
  \noindent \textbf{Mots-clés.} Théories algébriques, catégories supérieures,
  théorie de l'homotopie, théorie des types dépendants.
\end{abstract}

\pagenumbering{roman}
\tableofcontents
\pagenumbering{arabic}
\addtocontents{toc}{\protect\setcounter{tocdepth}{0}}
\chapter*{Introduction}

This dissertation is intended as a contribution to the universal algebra and
homotopy theory of ``dependently typed theories''.

\section*{Motivation and context}
This dissertation is an attempt to substantiate two theses.

\begin{thesis}
  \label{thes:dependency-cellularity} 
  Type dependency in logic corresponds to the cellularity inherent in
  higher-dimensional algebraic structures.
\end{thesis}
\begin{thesis}
  \label{thes:dep-typed-higher-alg-struc} 
  The correspondence between type dependency and cellularity is particularly
  well\=/suited to the description of \emph{homotopy-coherent} higher-dimensional
  algebraic structures in spaces (\emph{higher algebraic structures}).
\end{thesis}
\subsection*{Dependently typed theories and cellularity}
\cref{thes:dependency-cellularity} is best illustrated by way of an example.
Small categories are algebraic structures on (directed multi)graphs---namely,
there is a forgetful monadic functor $\Cat\to\psh{\GG_1}$ from the category of
small categories to the presheaf category $\psh{\GG_1}\eqdef\fun{\GG_1\op}\Set$
of graphs. Graphs have a cellular structure that is easy to visualise---their
vertices are cells of dimension $0$, their edges, cells of dimension $1$, and
there are no cells of higher dimension. Then the algebraic structure of a
category (such as the operation of composition of morphisms) can be defined
using the cellular structure of its underlying graph.

Let us describe this syntactically, step by step using dependent types, in such
a way that the graph-cellularity inherent in the type dependency is clear. The
monadic functor $\Cat\to\psh{\GG_1}$ underlies the fact that the theory of
categories is an extension of the theory of graphs. A graph $X$ is a diagram
$s,t\colon E_X\rightrightarrows V_X$ of sets, where $E_X$ is the set of edges of
$X$ and $V_X$ is the set of vertices. This can be represented syntactically by
the \defn{dependent type signature}
\[
  \bS_{\GG_1}\eqdef \left\{\emptyset\vdash V \ \ \mathrm{type} \quad,\quad x\of
    V,y\of V\vdash E(x,y) \ \ \mathrm{type}\right\}
\]
which \emph{displays} the type $E$ of edges as dependent on the context $(x\of
V,y\of V)$ (we write $(x,y\of V)$ for short) of a pair of variables of the type
$V$ of vertices. A \emph{model} of $\bS_{\GG_1}$ (in $\Set$) is the data of a
set $V$ and for every pair $(x,y)\in V\times V$, a set $E(x,y)$. Clearly, models
of $\bS_{\GG_1}$ are exactly graphs $E\to V\times V$. In other words, $\bS_{\GG_1}$ is
the \defn{dependently typed theory} of graphs.

The signature $\bS_{\GG_1}$ corresponds to the
category
\[
  \GG_1\eqdef
  \left\{
    \begin{tikzcd}
      D^0\ar[r,"s",shift left] \ar[r,shift right,"t"']
      &D^1
    \end{tikzcd}
  \right\}
\] by associating the type $V$ to $D^0$ and the type family $E$ to $D^1$.
The category $\GG_1$ is a \emph{direct} category---each representable in
$\psh{\GG_1}$ has a canonical \emph{boundary}, and every graph is a \emph{cell
  complex} constructed by ``gluing'' representables along boundaries. For
instance, via the correspondence $\bS_{\GG_1}\sim\GG_1$, the context $(x,y\of
V)$ can be understood as the graph with two vertices and no edges, that is the
boundary $\partial D^1$ of the representable $D^1\in\GG_1$. The boundary of the
representable $D^0$ is the empty graph, which corresponds to the empty context
$\emptyset$. Then the graph
\[
  \Gamma =\{ x\xto{f}y\xto{g}z \}
\]
with three vertices and two ``composable'' edges can be seen as the cell complex
\[
  \begin{tikzcd}
    \partial D^1\ar[r,"{(y,z)}"]\ar[d]
    & \{x\xto{f} y \quad z\}=\ft\Gamma 
    \ar[d]\\
    D^1\ar[r,"g"]
    &\Gamma\pomark
  \end{tikzcd}
\]
obtained by gluing the representable $D^1$ to the graph $\ft\Gamma$ along the
map $\partial D^1\to\ft\Gamma$ that takes the ``source'' and ``target'' vertices
of $D^1$ to $y$ and $z$ respectively. Via the correspondence
$\bS_{\GG_1}\sim\GG_1$, $\Gamma$ can be seen syntactically as the extension
\[
  \Gamma=(x,y,z\of V,f\of E(x,y),g\of E(y,z)) \tto (x,y,z\of V,f\of
  E(x,y))=\ft\Gamma
\]
of the context $\ft\Gamma$ with the variable $g\of E(y,z)$.

The \emph{dependently typed theory of small categories} introduces a
\defn{term signature} with two operations
\[
  \bF_{\Cat} \eqdef
  \left\{
    x\of V\vdash i(x)\type E(x,x)\quad,\quad x,y,z\of V,f\of E(x,y),g\of
    E(y,z)\vdash c(g,f)\type E(x,z) 
  \right\}
\]
over the signature $\bS_{\GG_1}$, as well as a set $\bE_{\Cat}$ of three
\defn{equations}
\begin{align*}
  x,y\of V,f\of E(x,y)
  &\vdash c(i(y),f) = f \type E(x,y)\\
  x,y\of V,f\of E(x,y)
  &\vdash c(f,i(x)) = f \type E(x,y)\\
  x,y,z,w\of V, f\of E(x,y), g\of E(y,z),
  h\of E(z,w)
  &\vdash c(h,c(g,f))= c(c(h,g),f) \type E(x,w)
\end{align*}
over the term signature $\bF_\Cat$. Then a model of the theory
$\bT_\Cat\eqdef(\bF_\Cat,\bE_\Cat)$ over the signature $\bS_{\GG_1}$ is a graph
with the functions of sets encoded by the operations of $\bF_\Cat$, satisfying
the equations in $\bE_\Cat$. It is easy to see that the models of $\bT_\Cat$ are
exactly small categories.

Similarly to small categories, a number of other familiar algebraic structures
can be \emph{classified} by dependently typed algebraic theories, each with its
own notion of cellularity. For instance, $\omega$\=/categories (respectively,
$n$\=/categories) are defined by a theory over the dependent type signature of
\emph{globes} (respectively, globes of dimension $\leq n$). Another example are
coloured planar operads, which are defined by a theory over the dependent
type signature of \emph{corollas} or \emph{elementary trees}.

This idea of associating type dependency to cellularity is not new---for
instance, it is clearly present in the theory of FOLDS of
\cite{makkai1995folds}, and in the groupoidal and globular $\omega$\=/groupoidal
structure of identity types of \cite{hofmann1998groupoid} and
\cite{lumsdaine2009weak,van2011types}. More recently, \cite[App.
A]{brunerie2016homotopy} gives a dependently typed theory of weak
$\omega$\=/groupoids, and \cite{finster2017type,benjaminthese2020} describe a
dependently typed algebraic theory of weak $\omega$\=/categories. Nevertheless,
a \emph{robust} definition and classification of dependently typed algebraic
theories and the corresponding algebraic structures defined by them has not yet
been laid down, and this is one of the primary goals of this dissertation.

In this dissertation, we give a number of definitions of what we mean by
\defn{dependently typed algebraic theory}, using the correspondence between
type dependency and cellularity. We prove a classification theorem
(\cref{thm:classification-dep-alg-theories}) showing that all these definitions
are equivalent. Syntactically, dependently typed algebraic theories correspond
to a \emph{strict} subclass of the \emph{generalised algebraic theories} of
\cite{cartmell1978generalised}.

We use the classification theorem to recognise many dependently typed algebraic
theories. In fact, it turns out that \emph{every} locally finitely presentable
category is the category of models in $\Set$ of some dependently typed algebraic theory
(\cref{thm:lfp-cat-models-of-C-sorted-theory}), therefore dependently typed
algebraic theories are just as expressive (with respect to their models in
$\Set$) as \emph{essentially algebraic} or \emph{finite limit} theories
(\cref{thm:morita-equiv-gat-eat-C-cxl-cats}).

\subsection*{Cellular multicategories}
Continuing with the example of the theory $\bT_\Cat$ of small categories, remark
that the contexts in every declaration of $\bF_\Cat$ and $\bE_\Cat$ (to the left
of the turnstile symbol ``$\vdash$'') are all graphs (they correspond to cell
complexes in $\psh{\GG_1}$). We can therefore see the declarations in $\bF_\Cat$
as \emph{operations} or multimorphisms of a ``cellular multicategory'' that take
as input a (finite) graph, and whose output ``cell'' is a representable. For
example, the operations $c(-,-)$ and $i(-)$ can be visualised as follows.
\[
  \begin{tikzpicture}
    \tikzset{scale=.8}
    \node [treenode] (b) at (0, 0) [label = left : $c$] {};
    \node (e1) at ($(b) + (-1, .7)$) [above] {\small $ x$};
    \node (e3) at ($(b) + (1, .7)$) [above] {\small $w$};
    \node (e2) at ($(e1)!.5!(e3)$) {\small $y$};
    \node (t) at ($(b) + (0, -.7)$) [below] {{\small$E(x,w)$}};
    \draw (b) -- (e1);
    \draw (b) -- (e3);
    \draw[->] (e1) -- (e2);
    \draw[->] (e2) -- (e3);
    \draw (b) -- (t);
  \end{tikzpicture}
  \qquad\qquad
  \begin{tikzpicture}
    \tikzset{scale=.8}
    \node [treenode] (b) at (0, 0) [label = left : $i$] {};
    \node (e1) at ($(b) + (0,.7)$) [above] {\small $x$};
    \node (t) at ($(b) + (0, -.7)$) [below] {{\small$E(x,x)$}};
    \draw (b) -- (e1);
    \draw (b) -- (t);
  \end{tikzpicture}
\]
The term $c(f,i(x))$ can be visualised as the tree of operations below.
\[
  \begin{tikzpicture}
    \tikzset{scale=.8}
    \node [treenode] (b) at (0, 0) [label = left : $c$] {};
    \node (e1) at ($(b) + (-1, .7)$) [above] {\small $ x$};
    \node (e3) at ($(b) + (1, .7)$) [above] {\small $w$};
    \node (e2) at ($(e1)!.5!(e3)$) {\small $x$};
    \node (t) at ($(b) + (0, -.7)$) [below] {{\small$E(x,w)$}};
    \node [treenode] (b1) at (-.5, 1.9) [label = left : $i$] {};
    \node (e11) at ($(b1) + (0,.7)$) [above] {\small $x$};
    \node (t1) at ($(b1) + (0, -.7)$) [below] {};
    \draw (b1) -- (e11);
    \draw (b1) -- (t1);
    \draw (b) -- (e1);
    \draw (b) -- (e3);
    \draw[->] (e1) -- (e2);
    \draw[->] (e2) -- (e3);
    \draw (b) -- (t);
  \end{tikzpicture}
\]
This ``cellular cartesian multicategory'' approach to dependently typed
algebraic theories strictly generalises the point of view of multisorted
algebraic theories as \emph{cartesian multicategories} (equivalently, categories
with finite products) due to Lawvere \cite{Lawvere869} and Bénabou
\cite{benabou1968structures}. In this framework, multisorted algebraic theories
are cellular cartesian multicategories whose cells are all of dimension $0$
(they are ``points'', namely they have no dependencies).

The point of view of dependently sorted/typed algebraic theories as cellular
cartesian multicategories appears in \cite[Sec. II]{fiore2008second} (somewhat
implicitly), where they are called \emph{$\Sigma_0$\=/models with substitution}.
More generally, the theory of \emph{monads with arities} and \emph{theories with
  arities} of \cite{weber2007familial,mellies2010segal,Berger2012} can be seen
as an abstract theory of ``multicategories'' with general arities. However, a
general concrete theory of cellular ``species'' and multicategories
corresponding to dependently typed algebraic theories is (to the best of my
knowledge) far from being established.

In this dissertation, we outline the point of view of dependently sorted/typed
algebraic theories as cellular cartesian multicategories
(\cref{sec:C-collections-theories}) but we do not set up a general theory of
cellular multicategories.\footnote{I think that such a general theory should
exist, and should subsume \emph{globular operads} as particular cases.}

\subsection*{Higher algebraic structures and theories}
\cref{thes:dep-typed-higher-alg-struc} has to do with algebraic structures up to
homotopy.
\defn{Higher algebraic structures} are the ``right'' notion of algebraic structures up
to homotopy in spaces. A general blueprint describing a kind of higher algebraic
structure (such as $(\infty,1)$\=/categories, \oo-operads, stacks, spectra,
$E_\infty$\=/spaces) is as a diagram of spaces, along with operations that are
subject to ``equations''. These equations are quotients of spaces, which are
infinite towers of homotopy\=/coherent data, and this renders the definition of
higher algebraic structures much more subtle than their $\Set$\=/valued
(``discrete'' or ``0\=/truncated'') versions.

For instance, composition of morphisms in $(\infty,1)$-categories
(``\oo-categories'' for short), is homotopy-coherently associative. This is
usually described in different ways, through different model categories, such as
the model structure for quasicategories on simplicial sets and the model
structure for complete Segal spaces on bisimplicial sets. Indeed, most of the
theory of \oo-categories has been formally developed only in quasicategories.
Nevertheless, quasicategories are not in any well-understood way a canonical
model for \oo-categories. The question of providing a good syntactic theory of
\oo\=/categories is still an open problem (and under active research, see
\cite{riehl2017type,finster2018towards,allioux2021types}).

A syntactic counterpart (to the predicate of equality for quotients of sets) for
quotients of spaces are the \emph{Martin-Löf identity types} of \emph{Homotopy
  Type Theory} (HoTT). This is substantiated by the fact that HoTT can be
interpreted in the \oo\=/topos of spaces \cite{kapulkin2012simplicial}, and
indeed in any \oo\=/topos \cite{shulman2019infty1toposes}, in such a way that
identity types correspond to the ``predicates of equality'' (iterated
diagonals). Therefore, we might expect that some higher algebraic counterpart,
that uses identity types, of the dependently typed algebraic theory $\bT_\Cat$
of categories is such a syntactic candidate of the theory of \oo\=/categories.
In general, we might expect a \defn{dependently typed higher algebraic theory}
to be an extension of HoTT with types, terms, and \emph{identity paths} instead
of equations in more or less the same way as the dependently typed algebraic
theories of \cref{thes:dependency-cellularity}.

\subsection*{Homotopical models of algebraic theories}
In this dissertation, we do not attempt to give a general theory of dependently
typed higher algebraic theories, which is a problem that requires a lot more
groundwork in order to be tackled.

We restrict ourselves to describing a theory of \defn{homotopy models} of the
dependently typed algebraic theories of \cref{thes:dependency-cellularity}, that
generalises certain aspects of the theory of homotopy models of multisorted
algebraic theories due to
\cite{Schwede2001stablehtpyAlgTheories,badzioch2002algtheories,Rezk2002simplicialAlgTheories,bergner2006rigidification}.
Even in this case, we are not able to show that all of the rigidification theory
of \emph{op. cit.} generalises (although we conjecture that it does, see
\cref{sec:htpical-models-C-cxl-cats}). Nevertheless, we posit this
generalisation as suggestive of the existence of a general theory of dependently
typed higher algebraic theories.

\subsection*{Idempotent opetopic theories}
Dependent types are sufficiently expressive to allow for the ``shapes'' of
operations (and equations) of an algebraic structure to be abstracted into the
type signature. For example, consider the operation (of composition of morphisms
in small categories)
\[
  x,y,z\of V,f\of E(x,y),g\of E(y,z)\vdash c(g,f)\type E(x,z)
\]
in the term signature $\bF_\Cat$. We can extend the type signature $\bS_{\GG_1}$
with the \emph{type declaration}
\[
  x,y,z\of V,f\of E(x,y),g\of E(y,z),h\of E(x,z)\vdash C(x,y,z,f,g,h)\ \
  \mathrm{type}
\]
which corresponds to the shape of the operation of composition $c(-,-)$. The
shape of the dependent type $C$ can be visualised as follows.
\begin{center}
  \begin{tikzpicture}
    \node (a) at (0, 0){$x$};
    \node (b) at (1.5, 0){$y$};
    \node (c) at (.75, 1){$z$};
    \draw[->] (a) -- (c) node[midway,left]{\small$f$};
    \draw[->] (c) -- (b) node[midway,right]{\small$g$};
    \draw[->] (a) -- (b) node[midway,below]{\small$h$};
    \draw (.75, .35) node {$\Downarrow$};
  \end{tikzpicture} 
\end{center}
If we write $\bS_{\GG_1}'$ for the extension of $\bS_{\GG_1}$ with the previous
type declaration, then $\bS_{\GG_1}'$ corresponds to the category
\[
  \GG_1'\eqdef
  \left\{
    \begin{tikzcd}
      D^0\ar[r,shift right] \ar[r,shift left]
      &D^1\ar[r]\ar[r,shift right] \ar[r,shift left]
      &C
    \end{tikzcd}
  \right\}
\]
and the category $\Cat$ has a monadic functor $\Cat\to\psh{\GG_1'}$ sending a
small category $A$ to a presheaf whose fibre over $C$ is the set of commutative
triangles in $A$. Thus the
operation of composition in a category can be ``opetopified'' into a
representable cell of the type signature. In fact, the theory $\bT_\Cat$ of
categories can be ``totally opetopified''---all its operations and equations can
be integrated into a type signature $\bS_{\OO_{\leq 3}}$ that is an extension of
the signature $\bS_{\GG_1}$. The signature $\bS_{\OO_{\leq 3}}$ corresponds to a
direct category $\OO_{\leq 3}$ (the category of \emph{opetopes} of dimension
$\leq 3$), such that $\GG_1\subto\OO_{\leq 3}$ is a full subcategory. Moreover,
there is a monadic functor $\Cat\to\psh{\OO_{\leq 3}}$ that is \emph{fully
  faithful}---we have therefore transformed all the ``structure'' of a small
category into ``properties'' of a presheaf on $\OO_{\leq 3}$. Finally, $\Cat$ is
the category of models of a finite limit sketch on the category $\OO_{\leq
  3}\op$. This new dependently typed algebraic theory of small categories over
the signature $\bS_{\OO_{\leq 3}}$ is an example of an \emph{idempotent opetopic
  theory}.

The process of opetopification was originally defined for coloured symmetric
operads in \cite{Baez1998}. A similar process is fundamental to the
$\deux{\mathit L}$\emph{\=/structures} of \cite{makkai1995folds}.

In this dissertation, we study a class of idempotent opetopic theories obtained
from the category $\OO$ of opetopes. We show that we can associate a good notion
of \emph{homotopy-coherent} model to each of these theories.

\subsection*{Locally presentable \texorpdfstring{\oo}{infinity}-categories}
The collection of all higher algebraic structures of a given class forms a
\emph{locally presentable \oo\=/category}. The theory of locally presentable
\oo\=/categories perfectly subsumes the theory of locally presentable
$1$\=/categories, and we conjecture that just as dependently typed algebraic
theories correspond to locally finitely presentable $1$\=/categories, so do the
conjectural dependently typed higher algebraic theories give rise
to\footnote{There is a subtlety here, in that dependently typed higher algebraic
  theories likely correspond to finitely complete (=``lex'') \oo\=/categories,
  and locally finitely presentable \oo\=/categories correspond to
  \emph{idempotent\=/complete} lex \oo\=/categories.} locally finitely
presentable \oo\=/categories . As a starting point, it is clear from the
definitions that the homotopy models of a dependently typed algebraic theory
form a locally finitely presentable \oo\=/category. The conjecture is also
supported by \cite{kapulkin2019internal}, where it is shown that finitely complete
\oo\=/categories correspond to \emph{tribes}, which are an abstraction of type
theories with Martin-Löf identity types.

In a separate but related part of this dissertation (\cref{part:part2}), we
study locally presentable \oo\=/categories, and develop a theory of
\emph{(pre\=/)modulators}, which shows that every accessible factorisation system on a
locally presentable \oo\=/category can be obtained by iterating a
\emph{plus\=/construction}. 

\section*{Organisation}
This dissertation is divided into two parts.

\medskip \cref{part:part1} contains the main content of this dissertation, and
is concerned with
\cref{thes:dependency-cellularity,thes:dep-typed-higher-alg-struc}.

\cref{chap:contextual-categories} develops the theory of dependently
typed/sorted algebraic theories. We begin by describing the correspondence
between cellularity and type dependency. The key definition is that of
\emph{\lfd~category}, which is an equivalent reformulation of the definition of
``simple category'' from \cite{makkai1995folds}. We use this reformulation to
give multiple definitions of the dependently typed algebraic theories that we
are interested in. The approach we take is via the formalism of
\emph{$\lfdex$\=/contextual categories} (\cref{def:C-cxl-cat}), which allows us
to avoid explicitly working with syntax. We prove a classification theorem
(\cref{thm:classification-dep-alg-theories}) showing that these definitions are
equivalent, and which strictly generalises Lawvere's classification theorem for
(multisorted) algebraic theories. We immediately use this result to detect
several dependently typed algebraic theories. Finally, we describe the syntactic
characterisation of dependently typed algebraic theories explicitly, and show
that they are a strict subclass of Cartmell's generalised algebraic theories
\cite{cartmell1978generalised,cartmell1986}. Other results of interest in this
chapter are the initiality of the free contextual category on a type signature
(\cref{prop:ctxl-functor-cxl-C-initial}) and the various properties that it
possesses
(\cref{prop:cxl-C-fin-limits,prop:cxl-C-pushouts-nonempty-coprods,prop:cxl-C-display-map-epi,prop:cxl-C-codescent}).

\cref{cha:models-of-C-contextual-categories} develops the theory of homotopy
models in spaces of dependently typed algebraic theories. We begin by
classifying the categories of $\Set$\=/models of dependently typed algebraic
theories as exactly the locally finitely presentable categories, allowing us to
conclude that the classes of dependently typed algebraic theories, essentially
algebraic theories, and generalised algebraic theories are all
\emph{Morita\=/equivalent}
(\cref{thm:lfp-cat-models-of-C-sorted-theory,thm:morita-equiv-gat-eat-C-cxl-cats}).
After recalling some elements of the theory of simplicial model categories, we
show a \emph{rigidification} theorem
(\cref{thm:flasque-C-spaces-rigidification}) for homotopy models in spaces of the initial
dependently typed algebraic theory on a \lfd~category $\lfdex$. Finally, we
prove the existence of a model structure for homotopy models in spaces of any
dependently typed algebraic theory (\cref{thm:partial-soln-rigid-conjecture}),
and conjecture a general rigidification theorem. Other results of interest in
this chapter are the construction of a \emph{flasque} intermediate global model
structure on simplicial presheaves on the initial $\lfdex$\=/contextual category
(\cref{thm:flasque-model-struc-existence}) and the description of
$\lfdex$\=/spaces as sheaves of \oo\=/groupoids on the initial
$\lfdex$\=/contextual category (\cref{prop:C-spaces-model-struct-equals-Cech}).

\cref{chap:opetopic-theories} studies a class of idempotent theories over the
\lfd~category $\OO$ of \emph{opetopes}. We begin by recalling the construction
of the category $\OO$ from \cite{Kock2010,hothanh18}. We define a family of
parametric right adjoint monads $\optPolyFun^n$ whose algebras are the
\emph{opetopic algebras} (\cref{def:opetopic-algebra}). We show that every
category $\oAlg$ of opetopic algebras admits a fully faithful, accessible
monadic \emph{opetopic nerve functor} $\oAlg\subto\pshO$ to the category of
opetopic sets, and that $\oAlg$ is the category of models of a finite projective
sketch on $\OO\op$ (\cref{th:nerve-theorem-O}), which implies that $\oAlg$ is
the category of models of an \emph{idempotent opetopic theory}
(\cref{thm:OAlg-models-of-idempotent-O-theory}). We show that particular cases
of $\oAlg$ are the categories $\Cat,\Opd,\Cmbd$ of small categories, coloured
planar operads, and coloured combinads of \cite{Loday2012a}, thus each of these
is the category of models of an idempotent opetopic theory. Finally, we show
that a technique due to \cite{horel2015model} allows us to define a model
structure for \emph{homotopy\=/coherent opetopic algebras} and prove it Quillen
equivalent to a model structure on \emph{simplicial opetopic algebras}
(\cref{thm:rigidification-Oalg}). In the particular case of $\Cat$ and $\Opd$,
we obtain Rezk's model structure for Segal spaces and the planar version of
Cisinski-Moerdijk's model structure for Segal dendroidal spaces.

\medskip \cref{part:part2} of this dissertation studies accessible orthogonal
factorisation systems in, and accessible localisations of, locally presentable
\oo\=/categories.

\cref{chap:fact-systems} describes the theory of (orthogonal) factorisation
systems in \oo\=/categories, using the
pushout-product/pullback-hom tensor/enrichment of arrow \oo\=/categories. This
formalism is used to recover several results about modalities and lex modalities
(\cref{prop:stable-lex-FS,prop:mod-lccc,prop:mod-slex,prop:lex-loc=lex-mod}).

\cref{chap:SOA} develops the theory of \emph{pre\=/modulators}. We begin by
showing that with a slight modification provided by the pushout-product and
pullback-hom, Kelly's \emph{small object argument} to construct accessible
orthogonal factorisation systems in locally presentable $1$\=/categories admits
a generalisation to locally presentable \oo\=/categories (\cref{thm:!SOA2}). We
show that in the particular case of sheafification associated to a Grothendieck
topology, Kelly's construction can be simplified to a \emph{plus\=/construction}
for presheaves (\cref{lem:plus-is-plus-for-topology}), and that for every
accessible factorisation system generated by a \emph{pre\=/modulator} on a
locally presentable \oo\=/category, Kelly's construction can be simplified to
the plus\=/construction (\cref{thm:plus-construction}). Moreover, without loss
of generality, every accessible factorisation system on a locally presentable
\oo\=/category is presented by a pre\=/modulator
(\cref{prop:modulator-completion}). We use the theory of pre\=/modulators to
define \emph{modulators} and \emph{lex modulators}, whose respective
plus\=/constructions generate modalities and lex modalities
(\cref{thm:+modality,thm:lex+construction}). We show that every accessible lex
localisation of an \oo\=/topos is generated by a lex modulator
(\cref{prop:lex-mod-4-loc-lex}), thus lex modulators are a good generalisation
of Grothendieck topologies to \oo\=/topoi. Finally, we show that the
plus\=/construction of a lex modulator on an \oo\=/topos converges in
$(n+2)$\=/steps on any $n$\=/truncated object (\cref{prop:truncated-sheaf}).

\medskip Parts of \cref{chap:contextual-categories} are a collaboration with
Peter LeFanu Lumsdaine. Much of~\cref{chap:opetopic-theories} is a collaboration
with Cédric Ho Thanh. \cref{chap:fact-systems,chap:SOA} are part of a
collaboration with Mathieu Anel. I am indebted to them, without
whom this story would be far more incomplete than it is.


\begin{otherlanguage}{french}
  
\chapter*{Introduction (français)}

Cette dissertation se veut une contribution à l'algèbre universelle et à la
théorie de l'homotopie des théories à types dépendants.

\section*{Contexte}

Cette dissertation cherche à justifier deux \emph{thèses}.

\begin{these}
  \label{these:dependence-cellularite} 
  La dépendance de types (ou sortes) en logique correspond à une
  \og cellularité \fg{} intrinsèque aux structures algébriques de dimension supérieure.
\end{these}
\begin{these}
  \label{these:dep-typed-higher-alg-struc}
  La correspondance entre les types dépendants et la cellularité est particulièrement
  bien adaptée à la description des structures algébriques de dimension
  supérieure, \emph{à homotopie cohérente près}, dans les espaces (les types
  d'homotopie).
\end{these}
\subsection*{Théories à types dépendants et cellularité}
Illustrons la thèse \ref{these:dependence-cellularite} par un exemple. Une
\emph{petite catégorie} est une structure algébrique avec un graphe sous-jacent.
(Nota bene~: pour nous, \og graphe\fg{} = \og graphe dirigé\fg{}.) Plus
précisément, il existe un foncteur monadique d'oubli de la catégorie
$\Cat$ des petites catégories vers la catégorie
$\psh{\GG_1}\eqdef\fun{\GG_1\op}\Set$ des graphes. Un graphe admet une structure
cellulaire évidente --- ses cellules de dimension $0$ sont ses n\oe uds, celles
de dimension $1$ sont ses arêtes, et aucun graphe n'a de cellule de dimension
$>1$. De là, la structure algébrique d'une petite catégorie (telle que son
opération de \emph{composition} d'arêtes) s'exprime en utilisant cette structure
cellulaire de son graphe sous-jacent.

Décrivons cela de manière syntaxique, en utilisant les types dépendants, de
façon à mettre en évidence la cellularité des graphes inhérente à
l'axiomatisation. L'existence du foncteur monadique $\Cat\to\psh{\GG_1}$
correspond au fait que la théorie des petites catégories est une extension de la
théorie des graphes. Un graphe $X$ est un diagramme $s,t\colon E_X\rightrightarrows
V_X$ d'ensembles, où $E_X$ est l'ensemble des arêtes de $X$, et $V_X$ son
ensemble de n\oe uds. De manière syntaxique, la théorie des graphes est
représentée par une \defn{signature de types dépendants}
\[
  \bS_{\GG_1}\eqdef \left\{\emptyset\vdash V \ \ \mathrm{type} \quad,\quad x\of
    V,y\of V\vdash E(x,y) \ \ \mathrm{type}\right\}
\]
exhibant le {type} $E$ des arêtes comme dépendant du {contexte} $(x\of
V,y\of V)$ (que l'on raccourcit en $(x,y\of V)$) de deux variables du type
$V$ des n\oe uds. Un \emph{modèle} de $\bS_{\GG_1}$ (dans la catégorie
des ensembles) est la donnée d'un ensemble $V$, et pour chaque paire d'élements
$(x,y)\in V\times V$, d'un ensemble $E(x,y)$. Il est clair que les modèles de
$\bS_{\GG_1}$ sont exactement les graphes~;
autrement dit, $\bS_{\GG_1}$ est la \defn{théorie à types dépendants} des graphes.

La signature $\bS_{\GG_1}$ correspond à la catégorie avec deux objets
et deux morphismes parallèles 
\[
  \GG_1\eqdef
  \left\{
    \begin{tikzcd}
      D^0\ar[r,"s",shift left] \ar[r,shift right,"t"']
      &D^1
    \end{tikzcd}
  \right\}
\] en associant le type $V$ à
l'objet $D^0$, et la famille de types $E$ à l'objet $D^1$. La catégorie $\GG_1$
est \emph{directe} --- chaque représentable dans $\psh{\GG_1}$ admet un
\emph{bord}, et chaque graphe est un \emph{complexe cellulaire} construit en 
recollant des préfaisceaux représentables le long de leurs bords. Par
exemple, via la correspondance $\bS_{\GG_1}\sim\GG_1$, le contexte $(x,y\of V)$
correspond au graphe avec deux n\oe uds et sans arêtes, à savoir le bord
$\partial D^1$ du préfaisceau représentable $D^1\in\GG_1$. Le bord du
représentable $D^0$ est le graphe vide, qui correspond au contexte vide
$\emptyset$. Un autre exemple est le graphe
\[
  \Gamma =\{ x\xto{f}y\xto{g}z \}
\]
avec trois n\oe uds et deux arêtes composables, qui peut être vu comme le complexe
cellulaire
\[
  \begin{tikzcd}
    \partial D^1\ar[r,"{(y,z)}"]\ar[d]
    & \{x\xto{f} y \quad z\}=\ft\Gamma 
    \ar[d]\\
    D^1\ar[r,"g"]
    &\Gamma\pomark
  \end{tikzcd}
\]
obtenu en recollant le représentable $D^1$ au graphe $\ft\Gamma$, le long de
l'application $\partial D^1\to\ft\Gamma$ envoyant la source  et le
but  de $D^1$ sur les n\oe uds $y$ et $z$ respectivement. En utilisant la
correspondance $\bS_{\GG_1}\sim\GG_1$, $\Gamma$ se voit de manière syntaxique
comme l'extension
\[
  \Gamma=(x,y,z\of V,f\of E(x,y),g\of E(y,z)) \tto (x,y,z\of V,f\of
  E(x,y))=\ft\Gamma
\]
du contexte $\ft\Gamma$ avec la variable $g\of E(y,z)$.

La \emph{théorie algébrique à types dépendants des petites catégories} introduit
d'abord une \defn{signature de termes} $\bF_{\Cat}$ sur la signature
de types $\bS_{\GG_1}$, avec deux opérations
\[
  \bF_{\Cat} \eqdef
  \left\{
    x\of V\vdash i(x)\type E(x,x)\quad,\quad x,y,z\of V,f\of E(x,y),g\of
    E(y,z)\vdash c(g,f)\type E(x,z) 
  \right\} ,
\]
puis un ensemble $\bE_{\Cat}$ de trois \defn{équations}
\begin{align*}
  x,y\of V,f\of E(x,y)
  &\vdash c(i(y),f) = f \type E(x,y)\\
  x,y\of V,f\of E(x,y)
  &\vdash c(f,i(x)) = f \type E(x,y)\\
  x,y,z,w\of V, f\of E(x,y), g\of E(y,z),
  h\of E(z,w)
  &\vdash c(h,c(g,f))= c(c(h,g),f) \type E(x,w)
\end{align*}
sur la signature de termes $\bF_\Cat$. Un modèle de la théorie
$\bT_\Cat\eqdef(\bF_\Cat,\bE_\Cat)$ sur la signature $\bS_{\GG_1}$ est un graphe
(un modèle de $\bS_{\GG_1}$)
muni de fonctions d'ensembles encodées par les opérations de $\bF_\Cat$, qui
satisfont aux équations dans $\bE_\Cat$. Il n'est pas difficile de voir que les
modèles de $\bT_\Cat$ sont exactement les petites catégories.

Tout comme les petites catégories, de nombreuses autres familles de structures
algébriques bien connues peuvent être \emph{classifiées} par des théories
algébriques à types dépendants. Par exemple, les $\omega$\=/catégories
(respectivement, les $n$\=/catégories) sont classifiées par une théorie sur la
signature de types dépendants des \emph{globes} (respectivement, des globes de
dimension $\leq n$). Un autre exemple est la famille des opérades planaires
colorées, classifiée par une théorie sur la signature de types
dépendants des \emph{corolles} ou
\emph{arbres élémentaires} planaires.

Cette idée d'associer la dépendance entre types à la cellularité intrinsèque aux
structures algébriques n'est pas nouvelle --- elle est clairement présente dans
la théorie de \og FOLDS \fg{} de \cite{makkai1995folds}, ainsi que dans la
structure groupoïdale et $\omega$\=/groupoïdale des types d'identité dans
\cite{hofmann1998groupoid} et \cite{lumsdaine2009weak,van2011types}. Plus
récemment, \cite[App. A]{brunerie2016homotopy} propose une théorie à types
dépendants classifiant les $\omega$\=/groupoïdes faibles, et
\cite{finster2017type,benjaminthese2020} travaillent avec une théorie algébrique
à types dépendants classifiant les $\omega$\=/catégories faibles. Néanmoins, il
reste à établir une définition {générale}, avec une algèbre universelle
associée, des théories algébriques à types dépendants et les structures
algébriques qu'elles classifient, et cela constitue le premier but de cette
dissertation.

Dans le chapitre \ref{chap:contextual-categories} de cette dissertation, nous
proposons plusieurs définitions de ce que nous entendons par une \defn{théorie
  algébrique à types dépendants}, en utilisant la correspondance entre les types
dépendants et la cellularité, et dans le théorème
\ref{thm:classification-dep-alg-theories} de classification, nous prouvons que
ces définitions sont équivalentes. Nous donnons dans la section
\ref{sec:C-ctxl-cats-syntax} une description syntaxique des théories algébriques
à types dépendants, qui forment une sous-classe \emph{stricte} des \og théories
algébriques généralisées \fg{} de Cartmell \cite{cartmell1978generalised}.

À l'aide du théorème de classification, nous reconnaissons immédiatement
(section \ref{sec:examples-C-cxl-cats}), plusieurs théories algébriques à types
dépendants. Plus loin (théorème \ref{thm:lfp-cat-models-of-C-sorted-theory}),
nous montrons que \emph{toute} catégorie localement finiment présentable est la
catégorie des modèles (dans les ensembles) d'une théorie algébrique à types
dépendants. Nous concluons que les théories algébriques à types dépendants ont
la même puissance expressive (quant à leurs modèles dans les ensembles) que les
théories \emph{essentiellement algébriques} ou les \emph{esquisses projectives}
(corollaire \ref{thm:morita-equiv-gat-eat-C-cxl-cats}).

\subsection*{Multicatégories cellulaires}
Conservons l'exemple de la théorie $\bT_\Cat$ des petites catégories,
et observons que les contextes (qui se trouvent à gauche du symbole ``turnstile''
$\vdash$) dans chaque déclaration de $\bF_\Cat$ et de $\bE_\Cat$ sont tous des
graphes (via la correspondance $\bS_{\GG_1}\sim\GG_1$), à savoir des complexes
cellulaires dans $\psh{\GG_1}$. Cela permet de voir chaque déclaration dans
$\bF_\Cat$ comme une \emph{opération} ou un \emph{multimorphisme} d'une \og
multicatégorie cellulaire \fg{}, prenant en entrée un
graphe (fini) et dont la \og cellule \fg{} de sortie est un représentable. Par
exemple, les opérations $c(-,-)$ et $i(-)$ peuvent être visualisées comme suit
\[
  \begin{tikzpicture}
    \tikzset{scale=.8}
    \node [treenode] (b) at (0, 0) [label = left : $c$] {};
    \node (e1) at ($(b) + (-1, .7)$) [above] {\small $ x$};
    \node (e3) at ($(b) + (1, .7)$) [above] {\small $w$};
    \node (e2) at ($(e1)!.5!(e3)$) {\small $y$};
    \node (t) at ($(b) + (0, -.7)$) [below] {{\small$E(x,w)$}};
    \draw (b) -- (e1);
    \draw (b) -- (e3);
    \draw[->] (e1) -- (e2);
    \draw[->] (e2) -- (e3);
    \draw (b) -- (t);
  \end{tikzpicture}
  \qquad\qquad
  \begin{tikzpicture}
    \tikzset{scale=.8}
    \node [treenode] (b) at (0, 0) [label = left : $i$] {};
    \node (e1) at ($(b) + (0,.7)$) [above] {\small $x$};
    \node (t) at ($(b) + (0, -.7)$) [below] {{\small$E(x,x)$}};
    \draw (b) -- (e1);
    \draw (b) -- (t);
  \end{tikzpicture}
\]
et le terme $c(f,i(x))$ correspond à l'arbre d'opérations suivant.
\[
  \begin{tikzpicture}
    \tikzset{scale=.8}
    \node [treenode] (b) at (0, 0) [label = left : $c$] {};
    \node (e1) at ($(b) + (-1, .7)$) [above] {\small $ x$};
    \node (e3) at ($(b) + (1, .7)$) [above] {\small $w$};
    \node (e2) at ($(e1)!.5!(e3)$) {\small $x$};
    \node (t) at ($(b) + (0, -.7)$) [below] {{\small$E(x,w)$}};
    \node [treenode] (b1) at (-.5, 1.9) [label = left : $i$] {};
    \node (e11) at ($(b1) + (0,.7)$) [above] {\small $x$};
    \node (t1) at ($(b1) + (0, -.7)$) [below] {};
    \draw (b1) -- (e11);
    \draw (b1) -- (t1);
    \draw (b) -- (e1);
    \draw (b) -- (e3);
    \draw[->] (e1) -- (e2);
    \draw[->] (e2) -- (e3);
    \draw (b) -- (t);
  \end{tikzpicture}
\]
Ce point de vue des théories algébriques à types dépendants comme \og
multicatégories cellulaires cartésiennes \fg{} est une généralisation stricte de
la description, due à Lawvere \cite{Lawvere869} et à Bénabou
\cite{benabou1968structures}, des théories algébriques \og ordinaires \fg{}
comme multicatégories cartésiennes (= catégories avec produits finis). Dans
cette généralisation, les théories algébriques ordinaires sont exactement les
multicatégories cellulaires cartésiennes dont toute cellule est de dimension $0$
(elles sont \og ponctuelles \fg{}, autrement elles n'ont aucune dépendance).

Ce point de vue des théories algébriques à types dépendants comme
multicatégories cellulaires apparaît (implicitement) dans
\cite[Sec. II]{fiore2008second}, sous l'appellation 
\emph{$\Sigma_0$\=/models with substitution}. Plus généralement, la théorie
des \emph{monades à arités} et des \emph{théories à arités} de
\cite{weber2007familial,mellies2010segal,Berger2012} peut se voir comme une
théorie abstraite des multicatégories avec arités quelconques. Cependant, une
théorie générale d'\og espèces cellulaires \fg{} et des multicatégories
correspondant aux théories algébriques à types dépendants est (à ma
connaissance) loin d'être établie.

Dans cette dissertation, nous soulignons (section
\ref{sec:C-collections-theories}) le point de vue des théories algébriques à
types dépendants comme multicatégories cellulaires cartésiennes mais nous
n'établissons pas une théorie générale des multicatégories
cellulaires.\footnote{Je pense qu'une telle théorie générale devrait exister,
  récupérant les \emph{opérades globulaires} comme cas particuliers.}

\subsection*{Structures algébriques supérieures et leurs théories}
La thèse \ref{these:dep-typed-higher-alg-struc} porte sur les structures
algébriques à homotopie près. Les \defn{structures algébriques supérieures}
(``higher algebraic structures'' en anglais) sont la \og bonne\fg{} notion de
structures algébriques dans les espaces. Une recette générale pour décrire une
structure algébrique supérieure (telle que les $(\infty,1)$\=/catégories, les
\oo-opérades, les champs, les spectres, les $A_\infty$\=/~et
$E_\infty$\=/espaces) consiste à donner un diagramme dans l'\oo\=/catégorie des
espaces, ainsi que des opérations satisfaisant des axiomes. Or, ces axiomes sont
encodés par des quotients d'espaces, qui sont des tours infinies de
\emph{cohérences} homotopiques, ce qui rend la définition de structures
algébriques supérieures beaucoup plus subtile que leurs versions à valeurs dans
$\Set$ (les structures \og discrètes\fg{} ou \og $0$\=/tronquées\fg{}).

Par exemple, la composition de morphismes dans une $(\infty,1)$-catégorie (une
\oo\=/catégorie), est associative à homotopie cohérente près. Cela se décrit
habituellement en utilisant des \emph{catégories de modèles}, tel qu'avec la
structure de modèles des quasicatégories sur la catégorie des ensembles
simpliciaux, ou la structure de modèles des espaces de Segal complets sur les
ensembles bisimpliciaux. En effet, la vaste majorité de la théorie des
\oo-catégories n'existe formellement que dans la structure de modèles des
quasicatégories --- cependant, les quasicatégories ne sont pas des modèles
canoniques des \oo-categories. Le problème de donner une théorie syntaxique des
\oo\=/catégories est toujours ouvert (voir
\cite{riehl2017type,finster2018towards,allioux2021types}).

Les \emph{types d'identité de Martin\=/Löf} dans la \emph{théorie des types
  homotopiques} (\emph{HoTT}) sont des généralisations syntaxiques, aux
quotients d'espaces, du prédicat d'égalité pour les quotients d'ensembles. Ceci
se justifie par le fait que HoTT s'interprète dans l'\oo\=/topos des espaces
\cite{kapulkin2012simplicial} (et même dans tout \oo\=/topos
\cite{shulman2019infty1toposes}), de manière à ce que les types d'identité
soient interprétés comme les prédicats d'égalité (c'est\=/à\=/dire les
diagonales itérées) dans l'\oo\=/topos. Nous pourrions donc imaginer qu'une
version homotopique, utilisant les types d'identité, de la théorie algébrique à
types dépendants $\bT_\Cat$ soit un candidat pour une théorie syntaxique des
\oo\=/catégories. Plus généralement, une \defn{théorie algébrique supérieure à
  types dépendants} serait une extension de HoTT avec des types, termes, et
\emph{preuves d'identité}, plus ou moins de la même façon que les théories
algébriques à types dépendants de la thèse \ref{these:dependence-cellularite}.

\subsection*{Modèles homotopiques des théories algébriques}
Dans cette dissertation, nous n'essayons pas de donner une théorie générale des
théories algébriques supérieures à types dépendants, ce qui nécessiterait plus
de travail en amont.

Nous nous restreignons à la description d'une théorie des \emph{modèles
  homotopiques} des théories algébriques à types dépendants, généralisant
certains aspects de la théorie des modèles homotopiques des théories algébriques
ordinaires due à
\cite{Schwede2001stablehtpyAlgTheories,badzioch2002algtheories,Rezk2002simplicialAlgTheories,bergner2006rigidification}.
Même dans ce cas, nous ne parvenons pas à montrer que toute la théorie de
rigidification des \emph{op.~cit.} se généralise à notre cadre (nous le
conjecturons tout de même, voir la section \ref{sec:htpical-models-C-cxl-cats}).
Néanmoins, nous posons cette généralisation partielle comme évocatrice de
l'existence d'une théorie générale des théories algébriques supérieures à types
dépendants.

\subsection*{Théories opétopiques idempotentes}
Les types dépendants sont suffisamment expressifs pour abstraire les \og
formes\fg{} des opérations (et des équations) d'une structure algébrique dans la
signature des types. Par exemple, considérons l'opération de composition de
morphismes
\[
  x,y,z\of V,f\of E(x,y),g\of E(y,z)\vdash c(g,f)\type E(x,z)
\]
dans la signature de termes $\bF_\Cat$ de la théorie $\bT_\Cat$ des petites
catégories. Nous pouvons étendre la signature de types $\bS_{\GG_1}$ avec une
\emph{déclaration de type}
\[
  x,y,z\of V,f\of E(x,y),g\of E(y,z),h\of E(x,z)\vdash C(x,y,z,f,g,h)\ \
  \mathrm{type}
\]
correspondant à la forme de l'opération de composition $c(-,-)$. La
forme du type dépendant $C$ peut être vue comme la cellule suivante.
\begin{center}
  \begin{tikzpicture}
    \node (a) at (0, 0){$x$};
    \node (b) at (1.5, 0){$y$};
    \node (c) at (.75, 1){$z$};
    \draw[->] (a) -- (c) node[midway,left]{\small$f$};
    \draw[->] (c) -- (b) node[midway,right]{\small$g$};
    \draw[->] (a) -- (b) node[midway,below]{\small$h$};
    \draw (.75, .35) node {$\Downarrow$};
  \end{tikzpicture} 
\end{center}
Si l'on écrit $\bS_{\GG_1}'$ l'extension de $\bS_{\GG_1}$ avec la déclaration de
type précédente, alors la signature de types $\bS_{\GG_1}'$ correspond à la
catégorie
\[
  \GG_1'\eqdef
  \left\{
    \begin{tikzcd}
      D^0\ar[r,shift right] \ar[r,shift left]
      &D^1\ar[r]\ar[r,shift right] \ar[r,shift left]
      &C
    \end{tikzcd}
  \right\},
\]
et la catégorie $\Cat$ admet un foncteur monadique $\Cat\to\psh{\GG_1'}$ qui
envoie une petite catégorie $A$ sur un préfaisceau dont la fibre au-dessus de
l'objet $C$ est l'ensemble des triangles commutatifs dans la catégorie $A$.
L'opération de composition des petites catégories peut ainsi être \og
opétopifiée\fg{} en une cellule représentable de la signature de types. De plus,
la théorie $\bT_\Cat$ peut être complètement opétopifiée en intégrant toutes ses
opérations et ses équations dans une signature de types $\bS_{\OO_{\leq 3}}$ qui
est une extension de la signature $\bS_{\GG_1}$. La signature $\bS_{\OO_{\leq
    3}}$ correspond à une catégorie directe $\OO_{\leq 3}$ (la catégorie des
\emph{opétopes} de dimension $\leq 3$), telle que $\GG_1\subto\OO_{\leq 3}$ en
est une sous-catégorie pleine. Par ailleurs, il existe un foncteur monadique
$\Cat\to\psh{\OO_{\leq 3}}$ qui est \emph{pleinement fidèle}~; ainsi, tout ce
qui est ostensiblement de la \og structure\fg{} d'une petite catégorie 
s'exprime par des \og propriétés\fg{} d'un préfaisceau sur $\OO_{\leq 3}$. Enfin,
$\Cat$ est la catégorie des modèles d'une \emph{esquisse projective finie} sur
la catégorie $\OO_{\leq 3}\op$. Cette nouvelle théorie algébrique à types
dépendants sur la signature $\bS_{\OO_{\leq 3}}$ (dont la catégorie des modèles
est toujours $\Cat$) est un exemple d'une \defn{théorie opétopique idempotente}.

Le processus d'\og opétopification\fg{} a été introduit pour les opérades
symétriques colorées dans \cite{Baez1998}. Un processus similaire est
fondamental dans la description des $\deux{\mathit L}$\emph{\=/structures} de
\cite{makkai1995folds}.

Dans cette dissertation, nous étudions une classe de théories opétopiques
idempotentes obtenues à partir de la catégorie $\OO$ des opétopes. Nous montrons
l'existence d'une bonne notion de modèle \emph{à homotopie cohérente près} pour
chacune de ces théories.

\subsection*{\texorpdfstring{\oo}{Infini}-catégories localement présentables}
La collection de toutes les structures algébriques supérieures d'une même
famille, et les morphismes entre elles, forment une \emph{\oo\=/catégorie
  localement présentable}. La théorie des \oo\=/catégories localement
présentables généralise parfaitement celle des $1$\=/catégories localement
présentables, et nous conjecturons que tout comme les théories algébriques à
types dépendants correspondent aux $1$\=/catégories localement présentables, une
bonne notion de théorie algébrique supérieure à types dépendants donnera lieu
aux \oo\=/catégories localement présentables. Il découle directement des
définitions que les modèles homotopiques d'une théorie algébrique à types
dépendants forment une \oo\=/catégorie localement présentable, ce qui peut
servir de point de départ. Cette conjecture est aussi soutenue par les
résultats de \cite{kapulkin2019internal}, où l'on montre que les
\oo\=/catégories avec limites finies correspondent aux \emph{tribus}, qui sont
une abstraction des théories à types dépendants avec types d'identité.

Dans la partie \ref{part:part2} de cette dissertation, nous étudions les
\oo\=/catégories localement présentables, et nous développons une théorie des
\emph{pré\=/modulateurs} permettant de construire tout système de factorisation
accessible sur une \oo\=/catégorie localement présentable par l'itération une
\og construction-plus\fg{}.

\section*{Organisation}
Cette dissertation s'organise en deux parties.

\medskip La partie \ref{part:part1} concerne les thèses
\ref{these:dependence-cellularite} et \ref{these:dep-typed-higher-alg-struc}, et
forme le contenu principal de cette dissertation.

Le chapitre \ref{chap:contextual-categories} développe la théorie des théories
algébriques à types/sortes dépendants, et commence avec une description
détaillée de la correspondance entre les types dépendants et une notion
abstraite de \og cellularité\fg{}. La definition clef est celle de
\emph{catégorie directe localement finie}, qui est une reformulation équivalente
de la définition de \emph{catégorie simple} de \cite{makkai1995folds}. Nous
utilisons cette reformulation pour donner plusieurs définitions des théories à
types dépendants qui nous intéressent. Notre approche consiste à introduire la
notion de \emph{$\lfdex$\=/catégorie contextuelle} (définition
\ref{def:C-cxl-cat}), outil algébrique nous permettant d'éviter des
raisonnements syntaxiques. Nous démontrons le théorème
\ref{thm:classification-dep-alg-theories} de classification qui généralise
strictement celle de Lawvere pour les théories algébriques ordinaires. Ce
résultat est alors utilisé dans la section \ref{sec:examples-C-cxl-cats} pour
reconnaître plusieurs théories algébriques à types dépendants. Enfin, nous
donnons une caractérisation syntaxique des théories algébriques à types
dépendants, et nous montrons qu'elles forment une sous-classe stricte des
théories algébriques généralisées de Cartmell
\cite{cartmell1978generalised,cartmell1986}. Un autre résultat important de ce
chapitre est la proposition \ref{prop:ctxl-functor-cxl-C-initial}, qui démontre
l'initialité de la catégorie contextuelle libre sur une signature de types,
en se basant sur les propositions \ref{prop:cxl-C-fin-limits} à
\ref{prop:cxl-C-codescent}.

Le chapitre \ref{cha:models-of-C-contextual-categories} développe la théorie des
modèles homotopiques, dans les espaces, des théories algébriques à types
dépendants. Nous commençons par une classification des catégories des modèles
dans les ensembles (les espaces $0$\=/tronqués) des théories algébriques à types
dépendants comme étant exactement les catégories localement finiment
présentable. On en déduit que les classes des théories algébriques à types
dépendants, des théories essentiellement algébriques, et des théories
algébriques généralisées sont toutes \emph{Morita\=/équivalentes} (théorème
\ref{thm:lfp-cat-models-of-C-sorted-theory} et corollaire
\ref{thm:morita-equiv-gat-eat-C-cxl-cats}). Après des rappels sur les catégories
de modèles simpliciales, nous démontrons le théorème
\ref{thm:flasque-C-spaces-rigidification} de \emph{rigidification} pour les
modèles homotopiques dans les espaces de la théorie algébrique à types
dépendants initiale sur une catégorie directe localement finie $\lfdex$. Enfin,
nous montrons l'existence d'une structure de modèles pour les modèles
homotopiques dans les espaces de toute théories algébrique à types dépendants
(théorème \ref{thm:partial-soln-rigid-conjecture}), et nous conjecturons
l'existence d'un théorème de rigidification général. Dans ce chapitre nous
obtenons aussi d'autres résultats~:
la construction d'une structure de modèles \emph{flasque}
intermédiaire sur la catégorie des préfaisceaux simpliciaux sur la
$\lfdex$\=/catégorie contextuelle initiale (théorème
\ref{thm:flasque-model-struc-existence}) et la description des
$\lfdex$\=/espaces comme faisceaux d'\oo\=/groupoïdes sur la
$\lfdex$\=/catégorie contextuelle initiale (théorème
\ref{prop:C-spaces-model-struct-equals-Cech}).

Le chapitre \ref{chap:opetopic-theories} étudie une classe de théories
idempotentes sur la catégorie directe localement finie $\OO$ des
\emph{opétopes}. Nous commençons en rappelant la construction de la catégorie
$\OO$ due à \cite{Kock2010,hothanh18}. Nous définissons une famille de monades
adjointes à droite paramétriques $\optPolyFun^n$ dont les algèbres sont les
\emph{algèbres opétopiques} (définition \ref{def:opetopic-algebra}). Nous
montrons que chaque catégorie $\oAlg$ d'algèbres opétopiques admet un foncteur
$\oAlg\subto\pshO$ vers la catégorie $\pshO$ des ensembles opétopiques, appelé
\emph{le foncteur du nerf opétopique}, qui est monadique, pleinement fidèle, et
accessible. De plus, $\oAlg$ est la catégorie des modèles d'une esquisse
projective finie sur la catégorie $\OO\op$ (théorème \ref{th:nerve-theorem-O}),
ce qui implique qu'elle est la catégorie des modèles d'une \emph{théorie
  opétopique idempotente} (théorème
\ref{thm:OAlg-models-of-idempotent-O-theory}). En particulier, la catégorie des
petites catégories, la catégorie des opérades colorées planaires, et la
catégorie des combinades colorées de Loday \cite{Loday2012a} sont des catégories
d'algèbres opétopiques --- elles sont donc chacune classifiée par une théorie
opétopique idempotente. Enfin, nous montrons qu'une technique due à Horel
\cite{horel2015model} s'applique aux algèbres opétopiques, nous permettant de
donner une définition d'une structure de modèles pour les \emph{algèbres
  opétopiques à homotopie cohérente près}, et de montrer qu'elle est équivalente au sens
de Quillen à une structure de modèles sur les \emph{algèbres opétopiques
  simpliciales} (théorème \ref{thm:rigidification-Oalg}). Dans les cas
particuliers de $\Cat$ et de $\Opd$, nous récupérons la structure de modèles de
Rezk pour les espaces de Segal et la version planaire de la structure de modèles
de Cisinski--Moerdijk pour les espaces de Segal dendroïdaux.

\medskip La partie \ref{part:part2} de cette dissertation étudie les
localisations accessibles et les systèmes de factorisation accessibles dans les
\oo\=/catégories localement présentables.

Le chapitre \ref{chap:fact-systems} décrit la théorie des systèmes de factorisation
(orthogonale) dans les \oo\=/catégories en utilisant le tenseur et
l'enrichissement fournis par les constructions de \og pushout-product\fg{} et de
\og pullback-hom\fg{} dans les \oo\=/catégories de flèches. Nous utilisons ce
formalisme pour récupérer plusieurs résultats sur les modalités et les modalités
exactes à gauche ou \og lex\fg{} (propositions \ref{prop:stable-lex-FS},
\ref{prop:mod-lccc}, \ref{prop:mod-slex}, et \ref{prop:lex-loc=lex-mod}).

Le chapitre \ref{chap:SOA} développe la théorie des \emph{pré\=/modulateurs},
commençant par une modification légère (utilisant le pushout-product et le
pullback-hom) de \emph{l'argument du petit objet} de Kelly pour les
$1$\=/catégories localement présentables, nous permettant de l'adapter à la
construction de systèmes de factorisation accessibles dans les \oo\=/catégories
localement présentables (théorème \ref{thm:!SOA2}). Nous montrons que dans le
cas particulier de la faisceautisation associée à une topologie de Grothendieck,
la construction de Kelly se simplifie en une \emph{construction-plus} sur les
préfaisceaux (lemme \ref{lem:plus-is-plus-for-topology}), et qu'il en est de
même pour tout système de factorisation accessible engendré par un
\emph{pré\=/modulateur} sur une \oo\=/catégorie localement présentable (théorème
\ref{thm:plus-construction}). De plus, tout système de factorisation accessible
sur une \oo\=/catégorie localement présentable est engendré par un
pré\=/modulateur (proposition \ref{prop:modulator-completion}). Nous utilisons
la théorie des pré\=/modulateurs pour définir les \emph{modulateurs} et les
\emph{modulateurs lex}, dont les constructions plus respectives engendrent les
modalités et les modalités exactes à gauche (théorèmes \ref{thm:+modality} et
\ref{thm:lex+construction}). Nous montrons que chaque localisation exacte à
gauche et accessible d'un \oo\=/topos est engendrée par un modulateur lex
(proposition \ref{prop:lex-mod-4-loc-lex}),
ainsi les modulateurs lex sont une bonne généralisation des topologies de
Grothendieck aux \oo\=/topoï. Enfin, nous montrons que la construction-plus
associée à un modulateur lex sur un \oo\=/topos converge en $(n+2)$ itérations
sur chaque objet $n$\=/tronqué de l'\oo\=/topos
(proposition \ref{prop:truncated-sheaf}).


\end{otherlanguage}




\addtocontents{toc}{\protect\setcounter{tocdepth}{2}}
\part{Dependently sorted algebraic theories}
\label{part:part1}

\chapter{Contextual categories as monoids in collections}
\label{chap:contextual-categories}

In this chapter, we substantiate \cref{thes:dependency-cellularity} of this
dissertation, by showing that type dependency in theories corresponds to a
notion of cellularity inherent in the corresponding algebraic structures. We do
so by describing a theory of \emph{dependently sorted/typed algebraic theories}, that
strictly generalise the \emph{(multisorted) algebraic theories} of
\cite{Lawvere869} and \cite{benabou1968structures}. We give a
classification theorem for dependently sorted algebraic theories
(\cref{thm:classification-dep-alg-theories}), and use it to detect many examples.
In~\cref{sec:C-ctxl-cats-syntax} we describe the class of dependently typed
syntactic theories \emph{à la} Martin-Löf that correspond to our notion of dependently
typed algebraic theory.

\section{Locally finite direct categories}
\label{sec:locally-finite-direct-cats}
\begin{para}
 \label{para:lfd-cat-intro} 
 A small category $C$ is \emph{finite} if its set $\ob {C\arr}$ of morphisms is
 a finite set. We say that $C$ is \defn{locally finite} if for all $c$ in $C$,
 the slice category $C\slice c$ is finite. Namely, in every cartesian square
 in $\Cat$
 of the
 following form, where $1$ is the terminal category and $\t\colon C\arr\to C$ is
 the codomain functor, the category $A$ is finite.
 \[
    \begin{tikzcd}[sep=scriptsize]
      A\pbmark \ar[r] \ar[d] & C\arr \ar[d,"\t"]\\
      1 \ar[r] &C
    \end{tikzcd}
  \]
  Every locally finite category has finite hom-sets---that is, it is enriched
  over the category $\fin$ of finite sets---but the converse is not necessarily
  true. A locally finite category is finite if and only if its set of objects is
  a finite set.
  
  Let $C$ be locally finite. For every presheaf $X$ in $\psh C\eqdef
  \fun{C\op}\Set$, the category of elements $C/X$ is locally finite, since a
  morphism $p\colon A\to C$ in $\Cat$ is a discrete fibration if and only if the
  following square in $\Cat$ is cartesian.
  \[
    \begin{tikzcd}[sep=scriptsize]
      A\arr \ar[r,"p\arr"] \ar[d,"\t"'] &{C}\arr \ar[d,"\t"]\\
      A \ar[r,"p"] & C
    \end{tikzcd}
  \]
  Let $C\in\Cat$ and $c\in C$. We define a \defn{saturated cover}
  of $c$ to
  be a (not necessarily full) subcategory $J$ of $C \slice c$ such that every
  non-identity morphism $c' \to c$ in $C$ factors through some $j \to c$ in $J$.
\end{para}


\begin{example}
  \label{ex:saturated-covers}
  For any $C\in\Cat$ and $c\in C$, the full subcategory $C \slice c ^- \subset
  C\slice c$ obtained by removing the identity morphism $1_c\colon c \to c$ (the
  terminal object of $C\slice c$), as well as the the full subcategory
  $\{1_c\}\subset C\slice c$ consisting of \emph{only} the identity morphism,
  are both saturated covers of $c$.
\end{example}

\begin{remark}
 \label{rem:sat-cover-alt-defn} 
  Let $C\in\Cat$ and $c\in C$. A subcategory $J$ of $C\slice c$ is a saturated
  cover if and only if the functor $\coprod_{j \to c \in J} C\slice j \to C\slice c$
  is surjective on objects and morphisms of $C\slice c ^-$.
\end{remark}



\begin{definition}
  \label{def:direct-cat}
  A small category $C$ is \defn{direct} if the binary relation on $\ob C$
  \[
    c < d \quad\Leftrightarrow\quad \text{ there exists a non-identity morphism
    } c \to d
  \]
  is well-founded, namely there are no infinite sequences $ \ldots < c _2 < c _1
  < c _0 $ in $C$.
\end{definition}

\begin{proposition}
  If $C$ is a direct category and $f\colon A\to C$ is a functor that reflects
  identity morphisms, then $A$ is a direct category.
\end{proposition}
\begin{proof}
  A functor reflects identities if and only if it preserves non-identity
  morphisms. Hence for $\ldots<a_1<a_0$ in $A$, we have
   $\ldots<fa_1<fa_0$ in $C$.
\end{proof}

\begin{corollary}
   \label{cor:presheaf-elts-direct} 
  If $C$ is direct and $X\in\psh C$ is a presheaf on $C$, then its category of
  elements $C/X$ is direct.
\end{corollary}
\begin{proof}
  Any discrete fibration in $\Cat$ reflects identity morphisms.
\end{proof}

\begin{notation}
  We write $\Dir$ for the full subcategory of $\Cat$ consisting of all small
  direct categories, and $\deux\Dir$ for the same full sub-$2$\=/category of
  $\deux\Cat$. We write $\Dir\lfin\subset \Dir$ and
  $\deux\Dir\lfin\subset\deux\Dir$ for the locally finite direct categories.
\end{notation}

\begin{remark}
  A direct category has no non-identity endomorphisms.
\end{remark}

\begin{remark}
 \label{rem:dir-cat-slice-join-1} 
 If $C$ is a direct category, then for every $c$ in $C$ the inclusion $C\slice
 c^-\subset C\slice c$ exhibits the slice category $C\slice c$ as the \emph{free
   $C\slice c^-$\=/cocone}. In other words, $C\slice c$ is obtained from $C\slice
 c^-$ by freely adding a terminal object, namely $C\slice c = C\slice c^-\star 1$
 (where $-\star-$ is the \emph{join} of categories
 \cite[\textsection{3.1}]{joyal2008theory}).
\end{remark}

\begin{lemma}
  \label{lem:fin-sat-cover-fin-slice}
  Let $C\in\Dir$ and $c\in C$. If $J$ is a finite saturated cover of $c$, and if
  for every $j \to c$ in $J$, $C\slice j$ is finite, then $C\slice c$ is finite.
\end{lemma}
\begin{proof}
  By~\cref{rem:sat-cover-alt-defn}, $\coprod_{j \to c \in J} C\slice j \to
  C\slice c$ is surjective on objects and morphisms of $C\slice c^-$, thus
  $C\slice c^-$ is finite. Hence $C\slice c = C\slice c^-\star 1$ is finite.
\end{proof}

\begin{remark}
  A small category $C$ is ``simple'' (one-way, skeletal, finitely branching) in
  the sense of Makkai's FOLDS \cite{makkai1995folds} if and only if $C\op$ is a
  \lfd~category.
\end{remark}

\begin{proposition}
  \label{prop:loc-finite-finite-cover}
  Let $C \in \Dir$. Then $C$ is locally finite if and only if every $c$ in $C$ has a
  finite saturated cover consisting of non-identity morphisms $j\to c$.
\end{proposition}
\begin{proof}
  The ``only if'' direction is obvious, since for all $c\in C$, if $C\slice c$ is finite
  then $C\slice c^-$ is such a finite saturated cover.
  Conversely, as $C$ is direct, we
  proceed by induction on $c\in\ob C$, using the well-founded relation $<$.
  In the base case, there are no non-identity morphisms $c' \to c$, and $C\slice
  c = 1$ (the terminal category). In the induction step, consider a finite
  saturated cover $J$ of $c$ consisting of non-identity morphisms. By the
  induction hypothesis, for each $f \colon j \to c$ in $J$, $C \slice j$ is
  finite, and we may use \cref{lem:fin-sat-cover-fin-slice} to conclude.
\end{proof}

\begin{example}
 \label{exa:lfd-cats} 
 Many categories of finite-dimensional, finite cellular ``shapes'' are \lfd, or have wide
 \lfd~subcategories that determine the shapes represented by the objects of such
 categories. We list a few examples (that they are so follows in each case from
 \cref{prop:loc-finite-finite-cover}).
 \begin{enumerate}
 \item
   \label{exa:lfd-cats-set}
   Any set $S$ (seen as a discrete category).
 \item
   \label{exa:lfd-cats-omega}
   The ordinal $\omega$ (seen as a totally ordered poset).
 \item
    \label{exa:lfd-cats-grph}
    The category $\GG_1 = \{0\rightrightarrows 1\}$ with two objects and two
    parallel non-identity arrows.
  \item
    \label{exa:lfd-cats-globe}
    The  category of \emph{globes} $\GG$ \cite[Def. 1.4.5]{Leinster2004}.
  \item
    \label{exa:lfd-cats-eltree}
    The category $\elTree_p$ of \emph{planar elementary trees} 
    \cite[2.4.4]{Kock2011} or \emph{planar corollas}.
  \item
    \label{exa:lfd-cats-opetopes}
    The category $\OO$ of \emph{opetopes} (\cref{def:o}).
  \item
    \label{exa:lfd-cats-reedy}
    Every Reedy category $R$ has a wide subcategory $R'$ that is direct. In
    many (if not most) well-known examples, $R'$ is also locally finite,
    such as:
    \begin{enumerate}
    \item
      \label{exa:lfd-cats-simplex}
      $R=\bDelta$, the \emph{simplex} category ($\bDelta'$ is called the
      \emph{semi-simplex} category),
    \item
      \label{exa:lfd-cats-dendrex}
      $R=\bOmega_p$, the \emph{planar dendroidal} category \cite[Def.
      2.2.1]{Moerdijk2010}, ($\bOmega_p'$ is called the category of \emph{planar
        semi-dendrices}),
    \item
      \label{exa:lfd-cats-theta}
      $R=\bTheta$, Joyal's cell category \cite{joyal1997disks},
    \end{enumerate}
    where in each case $R'$ is the wide subcategory of monomorphisms.
  \end{enumerate}
\end{example}





\begin{definition}
  \label{def:loc-fin-cat-dimension}
  Let $C$ be a \lfd~category. The \defn{dimension} $\height(c) \in \NN$ of each $c$ in
  $C$ is the length of the longest chain $c_n<\ldots < c_1 < c_0=c$. The
  {dimension}
  of $C$ is $\height(C) \eqdef \min\{n\leq \omega \mid \forall c\in
  C, n>\height(c)\}$. The dimension $\height(X)$ of a presheaf $X\in\psh C$ is
  the dimension of its category of elements $C/X$.
\end{definition}
\begin{remark}
  The empty category $\emptyset$ is the only \lfd~category of
  dimension $0$. There are of course \lfd~categories of dimension $\omega$, such
  as $\GG,\bDelta',\bTheta'$ (and $\omega$ itself, seen as a poset).
\end{remark}

\begin{remark}
  In \cite{makkai1995folds}, dimension is called \emph{height}. Our terminology
  uses the idea that objects of a \lfd~category represent finite-dimensional
  cells, morphisms represent inclusions of cells of lower dimension as
  sub-cells of higher dimension, and every cell only has finitely many sub-cells.
\end{remark}

\begin{lemma}
 \label{lem:lfd-cat-grading} 
  Every $C$ in $\Dir\lfin$ admits a grading $\emptyset = C_0 \subset C_1
  \ldots$ by its full \lfd~subcategories of finite dimension, such that $C
  = \bigcup_n C_n$.
\end{lemma}

\begin{lemma}
  \label{lem:product-loc-fin-cat}
  The full subcategory $\Dir\lfin\subset\Cat$ is closed under finite products and
  coproducts.
\end{lemma}

\begin{remark}
  \label{rem:loc-k-small-dir-cats}
  The cardinal $\aleph_0$ (which we will write as $\omega$) can be replaced by
  any regular cardinal $\kappa$, giving definitions of \emph{locally
    $\kappa$-small} and of \emph{locally $\kappa$-small direct} categories.
  Since every small category is locally $\kappa$-small for some regular
  $\kappa$, every $C\in\Dir$ is locally $\kappa$-small direct for some
  regular $\kappa$.
\end{remark}

We fix a \lfd~category $\lfdex$ for the rest of this chapter.

\begin{para}
  \label{para:boudaries-finite-cell-cplxs}
  The \defn{boundary} of an object $c$ in $\lfdex$, denoted $\partial
  c\in\psh\lfdex$, is the sub-representable presheaf $\partial c\subto c$ that
  is the colimit of the diagram $\lfdex \slice c ^- \subset \lfdex\slice c \to
  \lfdex \subto \psh \lfdex$ obtained by composing the domain functor
  $\lfdex\slice c^-\to \lfdex$ (see \cref{ex:saturated-covers}) with the Yoneda
  embedding. We define the set of \defn{boundary inclusions} of $\lfdex$ to be
  the following set of maps in $\psh\lfdex$.
  \[
    I_\lfdex\eqdef\left\{\delta_c\colon \partial c\subto c\mid
      c\in\lfdex\right\}
  \]
  As $\lfdex$ is locally finite, $\partial c$ is a finite colimit of
  representables for every $c$ in $\lfdex$, and is therefore a finitely
  presentable object of $\psh\lfdex$.
\end{para}


\begin{lemma}
  \label{lem:fin-colim-lfd-fin-prshf}
  Let $J\in\Cat$ be a finite category, let $f \colon J\to \lfdex$ be a functor,
  and let $X$ be the colimit of $(J\xto{f}\lfdex\subto\psh\lfdex)$. Then, for
  every $c$ in $\lfdex$, $X_c$ is a finite set. Moreover, there are only
  finitely many $c\in\lfdex$ such that $X_c \neq \emptyset$.
\end{lemma}
\begin{proof}
  For any $c$ in $\lfdex$, $X_c = \colim_{j\in J} \lfdex (c, fj)$. Hence the
  first claim follows since $J$ is finite and since $\lfdex$ has finite hom-sets
  (\cref{para:lfd-cat-intro}).
  For the second claim, remark that $\coprod_{j\in J}\lfdex\slice {fj}$ is
  finite (since $J$ is finite and $\lfdex$ is locally finite). Then, suppose
  that $c$ is not in the image of $\coprod_{j\in J}\lfdex\slice {fj} \to
  \lfdex$. Since $\lfdex$ is direct, $\lfdex (c, fj) = \emptyset$ for all $j$ in
  $ J$, and so $X_c=\emptyset$.
\end{proof}

\begin{proposition}
  \label{cor:lfd-fin-pres-fin-cat-of-elts}
  For all $X$ in $\psh\lfdex$, $X$ is finitely presentable if and only if its category of
  elements $\lfdex/X$ is finite.
\end{proposition}
\begin{proof}
  The ``if'' direction is obvious. Conversely, let $X$ be finitely presentable.
  Every finitely presentable presheaf can be written as a finite colimit of
  representables, hence by \cref{lem:fin-colim-lfd-fin-prshf},
  $\ob{\lfdex/X}=\coprod_{c\in\lfdex}X_c$ is a finite set. This implies that $\lfdex/X$ is finite (see
  \cref{para:lfd-cat-intro}).
\end{proof}

\begin{corollary}
 \label{cor:fin-pres-implies-fin-dim} 
 Any finitely presentable object in $\psh\lfdex$ is of finite dimension.
\end{corollary}

\begin{corollary}
  \label{cor:fin-C-pback-fin-lims}
  The category $\fin_\lfdex$ of finitely presentable objects of $\psh\lfdex$ has
  pullbacks and non-empty finite products. It has finite limits if and only if
  $\lfdex$ is finite.
\end{corollary}


\begin{definition}
  Let $\cC$ be a category with an initial object $\emptyset$, and let $S$ be a
  class of morphisms of $\cC$. A \defn{relative $S$-cell complex} is the
  (transfinite) composite of a sequence $X\colon \kappa\to\cC$ (for some ordinal
  $\kappa$ seen as a totally ordered poset) such that each morphism
  $X_\lambda\to X_{\lambda+1}$ is a pushout of some morphism in $S$. It is an
  \defn{$S$-cell complex} if $X_0=\emptyset$, and it is \defn{finite} if
  $\kappa$ is finite.
\end{definition}

\begin{proposition}
  \label{prop:mono-rel-cell-complex}
  Let $f \colon X \mono Y$ be a monomorphism in $\psh \lfdex$. Then $f$ is the
  transfinite composite of an $\omega$-long sequence $X_{(-)}:\omega \to
  \psh\lfdex$ where each morphism $X_i\to X_{i+1}$ is a pushout of $\coprod_{c
    \in S}(\partial c \subto c)$ for some family of objects of $\lfdex$.
  Moreover, if $\coprod_{c \in \lfdex}(Y_c - f_c(X_c))$ is finite, then $f$ is a
  finite relative $I_\lfdex$-cell complex.
\end{proposition}
\begin{proof}({\cite[II.3.8]{GabrielZisman1967}}).
  By induction, starting with $X=X_0$, we construct a sequence of inclusions
  $X_0 \xto{i_0} \ldots X_n\xto{i_n} \ldots$, along with canonical maps $g_n
  \colon X_n \to Y$, such that
  \begin{enumerate}
  \item $g_n = g_{n+1}i_n$,
  \item if $\height(c)> n$ then $(i_n)_c\colon
    (X_n)_c \cong (X_{n+1})_c$ is a bijection,
  \item and if $\height (c)<n$, then
  $(g_n)_c\colon(X_n)_c \cong Y_c$ is a bijection.
  \end{enumerate}
  Begin with $X_0 = X$ and $g_0 = f$.
  At stage $X \to \ldots \to X_k$, note that (trivially when $k=0$ and by the
  induction hypotheses when $k>0$) for all
  $c$ such that $\height(c) = k$, there is a bijection between the set $Y_c$ and
  the set of commutative squares $\psh\lfdex\arr (\delta_c,g_k)$.
  We define $i_k$ to be the pushout
  \[
    \begin{tikzcd}[sep=scriptsize]
      {\coprod_{c} \partial c} \ar[r] \ar[d,hook]
      &{X_k} \ar[d,"i_k"] \\
      {\coprod_{c} c} \ar[r] &{X_{k+1}\pomark }
    \end{tikzcd}
  \]
  where the coproduct is indexed over all $c$ such that $\height(c)=k$ and all
  elements of $Z_c\eqdef(Y_c - f_c(X_c))$.
  We define $g_{k+1}$ to be the canonical morphism $X_{k+1} \to Y$. Next, since
  pushouts are calculated pointwise, for $\height(c)=k$, we have $(g_{k+1})_c
  \colon (X_{k+1})_c = X_c \coprod Z_c \cong Y_c$, and for $\height(c)>k$ we
  have $(i_k)_c\colon (X_{k+1})_c = X_c \coprod \emptyset \cong X_c$.
  Clearly, $g_\omega \colon \colim_{i<\omega}X_i \to Y$ is an isomorphism.
  If $\coprod_{c\in \lfdex} Z_c$ is finite, then each $Z_c$ is finite and
  if $k=\max\{d(c) \mid Z_c \neq \emptyset\}$, then $g_{k+1}$ is an isomorphism.
\end{proof}

\begin{proposition}
   \label{cor:fin-cell-cplx-iff-fin-pres} 
    For all $X$ in $\psh \lfdex$,
    $X$ is finitely presentable
    if and only if
    there exists a finite $I_\lfdex$-cell complex $\emptyset \to X_1 \to
    \ldots \to X$.
\end{proposition}
\begin{proof}
  One direction follows from \cref{prop:mono-rel-cell-complex} and
  \cref{lem:fin-colim-lfd-fin-prshf}, and the other follows since every $\partial c$ is
  finitely presentable.
\end{proof}



\section{Cell complexes and cell contexts}

We come to one of the principal definitions of~\cref{chap:contextual-categories,cha:models-of-C-contextual-categories}.
\begin{definition}
  \label{def:Cell-C}
  The category $\cell \lfdex$ of \defn{finite $\lfdex$-sorted cell contexts} has as
  its set of objects a graded set $\ob{\cell\lfdex} \eqdef \coprod_{n\in\NN}
  (\cell\lfdex)_n$
  of finite $I_\lfdex$-cell complexes $\emptyset \to \ldots X_n$
  inductively defined as:
    \begin{itemize}
    \item $(\cell\lfdex)_0$ consists only of the empty presheaf $\emptyset$,
    \item for every $\emptyset \to \ldots X_n$ in $(\cell\lfdex)_n$, $c\in
      C$ and every span $c\xhookleftarrow{\delta_c}\partial c\xto{\gamma} X_n$ in
      $\psh\lfdex$, we make a \emph{choice} of pushout square
      \[
        \begin{tikzcd}[sep=scriptsize]
          \partial c \ar[r,"\gamma"] \ar[d,hook,"\delta_c"']
          &X_n \ar[d,"p_{n+1}"] \\
          c \ar[r,"\gamma.c"] &X_{n+1} \pomark
        \end{tikzcd}
      \]
      giving $\emptyset \to\ldots \to X_n\xto{p_{n+1}} X_{n+1}$ in $(\cell\lfdex)_{n+1}$.
  \end{itemize}
  The morphisms of $\cell\lfdex$ are defined by $\cell \lfdex (\emptyset \to
  \ldots X, \emptyset \to \ldots Y)\eqdef \psh \lfdex (X,Y)$.
\end{definition}

\begin{remark} 
  A similar construction (though not of the category $\cell\lfdex$) can be
  found in H. Gylterud's PhD thesis \cite[E:31-34]{gylterudunivalent}. It seems
  clear that their point of view is very close to the idea of dependently typed syntax as
  operations with cellular arities presented here.\footnote{I would like to thank H. Gylterud
    for pointing this out to me recently.} 
\end{remark}

\begin{proposition}
 \label{lem:cell-C-fin-pres} 
  The forgetful functor $\cell \lfdex \to \psh \lfdex$ is fully faithful, and
  its essential image is the full subcategory $\fin_{\lfdex}$ of finitely
  presentable presheaves on $\lfdex$.
\end{proposition}
\begin{proof}
  Immediate from~\cref{cor:fin-cell-cplx-iff-fin-pres}.
\end{proof}
\begin{corollary}
 \label{cor:cell-C-finite-cocompletion} 
  $\cell\lfdex$ is a completion of $\lfdex$ under finite colimits.
\end{corollary}

\begin{example}
  It is useful to consider a simple example to fix ideas. Let
  $\GG_1\eqdef\{s,t\colon D^0\rightrightarrows D^1\}$ be the category from
  \cref{exa:lfd-cats}\ref{exa:lfd-cats-grph}, so $\psh{\GG_1}$ is the category
  of (directed multi)graphs. Then consider the finite graph $G\in\fin_{\GG_1}$ below.
  \[
    \begin{tikzcd}
      x\ar[r,bend left,"f"] \ar[r,bend right,"g"'] &y
    \end{tikzcd}
  \] There are four distinct isomorphic objects of $\cell{\GG_1}$ that represent
  $G$ via the equivalence $\cell{\GG_1}\simeq \fin_{\GG_1}$ (even though there
  is only one nontrivial automorphism of $G$). This is because when building $G$
  as a cell-complex, we could add $x$ before $y$ or \emph{vice versa}. This is
  best seen syntactically, by considering the \emph{type signature} associated
  to $\GG_1$ (that we will encounter in \cref{sec:C-ctxl-cats-syntax}), which is
  given below.
  \[
    \vdash D^0 \qquad\qquad x\of D^0,y\of D^0\vdash D^1 (x,y)
  \]
  Then any object of $\cell{\GG_1}$ is a \emph{context} over this type
  signature. For instance, two distinct objects of $\cell{\GG_1}$ that both
  represent the previous graph $G$ are:\footnote{The reader familiar with type
    theory will recognise that these two contexts are {not} the same up to
    renaming of variables.}
  \[
    x\of D^0,y\of D^0,f\of D^1(x,y),g\of D^1(x,y) \qquad y\of D^0,x\of D^0,f\of
    D^1(x,y),g\of D^1(x,y).
  \]
  In each case, the cell complex has the form
  \[
    \emptyset\to D^0\to D^0\amalg D^0\to D^1 \to D^1\amalg_{D^0\amalg D^0} D^1 
  \]
  except that in one case, the ``source'' of the two ``arrows'' is added to the
  complex before the ``target'' (and \emph{vice versa} in the other case).
\end{example}

\begin{remark}
  \label{rem:functor-fin-to-cell}
  The construction of \cref{prop:mono-rel-cell-complex}, with a choice of
  ordering of coproducts at each step, lets us defines a fully faithful functor
  $\fin_\lfdex \subto \cell\lfdex$ exhibiting the equivalence of categories.
  Likewise, we can define a choice of inclusion $\lfdex\subto \cell\lfdex$. We
  will assume these choices of inclusion from now on.
\end{remark}

\begin{remark}
    \label{rem:replacing-fin-with-cell}
    The principal reasons for replacing $\fin_\lfdex$ with $\cell\lfdex$ are (i)
    the latter is a small category (while the former is only \emph{essentially}
    small, see \cref{prop:cell-S-is-ctxl-cat}), and (ii) the latter has a
    canonical structure of a \emph{co-contextual category} (while the former
    does not, see \cref{rem:ctx-cat-structure-disrespect-equiv-cats}).
\end{remark}

\begin{remark}
  Several of our constructions are an abstraction---into the language of cell
  complexes---of syntactic constructions that can be found in
  \cite[\textsection{4}]{makkai1995folds}. In particular, Makkai defines
  $\mathrm{Con}[\bK]$ to be the opposite of the syntactic category on a simple
  category $\mathbf K$, and this is nothing but our $\cell{\bK\op}$ (as we will
  show in the sequel). Moreover, Makkai shows the equivalence
  $\mathrm{Con}[\bK]\simeq \fin_{\bK\op}$ by hand.
\end{remark}

\begin{remark}
  \label{rem:loc-k-small-dir-cell}
  If $D$ is a locally $\kappa$-small direct category
  (\cref{rem:loc-k-small-dir-cats}), and if we replace ``finite'' with
  ``$\kappa$-small'' in the definition of $\cell D$, and take colimits at limit
  ordinals, then $\cell D$ is equivalent to the category of $\kappa$-presentable
  objects of $\psh D$.
\end{remark}

\begin{para}
  \cref{cor:cell-C-finite-cocompletion} tells us that if $\cC$ is any finitely
  cocomplete category, precomposition and pointwise left Kan extension along the
  inclusion $\lfdex\subto\cell\lfdex$ together define an equivalence of
  categories $\fun\lfdex \cC \simeq \fun{\cell\lfdex}\cC\rex$ (where
  $\fun{\cell\lfdex}\cC\rex$ is the category of finitely cocontinuous functors).

  However, the ``cellular'' definition of $\cell\lfdex$ allows for a more
  general equivalence of functor categories.
\end{para}

\begin{definition}
 \label{defn:fun-pres-I-cell-cplxes} 
  Let $\cC$ be a category with an initial object. We write
  $\fun{\fin_\lfdex}\cC_{I_\lfdex}$ for the category of functors
  $\fin_\lfdex\to\cC$ that preserve the initial object and pushouts of all
  boundary inclusions, namely that \defn{preserve (finite) $I_\lfdex$\=/cell
    complexes}. We define $\fun{\cell\lfdex}\cC_{I_\lfdex}$ using the
  equivalence $\cell\lfdex\simeq\fin_\lfdex$.
\end{definition}

\begin{proposition}
 \label{prop:rex-iff-cell-extendable} 
  Let $\cC$ be a category with an initial object. Then every $\tilde f
  \in\fun{\cell\lfdex}\cC_{I_\lfdex}$ is a pointwise left Kan extension along
  $\lfdex\subto\cell\lfdex$ of its restriction $f\colon \lfdex\to\cC$.
\end{proposition}
\begin{proof}
  Let $\tilde f\colon \cell\lfdex\to\cC$ preserve finite $I_\lfdex$\=/cell
  complexes. Then it is a pointwise left Kan extension along
  $\lfdex\subto\cell\lfdex$ if and only if for every $X = \emptyset\to\ldots
  X_n$ in $\cell\lfdex$, $\tilde f X$ is a colimit in $\cC$ of the composite
  \[
    \lfdex/X_n\to\lfdex\xto{f} \cC .
  \]
  Let $\emptyset\to\ldots X_n$ be in $\cell\lfdex$. We proceed by induction on
  the pair $(\height(X_n),n)$ of the dimension of $X_n$ and $n$. In the base
  case, $\tilde f \emptyset$ is an initial object of $\cC$ since $\tilde f$
  preserves $I_\lfdex$\=/cell complexes. In the induction step for
  $\emptyset\to\ldots X_{n+1}$, we know that for some $c$ in $\lfdex$,
  \[
    \lfdex/X_{n+1} =
    \lfdex\slice c \amalg_{\lfdex\slice c^-} \lfdex/X_n.
  \]
  Therefore $\height(\partial c)<\height(X_{n+1})$, $\height(X_n)\leq
  \height(X_{n+1})$, and of course $n<n+1$. Hence by the induction hypothesis,
  $\tilde f (\emptyset\to\ldots X_n)$ is a colimit of
  $\lfdex/X_n\to\lfdex\xto{f}\cC$ and $\tilde f\partial c$ is a colimit of
  $\lfdex\slice c^- \to \lfdex\xto{f}\cC$. Combining the two,\footnote{And since
    a colimit indexed by a colimit of a diagram $I\to\Cat$ is the colimit of
    colimits. The reader is also invited to calculate the particular kind of
    colimit in $\Cat$ that we use---namely a pushout of an inclusion $A\subto
    A\star 1$---which is a particularly easy case of the calculation in
    \cite[Prop. 5.2]{FritschLatch1979}.} the colimit of
  $\lfdex/X_{n+1}\to\lfdex\xto f\cC$ is
  \[
    (fc)\amalg_{\tilde f\partial c}\tilde f
    (\emptyset\to\ldots X_n).
  \]
  But this is $\tilde f (\emptyset\to\ldots X_{n+1})$
  since $\tilde f$ preserves $I_\lfdex$-cell complexes.
\end{proof}

\begin{corollary}
 \label{cor:cell-C-rex-iff-cell-extendable} 
 For $\cC$ with an initial object, we have
 $\fun{\cell\lfdex}\cC_{I_\lfdex}\simeq \fun{\cell\lfdex}{\cC}\rex$. In other
 words, a functor $\cell\lfdex\to\cC$ preserves finite colimits if and only if
 it preserves finite $I_\lfdex$\=/cell complexes.
\end{corollary}


\section{\texorpdfstring{$\lfdex$}{C}-collections and \texorpdfstring{$\lfdex$}{C}-sorted theories}
\label{sec:C-collections-theories}


\begin{definition}
 \label{def:C-collection} 
  We call the presheaf category $\fun{\cell\lfdex\times\lfdex\op}\Set$ the
  category of \defn{(cartesian) $\lfdex$-collections}. We write it as
  $\coll\lfdex$.
\end{definition}

\begin{remark}
 \label{rem:C-collections-are-term-signatures} 
  We see the objects of $\coll\lfdex$ as $\lfdex$-sorted \emph{signatures of
    function symbols}. If $F\in\coll\lfdex$, then for a \emph{context}
  $\Gamma\in\cell\lfdex$ and a \emph{sort} $c\in\lfdex$, we see any $t\in
  F(\Gamma,c)$ as a ``function symbol'' with input sort $\Gamma$ and output sort
  ``of shape $c$'', that we may write $\Gamma\vdash t\type c$. Functoriality in
  $\Gamma$ is \emph{change of input variables} for a function symbol.
  Functoriality in $c$ describes the \emph{dependencies} (which are therefore
  themselves function symbols) in the output sort of a function symbol.

  This point of view appears in \cite{fiore2008second}, where it is exploited to
  define higher order signatures. Our ``$\lfdex$\=/sorted theories'' are
  equivalent to the \emph{$\Sigma_0$\=/models with substitution} in
  \cite[II.3]{fiore2008second}.

  We will make the correspondence with syntax precise
  in~\cref{sec:C-ctxl-cats-syntax}.
\end{remark}

\begin{proposition}
  The category of $\lfdex$-collections is equivalent to the full subcategory of
  $[\psh\lfdex,\psh\lfdex]$ consisting of the finitary endofunctors (those that
  preserve finitely filtered colimits).
\end{proposition}
\begin{proof}
  This is just left Kan extension along $\cell\lfdex\subto\psh\lfdex$ due to the
  isomorphism $\coll\lfdex\cong [\cell\lfdex,\psh\lfdex]$ and since $\psh\lfdex
  \simeq \Ind{\cell\lfdex}$ is an ind-completion of $\cell\lfdex$
  (\cref{cor:cell-C-finite-cocompletion}).
\end{proof}
\begin{corollary}
  The full inclusion $\coll\lfdex\subto [\psh\lfdex,\psh\lfdex]$ is monoidal and
  exhibits $\coll\lfdex$ as a monoidal subcategory of $[\psh\lfdex,\psh\lfdex]$
  (for the composition of endofunctors).
\end{corollary}
\begin{proof}
  The identity functor is finitary, and the composite of two finitary functors
  is finitary.
\end{proof}

\begin{corollary}
  The category $\Mon{\coll\lfdex}$ of monoids in $\coll\lfdex$ and monoid
  morphisms is equivalent to the category of finitary monads (those whose
  endofunctor preserves filtered colimits) on $\psh\lfdex$ and monad morphisms.
\end{corollary}

\begin{remark}[{\cite[II.3]{fiore2008second}}]
  This \emph{substitution} monoidal product $-\circ-$ on $\coll\lfdex$ can be explicitly
  calculated by the coend formula for left Kan extension along
  $\cell\lfdex\subto\psh\lfdex$. For every $F,G$ in
  $\coll\lfdex$, we have
  \[
    G\circ F (X,c) = \int^{Y\in\cell\lfdex} G(Y,c)\times \psh\lfdex(Y,F(X,-)).
  \]
\end{remark}
\begin{remark}
  \label{rem:loc-k-small-dir-mnd}
  If $D$ is a locally $\kappa$-small direct category (\cref{rem:loc-k-small-dir-cats}), and
  if we replace ``finite'' with ``$\kappa$-small'' in the definitions of $\cell
  D$ and $\coll D$, then $\Mon{\coll D}$ is equivalent to the category of
  $\kappa$-accessible monads on $\psh D$.  
\end{remark}

\begin{remark}
  Following \cref{rem:C-collections-are-term-signatures}, we can see monoids in
  $\coll\lfdex$ as \emph{theories} over signatures. We will detail this
  in~\cref{sec:C-ctxl-cats-syntax}.
\end{remark}

\begin{para}[Cellular arities]
    \label{para:cellular-arities}
    Recall that a functor $f\colon
    \cA\to\cB$ is \emph{dense} if the identity functor $1_\cB$ is a (pointwise)
    left Kan extension of $f$ along $f$. That is, for every X in $\cB$, the
    canonical cocone from $(f/X\to \cB)$ to $X$ is a colimit-cocone.
    
    Let $A$ be a small category. Let $i\colon A\subto \cC$ be a fully faithful
    and dense functor to a category $\cC$.  The density of $i$ is equivalent to its
    \emph{nerve functor}, namely the functor $\nu_i\colon \cC\to\psh A$ that
    takes $X$ to $\cC(i-,X)$, being fully faithful.
  
  A \defn{monad with $A$ as arities} \cite[Def. 1.8]{Berger2012} is a monad
  $T\colon \cC\to\cC$ such that for every $X\in\cC$, the composite functor
  $\nu_iT$ takes the canonical colimit-cocone in $\cC$ from $(i/ X \to \cC)$ to
  $X$ (given by the density of $i$) to a colimit cocone in $\psh A$. In other
  words, this condition states that $\nu_iT$ is a pointwise left Kan extension
  of $\nu_iTi$ along $i$. Note that $\nu_i$ reflects colimits since it is fully
  faithful; hence if $T$ has arities $A$ then $T$ preserves the canonical
  colimit cocones. In other words, if $T$ has arities $A$ then $T$ is a
  pointwise left Kan extension of $Ti$ along $i$ (but the converse is not
  necessarily true). We write $\Mnd_A(\cC)$ for the category of monads with $A$
  as arities, whose objects are the monads on $\cC$ with $A$ as arities and
  whose morphisms are the morphisms of monads.

  A \defn{theory with $A$ as arities} \cite[Def. 3.1]{Berger2012} is a pair
  $(\Theta,j)$ of a category $\Theta$ and an identity-on-objects functor
  $j\colon A\to\Theta$ such that the monad\footnote{For any functor $f\colon
    A\to B$ between small categories, $f_!\dashv f^*\dashv f_*$ denote the left
    and right adjoints to the precomposition functor $f^*\colon\psh B\to \psh
    A$.} $j^*j_!\colon \psh A\to\psh A$ preserves the essential image of
  $\nu_i\colon\cC\subto\psh A$. A morphism of theories
  $(\Theta_1,j_1)\to(\Theta_2,j_2)$ is a functor $\theta\colon
  \Theta_1\to\Theta_2$ such that $\theta j_1=j_2$. We write $\Law_A(\cC)$ for
  the category of theories with $A$ as arities.

  Every $T$ in $\Mnd_A(\cC)$ gives a theory with $A$ as arities, by the
  (identity-on-objects, fully faithful) factorisation of the composite $Ti\colon
  A\to\cC$, and conversely, every $(\Theta,j)$ in $\Law_A(\cC)$ gives a monad with
  arities $A$ by restriction of the monad $j^*j_!$ to $\cC$.
\end{para}

\begin{remark}
  Monads with arities were studied in \cite{weber2007familial} and (``Lawvere'')
  theories with arities in \cite{mellies2010segal}. We refer to the detailed
  review in \cite{Berger2012}.
\end{remark}

\begin{example}
  The presheaf category $\psh{\GG_1}$
  (\cref{exa:lfd-cats}\ref{exa:lfd-cats-grph}) is the category of directed
  graphs. If $\bDelta_0$ is the full subcategory of all the non-empty finite
  filiform graphs (graphs of the form $x_0\to x_1\to\ldots \to x_n$), then the
  free-category monad (associated to the monadic forgetful functor
  $\Cat\to\psh{\GG_1}$) is a monad with $\bDelta_0$ as arities. Its theory with
  arities is just the identity-on-objects free-category functor
  $\bDelta_0\to\bDelta$ to the simplex category.
\end{example}

\begin{proposition}[{\cite[Thm 3.4]{Berger2012}}]
 \label{prop:mnd-theory-arity-equiv} 
  Let $A\in\Cat$ and let $i\colon A\subto \cC$ be a fully faithful dense
  functor. The mutually inverse constructions in~\cref{para:cellular-arities}
  define an equivalence of categories $\Mnd_A(\cC)\simeq \Law_A(\cC)$.
\end{proposition}

\begin{proposition}
  [{\cite[1.12]{Berger2012}}]
 \label{prop:mnd-arities-cell-fin-monads} 
  Consider the fully faithful functor $i\colon \cell\lfdex\subto\psh\lfdex$.
  Then $i$ is dense and the category $\Mnd_{\cell\lfdex}(\psh\lfdex)$ of monads
  with $\cell\lfdex$ as arities is precisely the
  category of finitary monads on $\psh\lfdex$ and monad morphisms.
\end{proposition}
\begin{proof}
  Density of $i$ follows from \cref{lem:cell-C-fin-pres}, since $\fin_\lfdex$ is
  dense in $\psh\lfdex$. Next, $T\colon \psh\lfdex\to\psh\lfdex$ is finitary if
  and only if it is the left Kan extension of $Ti$ along $i$. Finally, the nerve
  functor $\nu_i\colon \psh\lfdex \subto \psh{\cell\lfdex}$ preserves filtered
  colimits (since $i$ is an ind-completion), and for each $X\in\psh \lfdex$, the
  category $i/ X$ is filtered. So if $T\colon\psh\lfdex\to\psh\lfdex$ is
  finitary, then the composite $\nu_iT$ preserves the required filtered
  colimit-cocones, hence $T$ has arities $\cell\lfdex$.
\end{proof}

\begin{corollary}
\label{cor:equiv-cell-monoids-law-theories}  
  We have equivalences of categories
  \[
    \Mon{\coll\lfdex}\simeq \Mnd_{\cell\lfdex}(\psh\lfdex) \simeq
    \Law_{\cell\lfdex}(\psh\lfdex).
  \]
\end{corollary}

\begin{definition}
  \label{def:C-sorted-theory}
  A \defn{$\lfdex$-sorted theory} is a theory with $\cell\lfdex$ as arities. The
  category of $\lfdex$-sorted theories is $\Law_{\cell\lfdex}(\psh\lfdex)$. A
  $\lfdex$-sorted theory is \defn{idempotent} if its associated finitary monad
  on $\psh\lfdex$ is an idempotent monad.
\end{definition}

\begin{remark}
  When $\lfdex$ is a set (\cref{exa:lfd-cats}\ref{exa:lfd-cats-set}), we recover
  the well-known equivalence between finitary monads on $\psh \lfdex \simeq
  \Set/\lfdex$ and $\lfdex$-sorted Lawvere theories.
\end{remark}

\begin{remark}
  \label{rem:idempotent-C-sorted-th-refl-localisn}
  Idempotent $\lfdex$-sorted theories correspond exactly to
  ($\omega$\=/)accessible reflective localisations $\psh\lfdex\localisation\cC$
  in the sense of \cref{sec:localisations}.\footnote{Since locally presentable
    $1$-categories are a particular case of locally presentable \oo\=/categories.}
\end{remark}

\begin{proposition}
  [{\cite[3.5]{Berger2012}}]
  \label{prop:cell-th-w-arities}
  An identity-on-objects functor $j\colon \cell\lfdex\to\Theta$ is a
  $\lfdex$-sorted theory if and only if it preserves finite $I_\lfdex$-cell
  complexes.
\end{proposition}
\begin{proof}
  The nerve functor $\nu_i$ is just the subtopos inclusion\footnote{Recall that
    for any functor $f\colon A\to B$ between small categories, $f_*\colon\psh
    A\to \psh B$ is the nerve functor $\nu_{\nu_f}$ of the nerve functor $\nu_f$
    of $f$ \cite[Exp. I, Prop. 5.4(2)]{SGA4-1}, and that any one of $f,f_!,f_*$
    being fully faithful implies the other two being so as well (\cite[Exp. I,
    Prop. 5.6]{SGA4-1}).} $\iota_*\colon \psh\lfdex \subto \psh{\cell\lfdex}$
  associated to the full inclusion $\iota\colon\lfdex\subto \cell\lfdex$. The
  functor $\nu_i$ preserves filtered colimits, and for each $X\in\psh \lfdex$,
  the category $i/ X$ is filtered. Thus since $j^*j_!$ preserves colimits, $j$
  is a theory with $\cell\lfdex$ as arities if and only if the nerve functor
  $\nu_j\colon\Theta\to\psh{\cell\lfdex}$ factors through
  $\psh\lfdex\subto\psh{\cell\lfdex}$. But this is equivalent to $j$ preserving
  finite colimits, since $\cell\lfdex\simeq\fin_\lfdex$ (and so a presheaf is in
  the subtopos $\psh\lfdex\subto\psh{\cell\lfdex}$---namely, is a sheaf---if and
  only if it preserves finite limits). We conclude by
  \cref{cor:cell-C-rex-iff-cell-extendable}.
\end{proof}

\begin{corollary}
  \label{cor:io-ff-cell-th-with-arities}
  Let $\cell\lfdex\xto{j_\cC}\Theta_\cC\subto\cC$ be the (identity-on-objects,
  fully faithful) factorisation of a functor
  $f\colon{\cell\lfdex}\to\cC$ preserving finite $I_\lfdex$-cell
  complexes. Then $j_\cC$ is a theory with $\cell\lfdex$ as arities .
\end{corollary}
\begin{proof}
  Since $\Theta_\cC\subto\cC$ reflects colimits and $f$ preserves
  $I_\lfdex$-cell complexes, $j_\cC$ preserves $I_\lfdex$-cell complexes.
\end{proof}

\begin{para}
  [Cellular nerves of algebras]
  \label{para:cellular-nerves-algebras}
  We summarise some results from \cite{Berger2012}. Let
  $i\colon A\subto \cC$ be a fully faithful dense functor, and let $T\colon
  \cC\to\cC$ be a monad with $A$ as arities. Then if $j\colon A\to\Theta_T$ is
  its associated theory with $A$ as arities, the full inclusion $i_T\colon
  \Theta_T\subto T\alg$ is dense, and any $X\in\psh\Theta_T$ is in the essential
  image of $\nu_{i_T}\colon T\alg\subto\psh\Theta_T$ if and only if $j^*X$ is in
  the essential image of $\nu_i$. \cite[Thm 1.10]{Berger2012} (this is the
  \emph{nerve theorem}).
\end{para}

\begin{definition}
  \label{def:model-of-C-sorted-theory}
  Let $j\colon \cell\lfdex\to\Theta$ be a $\lfdex$-sorted theory. Its category
  $(\Theta,j)\Mod$ of \defn{models} is the full, wide subcategory of
  $\psh{\Theta}$ that is the essential image of the nerve functor
  $T_j\alg\subto\psh{\Theta}$ from the category of algebras of its associated
  finitary monad $T_j\colon\psh\lfdex\to\psh\lfdex$.
\end{definition}

\begin{lemma}
  \label{lem:model-C-sorted-theory-iff-preserves-cplxs}
  Let $j\colon \cell\lfdex\to\Theta$ be a $\lfdex$-sorted theory.
  Then $(\Theta,j)\Mod$ is the full subcategory of $\psh\Theta$ consisting of
  all $X$ such that $j^*X\colon\cell\lfdex\op\to\Set$ takes
  $\emptyset\in\cell\lfdex$ to $1\in\Set$ and pushouts of maps in
  $I_\lfdex$ to pullbacks in $\Set$.
\end{lemma}

\begin{remark}
  Most of~\cref{sec:C-collections-theories} holds when $\lfdex$ is a just a
  small category (and not necessarily \lfd~as well), by replacing $\cell\lfdex$
  with the category $\fin_\lfdex$ (and ``finite $I_\lfdex$\=/cell complexes''
  with ``finite colimits''). However, our goal is the classification of
  dependently sorted algebraic theories
  and~\cref{prop:cell-th-w-arities,cor:equiv-cell-monoids-law-theories} will
  play important rôles.
\end{remark}



\section{An equivalence between \texorpdfstring{$\lfdex$}{C}-sorted theories and
  \texorpdfstring{$\lfdex$}{C}-contextual categories}
\label{sec:C-contextual-cats}
We come at last to the goal of this chapter, namely the classification of
(certain) dependently sorted algebraic theories. To do so, we will incarnate
these theories as algebraic objects called \defn{$\lfdex$\=/contextual
  categories}. \emph{Contextual categories} were introduced in
\cite{cartmell1978generalised}, where Cartmell also gives an equivalence
\cite[\textsection{2.4}]{cartmell1978generalised} between the category of
contextual categories and the category of \defn{generalised algebraic theories
  (GATs)}. Hence, our dependently sorted algebraic theories will be examples of GATs,
but they will be a \emph{strict} subclass of all GATs. In other words,
$\lfdex$\=/contextual categories will be a strict subclass of all contextual
categories.

\begin{definition}[{\cite[\textsection{2.2}]{cartmell1978generalised}}, {\cite[Def.
    1.2]{streicher1991}}, {\cite[Def. 1.2.1]{kapulkin2012simplicial}}]
  \label{def:ctxl-category}
  A \defn{contextual category} is a small category $\sfD$ along with the data of:
  \begin{enumerate}
  \item
    \label[axiom]{def:ctxl-category:str:grading}
    a grading of its objects as $\ob \sfD = \coprod_{n\in\NN} \sfD_n$,
  \item
    \label[axiom]{def:ctxl-category:str:terminal}
    an object $1 \in \sfD_0$,
  \item
    \label[axiom]{def:ctxl-category:str:father}
    ``parent'' functions $\ft - _n \colon \sfD_{n+1} \to \sfD_n$ (we will
    usually suppress the subscript),
  \item
    \label[axiom]{def:ctxl-category:str:proj}
    for each $\Gamma \in \sfD_{n+1}$, a distinguished map $\pr \Gamma \colon \Gamma \to \ft
    \Gamma$,
  \item
    \label[axiom]{def:ctxl-category:str:connect-map}
    for each $\Gamma \in \sfD_{n+1}$ and $f \colon \Gamma' \to \ft\Gamma$, an
    object $f^*\Gamma$ together with a ``connecting map'' $\qpb f \Gamma \colon
    f^*\Gamma \to \Gamma$;
  \end{enumerate}
  such that:
  \begin{enumerate}[resume]
  \item
    \label[axiom]{def:ctxl-category:ax:terminal1}
    $1$ is the unique object in $\sfD_0$ ($\sfD_0=\{1\}$);
  \item
    \label[axiom]{def:ctxl-category:ax:terminal2}
    $1$ is a terminal object of $\sfD$;
  \item
    \label[axiom]{def:ctxl-category:ax:pullback}
    for each $\Gamma \in \sfD_{n+1}$, and $f\colon\Gamma' \to \ft\Gamma$, we
    have $\ft{f^*\Gamma} = \Gamma'$, and the square
    \[
      \begin{tikzcd}[sep=scriptsize]
        f^*\Gamma \pbmark \ar[r,"{\qpb f \Gamma}"] \ar[d,"\pr{f^*\Gamma}"']
        &\Gamma \ar[d,"\pr\Gamma"]
        \\
        \Gamma' \ar[r,"f"] &\ft\Gamma
      \end{tikzcd}
    \] commutes and is cartesian (called the \emph{canonical pullback} of
    $\Gamma$ along $f$); and
  \item
    \label[axiom]{def:ctxl-category:ax:functorial}
    canonical pullbacks are strictly functorial; that is, for $\Gamma \in
    \sfD_{n+1}$, we have $1_{\ft\Gamma}^* \Gamma = \Gamma$ and $\qpb
    {1_{\ft\Gamma}} \Gamma = 1_\Gamma$; and for $\Gamma \in \sfD_{n+1},
    f\colon\Gamma'\to\ft\Gamma$ and $g\colon\Gamma''\to \Gamma'$, we have
    $(fg)^*\Gamma = g^*f^*\Gamma$ and $\qpb {fg}\Gamma = (\qpb f \Gamma)
    (\qpb g {f^*\Gamma})$.
  \end{enumerate}
  A \defn{morphism} of contextual categories is a functor $\sfD_1\to\sfD_2$
  preserving the data of
  \cref{def:ctxl-category:str:grading,def:ctxl-category:str:terminal,def:ctxl-category:str:father,def:ctxl-category:str:proj,def:ctxl-category:str:connect-map}.
  The category of contextual categories is written $\CxlCat$.
\end{definition}


\begin{proposition}
  \label{prop:cell-S-is-ctxl-cat}
  Let $\lfdex$ be a \lfd~category. Then ${\cell \lfdex \op}$ admits a canonical
  structure of a contextual category.
\end{proposition}
\begin{proof}
  From \cref{rem:replacing-fin-with-cell}, ${\cell \lfdex \op}$ is small. Then,
  \begin{enumerate}\itemsep.25em
  \item[\ref{def:ctxl-category:str:grading}] We have the same grading $\ob
    {\cell\lfdex\op}\eqdef \coprod_n (\cell \lfdex)_n$ (see \cref{def:Cell-C}).
  \item[\ref{def:ctxl-category:str:terminal}] $\emptyset$ is
    \ref{def:ctxl-category:ax:terminal1} the unique element of $(\cell\lfdex)_0$
    and \ref{def:ctxl-category:ax:terminal2} is a terminal object of
    $\cell\lfdex\op$.
  \item[\ref{def:ctxl-category:str:father}] For every $X = \emptyset \to \ldots
    X_{n+1}$ in $(\cell \lfdex ) _{n+1}$, we define $\ft X\eqdef\emptyset \to
    \ldots X_n$, and \ref{def:ctxl-category:str:proj} $\pr X$ in $\relHom{\cell
      \lfdex \op}X{\ft X}$ as the chosen morphism $X_n\to X_{n+1}$ in
    $\relHom{\cell\lfdex}{\ft X}{X} 
    $.
  \item[\ref{def:ctxl-category:str:connect-map}] For $X = \emptyset \to \ldots
    X_{n+1}$ in $(\cell \lfdex ) _{n+1}$, the map $X_n\to X_{n+1}$ comes with a
    unique choice of pushout square (\cref{def:Cell-C}) as in the left square in
    the following diagram. Let $Y = \emptyset \to \ldots Y_m$ be in $(\cell
    \lfdex) _{m}$, and let $f \colon X_n \to Y$ be in $\relHom {\cell \lfdex}{\ft X}
    Y$.
    \[
      \begin{tikzcd}[sep=scriptsize]
        \partial c \ar[r,"\gamma"] \ar[d,hook]
        &X_n \ar[r,"f"] \ar[d]
        &Y_m \ar[d] \\
        c \ar[r,"\gamma.c"] \ar[rr,bend right=20,"f\gamma.c"']
        &X_{n+1}\pomark \ar[r,,"{\qpb f X}",dashed]
        &f^*X  
      \end{tikzcd}
    \]
     We define $\emptyset \to
    \ldots Y_m\to f^*X$ in $(\cell \lfdex)_{m+1}$ by choosing the outer pushout
    square in the previous diagram, and we define the morphism $\qpb f X $ in
    $\relHom{\cell\lfdex} X {f^*X}$ as the unique dotted arrow.

  \item[\ref{def:ctxl-category:ax:pullback}] In the previous diagram, the outer
    and left-hand squares are cocartesian, thus the square on the right is
    cocartesian in $\psh\lfdex$ and cartesian in $\cell\lfdex\op$.
  \item[\ref{def:ctxl-category:ax:functorial}] Given $X = \emptyset \to \ldots
    X_{n+1}$ in $(\cell \lfdex ) _{n+1}$ as above, $Y = \emptyset \to \ldots Y_m$ in
    $(\cell \lfdex) _{m}$, and $Z = \emptyset \to \ldots Z_l$ in $(\cell \lfdex)
    _{l}$, as well as maps $f\colon X_n\to Y_m$ and $g\colon Y_m\to Z_l$ in
    $\psh\lfdex$, then functoriality follows since $(gf)\gamma.c =
    g(f\gamma).c$.\qedhere
  \end{enumerate}
\end{proof}

\begin{definition}
 \label{def:initial-C-cxl-cat} 
  The \defn{initial $\lfdex$-contextual category} is defined to be
  $\cell\lfdex\op$ (\emph{qua} contextual category). We write it as $\cxl\lfdex$.
\end{definition}

\begin{remark}
  \label{rem:ctx-cat-structure-disrespect-equiv-cats}
  The structure of a contextual category cannot be transferred across an
  equivalence of categories. Hence, while $\cxl\lfdex$ is a contextual category,
  $\fin_\lfdex\op$ is not necessarily one\footnote{It is clearly not one (except
    when $\lfdex=\emptyset$) if our category $\fin$ of finite sets is only
    \emph{essentially} small. But this is not the only problem---a presheaf in
    $\fin_\lfdex$ has no ``canonical parent''.}. So contextual categories are
  category-theoretically ``evil'', but this is simply because they are the
  algebraic counterparts of syntactic objects.\footnote{Though we could
    reasonably object to the use of the word ``category'' in the name.}
\end{remark}

\begin{proposition}
 \label{prop:cxl-C-fin-limits} 
   $\cxl\lfdex$ has finite limits.
\end{proposition}
\begin{proof}
 $\cell\lfdex\simeq\fin_\lfdex$ is finitely cocomplete.
\end{proof}

\begin{proposition}
 \label{prop:cxl-C-pushouts-nonempty-coprods} 
  $\cxl\lfdex$ has pushouts and finite non-empty coproducts, and has
  finite colimits if and only if $\lfdex$ is finite.
\end{proposition}
\begin{proof}
  Follows from \cref{cor:fin-C-pback-fin-lims}.
\end{proof}

\begin{proposition}
 \label{prop:cxl-C-display-map-epi} 
  Every map $\pr\Gamma$ in $\cxl\lfdex$ as in \cref{def:ctxl-category:str:proj}
  of~\cref{def:ctxl-category} is an epimorphism.
\end{proposition}
\begin{proof}
  Since it is a monomorphism in $\cell\lfdex$.
\end{proof}

\begin{proposition}
 \label{prop:cxl-C-codescent} 
  Every canonical pullback square in $\cxl\lfdex$ as in
  \cref{def:ctxl-category:ax:pullback} of \cref{def:ctxl-category} has
  \emph{codescent}, namely:
  \begin{enumerate}
  \item it is also cocartesian,
  \item and in any commuting cube below, where the canonical pullback is the top face
    \[
      \begin{tikzcd}[sep=small]
      f^*X_{n+1} \ar[rr,"f.X"]\ar[rd] \ar[dd]
      && X_{n+1} \ar[rd] \ar[dd]\\
      &Y_m
      && X_n \ar[from=ll,"f"near start, crossing over]\ar[dd]\\
      \Xi_1 \ar[rd] \ar[rr]
      &&\Xi_2 \ar[rd]\\
      &\Xi_3 \ar[from = uu, crossing over] \ar[rr]
      &&\Xi_4 
    \end{tikzcd}
  \]
  if the front and right hand faces are cocartesian, then the bottom face is
  cartesian if and only if the back and left faces are cocartesian.
  \end{enumerate}
\end{proposition}
\begin{proof}
  Pushouts of monomorphisms in $\psh\lfdex$ are \emph{effective} (or ``van
  Kampen'') colimits. Hence the pushouts of boundary inclusions in $\cell\lfdex$
  satisfy descent (since by \cref{cor:fin-C-pback-fin-lims}, $\cell\lfdex$ has
  pullbacks).
\end{proof}

\begin{corollary}
  For every canonical pullback square in $\cxl\lfdex$ as in
  \cref{prop:cxl-C-codescent}, the pullback and pushout operations define an
  equivalence of categories
  \[
    {{\cxl\lfdex}}\coslice{f^*{X_{n+1}}}\simeq {\cxl\lfdex}\coslice {X_{n+1}}
    \times_{{\cxl\lfdex}\coslice{X_n}} {\cxl\lfdex}\coslice{Y_m}
  \]
  between the coslice category $\cxl\lfdex\coslice{f^*X_{n+1}}$ and the category
  of cocartesian natural transformations under the cospan $Y_m\xto{f}X_n\ot
  X_{n+1}$.
\end{corollary}

\begin{remark}
    \cref{prop:cxl-C-fin-limits,prop:cxl-C-pushouts-nonempty-coprods,prop:cxl-C-display-map-epi,prop:cxl-C-codescent}
    can be seen as consequences of the fact that the morphisms of the initial
    $\lfdex$\=/contextual category are all made up of only variables, since
    $\cxl\lfdex$ has no term-constructors
    (see~\cref{prop:type-sig-synt-cat-C-cxl-cat}).
\end{remark}

\begin{para}
  \cref{rem:functor-fin-to-cell} gives us a full inclusion
  $\lfdex\subto\cell\lfdex$ taking every $c$ in $\lfdex$ to a cell context
  $\emptyset\to\ldots\to \ft c\to c$ in $(\cell\lfdex)_{k_c}$, where $k_c =
  \sum_{d \in \lfdex} \abs{\relHom \lfdex d c}$ is the cardinality of the set of
  objects of $\lfdex\slice c$. Moreover, $\ft c \to c$ is the canonical morphism
  of colimits given by the inclusion $\lfdex\slice c^-\subset\lfdex\slice c$. We
  will show that the inclusion $\lfdex\op\subto \cxl\lfdex$ is universal among
  contextual categories for these properties.
\end{para}

\begin{definition}
  \label{def:contextual-functor-from-direct-cat}
  Let $\sfD$ be any contextual category. Then a functor $f \colon \lfdex\op \to
  \sfD$ is a \defn{contextual functor (under $\lfdex$)} if for all $c\in\lfdex$,
  \begin{enumerate}
  \item
    \label{def:contextual-functor-from-direct-cat:1}
    $fc$ is in $\sfD_{k_c}$, where $k_c = \sum_{d \in \lfdex} \abs{\relHom
      \lfdex d c} = \abs{\ob{\lfdex\slice c}}$,
  \item
    \label{def:contextual-functor-from-direct-cat:2}
    $\ft{fc}$ is a limit of $(\lfdex\slice c^-)\op\to\lfdex\op\xto{f}\sfD$,
    and $\pr{fc}\colon fc\to\ft{fc}$ is the canonical morphism of limits given
    by $\lfdex\slice c^-\subset \lfdex\slice c$.
  \end{enumerate}
\end{definition}

\begin{proposition}
  Let $\sfD$ be a contextual category and $f\colon\lfdex\op\to\sfD$ a contextual
  functor. Then $f$ has a pointwise right Kan extension along
  $\lfdex\op\subto\cxl\lfdex$. 
\end{proposition}
\begin{proof}
  Since $\sfD$ is a contextual category, $\sfD\op$ has an initial object and
  pushouts of $\pr{fc}\op\colon\ft{fc}\to fc$ for every $c\in\lfdex$. It thus
  has finite $f\op(I_\lfdex)$-cell complexes, where $f\op(I_\lfdex) =
  \{\pr{fc}\op\mid c\in\lfdex\}$. But then $f\op$ has a pointwise left Kan
  extension along $\lfdex\subto\cell\lfdex$, which is just the functor taking
  every $I_\lfdex$-cell complex to the corresponding $f\op(I_\lfdex)$-cell
  complex.
\end{proof}

\begin{proposition}
  [Initiality of $\lfdex\op\subto\cxl\lfdex$]
  \label{prop:ctxl-functor-cxl-C-initial}
  Let $\sfD$ be a contextual category.
  \begin{enumerate}
  \item
   \label{prop:ctxl-functor-cxl-C-initial:1} 
    Every morphism $\cxl\lfdex\to\sfD$ in $\CxlCat$ is the right Kan
    extension of some (essentially unique) contextual functor $\lfdex\op\to\sfD$.
  \item
    \label{prop:ctxl-functor-cxl-C-initial:2}
    Every contextual functor $f\colon\lfdex\op\to\sfD$ has a right Kan
    extension along $\lfdex\op\subto\cxl\lfdex$ that is a morphism of contextual
    categories.
  \end{enumerate}
\end{proposition}
\begin{proof}
  \ref{prop:ctxl-functor-cxl-C-initial:1} Any morphism $g\colon
  \cxl\lfdex\to\sfD$ in $\CxlCat$ preserves grading, thus its restriction along
  $\lfdex\op\subto\cxl\lfdex$ satisfies
  \cref{def:contextual-functor-from-direct-cat}\ref{def:contextual-functor-from-direct-cat:1}.
  Now $g$ also preserves canonical pullbacks, so $g\op$ preserves
  $I_\lfdex$-cell complexes and hence finite colimits
  (\cref{cor:cell-C-rex-iff-cell-extendable}).
  Thus $g$ is a right Kan extension of its restriction, and so its restriction
  necessarily satisfies
  \cref{def:contextual-functor-from-direct-cat}\ref{def:contextual-functor-from-direct-cat:2}.

  \ref{prop:ctxl-functor-cxl-C-initial:2} Let $f$ be a contextual functor.
  Define the functor $\tilde f\colon\cxl\lfdex\to\sfD$ as follows: let $\tilde
  f(\emptyset)\eqdef 1_\sfD$ and inductively (on the grading of $\cxl\lfdex$),
  let the chosen pushout squares in $\cell\lfdex$ be taken to the corresponding
  canonical pullbacks in $\sfD$. By \cref{cor:cell-C-rex-iff-cell-extendable},
  this defines a right Kan extension of $f$ along $\lfdex\op\subto\cxl\lfdex$.
  Further, the definition ensures that $\tilde f$ preserves the grading, parent
  maps and canonical pullbacks, thus is a morphism in $\CxlCat$.
\end{proof}


\begin{corollary}
  \label{cor:cxl-morphism-to-C-sorted-theory}
  Let $\cxl\lfdex\xto{j_f\op}\Theta_f\op\subto\sfD$ be the (identity-on-objects,
  fully faithful) factorisation of the underlying functor of a morphism $f\colon
  \cxl\lfdex\to\sfD$ of contextual categories. Then
  $j_f\colon\cell\lfdex\to\Theta_f$ is a $\lfdex$-sorted theory.
\end{corollary}
\begin{proof}
  Follows from~\cref{cor:io-ff-cell-th-with-arities} and the previous
  proposition.
\end{proof}

\begin{remark}
  The category $\CxlCat$ is not (in any good way) a $2$\=/category, which is why
  \cref{prop:ctxl-functor-cxl-C-initial} is not stated as an equivalence of
  functor categories. However $\CxlCat$ embeds via a left adjoint into various
  $2$\=/categories of models of dependent type theory (such as that of
  \emph{categories with attributes}). Since Kan extension along a fully faithful
  functor is a fully faithful functor between functor categories, an equivalence
  of (full) hom-categories in such a $2$\=/category follows from essential
  surjectivity, which is just the content of
  \cref{prop:ctxl-functor-cxl-C-initial}.
\end{remark}

\begin{remark}
  [Semi-simplicial, globular, opetopic types]
  \cref{def:contextual-functor-from-direct-cat} gives a very concrete
  description of what \emph{type-theoretic} $\lfdex$-objects in any contextual
  category $\sfD$ are. Hence, for example, a semi-simplicial type in any
  dependent type theory $\sfT$ \emph{à la} Martin-Löf is exactly a contextual
  functor ${\bDelta'}\op\to \sfD_\sfT$ to the syntactic contextual category
  $\sfD_\sfT$ of $\sfT$. A semi-simplicial type in a context $\Gamma$ is a
  contextual functor to the slice contextual category ${\sfD_\sfT}\slice
  \Gamma$.
\end{remark}



\begin{definition}
  \label{def:model-of-cxl-cat}
  For $\sfD$ a contextual category and $\cC$ a category, a \defn{$\sfD$\=/model
    in $\cC$} is a functor $\sfD\to\cC$ preserving the terminal object and
  taking canonical pullbacks (of \cref{def:ctxl-category:ax:pullback}
  of~\cref{def:ctxl-category}) to cartesian squares. The category $\sfD\Mod_\cC$
  of $\sfD$\=/models in $\cC$ is the full subcategory of $\fun\sfD\cC$
  consisting of the $\sfD$\=/models in $\cC$. We simply write $\sfD\Mod$ for the
  category of $\sfD$\=/models in $\Set$.
\end{definition}

\begin{proposition}
  We have canonical equivalences $\fun{\cell\lfdex}\cC_{I_\lfdex}\simeq
  \cxl\lfdex\Mod_\cC$ and $\psh\lfdex\simeq \cxl\lfdex\Mod$.
\end{proposition}
\begin{proof}
  Immediate from \cref{cor:cell-C-rex-iff-cell-extendable}.
\end{proof}

\begin{para}
 \label{para:cxl-cat-to-C-sorted-theory-models} 
  \cref{cor:cxl-morphism-to-C-sorted-theory} allows us to associate a
  $\lfdex$-sorted theory to every morphism $\cxl\lfdex\to\sfD$ in $\CxlCat$.
  However, it is entirely possible for two distinct morphisms
  $\cxl\lfdex\to\sfD_1$ and $\cxl\lfdex\to\sfD_2$ to give the same
  $\lfdex$-sorted theory.
  
  Moreover, if $f\colon \cxl\lfdex\to\sfD$ in $\CxlCat\coslice{\cell\lfdex}$ is any
  contextual category under $\cxl\lfdex$, and if
  $j_f\colon\cell\lfdex\to\Theta_f$ is its associated $\lfdex$-sorted theory,
  then every $X\in\sfD\Mod$ gives, by restriction along $\Theta_f\op\subto\sfD$,
  a model of the $\lfdex$-sorted theory $j_f$.
  We thus have a
  forgetful functor $\sfD\Mod\to (\Theta_f,j_f)\Mod$, which is not in general an
  equivalence of categories.  
\end{para}
\begin{definition}
 \label{def:C-cxl-cat} 
  Let $f\colon \cxl\lfdex\to\sfD$ be a morphism in $\CxlCat$, and
  let
  \[
    \begin{tikzcd}[sep=small]
      \cxl\lfdex \ar[r,"j_f\op"]
      &\Theta_f\op \ar[r,hook]
      &\sfD
    \end{tikzcd}
  \]
  be the (identity-on-objects, fully faithful) factorisation of its underlying
  functor. Then $f$ is a \defn{$\lfdex$-contextual category} if for every
  morphism $g\colon \cxl\lfdex\to\sfD'$ in $\CxlCat$ and every triangle
  $hj_f\op=g$ (in $\Cat$), there exists a unique morphism $\tilde h\colon
  \sfD\to\sfD'$ in $\CxlCat$ making the following diagram commute.
 \[
   \begin{tikzcd}[sep=scriptsize]
     \cxl\lfdex\ar[r,"j_f\op"]\ar[dr,"g"']
     &{\Theta_f\op} \ar[d,"h"] \ar[r,hook]
     &\sfD \ar[dl,dashed,"\exists!\tilde h"]\\
     &\sfD'
   \end{tikzcd}
 \]
 The category $\CxlCat_\lfdex$ of $\lfdex$-contextual categories is the full
 subcategory of the coslice category $\CxlCat\coslice{\cxl\lfdex}$ consisting of
 the $\lfdex$-contextual categories.
\end{definition}

\begin{remark}
    The identity functor $1_{\cxl\lfdex}$ is the initial object of
    $\CxlCat_\lfdex$ (see~\cref{def:initial-C-cxl-cat} and
    \cref{prop:ctxl-functor-cxl-C-initial}).
\end{remark}

\begin{remark}
   \label{rem:is-cxl-cat-C-cxl-cat} 
    An obvious question is whether every contextual category is a
    $\lfdex$\=/contextual category for some \lfd~category $\lfdex$. I do not
    know the answer to this question. Nevertheless,
    in~\cref{cha:models-of-C-contextual-categories}, we will see that every
    contextual category is \emph{Morita equivalent} to some
    $\lfdex$\=/contextual category, namely both have equivalent categories of
    $\Set$\=/models.
\end{remark}

\begin{lemma}
 \label{lem:C-cxl-cat-to-C-sorted-theory-fully-faithful} 
  The association $f\mapsto (j_f\colon \cell\lfdex\to \Theta_f)$ defines a fully
  faithful functor $\CxlCat_\lfdex\to \Law_{\cell\lfdex}(\psh\lfdex)$.
\end{lemma}
\begin{proof}
  Given a morphism $h$ in $\CxlCat_\lfdex$ as on the left below, the induced
  morphism $\bar h$ between (identity-on-objects, fully faithful) factorisations
  is in $\Law_{\cell\lfdex}(\psh\lfdex)$. Conversely, given a
  morphism $k$ in $\Law_{\cell\lfdex}$ as on the right below, there exists a
  unique morphism $\tilde k$ in $\CxlCat_\lfdex$. It is readily verified that the two
  constructions are mutually inverse.
  \[
    \begin{tikzcd}[sep=scriptsize]
      \cxl\lfdex\ar[r,"j_f\op"]\ar[dr,"j_g\op"']
     &{\Theta_f\op} \ar[d,dashed,"\bar h\op"] \ar[r,hook]
     &\sfD \ar[d,"h"]\\
     &{\Theta_g\op} \ar[r,hook]
     &\sfD'
   \end{tikzcd}
   \qquad
   \begin{tikzcd}[sep=scriptsize]
      \cxl\lfdex\ar[r,"j_f\op"]\ar[dr,"j_g\op"']
     &{\Theta_f\op} \ar[d,"k\op"] \ar[r,hook]
     &\sfD \ar[d,dashed,"\tilde k"]\\
     &{\Theta_g\op} \ar[r,hook]
     &\sfD'
   \end{tikzcd}
  \]
\end{proof}

\begin{proposition}
  \label{prop:C-cxl-cat-ess-surj}
  The fully faithful functor $\CxlCat_\lfdex\to \Law_{\cell\lfdex}(\psh\lfdex)$
  is essentially surjective.
\end{proposition}

We defer the proof of \cref{prop:C-cxl-cat-ess-surj} for the moment. We are now
able to state and prove the main theorem of \cref{chap:contextual-categories}.

\begin{theorem}
  [Classification of dependently sorted algebraic theories]
 \label{thm:classification-dep-alg-theories} 
  Let $\lfdex$ be a \lfd~category.
  The categories
  \begin{enumerate}
  \item
   \label{thm:classification-dep-alg-theories:monoids} 
    $\Mon{\coll\lfdex}$ of monoids in 
    cartesian $\lfdex$\=/collections,
  \item
    \label{thm:classification-dep-alg-theories:fin-monads}
    $\Mnd_{\cell\lfdex}(\psh\lfdex)$ of finitary monads on $\psh\lfdex$,
  \item
    \label{thm:classification-dep-alg-theories:C-sorted-theories}
    $\Law_{\cell\lfdex}(\psh\lfdex)$ of $\lfdex$\=/sorted theories,
  \item
    \label{thm:classification-dep-alg-theories:C-cxl-cats}
    and $\CxlCat_\lfdex$ of $\lfdex$\=/contextual categories,
  \end{enumerate}
  are equivalent.
\end{theorem}
\begin{proof}
  The equivalences $\ref{thm:classification-dep-alg-theories:monoids} \simeq
  \ref{thm:classification-dep-alg-theories:fin-monads}
  \simeq\ref{thm:classification-dep-alg-theories:C-sorted-theories}$ are
  from~\cref{cor:equiv-cell-monoids-law-theories}. The functor
  $\ref{thm:classification-dep-alg-theories:C-sorted-theories}\to
  \ref{thm:classification-dep-alg-theories:C-cxl-cats}$ is fully faithful by
  \cref{lem:C-cxl-cat-to-C-sorted-theory-fully-faithful} and essentially
  surjective by \cref{prop:C-cxl-cat-ess-surj}.
\end{proof}

The rest of \cref{sec:C-contextual-cats} is devoted to proving
\cref{prop:C-cxl-cat-ess-surj}. We fix a $\lfdex$\=/sorted theory
$j\colon\cell\lfdex\to\Theta$ throughout.

\begin{definition}
  We have fixed an inclusion $\lfdex\subto\cell\lfdex$. So for every $c$ in
  $\lfdex$, we have maps $j\delta_c\colon j\partial c\to jc$ in $\Theta$ (and so
  in $(\Theta,j)\Mod$). Let $jI_\lfdex$ be the set of these maps (between
  representables) in $(\Theta,j)\Mod$. We define the opposite category
  $\sfD_j\op$ of the \defn{contextual completion $\sfD_j$ of $(\Theta,j)$} as a category
  of finite $jI_\lfdex$\=/cell complexes in $(\Theta,j)\Mod$, proceeding as in
  \cref{def:Cell-C}. To begin, we set $(\sfD_j)_0\eqdef\{j\emptyset\}$. For any
  $j\emptyset\to X_1\to\ldots X_n$ in $(\sfD_j)_n$, and any span $jc\ot
  j\partial c\to X_n$ in $(\Theta,j)\Mod$, we make a choice of pushout in
  $(\Theta,j)\Mod$, giving $j\emptyset\to X_1\to\ldots X_n\to X_{n+1}$ in
  $(\sfD_j)_{n+1}$. Finally, we define the hom-sets as
  $\relHom{\sfD_j\op}{j\emptyset\to\ldots X}{j\emptyset\to\ldots Y}\eqdef
  \relHom{(\Theta,j)\Mod}XY$.
\end{definition}

\begin{proposition}
  $\sfD_j$ is a contextual category, and there is a full inclusion
  $\Theta\op\subto\sfD_j$ such that $\cxl\lfdex\to\Theta\op\subto\sfD_j$ is a
  morphism of contextual categories.
\end{proposition}
\begin{proof}
  It is readily verified that $\sfD_j$ is a contextual category (just as in the
  proof of \cref{prop:cell-S-is-ctxl-cat}). Since both $j$ and the composite
  $\cxl\lfdex\to\Theta\subto(\Theta,j)\Mod$ preserve $I_\lfdex$\=/cell
  complexes, there is a grading-preserving full embedding
  $\Theta\op\subto\sfD_j$, which suffices for
  $\cxl\lfdex\to\Theta\op\subto\sfD_j$ to be a morphism in $\CxlCat$.
\end{proof}

\begin{proposition}
  The forgetful functor $\sfD_j\Mod \to (\Theta,j)\Mod$
  (see \cref{para:cxl-cat-to-C-sorted-theory-models}) is an equivalence of
  categories.
\end{proposition}
\begin{proof}
  By definition of $\sfD_j$, the nerve functor associated to the inclusion
  $\sfD_j\op\subto (\Theta,j)\Mod$ factors through $\sfD_j\Mod$, and provides
  the required inverse.
\end{proof}

\begin{para}
  The objects of $\Theta$ are images of $I_\lfdex$\=/cell contexts. Hence, since
  both $j$ and the composite $\cell\lfdex\to\Theta\subto(\Theta,j)\Mod$ preserve
  $I_\lfdex$\=/cell complexes, they are $jI_\lfdex$\=/cell complexes in
  $(\Theta,j)\Mod$. However, $\Theta$ does not have pushouts of $j\delta_c$ for
  \emph{all} morphisms of models/algebras $j\partial c \to X$ in $\Theta$ (but
  only for those\footnote{These will be called
    \emph{free} morphisms.} in the image of $j$). An object in the contextual
  completion $\sfD_j$ is therefore \emph{not in general a free
    algebra/model}---namely, it is not necessarily in the image of
  $j_!\colon\psh\lfdex\to(\Theta,j)\Mod$.
\end{para}

\begin{example}
 \label{exa:pointed-cats-ctxl-cmpletion} 
  Consider the category $\Cat_*$ of pointed categories (the coslice category
  $\Cat\coslice 1$ under the terminal category $1$). Then $\Cat_*$ is the
  category of models of a $\GG_1$\=/sorted theory
  $j_{\Cat_*}\colon\cell{\GG_1}\to \Theta_{\Cat_*}$, since it is the category of
  algebras of a finitary monad on $\psh{\GG_1}$. It is easy to see that the
  pointed category $\{f\colon* \to x\}$, whose chosen object is $*$, is
  not free on any graph. But it is in the contextual completion $\sfD_{\Cat_*}$
  using the pushout below in $\Cat_*$, where the span is in
  $\Theta_{\Cat_*}$.
  \[
    \begin{tikzcd}[sep=scriptsize]
      j_{\Cat_*}\{y,x\}\ar[r,"{(*,x)}"] \ar[d,"j\delta"']
      &j_{\Cat_*}\{x\}\ar[d]\\
      j_{\Cat_*}\{f\colon y\to x\}\ar[r]
      &\{f\colon* \to x\}\pomark
    \end{tikzcd}
  \]
\end{example}

\begin{remark}
  When $\lfdex$ is a set, then every object in $\sfD_j\op$ is a coproduct of
  objects of the form $jc, c\in\lfdex$ and is thus a free algebra/model. Hence
  in this case (multisorted algebraic theories), $\Theta$ and $\sfD_j\op$
  coincide.
\end{remark}

\begin{definition}
  For every $c$ in $\lfdex$, we say that a morphism $f\colon jc\to Y$ in the
  image of $\sfD_j\op\subto (\Theta,j)\Mod$ is \defn{active} if there is no
  $j\emptyset\to Y_1\to\ldots Y_n = Y$ in $\sfD_j\op$ such that
  $f$ factors through
  $jc\to Y_k\to Y_{k+1}\to\ldots Y_n$ for some $k<n$. Hence
  any $jc\to j\emptyset$ is active, and for any $jc\to Y$, there is some
  $j\emptyset\to Y_1\to\ldots Y_n = Y$ in $\sfD_j\op$ and some $k\leq n$ such
  that $jc\to Y_k$ is active.


  A morphism in $\Theta$ is \defn{free} if it is in the image of $j$.
\end{definition}

\begin{definition}
  Let $j\emptyset\to\ldots Y_n$ be in $\sfD_j$. A \defn{free replacement} of
  $Y_n$ is a map $\bar Y\to Y_n$ from some object $\bar Y$ in $\Theta$, such
  that every map $\Gamma\to Y_n$ in $(\Theta,j)\Mod$ from any $\Gamma$ in
  $\Theta$ factors as $\Gamma\to\bar Y\to Y_n$.
\end{definition}

To understand the idea of the following construction-proposition, the
reader may skip to~\cref{exa:free-repl-explanation,rem:free-repl-explanation} to
get a feel for what's going on.

\begin{proposition}
  \label{prop:ctxl-cmpletion-has-free-replacements}
  Every object in $\sfD_j$ has a free replacement.
\end{proposition}
\begin{proof}
  For any $j\emptyset\to\ldots Y_n$ in
  $(\sfD_j)_n$, we will define a sequence of maps $\bar Y^k\to Y_k, 0\leq k\leq
  n$ from objects $\bar Y^k\in\Theta$, along with \emph{free} morphisms $\bar
  Y^k\to\bar Y^{k+1}$ making the following diagram commute.
  \[
    \begin{tikzcd}[sep=small]
      j\emptyset\ar[d,equals] \ar[r]
      &\bar Y^1\ar[r] \ar[d]
      &\ldots \ar[r]
      &\bar Y^n\ar[d]\\
      j\emptyset \ar[r]
      &Y_1\ar[r]
      &\ldots \ar[r]
      &Y_n
    \end{tikzcd}
  \]
  Every free morphism $\bar Y^k\to \bar Y^{k+1}$ will be the image by $j$ of a
  finite relative $I_\lfdex$\=/cell complex in $\cell\lfdex$. For every $k$,
 $\bar Y^k\to Y_k$ will be a free replacement. We proceed by induction on
 $j\emptyset\to\ldots Y_n$.
  \begin{enumerate}
  \item In the base case, we set $\bar{j\emptyset}\eqdef j\emptyset$ along with
    the identity morphism $1_{j\emptyset}$.
  \item In the induction step, we assume that we have reached $\bar Y^k\to Y_k$.
    Now, $Y_k\to Y_{k+1}$ is a pushout of the form below. By induction hypothesis, $f$
    factors as some morphism $\bar f\colon j\partial c\to \bar Y^k$ in $\Theta$
    followed by the free replacement $\bar Y^k\to Y_k$.
    \[
      \begin{tikzcd}[sep=small]
        j(\partial c)\ar[r,"f"] \ar[d] \ar[rr,bend left,"\bar f"]
        &Y_k\ar[d]
        &\bar Y^k\ar[l]\\
        jc\ar[r]
        &Y_{k+1}\pomark
      \end{tikzcd}
    \]
    Now $\bar f\colon j(\partial c)\to \bar Y^k$ is a morphism in $\Theta$, but
    it is not necessarily free. First, we write $\emptyset=(\partial c)_0\to
    (\partial c)_1\to\ldots (\partial c)_l= \partial c$ as an object of
    $\cell\lfdex$. We will define a sequence $\bar Y^k=\bar Y^k_0\to \bar
    Y^k_1\to\ldots\bar Y^k_l$ that is the image by $j$ of a finite relative
    $I_\lfdex$\=/cell complex, and such that for each $m\leq l$, $\bar Y^k\to\bar Y^k_{m}$
    has a retraction in $\Theta$. We will also define a free morphism $\bar
    f^m\colon j(\partial c)_m\to \bar Y^k_m$ making the square below
    commute, where $\bar Y^k_m\to \bar Y^k$ is the retraction.
    \[
      \begin{tikzcd}[sep=small]
        j(\partial c)_m \ar[d] \ar[r,"\bar f^m"]
        &\bar Y^k_m \ar[d,two heads]\\
        j(\partial c)\ar[r,"\bar f"]
        &\bar Y^k
      \end{tikzcd}
    \]
    We proceed by induction.
    \begin{enumerate}
    \item In the base case, we set $\bar Y^k_0\eqdef \bar Y^k$ and we set $\bar
      f^0\colon j\emptyset \to \bar Y^k_0$ which is clearly a free morphism.
    \item In the induction step, we assume we have reached the stage $m$. By
      induction hypothesis, $\bar f^m\colon j(\partial c)_m\to \bar Y^k_m$ is a
      free map. We begin a case analysis.
      \begin{enumerate}
      \item \textbf{If} there exists a \emph{free} morphism $\tilde f_{m}\colon
        j(\partial c)_{m+1}\to \bar Y^k_m$ making the triangles below commute,
        \[
          \begin{tikzcd}[sep=small]
            j(\partial c)_m\ar[r,"{\bar f^m}"] \ar[d]
            &\bar Y^k_m \\
            j(\partial c)_{m+1}
            \ar[ur,"{\tilde f_m}"']
          \end{tikzcd}
          \qquad\qquad
          \begin{tikzcd}[sep=small]
            j(\partial c)_{m+1} \ar[r,"{\tilde f_m}"] \ar[dr]
            &\bar Y^k_{m} \ar[d,two heads]\\
            &\bar Y^k
          \end{tikzcd}
        \]
        then we set $\bar Y^{k}_{m+1}\eqdef \bar Y^k_m$, and we set $\bar
        f^{m+1}\eqdef \tilde f_{m}$.
      \item \textbf{Else,} we define $\bar
        Y^k_{m+1}$ and $\bar f^{m+1}$ via the pushout of free maps below.
        \[
          \begin{tikzcd}[sep=small]
            j(\partial c)_m \ar[r,"{\bar f^m}"]\ar[d]
            &\bar Y^k_m \ar[d]\\
            (\partial c)_{m+1} \ar[r,"{\bar f^{m+1}}"']
            &{\bar Y^k_{m+1}}\pomark
          \end{tikzcd}
        \]
        Since we have a retraction $\bar Y^k_m\to\bar Y^k$ and a map $j(\partial
        c)_{m+1}\to\bar Y^k$, we obtain the
        retraction $\bar Y^k_{m+1}\to \bar Y^k$. We end our case analysis.
      \end{enumerate}
      Our inner induction is done.
    \end{enumerate}
    To continue with our outer induction, we consider the commutative diagram
    below in $(\Theta,j)\Mod$, using the free morphism $\bar f^l\colon
    j(\partial c)\to \bar Y^k_l$ and the retraction $\bar Y^k_l\to \bar Y^k$
    that we have just built.
    \begin{equation}
      \tag{$\star$}
      \label{eq:free-replacement}
      \begin{tikzcd}[sep=scriptsize]
        &\bar Y^k \ar[dr,hook]
        \\
        j(\partial c) \ar[rr,"{\bar f^l}",near start]\ar[dd] \ar[dr,"{f}"']
        \ar[ur,"{\bar f}"]
        &&\bar Y^k_l \ar[dl,"gr"] \ar[dd]\ar[ul,two heads,bend right,"r"']\\
        &Y_k\ar[from=uu,crossing over,"g"',near end]\\
        jc \ar[rr] \ar[dr]
        &&\bar Y^{k+1}\pomark \ar[dl,dashed]\\
        &Y_{k+1} \ar[from=uu,crossing over]
      \end{tikzcd}
    \end{equation}
    We define $\bar Y^{k+1}$ via the cocartesian back face, which is thus a
    cocartesian square of free morphisms. By definition of $\bar f^l$, we have
    $\bar f = r\bar f^l$, so we obtain the dotted arrow. Since the front left
    square is cocartesian by definition of $Y_{k+1}$, so is the intermediate
    (front right) square. It remains to be shown that $\bar Y^{k+1}\to Y_{k+1}$
    is a free replacement.

    We write $\rho:\psh\Theta\localisation (\Theta,j)\Mod:\nu$ for the reflective
    localisation associated to the fully faithful nerve functor $\nu$ of
    $\Theta\subto (\Theta,j)\Mod$. Consider the diagram below in $\psh\Theta$
    (where $\eta$ is the unit of the reflection).
    \begin{equation}
      \tag{$\star\star$}
      \label{eq:free-replacement-pushout}
      \begin{tikzcd}[sep=small]
        &\nu\bar Y^k \ar[dl,hook']\ar[dr]\\
        \nu\bar Y^k_l \ar[ur,two heads,bend left] \ar[rr] \ar[dd]
        &&\nu Y_k\ar[dd]\\
        &{}
        \\
        \nu\bar Y^{k+1} \ar[rr]
        &&\tilde Y_{k+1}\pomark \ar[r,"\eta"]
        &\nu Y_{k+1}
      \end{tikzcd}
    \end{equation}
    Let $j\Gamma$ be any representable in $\Theta$. By induction hypothesis, any
    map $j\Gamma\to \nu Y_k$ factors as $j\Gamma\to \nu\bar Y^k\to \nu Y_k$, and
    therefore as $j\Gamma\to \nu\bar Y^k_l\to \nu Y_k$. Since $j\Gamma$ is
    representable, and therefore \emph{tiny}, this implies that any map
    $j\Gamma\to \tilde Y_{k+1}$ factors as $j\Gamma\to \nu\bar Y^{k+1}\to \tilde
    Y_{k+1}$ . Therefore, the map $\nu\bar Y^{k+1}\to \tilde Y_{k+1}$ is an
    epimorphism. We can conclude if we show that in the previous diagram, the
    unit $\eta$ is an isomorphism, namely that $\tilde Y_{k+1}$ is in the
    essential image of $\nu$.

    Recall that this is the case if and only if the presheaf $j^*\tilde Y_{k+1}
    \in\psh{\cell\lfdex}$ preserves $I_\lfdex$\=/cell complexes. We proceed by
    induction on $\Gamma\in\cell\lfdex$. In the base case, since $\nu\bar
    Y^{k+1}\to\tilde Y_{k+1}$ is an epimorphism, we have $(j^*\tilde
    Y_{k+1})_\emptyset = 1$. In the induction step, let $\Gamma$ be defined by
    the pushout in $\psh\lfdex$ on the left below.
    \[
      \begin{tikzcd}[sep=small]
        \partial c'\ar[r]\ar[d,hook]
        &\ft\Gamma\ar[d]
        \\
        c'\ar[r]
        &\Gamma\pomark
      \end{tikzcd}
      \qquad\qquad
      \begin{tikzcd}[sep=small]
        j\partial c'\ar[r] \ar[d]
        &j\ft\Gamma\ar[d]\\
        jc'\ar[r]
        &\tilde Y_{k+1}
      \end{tikzcd}
    \]
    Consider a commutative square in $\psh\Theta$ as on the right above. We need
    to show that it factors through a unique morphism $j\Gamma\to\tilde
    Y_{k+1}$. For existence, we show that it factors as a commutative square
    that is a cocone from the cospan of free maps $jc'\ot j\partial c'\to
    j\ft\Gamma$ to $\nu\bar Y^{k+1}$, followed by the map $\nu\bar
    Y^{k+1}\to\tilde Y_{k+1}$. We factor each morphism to $\tilde Y_{k+1}$ in
    the right hand square above, as follows. Using \cref{eq:free-replacement}
    and the definition of $\bar Y^k_l$, if the map $j\ft\Gamma\to\tilde Y_{k+1}$
    factors through $\nu Y_k$, we may choose a map $j\ft\Gamma\to\bar Y^k_l$
    such that for every $c''\to\ft\Gamma$ in $\cell\lfdex$, there is a square
    below
    \[
      \begin{tikzcd}[sep=small]
        jc''\ar[r]\ar[d]
        &j\partial c \ar[d]
        \\
        j\ft\Gamma\ar[r]
        &\bar Y^k_l
      \end{tikzcd}
    \]
    where $jc''\to j\partial c$ is a free map. Doing the same for the map $c'\to
    \tilde Y_{k+1}$, then the square below commutes (since the map $j\partial
    c''\to\ft\Gamma$ is free).
    \[
      \begin{tikzcd}[sep=small]
        j\partial c'\ar[r] \ar[d]
        &j\ft\Gamma \ar[d]\\
        jc' \ar[r]
        &\nu\bar Y^{k+1}
      \end{tikzcd}
    \]
    Since $\nu\bar Y^{k+1}$ is in the image of $\nu$, we obtain the desired map
    $j\Gamma\to\nu\bar Y^{k+1}\to \tilde Y_{k+1}$. Uniqueness follows from the
    pushout in~\cref{eq:free-replacement-pushout}.
    Our induction is finished.\qedhere
  \end{enumerate}
\end{proof}

\begin{example}
 \label{exa:free-repl-explanation} 
  The construction in the proof of
  \cref{prop:ctxl-cmpletion-has-free-replacements} is best understood by looking
  at an example. We return to \cref{exa:pointed-cats-ctxl-cmpletion}, namely the
  $\GG_1$\=/sorted theory of pointed categories. The pointed category $\{f\colon
  * \to x\}$ is not free, and its free replacement provided by our
  construction is $j_{\Cat_*}\{f\colon y\to x\}$. Similarly, the free replacement of
  $\{f\colon *\to*\}$ is $j_{\Cat_*}\{f\colon x\to x\}$.
\end{example}

\begin{remark}
 \label{rem:free-repl-explanation} 
  In general, the free replacement
  of~\cref{prop:ctxl-cmpletion-has-free-replacements} consists in replacing all
  dependencies in any context $j\emptyset\to\ldots Y$ in $\sfD_j$ with
  variables, in a \emph{minimal} way.\footnote{I thank N. Jeannerod for a lively
    discussion on the syntactic version of the proof of
    \cref{prop:ctxl-cmpletion-has-free-replacements}.} 
\end{remark}

\noindent We can now prove our stated goal.
\begin{proof}(Of~\cref{prop:C-cxl-cat-ess-surj}.)

  We show that the contextual completion $\sfD_j$ is a $\lfdex$\=/contextual
  category.
  Consider a diagram of solid arrows, where $g$ is a morphism in $\CxlCat$.
  \[
    \begin{tikzcd}[sep=scriptsize]
      \cxl\lfdex\ar[r,"j\op"]\ar[dr,"g"']
      &{\Theta\op} \ar[d,"h"] \ar[r,hook]
      &\sfD_j \ar[dl,dashed,"\tilde h"]\\
      &\sfD'
    \end{tikzcd}
  \]
  We will construct the morphism $\tilde h$ in $\CxlCat$. Since
  $\Theta\subto\sfD_j\op$ is dense, it suffices to define $\tilde h$ on objects
  and on all maps $Y\to X$ in $\sfD_j$ such that $X$ is in $\Theta\op$. First
  for any $f\colon X\to X'$ in $ \Theta$, we let $\tilde h f \eqdef hf$. Then we
  proceed by induction. For every $\emptyset\to\ldots Y$ in $\sfD_j$, we
  use~\cref{prop:ctxl-cmpletion-has-free-replacements} to obtain a free
  replacement $\bar Y\to Y$. Then, we use~\cref{eq:free-replacement} and the
  induction hypothesis to define $\tilde h Y\to \tilde h \bar Y$ via the
  corresponding canonical pullback in $\sfD'$. Since $\bar Y\to Y$ is a free
  replacement, this suffices to define the image under $\tilde h$ of every map
  $Y\to X$ in $\sfD_j$, where $X$ is in $\Theta\op$. The uniqueness of $\tilde
  h$ follows from the definition of $\sfD_j$.
\end{proof}

\section{Examples of \texorpdfstring{$\lfdex$}{C}-contextual categories}
\label{sec:examples-C-cxl-cats}

We will immediately put~\cref{thm:classification-dep-alg-theories} to good use,
by using it to recognise several dependently sorted algebraic theories. For any
$\lfdex$\=/contextual category $\sfD$, its category of $\Set$\=/models
$\sfD\Mod$ is obviously locally finitely presentable. We will say that a locally
finitely presentable category $\cC$ is \defn{classified by $\sfD$} if it is
equivalent to $\sfD\Mod$.

\begin{example}
  [Semi-simplicial, globular, opetopic sets]
  \label{exa:C-sets}
  The identity monad on any \lfd~category $\lfdex$ is obviously finitary. Using
  \cref{exa:lfd-cats}, we therefore deduce the following $\lfdex$\=/contextual
  categories.
  \begin{enumerate}
  \item the $\GG$\=/contextual category $\cxl\GG$ classifying globular sets,
  \item the $\OO$\=/contextual category $\cxl\OO$ classifying opetopic sets,
  \item the $\bDelta'$\=/contextual category $\cxl{\bDelta'}$ classifying
    semi-simplicial sets.
  \end{enumerate}
\end{example}

\begin{example}
  [Simplicial, dendroidal sets] Let $R$ be a Reedy category such that its wide
  direct subcategory $R'$ is locally finite, and let $j\colon R'\to R$ be the
  wide inclusion. Then $j^*:\psh{R}\to \psh {R'}$ is monadic and $\omega$\=/accessible, and so
  we deduce the $R'$\=/contextual category classifying $\psh R$. Combined
  with~\cref{exa:lfd-cats}\ref{exa:lfd-cats-reedy}, this gives the dependently
  sorted algebraic theories classifying simplicial sets, dendroidal sets and
  $\Theta$\=/sets.
\end{example}

\begin{example}
  [Strict $\omega$\=/categories] Recall that we write $\GG$ for the category of
  globes. Let $\bTheta_0$ be the full subcategory of $\psh\GG$ called the
  \emph{globular site} \cite[Def. 1.5]{Berger2002cellular}. The Grothendieck
  topology on $\Theta_0$ is generated by the full inclusion $i\colon \GG\subto\Theta_0$
  of representables (thus $i_*\colon\psh\GG\subto\psh{\Theta_0}$ is the inclusion of
  sheaves into presheaves). Then the \emph{strict $\omega$\=/category} monad
  $T_\omega$ on $\psh\GG$ (the terminal globular operad) has $\bTheta_0$ as
  arities \cite[Example 4.18]{weber2007familial}. Moreover, every globular set
  in $\bTheta_0$ is finite. This implies that $T_\omega$ is finitary, thus the
  category $T_\omega\alg$ of strict $\omega$\=/categories is classified by a
  $\GG$\=/contextual category.
\end{example}

\begin{example}
  [Globular operads, weak $\omega$\=/categories] Let $\cE$ be a locally
  presentable, locally cartesian closed category. Let $F,G\colon \cD\to\cE$ be
  functors from any category $\cD$, and let $\alpha \colon G\to F$ be a
  \emph{cartesian} natural transformation.   
  Then $G$ preserves any colimit that $F$
  does.
  In effect, let $h\colon I\to \cD$ be a diagram
  with a colimit $X$ in $\cD$ that is preserved by $F$, namely $FX\cong \colim_i
  Fhi$. We have a cartesian transformation $\alpha_h \colon Gh\to
  Fh$.
  \[
    \begin{tikzcd}
      Ghi \pbmark \ar[r] \ar[d,"\alpha"']
      &GX \ar[d,"\alpha"]\\
      Fhi \ar[r]
      &FX = \colim_i Fhi
    \end{tikzcd}
  \]
  By universality of colimits in $\cE$, the map $\alpha_X\colon
  GX\to FX$ is the colimit of the maps $\alpha_{hi}\colon Ghi\to Fhi$, namely
  $GX\cong \colim_i Ghi$.

  A \emph{globular operad} is a cartesian monad $T$ on $\psh\GG$ equipped with a
  cartesian monad morphism $T\to T_\omega$ to the strict $\omega$\=/category
  monad \cite[Chs 4, 6]{Leinster2004}. Thus, since $T_\omega$ is finitary, so
  is $T$. Hence every globular operad has an associated $\GG$\=/contextual
  category, whose models are the algebras over the globular operad.

  Recall that a Batanin weak $\omega$\=/category is an algebra over a
  \emph{contractible} globular operad \cite[Def. 1.20]{Berger2002cellular}.
  \emph{A fortiori}, each contractible globular operad has an associated
  $\GG$\=/contextual category.

  The most important \emph{Grothendieck\=/Maltsiniotis coherators for weak
    $\omega$\=/categories} are all \emph{homogeneous globular theories}.
  Moreover, each corresponds to a contractible globular operad (this follows from
  \cite[Secs 6.6, 6.7]{ara2010infini}, subject to a conjecture proven in
  \cite{bourke2020iterated}). Therefore, each has a corresponding
  $\GG$\=/contextual category.

  Similarly, $n$\=/categories (with all varying degrees of
  strictness) are all classified by $\GG_{\leq n}$\=/contextual categories.
\end{example}

\begin{remark}
  An explicit syntactic presentation of a dependently typed algebraic theory for
  a particular Grothendieck\=/Maltsiniotis coherator is defined in
  \cite{finster2017type} and studied in \cite{benjaminthese2020}. Their
  syntactic presentation can be shown to be a syntactic presentation of a
  $\GG$\=/contextual category in the sense of
  \cref{def:dependently-typed-theory-syntactic}.
\end{remark}

\begin{example}
  [$\omega$\=/groupoids] Grothendieck $\omega$\=/groupoids are defined by
  \emph{coherators for $\omega$\=/groupoids}, which are certain
  identity-on-objects functors $\Theta_0\to C$ from the globular site
  \cite[1.5]{maltsiniotis2010grothendieck}. The category of
  \emph{$\omega$\=/groupoids of type $C$} is the full subcategory of $\psh C$ on
  all $C\op\to\Set$ such that the composite $\Theta_0\op\to C\op\to\Set$ is a
  sheaf (namely, is in the image of $\psh\GG\subto\psh{\Theta_0}$). Each
  coherator gives a monad $T_C$ on $\psh\GG$, such that $T_C\alg$ is the
  category of $\omega$\=/groupoids of type $C$. The monad $T_C$ is not a monad
  with $\Theta_0$ as arities, but it satisfies a nerve theorem (it is a
  \emph{nervous monad} in the terminology of \cite{BourkeGarner2019}) and is
  easily seen to be finitary. We therefore deduce a $\GG$\=/contextual category
  for each coherator for $\omega$\=/groupoids.
\end{example}

\begin{remark}
  An explicit syntactic presentation of a particular coherator for
  $\omega$\=/groupoids is given in \cite[App. A]{brunerie2016homotopy}. It is
  not clear whether their definition can be shown to be a syntactic presentation
  of a $\GG$\=/contextual category in the sense
  of~\cref{def:dependently-typed-theory-syntactic}. 
\end{remark}

\begin{example}
  [Planar coloured operads] The category $\Opd\pl$ of planar coloured operads in
  $\Set$ is the category of algebras of a finitary monad on $\psh{\elTree\pl}$
  (~\cref{exa:lfd-cats}\ref{exa:lfd-cats-eltree}). We deduce the
  $\elTree\pl$\=/contextual category classifying $\Opd\pl$.
\end{example}

\begin{remark}
  Coloured \emph{symmetric} $\Set$\=/operads are monadic over presheaves on the
  category $\elTree$ of (non\=/planar) corollas/elementary trees, which is a
  \emph{generalised} \lfd~category, namely $\elTree$ has non\=/trivial
  automorphisms. Syntactically, this corresponds to types of the signature
  associated to $\elTree$ having non\=/trivial auto-equivalences (loops in the
  universe). Since Homotopy Type Theory (HoTT) has identity types and a
  univalent universe, we conjecture that the $\elTree$\=/sorted theory of
  symmetric operads is an (elementary) example of a dependently typed higher
  algebraic theory (assuming any such gadget to be an extension of HoTT).

  Note that the category of symmetric $\Set$\=/operads is a locally finitely presentable
  category, therefore by~\cref{thm:lfp-cat-models-of-C-sorted-theory} it
  is classified by a $\lfdex$\=/contextual category over \emph{some}
  \lfd~category $\lfdex$.  
\end{remark}


\newcommand{\var}[1]{\mathrm{var}(#1)}

\newcommand{\emptyctx}{[]}

\newcommand{\Sig}{{\mathrm{\cS ig}}}

\section{Syntactic presentations of \texorpdfstring{$\lfdex$}{C}-contextual
  categories}
\label{sec:C-ctxl-cats-syntax}

Thus far, we have only worked with algebraic structures that we \emph{claim}
correspond to syntactic objects. By
\cite[\textsection{2.4}]{cartmell1978generalised} and
\cref{thm:classification-dep-alg-theories}, dependently sorted algebraic
theories in our sense correspond to certain (but not all) \emph{generalised
  algebraic theories} (GATs), which are syntactic objects. So we could start
from the definition of GATs, and carve out the precise subclass corresponding to
$\lfdex$-contextual categories using syntactic constraints. However, this would
mean casting aside much of the fruit of our previous labours, since
$\cell\lfdex$ and $\lfdex$\=/contextual categories are just syntax done up as
an algebraic gadget, and working with them is tantamount to working with syntax.
Therefore, we will define syntactic presentations of $\lfdex$-contextual
categories from the ground up, and \emph{observe} that they form a strict
subclass of GATs.

\begin{para}
  [Structural MLTT]
  \label{para:struc-rules-mltt}
  The basic structure common to all our dependently sorted theories is
  \defn{structural (Martin-Löf) dependent type theory (MLTT)}, which consists of the forms of
  judgment and ``structural rules'' of Martin-Löf Type Theory (introduced
  in\cite{MartinLof1975}). It can be found in \cite[Sec. 2]{Hofmann1997} and
  \cite[App. A.1]{kapulkin2012simplicial}; we will recall it briefly
  (following \cite{kapulkin2012simplicial}).

  To begin with, there are three classes of raw syntax: \emph{contexts},
  \emph{types} and \emph{terms},\footnote{These are sometimes \cite[Sec.
    2.3]{Hofmann1997} called \emph{pre}-contexts, -types, and -terms to
    emphasise that they have not been \emph{judged} to be so.} the latter containing
  an infinite set of (term) \emph{variables}; an element of each of these
  classes is a tree of symbols. This raw syntax is quotiented by
  \emph{alpha-equivalence} and the operation of \emph{(capture-free) substitution} is
  defined on it.\footnote{These involve the replacement of bound variables; they will
    not be very important for us.
  }  Next, the four \emph{judgment forms} are introduced.
  \[
    \Gamma\vdash A \qquad\Gamma\vdash A=A' \qquad\Gamma\vdash a\type A
    \qquad\Gamma\vdash a=a'\type A
  \]
  The first says that \emph{$A$ is a type in the context $\Gamma$}, the third
  that \emph{$a$ is a term of type $A$ in the context $\Gamma$}, and the second
  and fourth judge types and terms to be \emph{definitionally equal}. As in
  \cite[App. A]{kapulkin2012simplicial}, we take the judgment form for contexts to be
  derived from those of the previous forms. Namely, for $n$ in $\NN$, a list $\Gamma
  = (x_i\of A_i)_{i<n}$ is said to be a \emph{context} (written $\vdash \Gamma\ctx$) as an
  abbreviation of $\forall i<n, (x_j\of A_j)_{j<i}\vdash A_i$.

  We define the \emph{judgments} to be statements of the previous forms
  that can be derived (as conclusions of proof-trees) using the following 
  inference rules
  \[
    \MLVble{x_1\of A_1,\ldots,x_n\of A_n\ctx}
    {x_1\of A_1,\ldots,x_n\of A_n \vdash x_i\type A_i }
    \qquad\qquad
    \MLSubst{\Gamma\vdash a\type A}
    {\quad\Gamma,x\of A,\Delta\vdash \cJ}
    {\Gamma,\Delta[a/x]\vdash \cJ[a/x]}
  \]
  \[
    \MLWkg{\Gamma\vdash A}
    {\quad\Gamma,\Delta\vdash\cJ}
    {\Gamma,x\of A,\Delta\vdash\cJ}
  \]
  (where $[a/x]$ is the substitution of $x$ with $a$) as well as the following
  rules for definitional equality.
  \[
    \prftree{\Gamma\vdash A}
    {\Gamma\vdash A=A}
    \qquad\quad
    \prftree{\Gamma\vdash A=B}
    {\Gamma\vdash B=A}
    \qquad
    \prftree{\Gamma\vdash A=B}
    {\quad\Gamma\vdash B=C}
    {\Gamma\vdash A=C}
  \]
  \[
    \prftree{\Gamma\vdash a\type A}
    {\Gamma\vdash a=a\type A}
    \qquad
    \prftree{\Gamma\vdash a=b\type A}
    {\Gamma\vdash b=a\type A}
    \qquad
    \prftree{\Gamma\vdash a=b\type A}
    {\quad\Gamma\vdash b=c\type A}
    {\Gamma\vdash a=c\type A}
  \]
  \[
    \prftree{\Gamma\vdash a\type A} {\quad\Gamma\vdash A=B} {\Gamma\vdash a\type B}
    \qquad\qquad \prftree{\Gamma\vdash a=b\type A} {\quad\Gamma\vdash A=B}
    {\Gamma\vdash a=b\type B}
  \]
  The rules \textsf{Wkg} and \textsf{Subst} of MLTT are \emph{admissible}
  \cite[E2.7]{Hofmann1997}, meaning that every proof-tree that uses them can be
  replaced by one that does not. All the dependently sorted theories that we
  consider will be extensions of MLTT in which \textsf{Wkg} will remain
  admissible, but not necessarily \textsf{Subst}.

  Finally, we briefly recall the technique of \emph{structural induction}: when
  proving a statement about all judgments of a given form, we will use the
  well-founded partial order on proof\=/trees to reason by induction---treating
  each inference rule whose conclusion is of the desired form as a separate
  case. We will generally skip the cases of the rules for definitional equality,
  since these will be trivial.
\end{para}

\begin{definition}
    The following definitions make sense in any extension $\sfT$ of structural
    MLTT that we will consider.

    For any two contexts ${}\vdash \Gamma\ctx$ and ${}\vdash(x_i\of
    A_i)_{i<n}\ctx$ of $\sfT$, a \defn{context morphism} $\sigma\colon
    \Gamma\to(x_i\of A_i)_{i<n}$ is a list $\sigma=(t_i)_{i<n}$ of term
    judgments $\Gamma\vdash t_i \type A_i[x_j/t_j]_{j<i}$ in $\sfT$. The
    \defn{syntactic category of $\sfT$} has as its objects equivalence classes
    of context judgments (up to renaming of variables and definitional equality
    in $\sfT$), and as its morphisms, equivalence classes of context morphisms
    (up to renaming of variables and definitional equality in $\sfT$). The
    syntactic category of $\sfT$ is always a contextual category.
\end{definition}

\begin{para}
  A \emph{dependent type signature} will be a generalisation of the set of sorts of a
  multisorted algebraic theory. It will allow for sorts (called \emph{types}) that,
  in order to be well-defined, require a finite context of typed variables whose
  types are recursively well\=/defined. For instance, the sequence of statements
  \begin{align*}
    &\vdash A \\
    x\of A,y \of A
    &\vdash B \\
    x\of A,y \of A, z \of B(x,y)
    &\vdash C \\
    x \of A , y \of B(x,x) , z \of C(x,x,y)
    &\vdash D 
  \end{align*}
  is a type signature that defines $A$ to be a type with respect to the empty
  context, $B$ a type with respect to the context $x\of A ,y \of A$, and so on. This
  is captured by the following mutually inductive definitions.
\end{para}



\begin{definition}
  \label{def:type-signature}
  A \defn{(dependent) type signature} is a graded set $S = \coprod_{n \in \NN}
  S_n$ such that each $S_{j}$ is a set of pairwise free type declarations over the type signature
    $S_{<j} \eqdef \coprod_{i<j} S_i$.
\end{definition}

\begin{definition}
 \label{def:type-declaration} 
  A \defn{type declaration} over a type signature $S$ is an inference rule of the
  form below, where $\Gamma$ is a context of $S$ and the type symbol $A$ is \emph{fresh}, \ie
  it does not appear in $S$.
  \[
    \prfaxiom{\Gamma \vdash A}
  \]
  Two type declarations $s,s'$ over $S$ are 
    \emph{free} if $s$ is a type declaration over the signature
    $S \coprod \{ s' \} $ (hence, vice versa). 
\end{definition}

\begin{definition}
 \label{def:context-of-signature} 
  The \defn{type theory $\sfT_S$ associated to} a type signature $S$ is the
  type theory obtained by extending structural MLTT with the type symbols and
  inference rules in $S$ in the obvious well-defined manner. 

  A \defn{context} of a type signature $S$ is a context of the type theory
  $\sfT_S$.
\end{definition}

    


\begin{remark}
  \label{rem:signature-context-minimality}
  (Type declarations are minimally graded.)
  Without loss of generality, we may assume that the grading of a type signature $S =
  \coprod_n S_n$ satisfies the property that for every type declaration $\Gamma
  \vdash A$ in $S_{k}$, $\Gamma$ is not a context of any of the type signatures
  $S_{<j}$ for $j<k$.
\end{remark}

\begin{remark}
  \label{rem:signature-context-canonical-form}
  (Increasing contexts.)
  Without loss of generality, we may assume that for every type declaration
  $\Gamma \vdash A$ in $S$, if we have $\Gamma = x_1\of A_1 , \ldots, x_k\of
  A_k$ where each type $A_i$ is obtained from a type declaration in $S_{n_i}$,
  then for every $1\leq i \leq j \leq k$, we have $n_i \leq n_j$. 
\end{remark}

\begin{para}
  The type theory $\sfT_S$ associated to a signature $S$ satisfies a crucial
  property, which says that all of its terms are variables.
\end{para}

\begin{lemma}
  \label{lem:signature-term-judg-variable}
  Let $S$ be a type signature. Any term judgment of $\sfT_S$ is necessarily of
  the form $\Gamma \vdash x \type A$, where $x\of A $ is in $ \Gamma$.
\end{lemma}
\begin{proof}
  First, a straightforward structural induction proves that \textsf{Wkg} and
  \textsf{Subst} remain admissible in $\sfT_S$. Then
  the result follows by structural induction, since there is only one inference
  rule that concludes with a term judgment.
  \[
      \MLVble{x_1\of A_1,\ldots,x_n\of A_n\ctx}
      {x_1\of A_1,\ldots,x_n\of A_n \vdash x_i\type A_i }\qedhere
    \]
\end{proof}
\begin{corollary}
  \label{cor:signature-context-morphism-variables}
  Any context morphism $\sigma \colon \Gamma \to \Gamma' = x_1 \of A_1 , \ldots
  , x_k \of A_k$ in $\sfT_S$ is of the form $\sigma = (y_1,\ldots , y_k)$, where
  $y_i \of A_i(y_1,\ldots ,y_{i-1}) \in \Gamma$ for every $1\leq i\leq k$.
\end{corollary}

\begin{definition}
  Let $S$ be a type signature. The \defn{category $\lfdex_S$ associated to $S$} is
  defined as follows.
  \begin{enumerate}
  \item For each type declaration $\Gamma \vdash A$ in $S$, $(\Gamma,A)$ is an
    object of $\lfdex_S$.
  \item Morphisms $(\Gamma,A) \to (\Gamma',B)$ are context morphisms
    $(\Gamma',x\of B) \to (\Gamma,x\of A)$, and composition is that of context
    morphisms.
  \end{enumerate}
  In other words, $\lfdex_S\op$ is the full subcategory of the syntactic
  category of $\sfT_S$ consisting of the contexts $\Gamma,x\of A$ for each
  $\Gamma \vdash A$ in $S$.
\end{definition}


\begin{lemma}
  \label{lem:type-sig-is-dir-cat}
  For every type signature $S$, the category $\lfdex_S$ is direct.
\end{lemma}
\begin{proof}
  It suffices that for $\Gamma \vdash A $ in $ S_n$ and $\Gamma' \vdash B $ in $
  S_m$, if there exists a non-identity morphism $(\Gamma',B) \to (\Gamma,A)$,
  then $n>m$. By induction on $n$:
  \begin{enumerate}
  \item $(n=0)$ By~\cref{cor:signature-context-morphism-variables}, given a
    morphism $x\of A \to \Gamma',x\of B$, we necessarily have $A\equiv B$ is
    $x\of A \vdash x:A$. Hence there are no non-identity context morphisms from
    $x\of A$ to $\Gamma',y\of B$.
\item $(n=j+1)$ By \cref{lem:signature-term-judg-variable}, $\Gamma',y\of B$ is
    a context over $S_{<j} \coprod \{\Gamma \vdash A\}$. Hence $m\leq j+1$. If
    $m=j+1$, then pairwise freeness implies that $\Gamma \vdash A = \Gamma'
    \vdash B$ and \cref{lem:signature-term-judg-variable} implies that the only
    term judgment $\Gamma,x\of A \vdash t \type A$ is $\Gamma,x\of A \vdash
    x\type A$, hence the morphism $\Gamma,x\of A \to \Gamma x\of A$ is the
    identity.\qedhere
\end{enumerate}
\end{proof}

\begin{proposition}
  \label{prop:type-sig-is-lfd-cat}
  For every type signature $S$, the category $\lfdex_S$ is a locally finite
  direct category.
\end{proposition}
\begin{proof}
  We will use~\cref{lem:type-sig-is-dir-cat} and
  \cref{prop:loc-finite-finite-cover}. For every $(\Gamma, A)$ in $\lfdex_S$,
  let $J_{(\Gamma,A)}$ be the set of all non-identity morphisms $(\Delta,B) \to
  (\Gamma,A)$ in $\lfdex_S$ that have no non-trivial factorisation $(\Delta,B)
  \to(\Gamma',A') \to (\Gamma,A)$. It is easy to see that $J_{(\Gamma,A)}$ is a
  saturated cover of $(\Gamma,A)$. We will show that $J_{(\Gamma,A)}$ is finite.

  Now, $\Gamma$ is finite, and by
  \cref{cor:signature-context-morphism-variables}, every $\sigma\in
  J_{(\Gamma,A)}$ is a finite list of variables of $\Gamma$. If $J_{(\Gamma,A)}$
  is infinite, there must exist distinct elements $\sigma' = (x'_1,\ldots
  ,x'_k)\colon \Gamma,x\of A \to \Gamma', y\of A'$ and $\sigma =
  (x_1,\ldots,x_l) \colon \Gamma,x\of A \to \Delta , z\of B$ of $J_{(\Gamma,A)}$
  such that every element of $\sigma$ appears in $\sigma'$, \ie we have $x_i =
  x'_{m_i}$ for all $1\leq i\leq l$. But then the tuple $(y_{m_1},\ldots,
  y_{m_l})$ is a non-identity context morphism $\Gamma',y\of A' \to \Delta,z\of
  B$ that factors $\sigma$ non-trivially, which contradicts the definition of
  $J_{(\Gamma,A)}$.
\end{proof}

\begin{proposition}
 \label{prop:type-sig-synt-cat-C-cxl-cat} 
  Let $S$ be a type signature. Then the syntactic category of $\sfT_S$ is
  isomorphic to the contextual category $\cxl{\lfdex_S}$.
\end{proposition}
\begin{proof}
  A straightforward induction, sending every context to the corresponding finite
  $\lfdex_S$\=/cell context.
\end{proof}


\begin{proposition}
 \label{prop:lfd-cat-to-type-sig} 
  Let $\lfdex$ be a \lfd~category. Then there exists a type signature $S$ such
  that $\lfdex_S\cong \lfdex$.
\end{proposition}
\begin{proof}
  We define $S_n$ to be the set of objects of $\lfdex$ of dimension $n$. It
  suffices to prove that for every $n$, $S_{<n}$ is a type signature. We proceed
  by induction, recalling the grading $\lfdex=\bigcup_n\lfdex_n$ by dimension
  from~\cref{lem:lfd-cat-grading}. In the base case, $\emptyset$ is clearly a
  type signature. In the inductive step, for every $c\in
  S_n$,~\cref{prop:type-sig-synt-cat-C-cxl-cat} associates a context $\Delta_c$
  of $S_{<n}$ to the cell context $\partial c$ in $\cxl{\lfdex_n}$, and
  $\Delta_c\vdash c$ is the required type declaration.  
\end{proof}

\begin{remark}
  The previous results establish a correspondence between dependent type
  signatures and \lfd~categories. They are simply the reworking in our setup of
  cell complexes of the same results from \cite{makkai1995folds}. 
\end{remark}

\begin{para}
    [Term signatures]
    \label{para:term-signs}
    A \emph{(dependently typed) term signature} over a type signature $S$ will
    be a generalisation of the set of function symbols of a multisorted
    algebraic theory. A dependently typed term signature will be a graded set of
    generating function symbols, that are dependently sorted by the function
    symbols of lower grading.

    This is best seen with an example, before stating the general definition.
    Consider the type signature corresponding to the \lfd~category $\GG_1$.
    \[
        \vdash D^0 \qquad\qquad x\of D^0,y\of D^0\vdash D^1
    \]
    Then an example of a term signature over $\GG_1$ is the following list of
    judgments.
    \begin{align*}
    &\vdash \star\type D^0\\
    x\of D^0
    &\vdash t(x) \type D^1(x,\star)\\
    x\of D^0
    &\vdash i(x) \type D^1(x,x)\\
    x_1,x_2,x_3\of D^0,y\of D^1(x_1,x_2), z\of D^1(x_2,x_3)
    &\vdash c(z,y)\type D^1(x_1,x_3)
  \end{align*}
  Remark the following properties of every judgment in the previous list.
  \begin{enumerate}
  \item The context to the left of the turnstile (the symbol ``$\vdash$'') is a
    context of the type signature $\GG_1$---namely, it is well-defined
    independently of the term signature.
  \item The output type to the right of the turnstile is well-defined assuming
    the preceding sublist of term judgments (but is not independent of the term
    signature).
  \end{enumerate}
  The second property is a standard one satisfied by all generalised algebraic
  theories. However, the first property imposes a strong restriction on the
  expressive power of GATs. Nevertheless, we will require it of each of our
  term signatures.
\end{para}

\begin{definition}
  Let $S$ be a dependent type signature. An \defn{($S$\=/typed) term signature}
  is a graded set $F\eqdef \coprod_nF_n$ where $F_n$ is a set of term
  declarations of output dimension $n$ over the term signature $F_{<n}\eqdef
  \coprod_{k<n}F_k$.
\end{definition}
\begin{definition}
  A \defn{term declaration} over an $S$\=/typed term signature $F$ is an
  inference rule of the form
  \[
    \prfaxiom{\Delta \vdash t\type A[\sigma]}
  \]
  where $\Delta$ is a context of $S$, $(\Gamma\vdash A)$ a type declaration in $S$,
  $\sigma\colon \Delta \to \Gamma$ a context morphism of the type theory $\sfT_F$
  associated to $F$, and where $t$ is a
  fresh term symbol (namely, it does not appear in $F$). The \defn{output
    dimension} of the term declaration $\Delta\vdash t\of A[\sigma]$ is the
  dimension of $(\Gamma,A)$ in $\lfdex_S$.
\end{definition}


\begin{definition}
  The \defn{type theory $\sfT_{S,F}$ associated to} an $S$\=/typed term signature $F$ is the
  type theory obtained by extending $\sfT_S$ with the term symbols and
  inference rules in $F$ in the obvious well-defined manner.
\end{definition}

\begin{para}
  [Term signatures and collections]
  \label{para:term-sign-collections}
  Let $F$ be an $S$\=/typed term signature. We will write $\sfC_F$ for the
  syntactic category of $\sfT_{S,F}$. There is an evident functor $f\colon
  \cxl{\lfdex_S}\to \sfC_F$ exhibiting $\sfC_F$ as a $\lfdex_S$\=/contextual
  category. Recall that we have an inclusion $i\colon \lfdex_S\op \subto
  \cxl{\lfdex_S}$. Then the functor
  \[
    \relHom{\sfC_F}{f-}{fi-} : \cxl{\lfdex_S}\op\times \lfdex_S\op\tto \Set
  \] is a $\lfdex_S$\=/collection (\cref{def:C-collection}).

  Conversely, given a $\lfdex_S$\=/collection $X\in\coll{\lfdex_S}$, we obtain
  an $S$\=/typed term signature $F$ by defining $F_n\eqdef \coprod X(\Delta,
  (\Gamma,A))$, where the coproduct is over all $\Delta \in \cxl{\lfdex_S}$
  and all $(\Gamma\vdash A)$ in $S_n$.
\end{para}

\begin{remark}
    The $\lfdex_S$\=/collection $\relHom{\sfC_F}{f-}{fi-}$ constructed from the
    term signature $F$ is a monoid in $\coll{\lfdex_S}$ (since it is obtained
    from the $\lfdex_S$\=/contextual category $\sfC_F$). It is not in general
    \emph{free} on some underlying collection (equivalently, a free monad on a
    finitary endofunctor). This is because in the term signature $F$, the output
    type of a term declaration of dimension $n$ may contain \emph{terms} (and
    not just term declarations) in its dependencies. Nevertheless, the monoid
    $\relHom{\sfC_F}{f-}{fi-}$ in $\coll{\lfdex_S}$ satisfies a kind of
    ``polygraphic freeness'', since it is constructed by generators in every
    dimension.

    Therefore, the constructions in~\cref{para:term-sign-collections} are not
    mutually inverse---they do not establish a correspondence between term
    signatures and collections. However, the term signatures arising from
    collections are such that the output types of their term declarations
    (generators) contain only term declarations (generators of lower dimension).
    This is analogous to the characterisation of strict $\omega$\=/categories
    free on globular sets as those polygraphs such that the sources and targets
    of generators are themselves generators of lower dimension.
\end{remark}


\begin{para}
  [Theories]
  \label{para:dep-typed-theories}
  A \emph{(dependently typed) theory} over an $S$\=/typed term signature $F$
  will be the generalisation of the set of equations of a multisorted algebraic
  theory.

  Continuing the example of the $\GG_1$\=/typed term signature
  in~\cref{para:term-signs}, consider the following equations.
  \begin{align*}
    x,y\of D^0,z\of D^1(x,y)
    &\vdash c(t(y),z) = t(x) \type D^1(x,\star)\\
    x,y\of D^0,z\of D^1(x,y)
    &\vdash c(i(y),z) = z \type D^1(x,y)\\
    x,y\of D^0,z\of D^1(x,y)
    &\vdash c(z,i(x)) = z \type D^1(x,y)\\
    x_1,x_2,x_3,x_4\of D^0, f\of D^1(x_1,x_2),\\ g\of D^1(x_2,x_3),
    h\of D^1(x_3,x_4)
    &\vdash c(h,c(g,f))= c(c(h,g),f) \type D^1(x_1,x_4)
  \end{align*}
  It is quite clear that this equational theory corresponds to the finitary
  monad on $\psh{\GG_1}$ whose algebras are categories with a chosen terminal
  object. 
  Remark that the previous equations all satisfy the same property from
  \cref{para:term-signs}, namely:
  \begin{enumerate}
  \item
   \label{theories-crucial-property-contexts} 
    The context to the left of the turnstile is a
    context of the type signature $\GG_1$---namely, it is well-defined
    independently of the term signature and the other equations.
  \end{enumerate}
  The previous presentation of the theory of categories with a chosen terminal
  object is not necessarily the most obvious one---we might be tempted to
  replace the first equation with
  \[
    x\of D^0, y\of D^0(x,\star) \vdash t(x)=y \type D^1(x,\star).
  \]
  However, this equation no longer satisfies the property
  \ref{theories-crucial-property-contexts}, which we will require of all our
  theories.
\end{para}

\begin{definition}
   \label{def:dependently-typed-theory-syntactic} 
    Let $F$ be an $S$\=/typed signature. A \defn{(dependently typed) theory over
      $F$} is a graded set $E\eqdef\coprod_nE_n$, such that each $E_n$ is a set
    of identifications of output dimension $n$ in the theory $E_{<n}\eqdef
    \coprod_{i<n}E_i$ over $F$.
\end{definition}

\begin{definition}
    An \defn{identification} in a theory $E$ over $F$ is an inference rule
    \[
        \prfaxiom{\Delta \vdash t=u \type A[\sigma]}
    \]
    where $\Delta$ is a context of $S$, $(\Gamma\vdash A)$ is a type declaration
    in $S$, and where $\sigma\colon \Delta\to \Gamma$ is a context morphism and
    $(\Delta\vdash t\type A[\sigma])$, $(\Delta\vdash u\type A[\sigma])$ are
    term judgments of the type theory associated to $E$. Its \defn{output
      dimension} is the dimension of the object $(\Gamma,A)$ in $\lfdex_S$.
\end{definition}

\begin{definition}
    The \defn{type theory $\sfT_{S,F,E}$ associated to} a theory $E$ over $F$ is
    the type theory obtained by extending $\sfT_{S,F}$ with the inference rules
    in $E$ in the obvious well-defined manner.
\end{definition}

\begin{para}
    [Theories to monoids in $\coll{\lfdex_S}$]
   \label{para:theories-to-monoids} 
    The canonical functor $\cxl{\lfdex_S}\to\sfC_E$ exhibits the syntactic
    category $\sfC_E$ of $\sfT_{S,F,E}$ as a $\lfdex_S$\=/contextual category.
    Then, by~\cref{thm:classification-dep-alg-theories}, it is associated to a
    finitary monad $T_{S,F,E}$ on $\psh{\lfdex_S}$ whose algebras are the models
    in $\Set$ of the theory $E$.
\end{para}

\begin{para}
    [Free monoids in $\coll{\lfdex_S}$ to theories]
    \label{para:theory-of-free-monads}
    Let $T\in\coll{\lfdex_S}$ be a $\lfdex_S$\=/collection, and let $F_T$ be the
    $S$\=/typed term signature obtained using the construction
    from~\cref{para:term-sign-collections}. We define a dependently typed theory
    $E_T$ over $F_T$ whose associated $\lfdex_S$\=/contextual category
    $\sfC_{E_T}$ corresponds to the free monad $T^\infty$ on the finitary
    endofunctor $T\colon\psh{\lfdex_S}\to\psh{\lfdex_S}$.

    We proceed by induction on $n\in\NN$. For every $(\Gamma,A)\in\lfdex_S$ of
    dimension $n$, every $\sigma\colon \Delta'\to\Delta$ in $\cxl{\lfdex_S}$ and
    every $t\in T(\Delta,(\Gamma,A))$ (a term declaration of $F_T$), then
    $\Delta'\vdash t[\sigma] \type A[\rho\sigma]$ is a term judgment in the type
    theory $\sfT_{S,F_T}$ (where $\rho\colon \Delta\to \Gamma$ is the context
    morphism of $\sfT_{S,F_T}$ associated to $t$). Since $T$ is a
    $\lfdex_S$\=/collection, we have an element $T_\sigma(t)\in
    T(\Delta',(\Gamma,A))$. We introduce the identification $\Delta'\vdash
    t[\sigma] = T_\sigma(t)\type A[\rho\sigma]$ in the theory $(E_T)_{<n}$
    (where the composite context morphism $\rho\sigma$ in the type theory
    $\sfT_{S,F,(E_T)_{<n}}$ is well defined by induction hypothesis on $n$).

    Let $E_T$ be the set of all identifications $\Delta'\vdash t[\sigma] =
    T_\sigma(t)\type A[\rho\sigma]$, graded by output dimension. Then for any
    context $\Delta$ of $S$ and $(\Gamma,A)\in\lfdex_S$, a context morphism
    $\Delta\to(\Gamma,x\of A)$ in the type theory $\sfT_{S,F_T,E_T}$ is
    precisely an element in $T^\infty(\Delta,(\Gamma,A))$ of the free monoid
    $T^\infty\in\coll{\lfdex_S}$, and the $\lfdex_S$\=/contextual category
    $\sfC_{E_T}$ corresponds to the free monoid $T^\infty$.
\end{para}

\begin{para}
    [Monoids in $\coll{\lfdex_S}$ to theories]
    \label{para:theory-of-gen-monad}
    A finitary monad $T$ on $\psh{\lfdex_S}$ is a monoid in $\coll{\lfdex_S}$.
    Hence it is an algebra of the free-monoid monad on $\coll{\lfdex_S}$. Its
    algebra map is a natural transformation $\alpha\colon T^\infty\to T$. Let
    $F_T$ be the $S$\=/typed term signature associated to the collection $T$ via
    the construction in \cref{para:term-sign-collections}, and let $E_T$ be the
    theory of the free monad $T^\infty$ constructed in
    \cref{para:theory-of-free-monads}. We will add identifications to $E_T$ to
    obtain a theory $E^T$ over $F_T$ corresponding to the monad $T$.

    Let $(\Gamma,A)\in\lfdex_S$ of dimension $n$ and $\Delta$ in
    $\cxl{\lfdex_S}$. A context morphism $\sigma\colon\Delta\to\Gamma)$ and a
    term judgment $\Delta\vdash t\type A[\sigma]$ in the type theory
    $\sfT_{S,F_T,E_T}$ are precisely the data of an element $t\in
    T^\infty(\Delta,(\Gamma,A))$. Then $\alpha(t)$ is an element of
    $T(\Delta,(\Gamma,A))$, namely we have a term judgment $\Delta\vdash
    \alpha(t)\type A[\alpha(\sigma)]$ where $\alpha(\sigma)$ is the substitution
    obtained by applying $\alpha$ to the terms of $\sigma$. Then, by induction
    hypothesis on $n$, we have that $\sigma=\alpha(\sigma)$ as context morphisms
    of the type theory $\sfT_{S,F_T,E^T_{<n}}$. We introduce the identification
    $\Delta\vdash t = \alpha(t)\type A[\sigma]$ in the theory $E^T_{<n}$. Let
    $E^T$ be theory over $F_T$ that is the extension of $E_T$ with all these
    identifications, graded by dimension. The canonical morphism
    $\sfC_{E_T}\to\sfC_{E^T}$ of syntactic categories corresponds to the
    morphism $T^\infty\to T$ of monoids in $\coll{\lfdex_S}$.

    Along with \cref{para:theories-to-monoids}, this establishes a
    correspondence between theories over $S$\=/typed term signatures and monoids
    in $\coll{\lfdex_S}$.
\end{para}
\begin{remark}
    We have given a syntactic description
    (\cref{def:dependently-typed-theory-syntactic}) of the dependently sorted
    algebraic theories that correspond to the algebraic classification
    of~\cref{thm:classification-dep-alg-theories}. These are a strict subclass
    of Cartmell's generalised algebraic theories (GATs)---notably, they are
    subject to the restriction~\ref{theories-crucial-property-contexts}.
\end{remark}
 
\begin{remark}
    \label{rem:gats-vs-C-theories-expressivity}
    An obvious question is whether the class of all GATs is \emph{strictly} more
    expressive than the subclass of dependently sorted algebraic theories in our
    sense. We will partially answer this in the negative
    in~\cref{sec:morita-equiv-eats-C-theories}.
\end{remark}



\chapter{Models of \texorpdfstring{$\lfdex$}{C}-contextual categories}
\label{cha:models-of-C-contextual-categories}

In this chapter, we will develop the notion of a \emph{model in spaces} of the
dependently sorted algebraic theories classified
by~\cref{thm:classification-dep-alg-theories}. We will see that certain
techniques from the theory of homotopy models of multisorted algebraic theories
from
\cite{Schwede2001stablehtpyAlgTheories,badzioch2002algtheories,bergner2006rigidification,Rezk2002simplicialAlgTheories}
make sense for models in spaces of dependently sorted algebraic theories.

In \cref{sec:morita-equiv-eats-C-theories}, we classify the $1$\=/categories
that arise as categories of models in \emph{discrete spaces} (sets) of
dependently sorted algebraic theories. We will see that these are exactly the
locally finitely presentable $1$\=/categories
(\cref{thm:lfp-cat-models-of-C-sorted-theory}).

\section{Morita equivalence with essentially algebraic theories}
\label{sec:morita-equiv-eats-C-theories}
Our goal in this section is to give a partial answer to the questions raised
in~\cref{rem:is-cxl-cat-C-cxl-cat,rem:gats-vs-C-theories-expressivity}. We will
show that our dependently sorted algebraic theories are \emph{Morita equivalent}
to the (finitary) \emph{essentially algebraic theories} of Freyd
(\cite[\textsection{1}]{freyd1972aspects}, \cite[3.D]{Adamek1994a}).

We begin with a short review of the theory of accessible cocontinuous
localisations of locally presentable categories from \cite{gabrielulmer1971}.
The content of this review is well-known and is entirely subsumed
by~\cref{part:part2}. We have included it for convenience and to fix
terminology.

\begin{para}
    [Orthogonality]
    \label{para:GU-rappels-orthogonality}
    Let $\cC$ be a category, and $l,r \in \cC\arr$ be morphisms of $\cC$. We say
    that $l$ is \emph{left orthogonal to}\index{orthogonal} $r$ (equivalently,
    $r$ is \emph{right orthogonal} to $l$), written $l \perp
    r$\index{$\perp$|see {orthogonal}}, if for any solid commutative square as
    follows, there exists a \emph{unique} dashed arrow making the two triangles
    commute (the relation $\perp$ is also known as the \emph{unique lifting
      property}\index{unique lifting property|see {orthogonal}}).
    \[
        \begin{tikzcd}
            \cdot \ar[d, "l" left] \ar[r] & \cdot \ar[d, "r"] \\
            \cdot \ar[r] \ar[ur, dashed] & \cdot
        \end{tikzcd}
    \]
    If $\cC$ has a terminal object $1$, then for all $X \in \cC$, we write $l
    \perp X$ if $l$ is left orthogonal to the unique map $X \to 1$. Let $\cL$
    and $\cR$ be two classes of morphisms of $\cC$. We write $\cL \perp \cR$ if
    for all $l \in \cL$ and $r \in \cR$ we have $l \perp r$. The class of all
    morphisms $f$ such that $\cL \perp f$ (respectively, $f \perp \cR$) is
    denoted $\cL^\bot$ (respectively, ${}^\bot \cR$).
\end{para}


\begin{para}
    [Cocontinous localisations]
    \label{para:GU-rappels-localisations}
    Let $J$ be a class of morphisms of a cocomplete category $\cC$. Recall that
    the \emph{cocontinuous localisation of $\cC$ at $J$} is a cocontinuous
    functor $\gamma_J : \cC \to J\inv \cC$ such that $\gamma_J f$ is an
    isomorphism for every $f \in J$, and such that $\gamma_J$ is initial (in the
    $2$\=/category of cocontinuous functors with domain $\cC$) for this
    property. We say that $J$ has the
    \emph{``2~out~of~3''}\index{``2~out~of~3''} property when for every
    composite $gf$ of a pair of composable morphisms $f,g$ in $\cC$, if any two
    of $f$, $g$ and $gf$ are in $J$, then so is the third.
\end{para}

\begin{para}
    [Locally presentable categories]
    Assume now that $C$ is a small category, and that $J$ is a \emph{set}
    (rather than a proper class) of morphisms of $\psh C$, and consider the full
    subcategory $C_J \subto \psh C$ of all those presheaves $X \in \psh C$ such
    that $J \perp X$. A category is \emph{locally presentable}\index{locally
      presentable category} if and only if it is equivalent to one of the form
    $C_J$. The pair $({}^\bot (J^\bot), J^\bot)$ forms an \emph{orthogonal
      factorisation system}\index{orthogonal factorization system}, meaning that
    any morphism $f$ in $\psh C$ can be factored as $f = p i$, where $p \in
    J^\bot$ and $i \in {}^\bot (J^\bot)$. Applied to the map $X \to 1$ to the
    terminal presheaf, this factorisation provides a left adjoint (i.e. a
    reflection) $a_J : \psh C \to C_J$ to the full inclusion $C_J \subto \psh
    C$. Then $a_J$ is the cocontinuous localisation of $\psh C$ at $J$.
    Moreover, the full inclusion $C_J\subto \psh C$ preserves
    $\kappa$\=/filtered colimits for any regular cardinal $\kappa$ such that
    every map in $J$ is $\kappa$\=/presentable. Namely, the localisation is
    \emph{accessible}.
\end{para}

\begin{para}
    All of the previous paragraph still holds when $\psh C$ is replaced by a
    locally presentable category $\cE$ and $J$ is a set of morphisms of $\cE$.
    We call these cocontinuous accessible localisations $\cE\localisation\cE_J$
    of locally presentable categories \defn{Gabriel-Ulmer localisations}.
\end{para}
\begin{para}
    [Local isomorphisms] With $\cE$ and $J$ as in the previous paragraph, the
    class $W_J$ of \emph{$J$-local isomorphisms}\index{local isomorphism} is the
    class of all morphisms $f \in \cE\arr$ such that for all $X \in \cE_J$, $f
    \perp X$, that is, $W_J = {}^\bot \cE_J$. In general, ${}^\bot (J^\bot)\subset
    W_J$ is a strict inclusion (see~\cref{para:reflective-FS}). The class $W_J$
    is the smallest class containing $J$ that satisfies the ``2~out~of~3''
    property and that is closed under colimits in $\cE\arr$ \cite[Satz
    8.5]{gabrielulmer1971}. Thus the Gabriel-Ulmer localisation $a_J$ is also
    the ordinary localisation of $\cE$ at $W_J$.
\end{para}

\begin{para}
    A \emph{presentation} of a locally $\kappa$-presentable category $\cC$ is a
    small category $A$ and a $\kappa$-accessible fully faithful right adjoint
    $\cC\subto \psh A$. Equivalently, $\cC$ is the category of algebras of an
    idempotent $\kappa$-accessible monad on $\psh A$. The left adjoint is a
    Gabriel-Ulmer localisation of $\psh A$ at a set of $\kappa$\=/presentable
    maps.
\end{para}

We are now ready to tackle the main goal of this section.

\begin{para}
    [Morita equivalence] We will say that two classes of theories are
    \defn{Morita equivalent} if for every theory of each class, there exists a
    theory of the other class with the same category of $\Set$\=/models. It is
    well-known \cite[Thm 3.36]{Adamek1994a} that a category $\cC$ is locally
    finitely presentable if and only if it is the category of $\Set$\=/models of
    an \emph{essentially algebraic theory} (EAT). Since Morita equivalence is a
    semantic notion, any class of theories $\cT$ is Morita equivalent to EATs if
    and only if every category of $\Set$\=/models of a theory in the class $\cT$
    is locally finitely presentable, and every locally finitely presentable
    category is the category of $\Set$\=/models of a theory in the class $\cT$.
    For instance, \cite[\textsection{1.4}]{cartmell1978generalised}\cite[Sec.
    6]{cartmell1986} outlines a syntactic technique that transforms any EAT into
    a GAT with the same category of $\Set$\=/models. The category of
    $\Set$\=/models of any GAT is the category of $\Set$\=/models of its
    syntactic contextual category (\cref{def:model-of-cxl-cat}), and so is
    clearly locally finitely presentable (since it is the category of
    $\Set$\=/models of a finite-limit sketch). Thus the classes of EATs and GATs
    are Morita equivalent.

    Applying Cartmell's syntactic technique to an EAT does not necessarily
    produce a GAT that is a dependently sorted algebraic theory (in our sense).
    Indeed, the syntactic transformation generally produces a GAT that does not
    respect the constraint~\ref{theories-crucial-property-contexts}. In the rest
    of this section, we will use a different technique to show that every
    locally finitely presentable category is classified by a dependently sorted
    algebraic theory. \emph{A fortiori}, this implies that the classes of EATs,
    GATs and dependently sorted algebraic theories are all Morita equivalent.
\end{para}

\begin{remark}
    The technique that we will use is remarkable in its simplicity---it is the
    quickest proof that I know of the Morita equivalence of EATs and GATs.
    Moreover, it seems likely to be exportable to homotopical models of dependently
    sorted algebraic theories.
\end{remark}

\begin{para}
    Let $\bDelta'$ denote the category of semi-simplices, namely the non-full
    subcategory of $\Cat$ on the finite non-empty ordinals and
    \emph{monomorphisms}. Recall from~\cref{exa:lfd-cats} that $\bDelta'$ is
    locally finite and direct. The \emph{semi-simplicial nerve} of a small
    category $A$ is the image of $A$ under the nerve functor
    $\Cat\to\psh{\bDelta'}$ of the non-full inclusion $\bDelta' \to \Cat$. For
    any $A$ in $\Cat$, the comma category $\bDelta'/A$ is the category of
    elements of the presheaf on $\bDelta'$ that is the semi-simplicial nerve of
    $A$. There is also a functor $\tau_A\colon \bDelta'/A\to A$, natural in $A$,
    that sends $f\colon [n]\to A$ to $f(n)$, and $[m]\xto{\alpha}[n]\xto{f}A$ to
    $f(\alpha(m)\leq n)$.
\end{para}
\begin{lemma}
    \label{lem:semisimpl-nerve-lfd}
    For any $A$ in $\Cat$, $\bDelta'/A$ is locally finite and direct.
\end{lemma}
\begin{proof}
    By the previous discussion, $\bDelta'/A \to \bDelta'$ is a discrete
    fibration. Then the result follows from~\cref{cor:presheaf-elts-direct}
    and~\cref{para:lfd-cat-intro}.
\end{proof}

\begin{proposition}
    \label{prop:semsimp-nerve-fully-faithful-upper-star}
    For any $A$ in $\Cat$, the functor $\tau_A^*\colon \psh A \to
    \psh{\bDelta'/A}$ is fully faithful.
\end{proposition}
\begin{proof}
    This is \cite[Prop. 6.9]{cisinski2003imagesdirectes}, applied to the model
    structure on $\Set$ with cofibrations and fibrations all maps, and weak
    equivalences the isomorphisms.
\end{proof}

\begin{proposition}
    \label{thm:loc-pres-cat-lfd-presentation}
    Every locally $\kappa$\=/presentable category $\cC$ has a presentation
    $\cC\subto \psh D$ where $D$ is locally finite and direct.
\end{proposition}
\begin{proof}
    $\cC$ has a presentation $\cC\subto\psh A$ for some $A$ in $\Cat$. Since
    $\tau_A^*$ is both a left and a right adjoint, and is fully faithful, the
    composite $\cC\subto \psh A\subto \psh{\bDelta'/A}$ is a
    $\kappa$\=/accessible fully faithful right adjoint. We conclude using
    \cref{lem:semisimpl-nerve-lfd}.
\end{proof}
\begin{theorem}
 \label{thm:lfp-cat-models-of-C-sorted-theory} 
    Every locally finitely presentable category is the category of models of a
    $\lfdex$-contextual category for some locally finite direct category
    $\lfdex$.
\end{theorem}
\begin{proof}
    By~\cref{thm:loc-pres-cat-lfd-presentation}, $\cC$ is the category of
    algebras of an idempotent finitary monad on $\psh{\bDelta'/A}$. Thus $\cC$ is
    the category of models of an idempotent $\bDelta'/A$\=/sorted theory
    (see~\cref{def:C-sorted-theory,rem:idempotent-C-sorted-th-refl-localisn}).
\end{proof}

\begin{corollary}
 \label{thm:morita-equiv-gat-eat-C-cxl-cats} 
    The classes of
    \begin{enumerate}
    \item essentially algebraic theories,
    \item generalised algebraic theories,
    \item and dependently sorted algebraic theories
    \end{enumerate}
    are Morita equivalent.
\end{corollary}

\begin{remark}
    \cite[Prop. 6.9]{cisinski2003imagesdirectes} holds for derivators of all
    co/complete model categories $\cM$ and not just $\Set$. Thus we can expect
    homotopical versions of the previous results for locally presentable
    \oo-categories. 
\end{remark}


\section{Background on simplicial model categories}
\label{sec:backgr-simpl-model-cats}
We assume familiarity with the essentials of the theory of model categories. We
will write $(\cM,\rmW,\rmC,\rmF)$ for a model category $\cM$ along with its
classes of \emph{weak equivalences}, \emph{cofibrations}, and \emph{fibrations}
respectively. A \emph{trivial cofibration} (respectively, \emph{trivial
  fibration}) is a map in $\rmW\cap\rmC$ (respectively, in $\rmW\cap\rmF$). When
the model structure is clear from context, we will simply write $\cM$. A model
category $\cM$ is \emph{combinatorial} if its underlying category is locally
presentable and its model structure cofibrantly generated. We will give a quick
run-through of the (small fragment of a vast) theory of combinatorial simplicial
model categories that is germane to our situation.

We will call the presheaf category $\SSet\eqdef\psh\bDelta$ of simplicial sets,
equipped with the usual cofibrantly generated Kan-Quillen model structure, the
\defn{model category of spaces}. We will write $\Delta_n\in\SSet$ for the
representable $n$-simplex. We use the conventional sets of generating
cofibrations and trivial cofibrations called the \emph{boundary} and \emph{horn}
inclusions, and denoted $I$ and $J$ respectively.
\begin{align*}
  I&\eqdef\{\delta_n\colon\partial\Delta_n\mono \Delta_n\mid [n]\in\bDelta\}\\
  J&\eqdef\{h^k_n\colon\Lambda^k_n\mono \Delta_n\mid [n]\in\bDelta,0\leq k\leq n\}
\end{align*}

Since $[0]$ is a terminal object in $\bDelta$, the
reflective adjunction $e_\Delta:\bDelta\localisation 1:[0]$ in $\deux\Cat$ gives
a reflective adjunction $(e_\Delta)_!:\SSet\localisation \Set:e_\Delta^*$ of
presheaf categories. We call the fully faithful right adjoint
$e_\Delta^*\colon\Set\subto\SSet$ (that takes a set to a constant presheaf) the
\emph{inclusion of discrete spaces into spaces}. Since $e_\Delta^*$ is also a
left adjoint, it preserves all limits and colimits.


\begin{para}
  [Joyal-Tierney calculus]
  Let $\cC$ be a $\SSet$\=/enriched category. If $f\colon A\to B$ and $g\colon
  X\to Y$ are maps in $\cC$, their \emph{pullback-hom} $\pbh f g$ is defined to
  be the following cartesian gap map in $\SSet$ (where $\map fg$ is just the
  hom-simplicial-set of the arrow category $\cC\arr$).
  \[
    \begin{tikzcd}[sep=scriptsize]
      \map BX \ar[drr,bend left=15,shift left=1.5,"\map Bg"] \ar[ddr,bend
      right,"\map fX"']
      \ar[dr,dashed,"\pbh fg"]\\
      &\map fg \pbmark \ar[r,] \ar[d]
      &\map BY \ar[d,"\map fY"]\\
      &\map AX \ar[r,"\map Ag"'] &\map AY
    \end{tikzcd}
  \]
  If $f\colon K\to L$ is a map of simplicial sets, and $g\colon X\to Y$ is a map
  in $\cC$, the \emph{pushout-product} $f\pp g$ is defined to be the following
  cocartesian gap map in $\cC$ (should the various tensors and the pushout
  exist in $\cC$).
  \[
    \begin{tikzcd}[sep=scriptsize]
      K\otimes X \ar[r] \ar[d]
      &K\otimes Y \ar[d] \ar[ddr,bend left]\\
      L\otimes X \ar[r] \ar[drr,bend right=15]
      &L\otimes X \coprod_{K\otimes X} K\otimes Y \pomark \ar[dr,dashed,"f\pp g"]\\
      &&L\otimes Y
    \end{tikzcd}
  \]
  The presheaf category $\SSet$ is complete and cocomplete as an
  $\SSet$\=/enriched category---then it is readily seen that the pushout-product
  and pullback-hom on maps in $\SSet$ make
  $(\SSet\arr,\pp,\pbh--,\emptyset\to\Delta_0)$ a closed symmetric monoidal
  $\SSet$\=/enriched category.
  
  Then for an arbitrary $\SSet$\=/enriched category $\cC$, the pullback-hom
  functor $\pbh--\colon(\cC\arr)\op\times\cC\arr\to\SSet\arr$ defines an
  enrichment of $\cC\arr$ over $\SSet\arr$, and if $\cC$ has enough
  (weighted)\footnote{In an enriched setting, ``(co)limit(s)'' will always mean
    weighted (co)limit(s), and if ever we confound conical (co)limits with the
    underlying ordinary (co)limits, it is because they coincide.} colimits, then
  the $\SSet\arr$\=/enriched category $\cC\arr$ is tensored over $\SSet\arr$ via
  the pushout-product $\pp\colon \SSet\arr\times\cC\arr\to\cC\arr$. If $\cC$ has
  enough limits, then $\cC\arr$ is cotensored over $\SSet\arr$; with slight
  abuse, we will write the cotensor as $\pbh--\colon
  (\SSet\arr)\op\times\cC\arr\to\cC\arr$. For $f\colon K\to L$ in $\SSet$ and
  $g\colon X\to Y$ in $\cC$ we readily calculate the cotensor $\pbh fg$ as the cartesian gap
  map below in $\cC$.
  \[
    \begin{tikzcd}[sep=scriptsize]
       X^L \ar[drr,bend left=20,shift left=1.5," g^L"] \ar[ddr,bend
      right," X^f"']
      \ar[dr,dashed,"\pbh fg"]\\
      & X^K\times_{Y^K}Y^L \pbmark \ar[r] \ar[d]
      & Y^L \ar[d," Y^f"]\\
      & X^K \ar[r," g^K"'] & Y^K
    \end{tikzcd}
  \]
\end{para}

\begin{remark}
  The important fact to retain is that for any $f$ in $\SSet\arr$ and $g,h$ in
  $\cC\arr$, we have a natural isomorphism $\pbh{f\pp g}h\cong \pbh g{\pbh fh}$
  of morphisms of simplicial sets. The compatibility with the $\SSet$-enrichment
  says that we also have a natural isomorphism $\map {f\pp g}h\cong \map g{\pbh
    fh}$ of hom-simplicial sets. We will see the construction again in the context of
  \oo-categories in \cref{chap:fact-systems}.
\end{remark}

\begin{remark}
  Let $\cC$ be any $\SSet$\=/enriched category. Then for any two maps $f,g$ in its
  underlying ($\Set$\=/enriched) category $\cC_0$, $f$ has the left lifting
  property against $g$ (denoted ``$f\pitchfork g$'')
  if and only if the morphism $\pbh fg$ of simplicial sets is surjective on
  $0$-simplices.
\end{remark}
\noindent Let $\cC$ be a complete and cocomplete $\SSet$\=/enriched category, with
a model structure on its underlying category. Then $\cC$ is a \defn{simplicial
  model category} if any of the following equivalent conditions holds, for
every cofibration $f\colon K\to L$ in $\SSet$, cofibration $g\colon A\to B$ in
$\cC$ and fibration $h\colon X\to Y$ in $\cC$.
\begin{enumerate}
\item $f\pp g$ is a cofibration in $\cC$ that is trivial if either $f$ or $g$ is
  so.
\item $\pbh f h$ is a fibration in $\cC$ that is trivial if either $f$ or $h$ is
  so.
\item $\pbh gh$ is a fibration in $\SSet$ that is trivial if either $g$ or $h$
  is so.
\end{enumerate}
The model category $\SSet$ is the archetypical {simplicial model category}.

\begin{remark}
  The preceding yoga can be exported \emph{mutatis mutandis} to $\cV$\=/enriched
  model categories for suitable closed monoidal model categories $\cV$.
\end{remark}

\begin{para}
  [Simplicial presheaves]
  \label{para:simpl-presheaves}

  For any small category $C$, we write $\sPsh C$ for the functor category
  $[C\op,\SSet]$ which, as a category of diagrams in a topos, is enriched,
  tensored and cotensored over $\SSet$. For $X,Y\in\sPsh C$ and $K\in\SSet$, we
  will write $\map XY\in\SSet$ for the hom-space (given by the enrichment), and
  $X\otimes K,X^K\in\sPsh C$ for the tensor and cotensor respectively. We recall
  the formulas (where $K^L$ is the internal hom in $\SSet$)
  \[
    \map XY = \int_{c\in C} Y_c^{X_c} \quad,\quad (X\otimes K)_{c} = X_c \times K
    \quad,\quad (X^K)_c = X_c^K.  
  \]
  The inclusion of discrete spaces into spaces gives a reflective adjunction
  $\sPsh C\localisation \psh C$.
  The right adjoint $\psh C\subto\sPsh C$ is also a left adjoint, hence
  it preserves all limits and colimits. We write $r\colon C\subset \psh C\subto
  \sPsh C$ for the composite with the Yoneda embedding of $C$ (thus $r$
  preserves limits). The Yoneda embedding $C\times \bDelta\subset \sPsh C$ is
  just the functor $(c,[n])\mapsto rc\otimes\Delta_n$.

  We will call a map $f\colon X\to Y$ in $\sPsh C$ an \defn{objectwise} or
  \defn{global weak
    equivalence} if for every $c\in C$, $f_c\colon X_c\to Y_c$ is a weak equivalence
  in $\SSet$. We will call any model structure on $\sPsh C$ whose weak
  equivalences are the objectwise weak equivalences a \defn{global model
    structure}. There are always two canonical global model structures on $\sPsh
  C$:
  \begin{enumerate}
  \item the \defn{projective} model structure, in which a map $f\colon X\to Y$
    is a fibration if and only if for each $c$ in $C$, $f_c\colon X_c\to Y_c$ is
    a fibration (equivalently, $\pbh {\emptyset\to c}f$ is a fibration) in
    $\SSet$ ,
  \item and the \defn{injective} model structure, in which a map $f\colon X\to
    Y$ is a cofibration if and only if for each $c$ in $C$, $f_c\colon X_c\to
    Y_c$ is a cofibration in $\SSet$ (equivalently, $f$ is a monomorphism in
    $\sPsh C$).
  \end{enumerate}
  We write these model categories as $(\sPsh C\proj,\rmW,\rmC\proj,\rmF\proj)$
  and $(\sPsh C\inj,\rmW,\rmC\inj,\rmF\inj)$ respectively. Both are
  combinatorial simplicial model categories for the $\SSet$-enrichment and
  co-/tensor on $\sPsh C$ \cite[A.3.3.2, A.3.3.4]{LurieHT}, and the identity
  adjunction $1_{\sPsh C} : \sPsh C\proj\rightleftarrows\sPsh C\inj : 1_{\sPsh
    C}$ is a Quillen equivalence. For any functor $u\colon A\to B$ between small
  categories, the adjunction $u_!\dashv u^*$  (respectively, $u^*\dashv u_*$) between
  $\sPsh A$ and $\sPsh B$ is Quillen for the projective (respectively,
  injective) model structure.    

  The advantage of the projective model structure is that its fibrations and
  trivial fibrations are objectwise. Moreover, there exist easy-to-describe sets
  of generating cofibrations and trivial cofibrations---for instance the sets
  $I\proj$ and $J\proj$ below.
  \begin{align*}
    I\proj&\eqdef\{\delta_n\pp (\emptyset\to c)\mid c\in
            C,[n]\in\bDelta\}\\
    J\proj&\eqdef\{h^k_n\pp (\emptyset \to c)\mid c\in
            C,[n]\in\bDelta,0\leq k\leq n\}
  \end{align*}
  The inconvenience of the projective model structure is that, in many
  situations, it has too few cofibrations---all its cofibrations are of course
  monomorphisms, but few monomorphisms are projective cofibrations. Dually, the
  cofibrations of the injective model structure are just the monomorphisms, but
  the injective fibrations and trivial fibrations become difficult to detect.
  There exist a host \cite{isaksen2004flasque,jardine2006intermediate} of
  \emph{intermediate} global model structures on $\sPsh C$ whose cofibrations
  contain the projective cofibrations but are still monomorphisms, namely the
  identity adjunction forms a sequence of Quillen equivalences $\sPsh
  C\proj\rightleftarrows\sPsh C_{\text{int}}\rightleftarrows\sPsh C\inj$ (where
  $\sPsh C_{\text{int}}$ denotes the intermediate model structure). With
  additional structure on the category $C$ (such as being a \emph{Reedy}
  category), some intermediate model structures still allow for a simple
  characterisation of (trivial) fibrations.

\end{para}

\begin{para}
  [Left Bousfield localisation]
  \label{para:left-bousfield-localisation}
  The machinery of left Bousfield localisation of combinatorial model categories
  is the model-category-theoretic approach to the theory of Gabriel-Ulmer
  (=accessible cocontinuous) localisations of locally presentable
  \oo-categories, which constitute the \emph{objet d'étude} of our
  \cref{part:part2}.

  Let $\cM$ be a simplicial model category. Recall that for every $X,Y\in\cM$ we
  can define a Kan fibrant, functorial \emph{homotopy function complex} $\hmap
  XY\in\SSet$ such that if $X$ is cofibrant and $Y$ is fibrant, we have $\hmap
  XY \simeq \map XY$. Moreover, $\hmap --$ preserves weak equivalences and any
  $w\colon X\to X'$ in $\cM$ is a weak equivalence if and only if $\hmap{X'}-\to
  \hmap X-$ is an objectwise weak equivalence. Up to weak equivalence
  in $\SSet$, homotopy function complexes only depend on the class of weak
  equivalences of $\cM$.

  Let $S$ be a class of maps in $\cM$. Any object $Y\in\cM$ is
  \defn{$S$\=/local} if it is fibrant in $\cM$ and if for every $X\to X'$ in
  $S$, $\hmap {X'}Y\to\hmap XY$ is a weak equivalence in $\SSet$. A morphism $X\to Y$
  in $\cM$ is an \defn{$S$\=/local equivalence} if for every $S$\=/local $Z$,
  $\hmap YZ\to \hmap XZ$ is a weak equivalence in $\SSet$. The \defn{left
    Bousfield localisation of $\cM$ at $S$} (should it exist) is defined to be
  the model structure\footnote{A model structure is uniquely determined by two
    of the three classes $\rmW,\rmC,\rmF$.} on $\cM$ whose weak equivalences are
  the $S$\=/local equivalences and whose cofibrations are the cofibrations of
  $\cM$. We denote the model category of the left Bousfield localisation of
  $\cM$ at $S$ by $\cM_S^l$ (or $\cM^l$ if $S$ is clear from context).

  The fundamental theorem of left Bousfield localisations states that if $\cM$
  is combinatorial, left proper and cellular and if $S$ is a set (rather than a
  proper class) of maps of $\cM$, then the left Bousfield localisation $\cM_S^l$
  exists, and is a combinatorial left proper, cellular, simplicial model
  category \cite[Thm 4.1.1]{Hirschhorn2009}.
  
  Let $F:\cM\rightleftarrows\cN:G$ be a Quillen adjunction (respectively, a
  Quillen equivalence). Let $\LL FS$ be the
  class of all maps $\LL FX\to \LL FX'$ in $\cN$, where $X\to X'$ is in $S$ and
  $\LL F$ is the left derived functor of $F$. If the left
  Bousfield localisations of $\cM$ at $S$ and of $\cN$ at $\LL FS$ exist, then
  $F:\cM^l\rightleftarrows\cN^l:G$ is a Quillen adjunction (respectively, a
  Quillen equivalence) \cite[Thm 3.3.20]{Hirschhorn2009}.
\end{para}

\begin{para}[Transferred model structures]
  Let $F:\cM\rightleftarrows\cC:U$ be an adjunction, where
  $(\cM,\rmW,\rmC,\rmF)$ is a model category. We say that a model structure
  $(\cC,\rmW',\rmC',\rmF')$ is the
  \emph{right-transferred model structure on $\cC$ along} $U$ if we have
  $\rmW'=U\inv\rmW$ and $\rmF'=U\inv\rmF$ (namely, $U$ preserves and reflects
  weak equivalences and fibrations). If the right-transferred model structure
  exists, then $F\dashv U$ is a Quillen adjunction. We recall two 
  useful statements that provide the existence of a right-transferred model
  structure and detect if the adjunction is a Quillen equivalence.
\end{para}
\begin{proposition} 
  \label{prop:right-induced-model-structure}
  Let $F : \cM \rightleftarrows \cC : U$ be an adjunction between locally
  presentable categories, where $\cM$ is a combinatorial model category with
  sets $I$ of generating cofibrations and $J$ of generating trivial
  cofibrations. Suppose that:
  \begin{enumerate}
  \item $U$ preserves $\kappa$-filtered colimits (for some regular cardinal
    $\kappa\geq\omega$),
  \item $U$ sends pushouts of maps in $F J$ to weak equivalences.
  \end{enumerate}
  Then the right-transferred model structure on $\cC$ along $U$ exists and is
  cofibrantly generated, with $FI$ and $FJ$ as sets of generating cofibrations
  and trivial cofibrations respectively.
\end{proposition}
\begin{proof}
  We choose a regular cardinal $\lambda\geq\kappa$ such that every generating
  cofibration and trivial cofibration in $\cM$ has a $\lambda$-small domain.
  Then $F$ preserves $\lambda$-small objects and $U$ preserves
  $\lambda$-filtered colimits, and we use {\cite[Prop.
    11.3.2]{Hirschhorn2009}}.
\end{proof}

\begin{proposition}
 \label{prop:quillen-eq-iff-unit-weq} 
  Let $F : \cM \rightleftarrows \cC : U$ be a Quillen adjunction of model
  categories such that $U$ preserves and reflects weak equivalences. Then
  $F\dashv U$ is
  a Quillen equivalence if and only if for every cofibrant $X$ in $\cM$, the unit
  $\eta_X\colon X\to UFX$ is a weak equivalence.
\end{proposition}
\begin{proof}
  See for instance \cite[Lem. 3.3]{erdal2019model}.
\end{proof}

\begin{para}[Simplicial algebras and sheaves]
  \label{para:simpl-algebras-sheaves}
  For any category $\cC$, we write $\s\cC$ for the category $[\bDelta\op,\cC]$ of
  \emph{simplicial objects} in $\cC$. Let $C$ be a small category, and consider
  a Gabriel-Ulmer localisation $F:\psh C\localisation \cC:U$ of the presheaf
  category $\psh C$ at a \emph{set} $S$ of maps in $\psh C$. Then $\cC$ is a
  locally $\kappa$-presentable category, and the fully faithful right adjoint
  $U\colon\cC\subto \psh C$ preserves $\kappa$-filtered colimits (for some large
  enough $\kappa$). Taking simplicial objects, and since colimits in functor
  categories are calculated pointwise, we get a reflective adjunction $\tilde F:
  \sPsh{C} \localisation \s\cC:\tilde U$ where $\tilde U$ preserves
  $\kappa$-filtered colimits. Moreover, $\sPsh C\localisation\s\cC$ is the
  Gabriel-Ulmer localisation generated by the set of maps
  $S_\Delta\eqdef\{\Delta_n\otimes w \mid w\in S\}$ in $\sPsh C$,
  where $S$ is identified with its image under $\psh C\subto\sPsh C$. Then
  $\s\cC$ is locally $\kappa$-presentable as a 
  $\SSet$\=/enriched category, in particular it is enriched, tensored and
  cotensored over $\SSet$. The reflective adjunction is also simplicial, thus
  $\tilde F$ preserves tensors and $\tilde U$ preserves cotensors.

  Let $(C,K)$ be a site, namely a pair of a small category $C$ and a
  Grothendieck topology $K$ on $C$. Then the full subcategory $\Sh_K C\subto
  \psh C$ of sheaves of sets on $C$ is the Gabriel-Ulmer localisation of $\psh
  C$ at the set $S_K\eqdef\{(R\mono c)\in\psh C\arr \mid c\in C, R\in K\}$ of
  subrepresentables that are covering sieves of the topology. From the previous
  paragraph, we have a reflective simplicial adjunction $F:\sPsh
  C\localisation\s\Sh_K C:U$ between simplicial presheaves on $C$ and simplicial
  sheaves on $C$.

  Following \cite[App. A]{duggerHollanderIsaksen2004}, we call the left
  Bousfield localisation of $\sPsh C\inj$ at the set $S_K$ of covering sieves
  the \emph{injective \v{C}ech model structure} on $\sPsh C$ and write it as
  $(\sPsh C\inj^\cech,\rmW^\cech,\rmC\inj,\rmF\inj^\cech)$. The \emph{projective
    \v{C}ech model structure}
  is the left Bousfield localisation of $\sPsh C\proj$ at the same set $S_K$.
  The identity adjunction
  remains a Quillen equivalence.
\end{para}

\begin{remark}
  [{\cite[6.2.2.6, 6.5.4]{LurieHT}}]
 \label{rem:cech-simpl-sheaves-pres-topological-loc} 
  The \v Cech model structures on $\sPsh C$ present the
  \oo-topos of sheaves of \oo-groupoids on $(C,J)$, namely the topological
  left-exact localisation of the \oo-category $\Psh C$ (of presheaves of
  \oo-groupoids) generated by the set $S_K$ (the $(-1)$\=/truncated maps that are
  the covering sieves of the topology).
\end{remark}

\begin{proposition}
  \label{prop:right-transf-cech-model-struc-Q-equivalence}
  Let $(C,K)$ be a site. Then the right-transferred model structure along
  $U\colon \s\Sh_KC\subto\sPsh C\inj^\cech$ exists and the adjunction $F\dashv
  U$ is a Quillen equivalence.
\end{proposition}
\begin{proof}
  By \cite[Prop. A.2]{duggerHollanderIsaksen2004}, for every $X\in\sPsh C$, the
  unit $\eta_X\colon X\to UFX$ is a weak equivalence in $\sPsh
  C\inj^\cech$. Hence, by~\cref{prop:quillen-eq-iff-unit-weq}, it
  suffices to show that the right-transferred model structure exists. We will
  use \cref{prop:right-induced-model-structure}. Let
  $j\colon X\mono Y$ be a generating trivial cofibration in $\sPsh C\inj^\cech$
  and let $f\colon FX\to A$ be a map in $\s\Sh_KC$. Consider the diagram below
  in $\sPsh C$.
  \[
    \begin{tikzcd}
      X\ar[r,"\sim"',"\eta"] \ar[d,tail,"\sim","j"']
      &UFX \ar[r,"Uf"] \ar[d,"UFj"']
      &UA\ar[d] \ar[r,equals]
      &UA\ar[d]\\
      Y\ar[r,"\sim","\eta"']
      &UFY \ar[r]
      &B\pomark \ar[r,"\sim","\eta"']
      &UFB
    \end{tikzcd}
  \]
  Since sheafification is left-exact, the monad $UF$ preserves monomorphisms,
  namely maps in $\rmC\inj$. Hence $UFj$ is in $\rmC\inj$; by
  ``2~out~of~3'', it is a trivial cofibration in $\sPsh C\inj^\cech$, and so the
  pushout $UA\to B$ is as well. Thus the composite $UA\to UFB$ is a weak
  equivalence, and it is the image by $U$ of the pushout of $Fj$ along $f$.
\end{proof}

\begin{remark}
 \label{rem:fully-faithful-functor-sheaves-hypercompletion} 
  Let $i\colon{A}\subto B$ be a fully faithful functor between small categories.
  Then $\psh A\simeq \Sh_{K_A}B$ is the category of sheaves of sets on $B$ for
  the topology generated by the cocones in $B$ under ${A}$. By
  \cref{prop:right-transf-cech-model-struc-Q-equivalence}, the right-transferred
  model structure along $i_*\colon\sPsh A\simeq\s\Sh_{K_A}B\subto\sPsh
  B\inj^\cech$ exists and $i^*\dashv i_*$ is a Quillen equivalence. We will
  write $\sPsh A_\cech$ for this model category.

  The right-transferred model structure $\sPsh A_\cech$ is \emph{not} in general
  a global model structure on $\sPsh A$. This is because $\sPsh A_\cech$
  presents a \emph{topological} left-exact localisation of the presheaf
  \oo-topos $\Psh B$, and such a localisation need not in general give a
  hypercomplete \oo-topos. On the other hand, a global model structure on $\sPsh
  A$ presents the presheaf \oo-topos $\Psh{A}$, which is always
  hypercomplete---indeed, it is the hypercompletion of the \oo-topos presented
  by $\sPsh A_\cech$.\footnote{In summary, the equivalence $\Sh_{K_A}B\simeq
    \psh A$ of $1$-categories does not translate to an equivalence, but a
    cotopological localisation $\infty\Sh_{K_A}B\localisation\Psh A$ of
    \oo-topoi.}
\end{remark}


\section{\texorpdfstring{$\lfdex$}{C}-sorted spaces}
\label{sec:C-sorted-spaces}
\begin{center}
We fix a \lfd~category $\lfdex$ for the rest of this chapter.
\end{center}
Recall from \cref{rem:functor-fin-to-cell} that we have a canonical inclusion
$i:\lfdex\subto \cell\lfdex$. Thus we can see the category $\psh\lfdex$ of
presheaves of \emph{sets} as a category of sheaves of sets on $\cell\lfdex$, via
the left-exact localisation $i^*:\psh{\cell\lfdex}\localisation\psh\lfdex:i_*$.
Due to \cref{rem:fully-faithful-functor-sheaves-hypercompletion}, we may
legitimately question whether this continues to be true for presheaves of
\emph{spaces} on $\lfdex$---that is, whether the \oo-topos of sheaves of
\oo-groupoids on $\cell\lfdex$ models the hypercomplete presheaf \oo-topos
$\Psh\lfdex$. This turns out to be the case, as we will show in this section.

\begin{definition}
    \label{def:C-sorted-spaces}
    The simplicial model category of \defn{$\lfdex$-sorted spaces} is the
    simplicial presheaf category $\sPsh\lfdex$ equipped with the
    \emph{injective} global model structure.
\end{definition}

\begin{remark}
  When $\lfdex$ is a set~(\cref{exa:lfd-cats}\ref{exa:lfd-cats-set}), the
  projective and injective model structures on $\sPsh\lfdex$ coincide.
\end{remark}

\begin{para}
  Since $\lfdex$ is direct, $\lfdex\op$ is an \emph{inverse}
  category,\footnote{To be an inverse category can be defined---indirectly---as
    being the opposite of a direct category.} and as such, has the structure of
  a \emph{Reedy} category (we refer the reader to \cite{RiehlVerity2014reedy}
  for a survey). The general theory of Reedy model structures then ensures that
  the injective global model structure on $\sPsh\lfdex$ coincides with the
  \emph{Reedy} global model structure.\footnote{This phenomenon is more
    generally explained by the fact that $\lfdex$, as a direct category, is an
    \emph{elegant Reedy category} \cite[Sec. 3]{bergner2013reedy}.} The payoff
  is that the injective fibrations and trivial fibrations in $\sPsh\lfdex$ are
  easy to describe, and we have the simple generating sets of cofibrations and
  trivial cofibrations below (where $\delta_c\colon\partial c\subto c$ are the
  boundary inclusions of \cref{para:boudaries-finite-cell-cplxs}).
  \begin{align*}
    I_\lfdex&\eqdef\{\delta_n\pp \delta_c\mid c\in \lfdex,[n]\in\bDelta\}\\
    J_\lfdex&\eqdef\{h^k_n\pp \delta_c\mid c\in \lfdex,[n]\in\bDelta,0\leq k\leq n\}
  \end{align*}
  By a little pushout-product yoga, we see that any $X\in\sPsh\lfdex$ is
  injective (trivial) fibrant if and only if for every $c$ in $\lfdex$, the map
  $\pbh{\delta_c}{X\to 1}\colon X_c\to M_cX$ from $X_c$ to its \emph{matching
    object} $M_cX\eqdef \map{\partial c}X$ is a (trivial) fibration in $\SSet$.
  Likewise, a map $f\colon X\to Y$ is an injective (trivial) fibration if and
  only if the map $\pbh{\delta_c}f\colon X_c\to Y_c\times_{M_cY} M_cX$ is a
  (trivial) fibration in $\SSet$ for every $c$ in $\lfdex$.
\end{para}

\begin{remark}
  [Inverse Reedy categories and type theory] Categories of inverse diagrams
  valued in a homotopy-theoretic model of univalent type theory, seen as Reedy
  diagram categories, were one of the first \cite{shulman2015inverse} examples
  extending Voevodsky's model of univalent type theory in spaces to diagram
  categories intended to model diagrams in \oo-topoi. The Reedy fibrancy
  condition exactly models what it means for a context/object to have a
  canonical projection to its parent/matching object, reflected in how smoothly
  the construction ([\textsection{11}, \emph{op. cit.}]) goes through.
\end{remark}

\begin{para}
  We now show that we are able to give an easy proof of our claim, namely that
  $\Psh\lfdex\simeq\infty\Sh(\cell\lfdex)$. This is partly thanks to the Reedy
  model structure, but also because
  the presheaves $\partial c\in\psh\lfdex$ are in the image of the nerve functor
  $\cell\lfdex\subto\psh\lfdex$ associated to $i\colon\lfdex\subto\cell\lfdex$.
\end{para}

\begin{proposition}
  \label{prop:C-spaces-model-struct-equals-Cech}
  Let $i\colon\lfdex\subto\cell\lfdex$ denote the canonical inclusion. The
  right-transferred model structure $\sPsh\lfdex_\cech$ along
  $i_*\colon\sPsh\lfdex\subto \sPsh(\cell\lfdex)\inj^\cech$ (from
  \cref{rem:fully-faithful-functor-sheaves-hypercompletion}) is the same as the
  model structure for $\lfdex$-sorted spaces (the injective global model
  structure $\sPsh\lfdex\inj$).
\end{proposition}
\begin{proof}
  Any model structure is determined by its cofibrations and fibrant objects, or
  equivalently, by its trivial fibrations and fibrant objects
  \cite[E.1.10]{joyal2008theory}. For every $X\in\sPsh\lfdex$, $i_*X$ is a
  simplicial sheaf, thus it is fibrant in
  $\sPsh(\cell\lfdex)\inj^\cech$---namely injective \v Cech fibrant---if and
  only if it is injective fibrant. Next, $i_*\colon
  \sPsh\lfdex\inj\to\sPsh(\cell\lfdex)\inj$ is right Quillen, so it preserves
  injective trivial fibrations, and as a left Bousfield localisation,
  $\sPsh(\cell\lfdex)\inj^\cech$ and $\sPsh(\cell\lfdex)\inj$ have the same
  trivial fibrations. So it suffices to show that for every object $X$ and map
  $f$ in $\sPsh \lfdex$,
  \begin{enumerate}
  \item if $i_*X$ is injective fibrant, then $X$ is Reedy fibrant,
  \item and if $i_*f$ is an injective trivial fibration, then $f$ is a Reedy
    trivial fibration.
  \end{enumerate}
  But each lifting property follows by adjointness, since for every $c$ in
  $\lfdex$, the maps $\delta_c\colon \partial c\subto c$ are in the image of
  $\cell\lfdex\subto\psh\lfdex$, and so the maps $\delta_n\pp\delta_c$ in $\sPsh
  \lfdex$ are images of monomorphisms under
  $i^*\colon\sPsh(\cell\lfdex)\to\sPsh\lfdex$.\qedhere
\end{proof}
\begin{corollary}
    \label{cor:Cech-model-struc-on-cellC-spaces-models-C-spaces}
    The projective and injective \v Cech model structures on
    $\sPsh(\cell\lfdex)$, associated to the inclusion $\lfdex\subto
    \cell\lfdex$, both present the hypercomplete \oo-topos $\Psh\lfdex$. (In the
    language of \cite{duggerHollanderIsaksen2004}, any $X$ in
    $\sPsh(\cell\lfdex)$ that has \v Cech descent, has descent for hypercovers.)
\end{corollary}

\begin{remark}
    The previous result seems related to the fact that, by general
    considerations in \cite{Anel2021enveloping}, any direct category $C$ is
    ``good'', namely the \emph{enveloping \oo-topos} of the $1$-topos $\psh C$
    is the presheaf \oo-topos $\Psh C$.
\end{remark}

\begin{remark}
    In the case of discrete spaces (sets), we have that
    $\psh\lfdex\simeq\Ind{\cell\lfdex}$ is the ind-completion of the finite
    cocompletion $\cell\lfdex\simeq\fin_\lfdex$ of $\lfdex$ (for colimits of
    sets). In the case of spaces, $\cell\lfdex$ is of course not a finite
    cocompletion (for colimits of spaces).
\end{remark}


\newcommand{\theoryex}{\bT}
\section{Homotopical models of algebraic theories and rigidification}
\label{sec:htpy-models-alg-theories}
In this section, we recall some elements of the theory of models in spaces of
ordinary multisorted algebraic theories from
\cite{Schwede2001stablehtpyAlgTheories,badzioch2002algtheories,bergner2006rigidification,Rezk2002simplicialAlgTheories}
that we would like to generalise to dependently sorted algebraic theories.

We fix a set $I$ of sorts throughout this section.

\begin{para}
    [$I$\=/sorted spaces]
   \label{para:I-sorted-spaces} 
    Let $i\colon I\subto \cell I$ be the fully faithful
    inclusion. Remark that the global projective and injective model structures
    on the category $\sPsh I\cong\SSet^I$ coincide, namely, they are both the
    model structure on $\sPsh I$ for \emph{$I$\=/sorted spaces}
    (\cref{def:C-sorted-spaces}). Then this implies that
    \begin{enumerate}
    \item $i_!:\sPsh I\rightleftarrows\sPsh(\cell I)\proj:i^*$ is a Quillen
        adjunction,
    \item and $i^*:\sPsh(\cell I)\proj\rightleftarrows\sPsh I:i_*$ is a Quillen
        adjunction.
    \end{enumerate}
\end{para}

\begin{para}
    [Simplicial $\theoryex$\=/algebras]
    \label{para:simpl-T-algs}
    An $I$\=/sorted algebraic theory is an identity-on-objects, finite-product
    preserving functor $\cxl I\to \theoryex$ (equivalently, $\theoryex$ is an
    $I$\=/contextual category). The category $\s\theoryex\Mod$ of \emph{strict}
    models of $\theoryex$ in $\SSet$ is the category of
    finite-product-preserving functors $\theoryex\to \SSet$. Equivalently, it is
    the category of simplicial objects in the category $\theoryex\Mod$ of
    $\Set$\=/models of $\theoryex$.

    A technique due to Quillen \cite[II.4]{quillen1967} proves that
    the transferred model structure along the monadic functor $U\colon
    \s\theoryex\Mod\to \sPsh I$ exists. We call $\s\theoryex\Mod$ with the
    transferred model structure the model category of \defn{simplicial
      $\theoryex$\=/algebras}.

    Recall the simplicial adjunction
    $N:\sPsh(\theoryex\op)\localisation \s\theoryex\Mod:h$
    from~\cref{para:simpl-algebras-sheaves}. Then the model structure for
    simplicial $\theoryex$\=/algebras is also the transferred model structure
    along the right adjoint $N$ \cite[Thm 5.1]{badzioch2002algtheories}. The
    identity-on-objects functor $j\colon\cxl I\to\theoryex$ and the theory of
    monads with arities gives a diagram of adjunctions below such that $j_!i_*\cong
    Nj_!$ is an exact adjoint square.
    \begin{equation}
        \tag{$\star$}
        \label{eq:htpy-T-alg-exact-adj-square}
        \begin{tikzcd}
            \sPsh I\ar[r,shift left,"j_!"] \ar[d,shift left,hook',"i_*"]
            &\s\theoryex\Mod \ar[l,shift left,"j^*"]\ar[d,shift left,hook',"N"]\\
            \sPsh(\cell I)\ar[r,shift left,"j_!"] \ar[u,shift left,"i^*"]
            &\sPsh(\theoryex\op) \ar[l,shift left,"j^*"] \ar[u,shift left,"h"]
        \end{tikzcd}
    \end{equation}
    In the previous diagram, $j_!:\sPsh I\rightleftarrows\s\theoryex\Mod:j^*$ is
    both the restriction of $j_!\dashv j^*$ along the nerve functors $i_*$ and
    $N$ as well as the monadic adjunction of the monad $\theoryex$ on $\sPsh I$.
    It is readily verified that all the adjunctions are Quillen for the
    projective global model structures on $\sPsh(\cell I)$ and
    $\sPsh(\theoryex\op)$. Moreover, since $j$ is identity-on-objects,
    $j^*\colon \sPsh(\theoryex\op)\to \sPsh(\cell I)$ preserves and reflects
    projective fibrations and global weak equivalences. Thus every model
    structure on the right in each of the Quillen pairs in
    \eqref{eq:htpy-T-alg-exact-adj-square} is the transferred model structure
    along the right adjoint.
\end{para}

\begin{para}
    [Homotopy $\theoryex$\=/algebras]
    \label{para:htpy-T-algs}
    A simplicial $\theoryex$\=/algebra is a finite-product preserving functor
    $F\colon \theoryex\to\SSet$. In other words, for every
    $\Gamma=i_1\times\ldots\times i_k$ in $\theoryex$, we have $F\Gamma\cong
    Fi_1\times\ldots\times Fi_k$. An obvious way to define models of $\theoryex$
    in the \oo\=/category of spaces is as functors $F\colon \theoryex\to\SSet$
    that take finite products to \emph{homotopy limits}. Since finite products
    in simplicial sets are already homotopy limits (weak equivalences of $\SSet$ are
    stable under product), we can define a \emph{homotopy $\theoryex$\=/algebra}
    to be a functor $\theoryex\to\SSet$ such that for every
    $\Gamma\in\theoryex$, the canonical map $F\Gamma\to Fi_1\times\ldots\times
    Fi_k$ is a weak equivalence in $\SSet$.

    The maps $F\Gamma\to Fi_1\times \ldots\times Fi_k$, natural in $F$, are
    co-represented by maps $\mquote\Gamma\to\Gamma$ in $\Set^\theoryex$, where
    $\Gamma$ is the representable presheaf on $\Gamma\in\theoryex\op$, and
    $\mquote\Gamma$ is the corresponding coproduct of $i_1,\ldots,i_k$
    calculated in the presheaf category $\Set^\theoryex$. Let $S\eqdef
    \{\mquote\Gamma\to\Gamma\mid\Gamma\in\theoryex\}$. We write
    $\sPsh(\theoryex\op)\proj^l$ for the left Bousfield localisation of the
    projective model structure at the set of maps $S$, and call this the
    \emph{model category of homotopy $\theoryex$\=/algebras}.
\end{para}
\begin{para}
    [Rigidification] For every $\Gamma\in\theoryex$, $\mquote\Gamma$ is a
    cofibrant object of $\sPsh(\theoryex\op)\proj$ (since it is a coproduct of
    representables). Since $h\colon\sPsh(\theoryex\op)\to\s\theoryex\Mod$ sends
    every morphism in $S$ to an isomorphism, the adjunction
    $h:\sPsh(\theoryex\op)\proj^l\rightleftarrows\s\theoryex\Mod:N$ is thus
    Quillen (see \cref{para:left-bousfield-localisation}).

    The following result is due to Badzioch.
\end{para}

\begin{theorem}
    [Rigidification of homotopy $\theoryex$\=/algebras]
    \label{thm:badzioch-main}
    The Quillen adjunction
    $h:\sPsh(\theoryex\op)\proj^l\rightleftarrows\s\theoryex\Mod:N$ is a Quillen
    equivalence.
\end{theorem}
\begin{proof}
    The single-sorted case is \cite[Thm 6.4]{badzioch2002algtheories}, and the
    multisorted case (treated in \cite{bergner2006rigidification}) is exactly
    similar.
\end{proof}

\begin{remark}
   \label{rem:bousfield-exact-adj-square-I-theory} 
    By Bousfield localisation and \cref{thm:badzioch-main}, the exact
    adjoint square~\eqref{eq:htpy-T-alg-exact-adj-square} gives an exact adjoint
    square
    \[
        \begin{tikzcd}
            \sPsh I\ar[r,shift left,"j_!"] \ar[d,shift left,hook',"i_*"]
            &\s\theoryex\Mod \ar[l,shift left,"j^*"]\ar[d,shift left,hook',"N"]\\
            \sPsh(\cell I)\proj^l\ar[r,shift left,"j_!"] \ar[u,shift left,"i^*"]
            &\sPsh(\theoryex\op)\proj^l \ar[l,shift left,"j^*"] \ar[u,shift left,"h"]
        \end{tikzcd}
    \]
    where the vertical adjunctions are Quillen equivalences (the left vertical
    adjunction is the particular case of \cref{thm:badzioch-main} for the
    initial $I$\=/sorted theory).
\end{remark}

\section{A \emph{flasque} model structure for homotopy
  \texorpdfstring{$\lfdex$}{C}-sorted spaces}
\label{sec:flasq-model-structure}
In \cref{sec:C-sorted-spaces}, we saw that $\lfdex$\=/sorted spaces can be seen
as sheaves of spaces over $\cell\lfdex$, generalising the equivalence
$\psh\lfdex\simeq\Sh(\cell\lfdex)$. However, $\psh\lfdex$ is also the category
of models in $\Set$ of the initial $\lfdex$\=/sorted
theory/$\lfdex$\=/contextual category $\cxl\lfdex=\cell\lfdex\op$. From this
point of view, $\psh\lfdex$ is the category of functors $\cxl\lfdex\to\Set$ that
preserve canonical pullbacks (\cref{def:model-of-cxl-cat}), or in other words,
that take canonical pushouts of boundary inclusions $\partial c\subto c$ in
$\cell\lfdex$ to pullbacks in $\Set$. Now the sheaf condition is also a
limit-preservation condition, but the two are not the same. They coincide for
the initial $\lfdex$-contextual category, but the $\Set$-models of an arbitrary
$\lfdex$-contextual category $\cxl\lfdex\to\sfD$, despite being characterised by
preserving canonical pullbacks, are not a category of sheaves on $\sfD$.

We can therefore ask the question: what is an $\SSet$\=/model \emph{up to
  homotopy} of a $\lfdex$\=/sorted theory? In fact, since we have a candidate
(namely presheaves of spaces on $\lfdex$) for the models up to homotopy of the
initial $\lfdex$\=/sorted theory, we could start by restricting to this
degenerate case. A first answer could be to define such a homotopy model to be a
simplicial presheaf $F\colon \cell\lfdex\op\to\SSet$ such that for every
canonical pushout of a boundary inclusion in $\cell\lfdex$ as below,
\[
  \begin{tikzcd}[sep=scriptsize]
    \partial c \ar[r,"\gamma"] \ar[d,hook,"\delta_c"']
    &X_n \ar[d,"p_{n+1}"] \\
    c \ar[r,"\gamma.c"] &X_{n+1} \pomark
  \end{tikzcd}
\]
the map $F(X_{n+1})\to F(X_n)\times_{F(\partial c)}Fc$ is a weak equivalence.
But this does not define a model of $\cxl\lfdex$ in the \oo\=/category of
spaces, because such a model should take canonical pushouts to \emph{homotopy}
pullbacks (pullbacks in the \oo-category of spaces), and pullbacks in $\SSet$
are not necessarily homotopy pullbacks.

When $\lfdex$ is a set, we are precisely in the case of multisorted algebraic
theories, all the canonical pullbacks involved are \emph{products}, and products
of arbitrary simplicial sets \emph{are} homotopy limits. This is what motivates
the considerations in
\cite{Rezk2002simplicialAlgTheories,badzioch2002algtheories,bergner2006rigidification}
for simplicial sets, and in \cite{Schwede2001stablehtpyAlgTheories} for pointed
simplicial sets.\footnote{There are model categories that model the \oo-topos of
  spaces where weak equivalences are not stable under products (\cite[Sec.
  5]{maltsiniotis2009cubique} contains a very nice combinatorial example), where
  this approach would not work.}

Thus in the general case, the natural condition that emerges for the simplicial
presheaf $F$ to satisfy is for $F(X_{n+1})$ to be a homotopy pullback
$F(X_n)\hopb_{F(\partial c)}Fc$. In this section, we will first set up some
machinery to define homotopical models of $\cxl\lfdex$ as fibrant objects of a
model category, and then substantiate the previous point of view via a
rigidification theorem with respect to $\lfdex$\=/sorted spaces.
\begin{definition}
    A \defn{homotopical $\lfdex$\=/sorted space} is an $F$ in
    $\sPsh(\cell\lfdex)$ such that $F\emptyset\to 1$ is a weak equivalence in
    $\SSet$, and for every canonical pushout in $\cell\lfdex$ of the form
    \[
    \begin{tikzcd}[sep=scriptsize]
      \partial c \ar[r,"\gamma"] \ar[d,hook,"\delta_c"']
      &\Gamma_n \ar[d,"p_{n+1}"] \\
      c \ar[r,"\gamma.c"] &\Gamma_{n+1} \pomark
    \end{tikzcd}
  \] we have $F_{\Gamma_{n+1}}\simeq
  F_{\Gamma_n}\hopb_{F_{\partial c}}F_c$ in $\SSet$. A \defn{homotopical
    model} of a $\lfdex$\=/contextual category $\cxl\lfdex\to\sfD$ is a
  simplicial presheaf $\sfD\to\SSet$ such that its restriction to
  $\sPsh(\cell\lfdex)$ is a homotopical $\lfdex$\=/sorted space.
\end{definition}

\begin{para}
  [Flasque boundaries]
  \label{para:global-flasque-boundaries}
  Since we have a full inclusion $i\colon\lfdex\subto\cell\lfdex$, every $c$ in
  $\lfdex$ can be seen as an object of $\cell\lfdex$. We can therefore calculate
  the colimit of $\lfdex\slice
  c^-\to\lfdex\subto\cell\lfdex\subto\psh{\cell\lfdex}$ as the ``boundary''
  inclusion $\mquote{\delta_c}\colon\mquote{\partial c}\subto ic$ of the
  representable presheaf $ic$ (we will not explicitly write the inclusion $i$
  from now on, and write $c$ instead of $ic$, since the presence of quotes
  $\mquote{}$ necessarily means that we are in $\psh{\cell\lfdex}$ or
  $\sPsh(\cell\lfdex)$). Now, the presheaf $\partial c\in\psh\lfdex$ is
  representable in $\psh{\cell\lfdex}$, but $\partial c$ is \emph{not} the same
  as $\mquote{\partial c}$ in $\psh{\cell\lfdex}$---indeed, we have an inclusion
  $\mquote{\partial c}\subto\partial c$ that is just the counit on $\partial c$
  of the idempotent comonad $i_!i^*$ on $\psh{\cell\lfdex}$.

  Using the subrepresentables $\mquote{\delta_c}\colon\mquote{\partial c}\subto
  c$, we will define an intermediate global model structure, called the
  \emph{$\partial$\=/flasque model structure}\footnote{We use the word
    \emph{flasque} in analogy with ``flabby sheaves'' (as in
    \cite{isaksen2004flasque}).}, on $\sPsh(\cell\lfdex)$. We will set up our
  cofibrant generation in precisely the same way as in
  \cite{isaksen2004flasque}, except that in our case, the crux of the
  construction (\cite[Lem. 3.9]{isaksen2004flasque}) admits a much simpler
  proof---it follows automatically from the injective Reedy structure on
  $\sPsh\lfdex$. We begin by fixing sets of putative generating cofibrations and
  trivial cofibrations.
  \begin{align*}
    I_\partial&\eqdef I\proj\cup \{\delta_n\pp \mquote{\delta_c}\mid
                c\in \lfdex,[n]\in\bDelta\}\\
    J_\partial&\eqdef J\proj\cup \{h^k_n\pp \mquote{\delta_c}\mid
                c\in \lfdex,[n]\in\bDelta,0\leq k\leq n\}
  \end{align*}
\end{para}

\begin{remark}
  When $\lfdex$ is a set, $I_\partial=I\proj$ and $J_\partial=J\proj$.
\end{remark}

\begin{definition}[{\cite[Defs 3.3, 3.6]{isaksen2004flasque}}]
  A map $f$ in $\sPsh(\cell\lfdex)$ is
  \begin{enumerate}
  \item a \emph{$\partial$\=/flasque fibration} if it is a
    $J_\partial$-injective\footnote{Has the right lifting property against every
      map in $J_\partial$.},
  \item a \emph{$\partial$\=/flasque cofibration} if it has the left lifting
    property against all objectwise acyclic $\partial$\=/flasque fibrations,
  \end{enumerate}
  and in either case, is \emph{objectwise acyclic} if it is also a global weak
  equivalence.
\end{definition}

\begin{lemma}[{\cite[Lem. 3.8]{isaksen2004flasque}}]
 \label{lem:flasque-struc-is-intermediate} 
  \begin{enumerate}
  \item A projective cofibration is a $\partial$\=/flasque cofibration, and a
    $\partial$\=/flasque fibration is a projective fibration.
  \item An injective fibration is a $\partial$\=/flasque fibration, and a
    $\partial$\=/flasque cofibration is an injective cofibration (a monomorphism).
  \end{enumerate}
\end{lemma}
\begin{proof}
  \begin{enumerate}
  \item Since $J\proj\subset J_\partial$, it follows that (objectwise acyclic)
    $\partial$\=/flasque fibrations are projective (trivial) fibrations. So
    projective cofibrations are $\partial$\=/flasque cofibrations.
  \item By pushout-product yoga, every map in $J_\partial$ is an injective
    trivial cofibration. So every injective (trivial) fibration is an (objectwise
    acyclic) $\partial$\=/flasque fibration, so every $\partial$\=/flasque cofibration is
    an injective cofibration.\qedhere
  \end{enumerate}
\end{proof}

\begin{lemma}[{\cite[Lem. 3.11]{isaksen2004flasque}}]
 \label{lem:J-cof-is-obj-acyclic} 
  If a map $f$ is a $J_\partial$\=/cofibration\footnote{A retract of a relative
    $J_\partial$\=/cell-complex.}, then it is an objectwise acyclic
  $\partial$\=/flasque cofibration.
\end{lemma}
\begin{proof}
    Since $f$ has the left lifting property against all $\partial$\=/flasque fibrations, it
    is a $\partial$\=/flasque cofibration. Since every map in $J_\partial$ is an injective
    trivial cofibration, so is every $J_\partial$-cofibration, thus $f$ is an
    objectwise weak equivalence.
\end{proof}

\begin{lemma}
 \label{lem:flasque-fibration-inv-image-reedy} 
    If $f\colon X\to Y$ is a $\partial$\=/flasque fibration, then $i^*f$ is a fibration 
    in $\sPsh\lfdex\inj$ (a Reedy fibration).
\end{lemma}
\begin{proof}
  By adjointness, $\pbh{\mquote{\delta_c}}f\cong \pbh{\delta_c}{i^*f}$.
\end{proof}

\begin{lemma}
  [{\cite[Lem. 3.9]{isaksen2004flasque}}]
 \label{lem:obj-acycl-flasque-fib-iff} 
  A map $f\colon X\to Y$ in $\sPsh(\cell\lfdex)$ is an objectwise acyclic
  $\partial$\=/flasque fibration if and only if it is an
  $I_\partial$\=/injective ($f$is a projective trivial fibration and for every
  $c\in\lfdex$ the map $\pbh{\mquote{\delta_c}}f$ is a trivial fibration in
  $\SSet$).
\end{lemma}
\begin{proof}
  The ``if'' direction is clear. For the ``only if'', let $f$ be an objectwise
  acyclic $\partial$\=/flasque fibration. By
  \cref{lem:flasque-struc-is-intermediate}, it is a projective trivial
  fibration. So we need to show that for every $c\in\lfdex$,
  $\pbh{\mquote{\delta_c}}f$ is a trivial fibration. But $i^*$ preserves global
  weak equivalences and by the previous lemma, $i^*f$ is a Reedy fibration. So
  $i^*f$ is a Reedy trivial fibration and by adjointness we are done.
\end{proof}

\begin{theorem}
  \label{thm:flasque-model-struc-existence}
  There is a cofibrantly generated, intermediate, global model structure on
  $\sPsh(\cell\lfdex)$ whose:
  \begin{enumerate}
  \item weak equivalences are the global weak equivalences,
  \item cofibrations are the $\partial$\=/flasque cofibrations,
  \item and whose fibrations are the $\partial$\=/flasque fibrations.
  \end{enumerate}
\end{theorem}
\begin{proof}
  We use the standard Kan recognition theorem \cite[Thm 11.3.1]{Hirschhorn2009},
  using \cref{lem:J-cof-is-obj-acyclic,lem:obj-acycl-flasque-fib-iff} and since
  $I_\partial$ and $J_\partial$ clearly permit the small object argument. Then
  by \cref{lem:flasque-struc-is-intermediate}, the model structure is
  intermediate.
\end{proof}

\begin{definition}
  The \defn{$\partial$\=/ flasque model structure} on $\sPsh(\cell\lfdex)$is the
  global model structure defined in~\cref{thm:flasque-model-struc-existence}.
  We write it as $(\sPsh(\cell\lfdex)_\partial, \rmW,
  \rmC_\partial,\rmF_\partial)$.
\end{definition}

\begin{proposition}
  $\sPsh(\cell\lfdex)_\partial$ is a proper simplicial model category.
\end{proposition}
\begin{proof}
  Any intermediate global model structure is proper. If $f$ is a cofibration in
  $\SSet$ and $h$ a $\partial$\=/flasque fibration, then we have an isomorphism
  $\pbh{\mquote{\delta_c}}{\pbh fh}\cong \pbh f {\pbh{\mquote{\delta_c}}h}$ in
  $\SSet\arr$, and since $\SSet$ is a simplicial model category, the map $\pbh f
  {\pbh{\mquote{\delta_c}}h}$ is a fibration which is trivial if either $f$ or
  $h$ is so.
\end{proof}

We finally come to the first reason for introducing this model structure, namely
that we recover a strict generalisation of the properties from
\cref{para:I-sorted-spaces}.

\begin{theorem}
 \label{thm:flasque-model-struc-generalises-projective}
  Let $i\colon \lfdex\subto \cell\lfdex$ denote the canonical inclusion. Then,
  \begin{enumerate}
  \item
    \label{thm:flasque-model-struc-generalises-projective:1}
    $i_!:\sPsh\lfdex\inj\rightleftarrows\sPsh(\cell\lfdex)_\partial:i^*$ is
    a Quillen adjunction,
  \item
    \label{thm:flasque-model-struc-generalises-projective:2}
    and $i^*:\sPsh(\cell\lfdex)_\partial\rightleftarrows\sPsh\lfdex\inj:i_*$
    is a Quillen adjunction.
  \end{enumerate}
\end{theorem}
\begin{proof}
  Since $i^*$ preserves global weak equivalences and monomorphisms, and sends
  $\partial$\=/flasque fibrations to fibrations in $\sPsh\lfdex\inj$
  (\cref{lem:flasque-fibration-inv-image-reedy}), it is both left and right
  Quillen.
\end{proof}

  

\begin{proposition}
 \label{lem:reedy-fibrant-is-htpy-coherent-C-space} 
  Let $X$ be fibrant in $\sPsh\lfdex\inj$. Then $i_*X$ is a homotopical
  $\lfdex$\=/sorted space.
\end{proposition}
\begin{proof}
  It is easy to see that $i_*X\colon\cell\lfdex\op\to\SSet$ takes canonical
  pushouts of every $\partial c\subto c$ to strict pullbacks. Since $X$ is Reedy
  fibrant, $X_c\to (i_*X)_{\partial c}$ is a fibration, so these pullbacks are
  homotopy pullbacks.
\end{proof}

\begin{remark}
  What the proof of \cref{lem:reedy-fibrant-is-htpy-coherent-C-space} really
  says is that Reedy fibrant objects in $\sPsh\lfdex$, \emph{qua} homotopical
  models in $\SSet$ of the initial $\lfdex$-contextual category are eminently
  \emph{type-theoretic} in nature, since not only do they take canonical
  pullbacks to homotopy pullbacks (as functors $\cxl\lfdex\to\SSet$), but they
  also take dependent projections to Kan fibrations.
\end{remark}

\begin{para}
   \label{para:flasque-local-eqs} 
  Recall that every object of $\cell\lfdex$ is a finite $I_\lfdex$\=/cell
  complex $\emptyset\to\ldots \Gamma_n$. For each such cell complex (that we
  write $\Gamma_n$ with slight abuse), we will define a subrepresentable
  presheaf $\mquote{\Gamma_n}\subto\Gamma_n$ in $\psh{\cell\lfdex}$. Proceeding
  by induction on the dimension and length of $\Gamma_n$ (as in the proof of
  \cref{prop:rex-iff-cell-extendable}), in the base case, the subrepresentable
  presheaf of the initial object $\emptyset\in\cell\lfdex$ is defined to be the
  empty presheaf (the initial object of $\psh{\cell\lfdex}$)
  $\mquote\emptyset\subto \emptyset$. Let $\Gamma_{n+1}$ be an object in
  $\cell\lfdex$ of the form given by the following left-hand square in
  $\cell\lfdex$. We define the subrepresentable
  $\mquote{\Gamma_{n+1}}\subto\Gamma_{n+1}$ by the cocartesian back face in the
  following right-hand cube in $\psh{\cell\lfdex}$ (the front face is only
  cartesian but not cocartesian since $\cell\lfdex\subto \psh{\cell\lfdex}$
  preserves limits but not colimits).
  \[
    \begin{tikzcd}[sep=scriptsize]
      \partial c\pbmark \ar[r] \ar[d,hook,"\delta_c"']
      &\Gamma_n \ar[d,hook] \\
      c \ar[r] &\Gamma_{n+1} \pomark
    \end{tikzcd}
    \qquad\qquad
    \begin{tikzcd}[row sep=tiny,column sep=small]
      \mquote{\partial c} \ar[rr]\ar[rd,hook] \ar[dd,hook]
      && \mquote{\Gamma_n} \ar[rd,hook] \ar[dd,hook]\\
      &\partial c
      && \Gamma_n \ar[from=ll, crossing over]\ar[dd,hook]\\
      c \ar[rd,equals] \ar[rr]
      && \mquote{\Gamma_{n+1}}\pomark \ar[rd,hook]\\
      &c \ar[from = uu,hook, crossing over] \ar[rr]
      && \Gamma_{n+1}
    \end{tikzcd}
  \]
\end{para}

\begin{remark}
 \label{rem:desc-flasque-spine-counit} 
 For any $c\in\lfdex$, we have $\mquote c=c$, and the previous definition of
 $\mquote{\partial c}$ coincides with the one in
 \cref{para:global-flasque-boundaries}. In fact, we could just as easily have
 defined $\mquote\Gamma\subto\Gamma$ as simply being the counit
 $i_!i^*\Gamma\subto\Gamma$, but inductively defining them as above
 ``automatises'' many of the following proofs.
\end{remark}

\begin{proposition}
  \label{prop:flasque-objs-are-cofibrant}
  For every $\Gamma\in\cell\lfdex$, \emph{$\mquote\Gamma$} is cofibrant in
  $\sPsh(\cell\lfdex)_\partial$.
\end{proposition}
\begin{proof}
  They are all $I_\partial$\=/cell complexes, since $\mquote\emptyset$ is the
  initial object and since $\mquote{\partial c}\subto c$ are generating
  $\partial$\=/flasque cofibrations.
\end{proof}

\begin{remark}
  It is important to note that the $\mquote\Gamma$ are \emph{not} projective
  cofibrant in $\sPsh(\cell\lfdex)$---unless $\lfdex$ is a set, which is just
  the case of multisorted Lawvere theories, in which case
  $\sPsh(\cell\lfdex)_\partial=\sPsh(\cell\lfdex)\proj$ and the $\mquote\Gamma$
  are all coproducts of objects in
  $\lfdex\subto\cell\lfdex\subto\psh{\cell\lfdex}$ (and thus coproducts of
  representables).\footnote{This cofibrancy plays a rôle in the rigidification
    of \cite{badzioch2002algtheories} (and the multisorted case in
    \cite{bergner2006rigidification}), and this was one of the motivations for
    defining the $\partial$\=/flasque model structure.}
\end{remark}

\begin{proposition}
 \label{prop:flasque-fib-local-is-htpy-coherent} 
  Let $X$ be a fibrant object in $\sPsh(\cell\lfdex)_\partial$, such that for
  every $\Gamma$ in $\cell\lfdex$, the map \emph{$X_\Gamma\to
    \map{\mquote\Gamma}X$} is a weak equivalence in $\SSet$. Then $X$ is a
  homotopical $\lfdex$\=/sorted space.
\end{proposition}
\begin{proof}
    First, $\mquote{\emptyset}$ is initial so by hypothesis $X_\emptyset\to 1$
    is an equivalence in $\SSet$. Next, we have the cube below in $\SSet$, whose
    front face is cartesian.
  \[
\begin{tikzcd}[row sep=small,column sep=tiny]
      X_{\Gamma_{n+1}} \ar[rr]\ar[rd] \ar[dd]
      && X_c \ar[rd,equals] \ar[dd]\\
      &\map{\mquote{\Gamma_{n+1}}}X\pbmark 
      && X_c \ar[from=ll, crossing over]\ar[dd]\\
      X_{\Gamma_n} \ar[rd] \ar[rr]
      && X_{\partial c} \ar[rd]\\
      &\map{\mquote{\Gamma_n}}X \ar[from = uu, crossing over] \ar[rr]
      && \map{\mquote{\partial c}}X
    \end{tikzcd}
  \]
  Since $X$ is $\partial$\=/flasque fibrant, $X_c\to\map{\mquote{\partial c}}X$
  is a fibration, so the front face is a homotopy pullback. By hypothesis,
  the intervening arrows are all weak equivalences, so the back face is a
  homotopy pullback (\cite[Prop. 13.3.13]{Hirschhorn2009}).
\end{proof}

\begin{proposition}
 \label{prop:ess-img-C-spaces-is-strict-local} 
  Any $X\in\sPsh(\cell\lfdex)$ is in the essential image of
  $i_*\colon\sPsh\lfdex\subto\sPsh(\cell\lfdex)$ if and only if for every
  $\Gamma$ in $\cell\lfdex$, \emph{$X_\Gamma\to\map{\mquote\Gamma}X$} is an
  isomorphism in $\SSet$.
\end{proposition}
\begin{proof}
  $X$ is in the essential image of $i_*$ if and only if
  $X_\Gamma\to\map{i_!i^*\Gamma}X$ is an isomorphism. By
  \cref{rem:desc-flasque-spine-counit}, this is what is desired.
\end{proof}

\begin{para}
   \label{para:flasque-bousfield-generators} 
  \cref{prop:flasque-fib-local-is-htpy-coherent,prop:ess-img-C-spaces-is-strict-local}
  tell us that the fibrant objects in $\sPsh(\cell\lfdex)_\partial$ that see the
  maps $\mquote\Gamma\subto\Gamma$ as weak equivalences are good candidates
  for the fibrant objects of a model structure for homotopical
  $\lfdex$-sorted spaces.

  \cref{prop:flasque-objs-are-cofibrant} tells us that we are particularly
  well-placed to perform a left Bousfield localisation at just this set of maps,
  since their domains and codomains are $\partial$\=/flasque cofibrant objects.
  So let us fix a name for this set of maps.
  \[
    S_\partial\eqdef
    \{s_\Gamma\colon\mquote\Gamma\subto\Gamma\mid\Gamma\in\cell\lfdex\}
  \]
\end{para}

\begin{definition}
  \label{defn:htpy-C-space}
  A \defn{homotopy $\lfdex$\=/sorted space} is an $S_\partial$-local object,
  namely a $\partial$\=/flasque fibrant object $X$ such that for every
  $\mquote{\Gamma}\subto\Gamma$ in $S_\partial$, the induced morphism $\hmap
  \Gamma X \to\hmap {\mquote{\Gamma}} X$ of homotopy function complexes is a weak
  equivalence in $\SSet$.
\end{definition}

\begin{definition}
  \label{defn:local-flasque-model-structure-C-spaces}
  The \defn{model structure for homotopy
    $\lfdex$\=/sorted spaces} is the left Bousfield localisation of
  $\sPsh(\cell\lfdex)_\partial$ at the set $S_\partial$. We write the model
  structure as
  $(\sPsh(\cell\lfdex)_\partial^l,\rmW^l,\rmC_\partial^l,\rmF_\partial^l)$.
\end{definition}

\begin{proposition}
  Every homotopy $\lfdex$\=/sorted space is a homotopical
  $\lfdex$\=/sorted space.
\end{proposition}
\begin{proof}
  By \cite[Prop. 3.4.1]{Hirschhorn2009}, $S_\partial$\=/local objects and
  fibrant objects of $\sPsh(\cell\lfdex)_\partial^l$ coincide. Since every map
  in $S_\partial$ is between $\partial$\=/flasque cofibrant objects, a
  $\partial$\=/flasque fibrant object $X$ is $S_\partial$\=/local if and only if
  $X_\Gamma\to\map{\mquote{\Gamma}}X$ is a weak equivalence for all $\Gamma$ in
  $\cell\lfdex$. We conclude by \cref{prop:flasque-fib-local-is-htpy-coherent}.
\end{proof}

\begin{proposition}
 \label{prop:local-flasque-to-C-spaces-Q-adj} 
  The adjunction
  $i^*:\sPsh(\cell\lfdex)_\partial^l\localisation\sPsh\lfdex\inj:i_*$ is
  Quillen.
\end{proposition}
\begin{proof}
  The maps in $S_\partial$ are between $\partial$\=/flasque cofibrant objects
  and are sent to isomorphisms by $i^*$, so this follows from
  \cref{thm:flasque-model-struc-generalises-projective}\ref{thm:flasque-model-struc-generalises-projective:2}
  and \cite[Thm 3.3.20]{Hirschhorn2009}.
\end{proof}

We would now like to show that it is a Quillen equivalence. It would be nice to
be able to generalise the technique from \cite{badzioch2002algtheories} for the
rigidification of ordinary algebraic theories. Unfortunately, this technique
relies fundamentally on the following statement, that I am unable to prove. In
fact, if we can, then the proof from \emph{op. cit.} works word-for-word.

\begin{conjecture}
 \label{conj-i-lower-star-pres-equivs} 
  Let $f$ be a global weak equivalence in $\sPsh\lfdex$. Then $i_*f$ is an
  $S_\partial$\=/local equivalence.
\end{conjecture}

\begin{remark}
  When $\lfdex$ is a set, the previous claim trivialises, since $i_*$ sends
  global weak equivalences to \emph{global} weak equivalences (and this is once
  again because weak equivalences are stable under products in $\SSet$). This is
  false for an arbitrary \lfd~category $\lfdex$, and it is easy to find a simple
  counter-example. Let $\lfdex=\{0\to 1\}$ be the ``walking arrow'' category.
  Recall that this corresponds to the type signature below.
  \[
    \vdash T_0 \qquad\qquad x\of T_0\vdash T_1
  \]
  Then consider the object $\Gamma$ of $\cell\lfdex$ that corresponds to the
  cell complex below, namely the context $(x\of T_0,y\of T_1(x),z\of T_1(x))$.
  \[
    \begin{tikzcd}[sep=small]
      0 \ar[r] \ar[d]
      &1\ar[d]\\
      1\ar[r] &\Gamma\pomark
    \end{tikzcd}
  \]
  Therefore, for any $f\in\sPsh\lfdex=\SSet\arr$, we have $(i_*f)_\Gamma =
  f_1\times_{f_0}f_1$. Consider a (trivial cofibration, fibration)
  factorisation in $\SSet$
  \[
    \begin{tikzcd}
      \Delta_0\ar[r,tail,"\sim"]\ar[rr,bend right,"f"]
      &E \ar[r,two heads,"g"]&\partial\Delta_2  
    \end{tikzcd}
  \]
  of an inclusion $f$ of a $0$\=/simplex into the boundary of the $2$\=/simplex.
  Then $f\to g$ is an objectwise weak equivalence in $\SSet\arr$. But $f_\Gamma$
  and $g_\Gamma$ are not equivalent.
\end{remark}

Luckily, we can make do without~\cref{conj-i-lower-star-pres-equivs}.

\begin{theorem}[Rigidification for $\lfdex$\=/sorted spaces]
  \label{thm:flasque-C-spaces-rigidification}
  The Quillen adjunction
  \[
    i^*:\sPsh(\cell\lfdex)_\partial^l\localisation\sPsh\lfdex\inj:i_*
  \]
  is a Quillen equivalence between the model category
  $\sPsh(\cell\lfdex)_\partial^l$ for homotopy $\lfdex$\=/sorted spaces and the
  model category $\sPsh\lfdex\inj$ for $\lfdex$\=/sorted spaces.
\end{theorem}
\begin{proof}
  We will show that it is a homotopy reflection and coreflection in the sense of
  \cite[E.2.17, E.2.24]{joyal2008theory}. Let $X$ be a fibrant-cofibrant object
  in $\sPsh(\cell\lfdex)_\partial^l$, then
  by~\cref{lem:flasque-fibration-inv-image-reedy}, $i^*X$ is Reedy fibrant, and
  by~\cref{lem:reedy-fibrant-is-htpy-coherent-C-space}, $i_*i^*X$ is
  $S_\partial$\=/local. So $X\to i_*i^*X$ is an $S_\partial$\=/local equivalence
  if and only if it is a global weak equivalence. We conclude since
  $(i_*i^*X)_\Gamma = \map{\mquote\Gamma}X$, and since $X$ is
  $S_\partial$\=/local. Hence the Quillen adjunction is a homotopy coreflection.
  To show it is a homotopy reflection, let $X$ be Reedy fibrant. Let $X'\to
  i_*X$ be a $\partial$\=/flasque cofibrant replacement. Since $i^*$ preserves global weak
  equivalences, $i^*X'\to i^*i_*X\cong X$ is a weak equivalence.\footnote{I
    thank P. LeFanu Lumsdaine for pointing out that the direct proof works.}
\end{proof}


\begin{theorem}
  \label{thm:id-Quillen-equiv-betw-flasque-local-and-Cech}
  The identity adjunction
  $\sPsh(\cell\lfdex)_\partial^l\rightleftarrows\sPsh(\cell\lfdex)\inj^\cech$ is
  a Quillen equivalence between the local flasque model structure and the
  injective \v Cech model structure.
\end{theorem}
\begin{proof}
  We have the commuting triangle of adjunctions below.
  \[
    \begin{tikzcd}[sep=small]
      &\sPsh\lfdex\inj \ar[dl,phantom,"{\scriptstyle\dashv}"description]
      \ar[dl,shift left=1ex,hook'] \ar[from=dl,shift left=1ex]
      \ar[dr,phantom,"{\scriptstyle\vdash}"description] \ar[dr,shift
      right=1ex,hook] \ar[from=dr,shift right=1ex]
      \\
      \sPsh(\cell\lfdex)_\partial^l
      \ar[rr,phantom,"{\scriptstyle\bot}"description] \ar[rr,shift left=.9ex]
      \ar[from=rr,shift left=.9ex] && \sPsh(\cell\lfdex)\inj^\cech
    \end{tikzcd}
  \]
  By ``2~out~of~3'' for Quillen equivalences, it suffices to show that the
  identity adjunction is Quillen. The identity functor takes
  $\partial$\=/flasque cofibrations to injective cofibrations and injective
  fibrant objects to $\partial$\=/flasque fibrant objects. Every
  subrepresentable $\mquote{\Gamma}\cong i_!i^*\Gamma\subto \Gamma$ is a
  covering sieve in the topology induced by $i\colon \lfdex\subto\cell\lfdex$,
  so every injective \v Cech fibrant object is an $S_\partial$\=/local object,
  so every $S_\partial$\=/local equivalence is a weak equivalence of
  $\sPsh(\cell\lfdex)\inj^\cech$. Thus the identity functor
  $\sPsh(\cell\lfdex)_\partial^l\to\sPsh(\cell\lfdex)\inj^\cech$ is left
  Quillen.\footnote{I thank U.P. for their birthday, the 13th of July.}
\end{proof}


\section{Homotopical models of \texorpdfstring{$\lfdex$}{C}\=/contextual
  categories}
\label{sec:htpical-models-C-cxl-cats}

In \cref{sec:C-sorted-spaces,sec:flasq-model-structure}, we have seen two ways
of describing homotopical models of the \emph{initial} $\lfdex$\=/contextual
category, as sheaves of spaces on $\cell\lfdex$
(\cref{prop:C-spaces-model-struct-equals-Cech,cor:Cech-model-struc-on-cellC-spaces-models-C-spaces})
and as homotopy-coherent models of the $\lfdex$\=/contextual category
$\cxl\lfdex$ (\cref{thm:flasque-C-spaces-rigidification}). We have also seen
that the two descriptions are equivalent
(\cref{thm:id-Quillen-equiv-betw-flasque-local-and-Cech}). These results
(especially \cref{thm:flasque-C-spaces-rigidification}) can be seen as
\emph{rigidification theorems} for the initial $\lfdex$\=/sorted theory.

We would like to be able to extend these results to other $\lfdex$\=/sorted
theories, such as the theories of small categories, $\omega$\=/categories
(strict and weak) and planar coloured operads. In particular, for a suitable
(perhaps for every?) $\lfdex$\=/contextual category $\cxl\lfdex\to\sfD$, we
would like to obtain a rigidification theorem in the form of a Quillen
equivalence generalising that of \cref{thm:badzioch-main}, between a suitable
local model structure for \emph{homotopy $\sfD$\=/algebras} and a model
structure on \emph{simplicial $\sfD$\=/algebras} (equivalently, models of
$\sfD$ in $\SSet$).
In this section, we will outline a tentative approach to solving this problem.
We conjecture that this approach allows for a solution, but there is an
obstruction that we are not yet able to surmount.

We fix a $\lfdex$\=/contextual category $\cxl\lfdex\to\sfD$ throughout this
section. We write the associated $\lfdex$\=/sorted theory as
$j\colon\cell\lfdex\to \Theta$, and the associated finitary monad on
$\psh\lfdex$ as $T$. We begin by showing that a \emph{$\partial$\=/flasque}
model structure on the category $\sPsh\Theta$ exists and is right-transferred
along the functor $j^*\colon\sPsh\Theta\to\sPsh(\cell\lfdex)_\partial$.

    

\begin{proposition}
   \label{prop:flasque-model-struc-Theta} 
    The right-transferred model structure on $\sPsh\Theta$ along the functor
    $j^*\colon\sPsh\Theta\to\sPsh(\cell\lfdex)_\partial$ exists, where
    $\sPsh(\cell\lfdex)_\partial$ is the $\partial$\=/flasque model structure of
    \cref{thm:flasque-model-struc-existence}.
\end{proposition}
\begin{proof}
    Since the functor $j$ is identity-on-objects, the image under $j_!$ of the
    sets $I\proj$ and $J\proj$ of generating projective (trivial) cofibrations
    of $\sPsh(\cell\lfdex)$ are the sets $I\proj$ and $J\proj$ of generating
    projective (trivial) cofibrations of $\sPsh{\Theta}$.
    
    For $c$ in $\lfdex$, recall the definition of $\mquote{\partial c}$ and
    $\delta_c\colon\mquote{\partial c}\subto c$ in $\psh{\cell\lfdex}$
    from~\cref{para:global-flasque-boundaries}. Then the map $j_!\delta_c\colon
    j_!\mquote{\partial c}\subto jc$ in $\psh\Theta$ is the sub-representable of
    $jc\in\Theta$ consisting of all maps $j\Gamma\to jc$ in $\Theta$ that factor
    through the image by $j$ of some map $c'\to c$ in $\lfdex$. Consider the set
    $j_!J_\partial= J\proj\cup \{h^k_n\pp j_!\mquote{\delta_c}\mid c\in
    \lfdex,[n]\in\bDelta,0\leq k\leq n\}$, that is the image under $j_!$ of the
    set $J_\partial$ of generating $\partial$\=/flasque trivial cofibrations
    (\cref{para:global-flasque-boundaries}).
    By~\cref{prop:right-induced-model-structure}, it suffices to show that any
    pushout of a map in $j_!J_\partial$ is taken by $j^*$ to a global weak
    equivalence in $\sPsh(\cell\lfdex)$. But since $j^*$ preserves colimits and
    monomorphisms, $j^*j_!\mquote{\delta_c}$ is an injective cofibration and
    $h^k_n\pp j^*j_!\mquote{\delta_c}$ is an injective trivial cofibration.
\end{proof}

\begin{definition}
    We call the model structure of \cref{prop:flasque-model-struc-Theta} the
    \defn{$\partial$\=/flasque model structure} on $\sPsh\Theta$, and write this
    model category as $\sPsh\Theta_\partial$.
\end{definition}


\begin{remark}
    Let $X$ be in the essential image of $\s\sfD\Mod\subto\sPsh\Theta$. Then
    $j^*X$ is in the essential image of $\sPsh\lfdex\subto\sPsh(\cell\lfdex)$.
    Then $X$ is fibrant in $\sPsh\Theta_\partial$ if and only if $j^*X$ is Reedy
    fibrant in $\sPsh\lfdex$.
\end{remark}

\begin{para}
    For every $\Gamma$ in $\cell\lfdex$, recall the definition of
    $\mquote\Gamma$ and the set $S_\partial =
    \{\mquote\Gamma\subto\Gamma\mid\Gamma\in\cell\lfdex\}$
    from~\cref{para:flasque-local-eqs}
    and~\cref{para:flasque-bousfield-generators}. Since $j_!$ is left Quillen
    and since $S_\partial$ is a set of maps between $\partial$\=/flasque
    cofibrant objects, the image $j_!S_\partial$ is also a set of maps between
    cofibrant objects of $\sPsh\Theta_\partial$.
\end{para}

\begin{definition}
  The \defn{model category of homotopy $\sfD$\=/algebras} is the left Bousfield
  localisation of the model category $\sPsh\Theta_\partial$ at the set
  of maps $j_!S_\partial$. We write it as $\sPsh\Theta_\partial^l$.
\end{definition}

\begin{remark}
  Any $X$ in $\sPsh\Theta$ is fibrant in $\sPsh\Theta_\partial^l$ if and only
  if $j^*X$ is fibrant in $\sPsh(\cell\lfdex)$. If $X$ is in the essential
  image of $\s\sfD\Mod$ then it is fibrant in $\sPsh\Theta_\partial^l$ if and
  only if it is fibrant in $\sPsh\Theta_\partial$.
\end{remark}

\begin{theorem}
  \label{thm:partial-soln-rigid-conjecture} 
  In the exact adjoint square associated to the $\lfdex$\=/sorted theory $j$
  \begin{equation}
    \tag{$\star$}
    \label{eq:htpy-T-alg-exact-adj-square}
    \begin{tikzcd}
      \sPsh \lfdex\inj\ar[r,shift left,"j_!"] \ar[d,shift
      left,hook',"i_*"]
      &\s\sfD\Mod \ar[l,shift left,"j^*"]\ar[d,shift left,hook',"N"]\\
      \sPsh(\cell\lfdex)_\partial^l\ar[r,shift left,"j_!"] \ar[u,shift
      left,"i^*"] &\sPsh\Theta_\partial^l \ar[l,shift left,"j^*"]
      \ar[u,shift left,"h"]
    \end{tikzcd}
  \end{equation}
  the lower horizontal adjunction is a Quillen adjunction and the left
  vertical adjunction is a Quillen equivalence.
\end{theorem}
\begin{proof}
  This follows from~\cref{thm:flasque-C-spaces-rigidification} and from the
  definition of $\sPsh\Theta_\partial^l$.
\end{proof}

\begin{para}
  [Rigidification of homotopy $\sfD$\=/algebras]
  \label{para:rigidification-htpy-D-algs-conjecture} 
  \cref{thm:partial-soln-rigid-conjecture} is a partial generalisation of
  \cref{rem:bousfield-exact-adj-square-I-theory}. In order to obtain a full
  rigidification theory for $\lfdex$\=/contextual categories, we would like to
  transfer the model structure $\sPsh\lfdex\inj$ for $\lfdex$\=/sorted spaces
  along the right adjoint $\s\sfD\Mod\to\sPsh\lfdex$.

  When $\lfdex$ is a \emph{set}, then $\sPsh\lfdex\inj=\sPsh\lfdex\proj$ and
  this is just the technique from \cite[II.4]{quillen1967} (reprised in
  \cite[Thm 3.1]{Schwede2001stablehtpyAlgTheories}, \cite[Sec.
  7]{Rezk2002simplicialAlgTheories}) of constructing a model structure on the
  category of simplicial objects in models of a multisorted algebraic theory.
  Quillen's technique relies on the existence of a fibrant replacement functor
  that passes from $\sPsh\lfdex\inj$ to $\s\sfD\Mod$. When
  $\sPsh\lfdex\inj=\sPsh\lfdex\proj$, we can choose any
  finite-limit-preserving (finite products suffice) fibrant replacement
  functor on objects of $\SSet$, such as Kan's $\mathrm{Ex}^\infty$.

  However, in the general case where $\lfdex$ is any \lfd~category, we need a
  \emph{Reedy} fibrant replacement functor on $\sPsh\lfdex\inj$ that passes to
  $\s\sfD\Mod$. This seems to be a little difficult, and I have not yet found
  a way out of this impasse.
\end{para}

\chapter{Opetopic theories}
\label{chap:opetopic-theories}
In this chapter, we will describe certain familiar algebraic structures (small
categories and coloured planar operads among them) as $\Set$\=/models of
\emph{idempotent} $\OO$\=/sorted theories, where $\OO$ is the \lfd~category of
\defn{opetopes}. We thus transform what is ostensibly part of the ``structure''
of a category or operad (like its operation of composition) into ``properties''
on some $\OO$-sorted set (an \emph{opetopic set}).

Opetopes (\emph{ope}ration poly\emph{topes}) are certain cellular shapes that
are ``tree-like'' in every dimension---namely, every opetope can be seen both as
an operation (an elementary planar tree or corolla) in its dimension, and as a
planar tree of opetopes of lower dimension. For instance, here is an opetope $\mu$ of
dimension $2$.
\[
  \tikzinput{opetope-composition}{comp4}
  \qquad\qquad
  \begin{tikzpicture}
    \tikzset{scale=.8}
    \node [treenode] (b) at (0, 0) [label = left : $\mu$] {};
    \node (e1) at ($(b) + (-1, .7)$) [above] {\small $ f$};
    \node (e4) at ($(b) + (1, .7)$) [above] {\small $i$};
    \node (e2) at ($(e1)!.35!(e4)$) {\small $g$};
    \node (e3) at ($(e1)!.7!(e4)$) {\small $h$};
    \node (t) at ($(b) + (0, -.7)$) [below] {\small $k$};
    \draw (b) -- (e1);
    \draw (b) -- (e2);
    \draw (b) -- (e3);
    \draw (b) -- (e4);
    \draw (b) -- (t);
  \end{tikzpicture}
\]
It can be seen as a filiform tree $a\to b\to c\to d \to e$ of opetopes of
dimension $1$, as well as a planar corolla (a rooted tree with one node). It
encodes the shape of an operation of composition of morphisms in a
category, that takes as input a filiform tree and returns as output a filiform
tree of length $1$.
Hence, $\mu$ can be seen as the \emph{compositor} of $f$, $g$, $h$, and $i$, and
outputs the actual composition $k$.

A second example, one dimension higher, is that of a planar coloured
$\Set$-operad $\cP$ (a.k.a. a nonsymmetric multicategory), whose compositors
have planar trees of composable multimorphisms of $\cP$ as arities. 
\[
    \tikzinput{opetope-3-graphical}{ex1'.labels}
\]
Here, $A$ is the compositor of the pasting of $\alpha$, $\beta$, $\gamma$, and
$\delta$ as on the left, and points towards the actual composition $\alpha 
(\beta, \gamma (\delta))$.

Heuristically extending this pattern, one infers that such an algebraic
structure one dimension above that of planar coloured operads should have an
operation of composition whose arities are suitably ``planar'' trees of
``operations'' whose inputs are \emph{planar trees}. Indeed, such algebraic
structures are precisely the (coloured) \emph{$PT$\=/combinads in $\Set$}
(combinads over the combinatorial pattern of planar trees) of Loday
\cite{Loday2012a}.

The goal of this chapter is to give a precise definition of this hierarchy of
algebraic structures (that we will call \emph{opetopic algebras}), and to show
that they are naturally encoded by idempotent dependently sorted algebraic
theories in the sense of \cref{chap:contextual-categories}.
We end by giving a partial answer to the conjecture
of~\cref{para:rigidification-htpy-D-algs-conjecture}. Namely, we show that there
exists a model structure on the category of simplicial opetopic algebras, that
has an associated rigidification theorem. 
\begin{notation}
    [Change of notation] We will no longer use $\omega$ to denote the first
    infinite ordinal/cardinal, since we will require it in our notation for
    opetopes. Therefore, for this chapter \emph{only}, we use
    $\infty$ to denote the first infinite ordinal.
\end{notation}
\section{Polynomial endofunctors and monads}
\label{sec:poly-functors-monads}
In this section, we recall some of the theory of polynomial endofunctors and
monads in $\Set$ that we will require to define the category of opetopes. 
\begin{para}
  \label{para:polynomial-functor}
  A \defn{polynomial (endo)functor} $P$\index{polynomial endofunctor} is a
  diagram of the form below in $\Set$.
  \[
    \polynomialfunctor{I}{E}{B}{I}{s}{p}{t}
  \]
  We use the following terminology for $P$. The
  elements of $B$ are its \emph{nodes} or
  \emph{operations}, the elements of $I$ are its
  \emph{colours} or \emph{sorts}, and for every node $b$, the
  elements of the fibre $E_b \eqdef p^{-1} (b)$ are the
  \emph{inputs} of $b$. For every
  input $e$ of a node $b$, we denote its colour by $s_e (b) \eqdef s (e)$.  
  \[
    \tikzinput[.8]{polynomial-functors}{operation}
  \]
  $P$ is \defn{finitary} if the fibres of $p\colon E \to B$ are finite sets. All
  polynomial functors we consider will be finitary. A morphism $P\to P'$ of polynomial
  functors is a diagram 
  \[
    \begin{tikzcd}[sep=scriptsize]
      I \ar[d, "f_0"'] & E \pullbackcorner \ar[r, "p"] \ar[d, "f_2" left] \ar[l,
      "s" above] & B \ar[r, "t"] \ar[d, "f_1" left] &
      I \ar[d, "f_0" left] \\
      I' & E' \ar[r, "p'"] \ar[l, "s'" above] & B' \ar[r, "t"] & I'
    \end{tikzcd}
  \]
  where the middle square is cartesian (a pullback square). If $P$ and $P'$ have
  the same set $I$ of colours, then a morphism from $P$ to $P'$ \emph{over $I$}
  is a commutative diagram as above, but where $f_0$ is the identity function.
  Let $\PolyEnd$\index{$\PolyEnd$} denote the category of polynomial functors,
  and $\PolyEnd (I)$\index{$\PolyEnd (I)$} the category of polynomial functors
  over $I$.
\end{para}

\begin{remark}
  \label{rem:polynomial-functor-functor}
  Every $P$ in $\PolyEnd$ has an associated composite endofunctor
  \[
    \begin{tikzcd}[sep=scriptsize]
      P\colon \Set / I
      \ar[r,"s^*"] &\Set / E
      \ar[r,"p_*"] &\Set / B
      \ar[r,"t_!"] &\Set / I
    \end{tikzcd}
  \]
  Explicitly, for $X = (X_i \mid i \in I) \in \Set / I$, we calculate $P(X)$ as
  the ``polynomial''
  \begin{equation}
    \label{para:polynomial-functor:eval}
    P (X)_j = \sum_{b \in B_j} \prod_{e \in E_b} X_{s(e)}, \quad j \in I
  \end{equation}
  where $B_j \eqdef t^{-1} (j)$. Visually, elements of $P (X)_j$ are nodes $b
  \in B$ such that $t b = j$, and whose inputs are decorated by elements of
  $(X_i \mid i \in I)$ in a manner compatible with their colours. Graphically,
  an element of $P(X)_i$ can be represented as
  \[
    \tikzinput{polynomial-functors}{evaluation}
  \]
  with $t(b) = i$ and $x_j \in X_{s_{e_j} b}$ for $1 \leq j \leq k$. Moreover,
  the endofunctor $P: \Set / I \to \Set / I$ preserves connected limits: $s^*$
  and $p_*$ preserve all limits (as right adjoints), and $t_!$ preserves and
  reflects connected limits. It is finitary (preserves filtered colimits) if and
  only if the polynomial functor is finitary (\cref{para:polynomial-functor}).

  This gives a fully faithful functor $\PolyEnd(I)\to [\Set/I,\Set/I]_\cart$,
  the latter being the category of endofunctors of $\Set / I$ and cartesian
  natural transformations\footnote{A natural transformation is
    \emph{cartesian}\index{cartesian natural transformation} if all its
    naturality squares are cartesian.}. In fact, the image of this full
  embedding consists precisely of those endofunctors that preserve connected
  limits \cite[1.18]{Gambino2013}. The composition of endofunctors gives
  $[\Set/I,\Set/I]_\cart$ the structure of a monoidal category, and
  $\PolyEnd(I)$ is stable under this monoidal product \cite[1.12]{Gambino2013}.
  The identity polynomial functor $\polynomialinline{I}{I}{I}{I}$ is associated
  to the identity endofunctor; thus $\PolyEnd(I)$ is a monoidal subcategory of
  $[\Set/I,\Set/I]_\cart$.
\end{remark}
\begin{remark}
  We refer the reader to \cite{Kock2011, Gambino2013} for a detailed account of
  the theory of polynomial functors.
\end{remark}

\begin{definition}
  [{\cite[1.0.3]{Kock2011}}]
  \label{def:polynomial-tree}
  A polynomial functor $T$ given by
  \[
    \polynomialfunctor{T_0}{T_2}{T_1}{T_0}{s}{p}{t}
  \]
  is a \defn{(polynomial) tree}\index{polynomial tree} if
  \begin{enumerate}
  \item the sets $T_0$, $T_1$ and $T_2$ are finite (in particular, each node has
    finitely many inputs);
  \item $t$ is injective;
  \item $s$ is injective, and the complement $T_0 - \im s$
    has a single element, the \emph{root}\index{root} of $T$;
  \item let $T_0 = T_2 + \{ r \}$, with $r$ the root, and define the
    \emph{walk-to-root}\index{walk-to-root function} function $\sigma$ by
    $\sigma (r) = r$, and otherwise $\sigma (e) = t p (e)$; then we ask that for
    all $x \in T_0$, there exists $k \in \NN$ such that $\sigma^k (x) = r$.
  \end{enumerate}
  We call the colours of a tree its \emph{edges} and the inputs of a node the
  \emph{input edges}\index{input edge} of that node. Let $\Tree$\index{$\Tree$}
  (denoted \textbf{\emph{TEmb}} in \cite{Kock2011}) be the full subcategory
  of $\PolyEnd$ whose objects are trees. An \defn{elementary
    tree}\index{elementary polynomial tree} is a tree with at most one node. Let
  $\elTree$\index{$\elTree$} be the full subcategory of $\Tree$ spanned by
  elementary trees.
\end{definition}

\begin{remark}
  $\Tree$ is the category of \emph{symmetric} or \emph{non-planar} trees (the
  automorphism group of a tree is in general non-trivial) and its morphisms
  correspond to inclusions of non-planar subtrees. The category $\Omega$ of
  \emph{dendrices} is the image of $\Tree$ under the ``free coloured symmetric
  operad'' functor \cite[1.3]{Kock2011}.
\end{remark}

\begin{definition}
  [$P$-tree]
  \label{def:p-tree}
  For $P $ in $\PolyEnd$, the over-category $\Tree / P$ is the category
  $\tree P$\index{$\tree P$} of \defn{$P$-trees}\index{$P$-tree}. The
  fundamental difference between $\Tree$ and any $\tree P$ is that the latter is
  always rigid i.e. it has no non-trivial automorphisms \cite[Prop.
  1.2.3]{Kock2011}. In particular, this implies that $\PolyEnd$ does not have a
  terminal object.
\end{definition}

\begin{notation}
  \label{not:p-tree}
  Every $P$-tree $T \in \tree P$ corresponds to a morphism from a tree (which we
  denote by $\underlyingtree{T}$\index{$\underlyingtree{T}$}) to $P$, so that $T
  : \underlyingtree{T} \to P$. We point out that $\underlyingtree{T}_1$ is the
  set of nodes of $\underlyingtree{T}$, while $T_1 : \underlyingtree{T}_1 \to
  P_1$ is a decoration of the nodes of $\underlyingtree{T}$ by nodes of $P$, and
  likewise for edges.
  \end{notation}

\begin{definition}
  [{\cite[2.1.1]{Kock2011}}]
  \label{def:category-of-elements}
  For $P \in \PolyEnd$, its \defn{category of elements}\index{category of
    elements of a polynomial functor}\footnote{Not to be confused with the
    category of elements of a presheaf.} $\elt P$ is the over-category $\elTree/P$.
  For $P$ as in \cref{para:polynomial-functor}, the set of objects
  of $\elt P$ is $I + B$, and for each $b \in B$, there is a morphism $\tgt :
  t(b) \to b$\index{$\tgt$}, and a morphism $\src_e : s_e(b) \to
  b$\index{$\src_e$} for each $e \in E_b$. 
\end{definition}

\begin{lemma}
 \label{lem:elt-P-lfd} 
  For every $P\in\PolyEnd$, $\elt P$ is a \lfd~category.
\end{lemma}
\begin{proof}
  Immediate, since there is no non-trivial composition of arrows.
\end{proof}

\begin{proposition}
  [{\cite[Prop. 2.1.3]{Kock2011}}]
  \label{prop:category-of-elements-slice}
  There is an equivalence of categories between the presheaf category $\psh{\elt
    P}$ and the slice category $\PolyEnd\slice P$.
\end{proposition}

\begin{notation}[Addresses]
  \label{not:polynomial-addresses}
  
  Let $T \in \Tree$ with $\sigma$ its walk-to-root
  function (\cref{def:polynomial-tree}). We define the
  \defn{address}\index{address} function $\addr$\index{$\addr$|see {address}} on
  edges inductively as follows:
  \begin{enumerate}
  \item if $r$ is the root edge, let $\addr r \eqdef []$,
  \item if $e \in T_0 - \{ r \}$ and if $\addr \sigma (e) = [x]$, define $\addr
    e \eqdef [x e]$.
  \end{enumerate}
  The address of a node $b \in T_1$ is defined as $\addr b \eqdef \addr t (b)$,
  namely that of its target.
  Note that this function is injective since $t$ is. Let
  $T^\nodesymbol$\index{$T^\nodesymbol$|see {node address}} denote its image,
  the set of \emph{node addresses}\index{node address} of $T$, and let
  $T^\leafsymbol$\index{$T^\leafsymbol$|see {leaf address}} be the set of
  addresses of leaf edges\index{leaf address}, i.e. those not in the image of
  $t$.

  Assume now that $T : \underlyingtree{T} \to P$ is a $P$-tree. If $b \in
  \underlyingtree{T}_1$ has address $\addr b = [p]$, write $\src_{[p]} T \eqdef
  T_1 (b)$. For convenience, we let $T^\nodesymbol \eqdef
  \underlyingtree{T}^\nodesymbol$, and $T^\leafsymbol \eqdef
  \underlyingtree{T}^\leafsymbol$.
\end{notation}

\begin{remark}
  The formalism of addresses is purely bookkeeping syntax for the operations of
  grafting and substitution on trees. The syntax of addresses will extend to the
  category $\OO$ of opetopes and will allow us to give a precise description of
  the composition of arrows in $\OO$ (see \cref{def:o}) as well as certain
  constructions in $\psh\OO$.
\end{remark}

\begin{notation}
  \label{not:marked-tree}
  We denote by $\treewithleaf P$\index{$\treewithleaf P$} the set of $P$-trees
  with a marked leaf, namely a tree along with the address of one of its leaves.
  Similarly, we denote by $\treewithnode P$\index{$\treewithnode P$} the set of
  $P$-trees with a marked node.
\end{notation}

        

\begin{para}[Grafting]
  \label{para:grafting}
  Let $P$ be as in \cref{para:polynomial-functor}. For $i \in I$, define
  $\itree{i} \in \tree P$ as having underlying tree $\polynomialinline
  {\{i\}}{\emptyset}{\emptyset}{\{i \}}$, along with the obvious morphism to $P$
  picking out $i \in I$. Then $\itree i$ is the $P$-tree with no nodes and a
  unique edge decorated by $i$. Define $\ytree{b} \in \tree P$, the
  \emph{corolla}\index{corolla} at $b$, as having underlying tree
  \[
    \polynomialfunctor {s(E_b) + \{*\}}{E_b}{\{b \}}{s(E_b) + \{*\},} {s}{}{}
  \]
  where the rightmost map sends $b$ to $*$, and where the morphism $\ytree{b}
  \to P$ is the identity on $s (E_b) \subseteq I$, maps $*$ to $t (b) \in I$, is
  the identity on $E_b \subseteq E$, and maps $b$ to $b \in B$. Then $\ytree b$
  is the $P$-tree with a single node (a corolla) decorated by $b$. Observe that
  for $T \in \tree P$, giving a morphism $\itree{i} \to T$ is equivalent to
  specifying the address $[p]$ of an edge of $T$ decorated by $i$. Likewise,
  morphisms of the form $\ytree{b} \to T$ are in bijection with addresses of
  nodes of $T$ decorated by $b$.

  For $S, T \in \tree P$, $[l] \in S^\leafsymbol$ such that the leaf of $S$ at
  $[l]$ and the root edge of $T$ are decorated by the same $i \in I$, define the
  \defn{grafting}\index{grafting} $S \graft_{[l]} T$ of $S$ and $T$ on $[l]$ by
  the following pushout in $\tree P$ (\cite[Prop. 1.1.19]{Kock2011}).
  \[
    \pushoutdiagram {\itree{i}}{T}{S}{S \graft_{[l]} T} {[]}{[l]}{}{}
  \]
    If $S$ (resp. $T$) is a trivial tree, then $S \graft_{[l]} T = T$
    (resp. $S$). We assume, by convention, that the grafting operator $\circ$
    associates to the right.
\end{para}

\begin{proposition}
    [{\cite[Prop. 1.1.21]{Kock2011}}]
    \label{prop:polynomial-functor:trees-are-graftings}
    Every $P$-tree is either of the form $\itree{i}$, for some $i \in I$, or
    obtained by iterated graftings of corollas (i.e. $P$-trees of the form
    $\ytree{b}$ for $b \in B$).
\end{proposition}

\begin{notation}
    [Total grafting]
    \label{not:total-grafting}
    Let $T, U_1, \ldots, U_k \in \tree P$, where $T^\leafsymbol = \left\{ [l_1],
    \ldots, [l_k] \right\}$, and assume the grafting $T \graft_{[l_i]} U_i$ is
    defined for all $i$. Then the \emph{total grafting} of $U_1,\ldots,U_k$ onto
    the leaves of $T$ will be denoted
    concisely by
    \begin{equation}
        \label{eq:big-grafting}
        T \biggraft_{[l_i]} U_i
        =
        ( \cdots (T \graft_{[l_1]} U_1) \graft_{[l_2]} U_2 \cdots) \graft_{[l_k]} U_k .
    \end{equation}\index{$T \biggraft_{[l_i]} U_i$}
    Clearly, the result does not depend on the order in which
    graftings are performed.
\end{notation}

\begin{para}
    [Substitution]
    \label{para:substitution}
    Let $[p] \in T^\nodesymbol$ and $b = \src_{[p]} T$. Then $T$ can be
    decomposed as
    \begin{equation}
        \label{eq:node-decomposition}
        T = A \graft_{[p]} \ytree{b} \biggraft_{[e_i]} B_i ,
    \end{equation}
    where $E_b = \{e_1, \ldots, e_k\}$, and $A, B_1, \ldots, B_k \in \tree P$
    and the corolla $\ytree b$ at the node $b$ at address $[p]$ is grafted onto
    a leaf of $A$. For $U$ a $P$-tree with a bijection $\readdress :
    U^\leafsymbol \to E_b$ over $I$, we define the
    \emph{substitution}\index{substitution} $T \subst_{[p]} U$\index{$T
      \subst_{[p]} U$|see {substitution}} as
    \begin{equation}
        \label{eq:subst}
        T \subst_{[p]} U
        \eqdef
        A \graft_{[p]} U \biggraft_{\readdress^{-1} e_i} B_i .
    \end{equation}
    In other words, the node at address $[p]$ in $T$ has been replaced by $U$,
    and the map $\readdress$ provides ``rewiring instructions'' to connect the
    leaves of $U$ to the rest of $T$.
\end{para}


\begin{para}
  [Polynomial monads]
  \label{para:polynomial:monads}
  A \defn{polynomial monad over $I$}\index{polynomial monad} is a monoid in
  $\PolyEnd(I)$ (thus necessarily a cartesian monad on $\Set/I$).\footnote{A
    monad on $C$ is \emph{cartesian}\index{cartesian monad} if it is a monoid in
    $\fun CC_\cart$.} We denote $\PolyMnd(I)$ the category of monoids in
  $\PolyEnd(I)$.

  If $M$ is in $\PolyMnd(I)$ and $M'$ is in $\PolyMnd(J)$ then a morphism of
  polynomial monads $M\to M'$ is a morphism of polynomial functors that respects
  the monoid structure. The category of all polynomial monads is denoted
  $\PolyMnd$.
\end{para}
\begin{definition}[{\cite[1.2.7]{Kock2011}}]
  \label{def:free-polymonad}
  Given a polynomial endofunctor $P$, we define the \defn{free monad on $P$},
  denoted $P^\star$, as
  \[
    \polynomialfunctor {I}{\treewithleaf P}{\tree P}{I} {s}{p}{t}
  \]
  where $s$ maps a $P$-tree with a marked leaf to the decoration of that leaf,
  $p$ forgets the marking, and $t$ maps a tree to the decoration of its root.
  Remark that for $T \in \tree P$ we have $p^{-1} T = T^\leafsymbol$.
\end{definition}

\begin{proposition}
  [{\cite[Prop. 1.2.8]{Kock2011}}]
  \label{prop:polynomial-monads-star-algebras}
    The forgetful functor from $\PolyMnd (I)$ to $\PolyEnd (I)$ is a monadic right
    adjoint, and its left adjoint is given by $(-)^\star$.
\end{proposition}


\begin{para}
  [Readdressing]
  \label{def:readdressing}
  Let $M$ be a polynomial monad as on the left below. By
  \cref{prop:polynomial-monads-star-algebras}, $M$ is an algebra for the
  free-monoid monad $(-)^\star$, and we will write its structure map $M^\star
  \to M$ as on the right.
  \[
    \polynomialfunctor{I}{E}{B}{I,}{s}{p}{t} \qquad
    \begin{tikzcd}[sep=scriptsize]
      I\ar[d, equal]
      &\treewithleaf M \pullbackcorner
      \ar[l] \ar[d, "\readdress"'] \ar[r]
      &\tree M \ar[r] \ar[d, "\tgt"]
      &I\ar[d, equal]\\
      I
      &E \ar[l] \ar[r]
      &B \ar[r]
      &I
    \end{tikzcd}
  \]
  We call $\readdress_T : T^\leafsymbol \cong E_{\tgt T}$ the
  \emph{readdressing}\index{readdressing} function of $T$, and $\tgt T \in B$ is
  called the \emph{target} of $T$. If we think of any $b \in B$ as the corolla
  $\ytree{b}$, then the target map $\tgt$ ``contracts'' a tree to a corolla, and
  since the middle square is a pullback, the number of leaves is preserved. The
  map $\readdress_T$ establishes a coherent correspondence between the set
  $T^\leafsymbol$ of leaf addresses of a tree $T$ and the set $E_{\tgt T}$ of
  inputs of $\tgt T$.
\end{para}



\begin{definition}
  \label{def:baez-dolan-construction}
  For a polynomial monad $M$, we define
  its \defn{Baez--Dolan $(-)^+$ construction}\index{Baez--Dolan $(-)^+$ construction} $M^+$ to
  be
  \[
    \polynomialfunctor{B}{\treewithnode M}{\tree M}{B,}{\src}{p}{\tgt}
  \]
  where $\src$ maps an $M$-tree with a marked node to the label of that node,
  $p$ forgets the marking, and $\tgt$ is the target map. If $T \in \tree M$,
  remark that $p^{-1} T = T^\nodesymbol$ is the set of node addresses of $T$. If
  $[p] \in T^\nodesymbol$, then $\src ([p]) \eqdef \src_{[p]} T$.
\end{definition}

\begin{proposition} [{\cite[3.2]{Kock2010}}]
    \label{prop:polynomial-functor:+:is-monad}
    The polynomial functor $M^+$ has a canonical structure of a polynomial
    monad.
\end{proposition}

\begin{remark}
  The $(-)^+$ construction is an endofunctor on $\PolyMnd$. If we begin with a
  polynomial monad $M$, then the colours of $M^+$ are the operations of $M$. The
  operations of $M^+$, along with their output colour, are given by the monad
  multiplication of $M$: they are the \emph{relations} of $M$, i.e. the
  reductions of $M$\=/trees to operations of $M$. The monad multiplication on
  $M^+$ is given as follows: the reduction of a tree of $M^+$ to an operation of
  $M^+$ (which is an $M$\=/tree) is obtained by \emph{substituting} $M$\=/trees
  into nodes of $M$\=/trees.
\end{remark}


\begin{proposition}
	\label{prop:algebraic-realization:equivalence-P+-algebras}
  For $M$ a polynomial monad, there is an equivalence of categories between the
  category of $M^+$-algebras and the slice category $\PolyMnd(I)\slice M$.
\end{proposition}
\begin{proof}
  For an $M^+$-algebra $M^+X \xrightarrow{x} X$ in $\Set/B$, let $\Phi X$ in
  $ \PolyEnd(I)\slice M$ be
  \[
    \begin{tikzcd}[sep=scriptsize]
      I \ar[d, equal] & E_X \pullbackcorner \ar[l] \ar[d] \ar[r] & X \ar[r]
      \ar[d] & I \ar[d, equal]
      \\
      I & E \ar[l] \ar[r] & B \ar[r] & I.
    \end{tikzcd}
  \]
  There is an evident bijection $\tree \Phi X \cong M^+X$ in $\Set/I$, and the
  structure map $x$ extends by pullback along $E_X \to X$ to a map $(\Phi
  X)^\star \to \Phi X$ in $\PolyEnd(I)$. It is easy to verify that this
  determines a $(-)^\star$-algebra structure on $\Phi X$, and that the map $\Phi
  X \to M$ in $\PolyEnd(I)$ is a morphism of $(-)^\star$-algebras. Conversely,
  given an $N \in \PolyMnd(I)\slice M$ whose underlying polynomial functor is
  $\polynomialinline {I} {E'} {B'} {I}$, then the bijection $\tree N \cong M^+
  B'$ in $\Set/I$ and the $(-)^\star$-algebra map $N^\star \to N$ provide a map
  $M^+ B' \xrightarrow{\Psi N} B'$ in $\Set/I$. It is easy to verify that $\Psi
  N$ is the structure map of a $M^+$-algebra and that the  constructions $\Phi$
  and $\Psi$ are functorial and mutually inverse.
\end{proof}



\section{The category \texorpdfstring{$\OO$}{O} of opetopes}
\label{sec:opetopes}

In this section, we will use the formalism of polynomial functors and polynomial
monads of \cref{sec:poly-functors-monads} to define the \lfd~category $\OO$ of
opetopes. As we will see, the category $\OO$ is a \lfd~category. Our
construction of opetopes is precisely that of \cite{Kock2010}, and by \cite[Thm
3.16]{Kock2010}, also that of \cite{Leinster2004}, and by \cite[Cor.
2.6]{Cheng2004a}, also that of \cite{Cheng2003}. The definition of the category
of opetopes that we use first appears in \cite{hothanh18}, and this section is
largely drawn from \cite{HTLS2019opetopic}.

\begin{para}[Defining opetopes]
  \label{para:opetopes-definition}
  Let $\optPolyFun^0$ be the identity polynomial monad, as on the left
  below, and let $\optPolyFun^{n+1} \eqdef (\optPolyFun^{n})^+$ for all
  $n\in\NN$. Write $\optPolyFun^n$ as on the right below.
  \[
    \polynomialfunctor {\{* \}}{\{* \}}{\{* \}}{\{* \},} {}{}{} \qquad
    \polynomialfunctor {\OO_n}{E_{n+1}}{\OO_{n+1}}{\OO_n .} {\src}{p}{\tgt}
  \]
  The set of \defn{$n$-dimensional opetopes} (or simply \defn{$n$-opetopes}) is
  $\OO_n$. Thus for $n\geq 0$, an $(n+2)$-opetope $\omega$ is a
  $\optPolyFun^n$\=/tree. It is called \emph{degenerate}\index{degenerate
    opetope} if it is a $\optPolyFun^{n}$-tree with no nodes (thus
  $\omega=\itree\phi$ for some $\phi\in\OO_{n}$); it is \emph{non
    degenerate}\index{non-degenerate opetope} otherwise. The readdressing
  function (see \cref{def:readdressing}) gives a bijection $\readdress_\omega \colon
  \omega^\leafsymbol \to (\tgt \omega)^\nodesymbol$ between the leaves of
  $\omega$ and the nodes of $\tgt \omega$, preserving the decoration by
  $n$-opetopes.
\end{para}

\begin{example}
 \label{exa:opetopes} 
  \begin{enumerate}
  \item The unique $0$-opetope is denoted $\optZero$\index{$\optZero$|see {point
        (opetope)}} and called the \emph{point}\index{point (opetope)}.
  \item The unique $1$-opetope is denoted $\optOne$\index{$\optOne$|see {arrow
        (opetope)}} and called the \emph{arrow}\index{arrow (opetope)}. We will
    also denote it by $\ytree{\optZero}$ (strictly speaking, it is not a
    tree, but it is the unary operation of $\optPolyFun^0$ and so can
    be seen as a corolla with a single input).
      
  \item Any $\omega\in \OO_{n+2}$ is a $\optPolyFun^{n}$-tree, namely a tree whose
    nodes are labeled with $(n+1)$-opetopes, and edges are labeled with
    $n$-opetopes. In particular, $2$-opetopes are $\optPolyFun^0$-trees, which
    are the
    filiform or linear trees, and thus in bijection with $\NN$. We will refer to
    them as \emph{opetopic integers}\index{opetopic integers}, and write
    $\optInt{n}$ for the $2$-opetope having exactly $n$ nodes.
  \end{enumerate}
\end{example}

\begin{para}
  [Higher dimensional addresses]
  \label{para:higher-address}
  By definition, an opetope $\omega$ of dimension $n \geq 2$ is a
  $\optPolyFun^{n-2}$-tree, thus the formalism of tree addresses
  (\cref{not:polynomial-addresses}) can be applied to designate nodes of
  $\omega$. We will iterate this to give
  \defn{higher dimensional addresses}\index{higher dimensional address}, which
  are an efficient syntax for the morphisms of the category of opetopes.
  Start by defining the sets $\AA_n$ of \emph{$n$-addresses} as follows:
  \[
    \AA_0 \eqdef \left\{ * \right\} , \qquad \AA_{n+1} \eqdef \operatorname{lists} (\AA_n) ,
  \]
  where $\operatorname{lists} (X)$ is the set of finite lists (or words) on the
  alphabet $X$.
  Explicitly, the unique $0$-address is $*$ (also written $[]$ by convention),
  while an $(n+1)$-address is a sequence of $n$-addresses. Such sequences are
  enclosed by brackets. The address $[]$, namely the empty
  word, is in $\AA_n$ for all $n \geq 0$. 
  Here are examples:
  \[
    [] \in \AA_1, \qquad [*** \,*] \in \AA_1, \qquad [[][*][]] \in \AA_2, \qquad
    [[[[*]]]] \in \AA_4 .
  \]
  For any opetope $\omega$, its nodes can be specified uniquely using higher
  addresses. We will say that the $0$-opetope $\optZero$ has no nodes, and that
  the address of the unique node\footnote{Recall that although $\optOne$ is not
    a tree, with slight abuse we see it as a unary corolla (\cref{exa:opetopes}).} of the
  $1$-opetope $\optOne$ is $*\in\AA_0$.

  Let $n \geq 2$, and assume by induction that that for all $0<k < n$ and any
  $k$-opetope $\psi$, the nodes of $\psi$ are assigned $(k-1)$-addresses, namely
  that we have an injective map $\addr : \psi^\nodesymbol \to \AA_{k-1}$. Now,
  an $n$-opetope $\omega$ is a $\optPolyFun^{n-2}$-tree $\omega\colon
  \underlyingtree{\omega} \to \optPolyFun^{n-2}$. Write its underlying tree
  $\underlyingtree{\omega}$ as
  $\polynomialinline{I_\omega}{E_\omega}{B_\omega}{I_\omega}$. A node $b \in
  B_\omega$ has an address $\addr b$, which is a list of edges of
  $\underlyingtree\omega$ describing the path from the root to $b$. Write this
  address as $[e_1 \cdots e_m]$, where $e_1, \ldots, e_m \in E_\omega$. The edge
  $e_1$ is an input edge of the root node $b_1$ of $\underlyingtree\omega$, and
  so it corresponds to a node $c_1$ of the $(n-1)$-opetope $\psi_1$ decorating
  $b_1$. By induction, $c_1$ has a higher address $[q_1]\in\AA_{n-2}$. Likewise,
  $e_2$ is an input edge of a node $b_2$ whose target is $e_1$. If $\psi_2$ is
  the $(n-1)$-opetope decorating $b_2$, then $e_2$ corresponds to a node $c_2$
  of $\psi_2$, that has an $(n-2)$-address $[q_2] \eqdef \addr b_2 \in
  \AA_{n-2}$. Repeating the argument, each $e_i$ in the list $[e_1 \cdots e_m]$
  gives rise to an $(n-2)$-address $[q_i]$, and so we assign the higher address
  $[[q_1]\cdots[q_m]]$ to $b$.

  The crux of the construction is that it allows us to specify, explicitly, every
  sub-opetope of an opetope (of strictly higher dimension) via its (the sub-opetope's)
  address. Henceforth, we use $\omega^\nodesymbol$ to denote the set of
  \emph{higher} addresses of the nodes of $\omega$, and likewise for
  $\omega^\leafsymbol$. Moreover, if $[p] \in \omega^\nodesymbol$ is a node
  higher address of $\omega$, then $\src_{[p]} \omega$ will now denote the
  $(n-1)$-opetope that is the decoration of the node at $[p]$. Let $[l] = [p[q]]
  \in \AA_{n-1}$ be an address such that $[p] \in \omega^\nodesymbol$ and
  $[q]\in(\src_{[p]}\omega)^\nodesymbol$. Then as a shorthand, we write the
  corresponding edge address of $\omega$ as
  \begin{equation}
    \label{eq:preliminaries:edg}
    \edg_{[l]} \omega \eqdef \src_{[q]} \src_{[p]} \omega .
  \end{equation}
\end{para}

\begin{example}
  \label{ex:higher-addresses:3}
  Consider the $2$-opetope on the left, called $\optInt{3}$:
  \[
    \tikzinput{opetope-2-graphical}{3} \qquad\qquad
    \tikzinput[.8]{opetope-2-tree}{3,3.address}
  \]
  Its underlying pasting diagram consists of $3$ arrows $\optOne$ grafted
  linearly. Since the only node address of $\optOne$ is $* \in \AA_0$, the
  underlying tree of $\optInt{3}$ can be depicted as on the right. On the left
  of this tree are the decorations: nodes are decorated with $\optOne \in
  \OO_1$, while the edges are decorated with $\optZero \in \OO_0$. For each node
  in the tree, the set of input edges of that node is in bijective
  correspondence with the node addresses of the decorating opetope, and this
  address is written on the right of each edge. In this low dimensional example,
  these addresses can only be $*$. Finally, on the right of each node of the
  tree is its $1$-address, which is just a sequence of $0$-addresses giving
  ``walking instructions'' to get from the root to that node.

  The $2$-opetope $\optInt{3}$ can then be seen as a corolla in some $3$-opetope
  as follows:
  \[
    \tikzinput{opetope-2-tree}{3.+,3.+.address}
  \]
  As previously mentioned, the set of input edges is in bijective correspondence
  with the set of node addresses of $\optInt{3}$. Here is now an example of a
  $3$-opetope, with its annotated underlying tree on the right (the $2$-opetopes
  $\optInt{1}$ and $\optInt{2}$ are analogous to $\optInt{3}$):
  \[
    \tikzinput{opetope-3-graphical}{ex1}
    \tikzinput{opetope-3-tree}{ex1,ex1.address}
  \]
\end{example}

\begin{remark}
  \label{prop:opetope-target}
  Let $\omega \in \OO_n$ with $n \geq 2$. We describe its target $\tgt \omega
  \in \OO_{n-1}$ and readdressing function $\readdress_\omega \colon
  \omega^\leafsymbol \to (\tgt \omega)^\nodesymbol$ in terms of higher
  addresses.
  \begin{enumerate}
  \item If $\omega$ is degenerate, namely $\omega = \itree{\phi}$ for some $\phi
    \in \OO_{n-2}$, then $\tgt \omega = \ytree{\phi}$ and $\readdress_\omega :
    \omega^\leafsymbol = \left\{ [] \right\} \to \ytree{\phi}^\nodesymbol =
    \left\{ [] \right\}$ obviously maps the unique edge $[]$ to the unique node
    $[]$.

  \item If $\omega = \ytree{\psi}$ (a corolla) for some $\psi \in \OO_{n-1}$,
    then $\tgt \omega = \psi$. Further, $\omega^\leafsymbol = \left\{ [[q]] \mid
      [q] \in \psi^\nodesymbol \right\}$, and $\readdress_\omega$ maps $[[q]]$
    to $[q]$.

  \item Otherwise, $\omega$ decomposes as the grafting $\omega = \nu
    \graft_{[l]} \ytree{\psi}$ of a corolla onto a tree, for some $\nu \in
    \OO_n$, $\psi \in \OO_{n-1}$, and $[l] \in \nu^\leafsymbol$, and its target
    is the substitution
    \[
      \tgt \omega = (\tgt \nu) \subst_{\readdress_\nu [l]} \psi .
    \]
    The readdressing function $\readdress_\omega : \omega^\leafsymbol \to (\tgt
    \omega)^\nodesymbol$ is given as follows. Let $[j] \in \omega^\leafsymbol$.
    \begin{enumerate}
    \item If $[l] \sqsubseteq [j]$ ($[l]$ is a prefix of $j$), then $[j] =
      [l[q]]$ for some $[q] \in \psi^\nodesymbol$, and $\readdress_\omega [l[q]]
      = (\readdress_\nu [l]) \cdot [q]$ (where $\cdot$ is concatenation).
    \item If $[l] \not\sqsubseteq [j]$, then $[j] \in \nu^\leafsymbol$. Assume
      $\readdress_\nu [l] \sqsubseteq \readdress_\nu [j]$. Then $\readdress_\nu
      [j] = (\readdress_\nu [l]) \cdot [[q]] \cdot [a]$, for some $[q] \in
      (\src_{\readdress_\nu [l]} \tgt \nu)^\nodesymbol = (\tgt
      \psi)^\nodesymbol$, and let $\readdress_\omega [j] = (\readdress_\nu [l])
      \cdot (\readdress_\psi^{-1} [q]) \cdot [a]$.
    \item If $\readdress_\nu [l] \not\sqsubseteq \readdress_\nu [j]$, then
      $\readdress_\omega [j] = \readdress_\nu [j]$.
    \end{enumerate}
  \end{enumerate}
\end{remark}

\begin{para} [The category of opetopes]
  \label{sec:opetopes:o}
  We are now equipped to give a definition of the category $\OO$ of opetopes in
  terms of generators and relations. While the category of opetopes was first
  introduced in \cite{Cheng2003}, the definition we use is due to
  \cite{hothanh18}.
\end{para}
\begin{proposition}
  [Opetopic identities, {\cite[Thm 4.1]{hothanh18}}]
  \label{th:opetopic-identities}
  Let $\omega \in \OO_n$ with $n \geq 2$.
  \begin{enumerate}
  \item (Inner edge) For $[p[q]] \in \omega^\nodesymbol$ (forcing $\omega$ to be
    non-degenerate), we have $\tgt \src_{[p[q]]} \omega = \src_{[q]} \src_{[p]}
    \omega$.
  \item (Globularity 1) If $\omega$ is non-degenerate, we have $\tgt \src_{[]}
    \omega = \tgt \tgt \omega$.
  \item (Globularity 2) If $\omega$ is non-degenerate, and $[p[q]] \in
    \omega^\leafsymbol$, we have $\src_{[q]} \src_{[p]} \omega =
    \src_{\readdress_\omega [p[q]]} \tgt \omega$.
  \item (Degeneracy) If $\omega$ is degenerate, we have $\src_{[]} \tgt \omega =
    \tgt \tgt \omega$.
  \end{enumerate}
\end{proposition}

\begin{definition}
  [{\cite[4.2]{hothanh18}, \cite[Def. 3.3.2]{HTLS2019opetopic}}]
  \label{def:o}
  With these identities at hand, we define the \defn{category $\OO$ of opetopes}
  by generators and relations as follows.
  \begin{enumerate}
  \item (Objects) We set $\ob \OO = \sum_{n \in \NN} \OO_n$.
  \item (Generating morphisms) Let $\omega \in \OO_n$ with $n \geq 1$. We define
    a morphism $\tgt\colon\tgt \omega\to \omega$, called the \emph{target
      map}\index{target map}. If $[p] \in \omega^\nodesymbol$, we define a
    morphism $\src_{[p]}\colon\src_{[p]} \omega \to \omega$, called a
    \emph{source map}\index{source map}. A \emph{face map}\index{face map} is
    either a source or the target map.
  \item (Relations) We impose that the following squares, that are well defined
    thanks to \cref{th:opetopic-identities}, commute. Let $\omega \in \OO_n$
    with $n \geq 2$.
    \begin{enumerate}
    \item \condition{Inner}\index{inner@\condition{Inner}} for $[p[q]] \in
      \omega^\nodesymbol$ (forcing $\omega$ to be non-degenerate):
      \[
        \squarediagram {\src_{[q]} \src_{[p]} \omega}{\src_{[p]} \omega}
        {\src_{[p[q]]}\omega}{\omega}
        {\src_{[q]}}{\tgt}{\src_{[p]}}{\src_{[p[q]]}}
      \]
    \item \condition{Glob1}\index{glob1@\condition{Glob1}} for every $\omega$
      that is non-degenerate:
      \[
        \squarediagram {\tgt \tgt \omega}{\tgt \omega}{\src_{[]} \omega}{\omega
          .} {\tgt}{\tgt}{\tgt}{\src_{[]}}
      \]
    \item \condition{Glob2}\index{glob2@\condition{Glob2}} if $\omega$ is
      non-degenerate, then for $[p [q]] \in \omega^\leafsymbol$:
      \[
        \squarediagram {\src_{\readdress_\omega [p [q]]} \tgt \omega}{\tgt
          \omega} {\src_{[p]} \omega}{\omega .} {\src_{\readdress_\omega
            [p[q]]}}{\src_{[q]}}{\tgt} {\src_{[p]}}
      \]
    \item \condition{Degen}\index{degen@\condition{Degen}} if $\omega$ is
      degenerate:
      \[
        \squarediagram {\tgt \tgt \omega}{\tgt \omega}{\tgt \omega}{\omega .}
        {\tgt}{\src_{[]}}{\tgt}{\tgt}
      \]
    \end{enumerate}
  \end{enumerate}
\end{definition}

\begin{remark}
  \cite[4.2]{hothanh18} has pictorial examples of these relations.
\end{remark}

\begin{proposition}
  \label{lem:O-lfd-cat}
  $\OO$ is a \lfd~category.
\end{proposition}
\begin{proof}
  $\OO$ is direct since every generating morphism strictly increases dimension.
  By~\cref{prop:loc-finite-finite-cover}, it is locally finite since for every
  $\omega\in\OO$, the source maps and target map of $\omega$ form a finite saturated cover
  of non-identity morphisms.
\end{proof}

\begin{remark}
  We have $\omega\in\OO_n$ if and only if $\height(\omega)=n$
  (see \cref{def:loc-fin-cat-dimension}).
\end{remark}

\begin{notation}
  \label{not:o}
  For $n \in \NN$, let $\OO_{\leq n}$ be the full subcategory of $\OO$ of
  opetopes of dimension $\leq n$. The full subcategories $\OO_{< n}$, $\OO_{\geq
    n}$, $\OO_{> n}$ are defined similarly.
\end{notation}
  
\begin{para}[Opetopic sets]
  \label{sec:opetopes:opetopic-sets}
  The category of \defn{opetopic sets} is the presheaf category $\pshO$.
  Following \cref{sec:locally-finite-direct-cats}, for $X \in \pshO$ and $\omega
  \in \OO$, we will refer to the elements of the set $X_\omega$ as the
  \emph{cells}\index{cell} of $X$ of \emph{shape} $\omega$. The representable
  presheaf on $\omega \in \OO$ is denoted $O[\omega]$. 
\end{para}

\begin{definition}
 \label{def:opetopic-spine-boundary} 
  As
  in~\cref{para:boudaries-finite-cell-cplxs}, the \defn{boundary} $\partial
  O[\omega]$ of $\omega$ is the subpresheaf of $O[\omega]$ not containing the
  identity morphism, we write $\delta_\omega\colon \partial O[\omega]\subto
  O[\omega]$ for the boundary inclusion, and the set of all boundary inclusions is
  denoted $\sfB$.

  The \defn{spine}\index{spine} $S [\omega]$ of $\omega$ is the maximal
  subpresheaf of $\partial O [\omega]$ not containing the cell $\tgt\colon \tgt
  \omega\to \omega$, and we write $S_\omega \colon S [\omega] \subto O [\omega]$
  for the \emph{spine inclusion} of $\omega$. The set of spine inclusions is
  denoted $\sfS$.
\end{definition}

\begin{notation}
 \label{notn:empty-incl-O} 
  For every $\omega \in \OO$ we note $\mathsf{o}_\omega\colon \emptyset\to
  O[\omega]$ and $\sfO\eqdef\{\mathsf{o}_\omega\mid \omega \in\OO\}$.
\end{notation}

\begin{lemma}
  \label{lemma:opetopes-technical:boundary-pushout}
  For $\omega$ in $\OO$ with $\height(\omega) \geq 1$, the following square in
  $\pshO$ is cartesian and cocartesian\footnote{Recall that this follows from
    van-Kampen-ness of pushouts of monomorphisms
    in a topos.}, where all arrows are canonical
  inclusions.
  \[
    \squarediagram {\partial O [\tgt \omega]}{S [\omega]}{O [\tgt \omega]}
    {\partial O [\omega]} {}{\delta_{\tgt\omega}}{S_\omega}{}
  \]
\end{lemma}
\begin{proof}
  Seen as subpresheaves of $O[\omega]$, the intersection of $S[\omega]$ and
  $O[\tgt\omega]$ is precisely $\partial O[\tgt\omega]$ (namely, the inclusions
  of the leaves and the root edge of $\omega$). The union of $S[\omega]$ and
  $O[\tgt\omega]$ is $\partial O[\omega]$ (by definition of $S[\omega]$).
\end{proof}

\begin{lemma}
  \label{lemma:opetopes-technical:spine-pushout}
  Let $n \geq 1$, $\omega \in \OO_n$, $[l] \in \omega^\leafsymbol$, and $\psi
  \in \OO_{n-1}$ be such that $\edg_{[l]} \omega = \tgt \psi$, so that the
  grafting $\omega \graft_{[l]} \ytree{\psi}$ is in $\OO_n$. Then the following
  square of inclusions is cocartesian and cartesian in $\pshO$.
  \[
    \squarediagram {O [\edg_{[l]} \omega]}{O [\psi]}{S [\omega]} {S [ \omega
      \graft_{[l]} \ytree{\psi}]} {\tgt}{\src_{[l]}}{\src_{[l]}}{}
  \]
\end{lemma}
\begin{proof}
  Exactly similar to the previous lemma.
\end{proof}

\begin{notation}
  Let $F : \OO \to \ob{\cC\arr}$ be a function that maps opetopes to morphisms
  in some category $\cC$, and $\sfM$ the set of maps defined by $\sfM \eqdef
  \left \{ F (\omega) \mid \omega \in \OO \right \}$. Then for $n \in \NN$, we
  define $\sfM_{\geq n} \eqdef \left \{ F (\omega) \mid \omega \in \OO_{\geq n}
  \right \}$, and similarly for $\sfM_{> n}$, $\sfM_{\leq n}$, $\sfM_{< n}$, and
  $\sfM_{= n}$. For convenience, the latter is abbreviated $\sfM_n$. If $m \leq
  n$, we also let $\sfM_{m, n} \eqdef \sfM_{\geq m} \cap \sfM_{\leq n}$. By
  convention, $\sfM_{\leq n} = \emptyset$ if $n < 0$. For example, $\sfS_{\geq 2} =
  \left \{ S_\omega \mid \omega \in \OO_{\geq 2} \right \}$, and $\sfS_{n, n+1} =
  \sfS_n \cup \sfS_{n+1}$.
\end{notation}

\begin{proposition}
  \label{lemma:opetopes-technical:h-lift}
  Let $X \in \pshO$ be an opetopic set.
  \begin{enumerate}
  \item If $\sfS_{n,n+1} \perp X$, then $\sfB_{n+1} \perp X$. Equivalently,
    $\sfB_{n+1}$ is a subset of the $\sfS_{n,n+1}$\=/local isomorphisms
    ($\sfB_{n+1}\subset ({}^\bot(\sfS_{n,n+1}^\bot))^{3/2}$ in the notation of
    \cref{sec:localisations}). Thus every morphism in $\sfB_{\geq n+1}$ is an
    $\sfS_{\geq n}$-local isomorphism.
  \item If $\sfS_{n, n+1} \perp X$ and $\sfB_{n+2} \perp X$, then $\sfS_{n+2} \perp
    X$. Thus if $\sfS_{n, n+1} \perp X$ and $\sfB_{\geq n+2} \perp X$, then
    $\sfS_{\geq n} \perp X$.
  \end{enumerate}
\end{proposition}
\begin{proof}
  \begin{enumerate}
  \item Let $\omega \in \OO_{n+1}$. Note that the following triangle commutes.
    \[
      \triangleURdiagram {S [\omega]}{\partial O [\omega]}{O [\omega] .}
      {S_\omega}{i}{\delta_\omega}
    \]
    By ``2~out~of~3'' for $\sfS_{n,n+1}$-local isomorphisms, it is enough to show
    that $i$ is one. Suppose that $\sfS_{n,n+1} \perp X$. Take a morphism $f : S
    [\omega] \to X$. The existence of a lift $\partial O [\omega] \to X$ follows
    from the existence of a lift $O [\omega] \to X$, since $S_\omega \perp X$.
    For uniqueness, consider two lifts $g, h : \partial O [\omega] \to X$ of
    $f$. By \cref{lemma:opetopes-technical:boundary-pushout}, in order to show
    that they are equal, it suffices to show that they coincide on $O [\tgt
    \omega]$, since they coincide on $S [\omega]$ (as they extend $f$). But
    since they coincide on $S [\omega]$, they must coincide on the subpresheaf
    $S [\tgt \omega] \subseteq S [\omega]$. Since $\sfS_n \perp X$, $g$ and $h$
    coincide on $O [\tgt \omega]$, and are thus equal.

  \item Let $\omega \in \OO_{n+2}$ and $f : S [\omega] \to X$. By assumption,
    the restriction $f|_{S [\tgt \omega]}$ of $f$ to $S [\tgt \omega]$ extends
    to a unique $g : O [\tgt \omega] \to X$. We now show that the following
    square commutes:
    \[
      \squarediagram {\partial O [\tgt \omega]}{S [\omega]}{O [\tgt \omega]}{X
        .} {}{}{f}{g}
    \]
    By \cref{lemma:opetopes-technical:boundary-pushout}, it suffices to show
    that $f$ and $g$ coincide on $S [\tgt \omega]$ and on $O [\tgt \tgt
    \omega]$. The former is tautological, and the latter follows from the
    hypothesis that $S_{\tgt \tgt \omega} \perp X$ and that $f$ and $g$ coincide
    on $S [\tgt \tgt \omega] \subseteq S [\tgt \omega]$. Therefore, the square
    above commutes, and by \cref{lemma:opetopes-technical:boundary-pushout}
    again, $f$ and $g$ extend to a morphism $h : \partial O [\omega] \to X$,
    which in turn extends to a morphism $i : O [\omega] \to X$, since by
    assumption $\sfB_{n+2} \perp X$.

    For uniqueness, consider two lifts $i, i' : O [\omega] \to X$ of $f$. By
    \cref{lemma:opetopes-technical:boundary-pushout}, they are equal if and only
    if their restriction $g, g' : O [\tgt \omega] \to X$ are equal. Since $g|_{S
      [\tgt \omega]} = f|_{S [\tgt \omega]} = g'|_{S [\tgt \omega]}$, and since
    by assumption $\sfS_{n+1} \perp X$, we have $g = g'$, and thus $i = i'$.
    \qedhere
  \end{enumerate}
\end{proof}

\begin{corollary}
  \label{coroll:opetopes-technical:h-lift}
  Let $X$ be an opetopic set such that $\sfS_{n, n+1} \perp X$. Then $\sfS_{\geq n}
  \perp X$ if and only if $\sfB_{\geq n+2} \perp X$.
\end{corollary}

\begin{proposition}
  \label{lemma:opetopes-technical:spines}
  Let $n \in \NN$, and $\omega \in \OO_{n+2}$. Then the inclusion $S[\tgt\omega]
  \subto S[\omega]$ is a relative $\sfS_{n+1}$-cell complex.
\end{proposition}
\begin{proof}
  We show that the morphism $S[\tgt\omega] \subto S[\omega]$ is a composite of
  pushouts of elements of $\sfS_{n+1}$. If $\omega$ is degenerate, say $\omega =
  \itree{\phi}$ for some $\phi \in \OO_n$, then $S [\tgt \omega] = S
  [\ytree{\phi}] = O [\phi] = S [\omega]$, so the result trivially holds.

  Assume that $\omega$ is not degenerate, let $X^{(0)} \eqdef S[\tgt\omega]$,
  and $[p_1] \succ \cdots \succ [p_k]$ be the node addresses of $\omega$, sorted
  in reverse lexicographical order. By induction, assume that $X^{(i-1)}$ is a
  subpresheaf of $S [\omega]$ containing the $(n+1)$-cells $\src_{[p_1]} \omega,
  \ldots, \src_{[p_{i-1}]} \omega \in S [\omega]$. Clearly, this holds when $i =
  1$, as $S [\tgt \omega]$ does not contain any $(n+1)$-cell.

  Take $[q] \in (\src_{[p_i]} \omega)^\nodesymbol$. By induction, and since
  $[p_i [q]] \succ [p_i]$, the $(n+1)$-cell $\src_{[p_i[q]]} \omega$ is in
  $X^{(i-1)}$. Further, the $n$-cell $\src_{[q]} \src_{[p_i]} \omega$ is present
  in $X^{(i-1)}$, since by \condition{Inner}, $\src_{[q]} \src_{[p_i]} \omega =
  \tgt \src_{[p_i[q]]} \omega$. Therefore, we have an inclusion $u_i : S
  [\src_{[p_i]} \omega] \to X^{(i-1)}$ mapping $\src_{[q]} \src_{[p_i]} \omega$
  to $\src_{[q]} \src_{[p_i]} \omega$, and let $X^{(i)}$ be the pushout
  \[
    \pushoutdiagram {S [\src_{[p_i]} \omega]}{X^{(i-1)}}{O[\src_{[p_i]} \omega]}
    {X^{(i)}} {u_i}{S_{\src_{[p_i]} \omega}}{}{}{}
  \]
  Clearly, $X^{(i)}$ is a subpresheaf of $S [\omega]$ containing the
  $(n+1)$-cell $\src_{[p_j]} \omega$ for $1 \leq j \leq i$, and the induction
  hypothesis is satisfied.

  Finally, $X^{(k)} \subseteq S [\omega]$ contains all the $(n+1)$-cells of $S
  [\omega]$, whence $X^{(k)} = S [\omega]$. By construction, the chain of
  inclusions $S [\tgt \omega] = X^{(0)} \subto X^{(1)} \subto \cdots \subto
  X^{(k)} = S [\omega]$ is a relative $\sfS_{n+1}$-cell complex.
\end{proof}

\begin{corollary}
  \label{lemma:comparison:tgt-spine-local}
  Let $n \in \NN$, and $\omega \in \OO_{n+2}$. Then the target map $\tgt \omega
  \to \omega$ of $\omega$ is an $\sfS_{n+1, n+2}$-local isomorphism.
\end{corollary}
\begin{proof}
  In the square below
  \[
    \squarediagram {S [\tgt \omega]} {O [\tgt \omega]} {S [\omega]} {O [\omega]}
    {S_{\tgt \omega}} {r} {\tgt} {S_\omega}
  \]
  the map $r$ is an $\sfS_{n+1}$-local isomorphism by
  \cref{lemma:opetopes-technical:spines}, and the horizontal maps are in
  $\sfS_{n+1, n+2}$. The result follows by ``2~out~of~3''. \qedhere
\end{proof}

\begin{corollary}
  \label{coroll:y-i:spine-local}
  Let $\omega \in \OO_n$.
  \begin{enumerate}
  \item $\tgt \tgt = \src_{[]} \tgt : \omega \to \itree{\omega}$ is in
    $\sfS_{n+2}$.
  \item The morphisms $\src_{[]}, \tgt : \omega \to \ytree{\omega}$ are $\sfS_{n+1,
      n+2}$-local isomorphisms.
  \end{enumerate}
\end{corollary}
\begin{proof}
  \begin{enumerate}
  \item The map $\tgt \tgt = \src_{[]} \tgt : \omega \to \itree{\omega}$ is
    precisely the spine inclusion $S_{\itree{\omega}}$ of the degenerate
    $(n+2)$-opetope $\itree{\omega}$.
  \item The source map $\src_{[]} : \omega \to \ytree{\omega}$ is precisely the
    spine inclusion $S_{\ytree{\omega}}$ of the $(n+1)$-opetope
    $\ytree{\omega}$. The target map $\tgt : \omega \to \ytree{\omega}$ is the
    morphism $\tgt : \tgt\tgt\itree{\omega} \to \tgt\itree{\omega}$ and is the
    vertical arrow in the diagram below.
    \[
      \triangleDLdiagram {\omega = S [\itree{\omega}]}{\ytree{\omega} = \tgt
        \itree{\omega}} {\itree{\omega}.} {\tgt}{S_{\itree{\omega}}}{\tgt}
    \]
    The horizontal arrow is an $\sfS_{n+1, n+2}$-local isomorphism by
    \cref{lemma:comparison:tgt-spine-local} and the diagonal arrow is in
    $\sfS_{n+2}$ by point (1). The result follows by ``2~out~of~3''.
    \qedhere
  \end{enumerate}
\end{proof}



\begin{notation}
  \label{not:subcategories-of-o}
  Let $m \in \NN$ and $n \in \NN \cup \{ \infty \}$ be such that $m \leq n$, and
  let $\OO_{m, n}$ be the full subcategory of $\OO$ consisting of opetopes
  $\omega$ of dimension $m \leq \height(\omega) \leq n$. Note that $\OO_{m,
    \infty} = \OO_{\geq m}$.
\end{notation}

\begin{definition}
  \label{def:truncation}
    
  The restriction functor $(-)_{m, n}\eqdef (\iota^{\geq m})^* : \psh{\OO_{\geq
      m}} \to \psh{\OO_{m, n}}$ associated to the inclusion $\iota^{\geq m} :
  \OO_{m, n} \to \OO_{\geq m}$ is called \defn{truncation}\index{truncation}. We
  write the left adjoint as $\iota^{\geq m}_!$ and the right adjoint as
  $\iota^{\geq m}_*$. Note that for $X \in \psh{\OO_{m, n}}$, the presheaf
  $\iota^{\geq m}_! X$ is the ``extension by $0$'', i.e. $(\iota^{\geq m}_!
  X)_{m, n} = X$, and $(\iota^{\geq m}_! X)_\psi = \emptyset$ for all $\psi \in
  \OO_{> n}$. On the other hand, $\iota^{\geq m}_* X$ is the ``canonical
  extension'' of $X$ into a presheaf over $\OO_{\geq m}$: we have $(\iota^{\geq
    m}_* X)_{m, n} = X$, and $\sfB_{> n} \perp \iota^{\geq m}_* X$, which
  uniquely determines $\iota^{\geq m}_* X$.

  Likewise, the inclusion $\iota^{\leq n} : \OO_{m, n} \to \OO_{\leq n}$ induces
  a restriction functor $\psh{\OO_{\leq n}} \to \psh{\OO_{m, n}}$, also denoted
  by $(-)_{m, n}$ and again called \emph{truncation}\index{truncation}, that has
  both a left adjoint $\iota^{\leq n}_!$ and a right adjoint $\iota^{\leq n}_*$.
  Explicitly, for $X \in \psh{\OO_{m, n}}$, the presheaf $\iota^{\leq n}_! X$ is
  the ``canonical extension'' of $X$ into a presheaf over $\OO_{\leq n}$:
  \[
    \iota^{\leq n}_! X = \colim_{O[\psi]_{m, n} \rightarrow X} O [\psi] .
  \]
  On the other hand, $\iota^{\leq n}_* X$ is the ``terminal extension'' of $X$
  in that $(\iota^{\leq n}_* X)_{m, n} = X$, and $(\iota^{\leq n}_* X)_\psi$ is
  a singleton, for all $\psi \in \OO_{< m}$.


  For $n < \infty$, we write $(-)_{\leq n}$ for $(-)_{0, n} : \psh{\OO_{\geq 0}}
  = \pshO \to \psh{\OO_{0, n}} = \psh{\OO_{\leq n}}$, and let $(-)_{< n} =
  (-)_{\leq n-1}$ if $n \geq 0$. Similarly, we note $(-)_{m, \infty} :
  \psh{\OO_{\leq \infty}} = \pshO \to \psh{\OO_{m, \infty}} = \psh{\OO_{\geq
      m}}$ by $(-)_{\geq m}$, and let $(-)_{> m} = (-)_{\geq m+1}$.
\end{definition}

\begin{proposition}
  \label{prop:properties-of-iotas}
  \begin{enumerate}
  \item The functors $\iota^{\geq m}_!$, $\iota^{\geq m}_*$, $\iota^{\leq n}_!$,
    and $\iota^{\leq n}_*$ are fully faithful.
  \item A presheaf $X \in \psh{\OO_{\geq m}}$ is in the essential image of
    $\iota^{\geq m}_!$ if and only if $X_{> n} = \emptyset$.
  \item A presheaf $X \in \psh{\OO_{\geq m}}$ is in the essential image of
    $\iota^{\geq m}_*$ if and only if for all $\omega \in \OO_{> n}$ we have
    $(\delta_\omega)_{\geq m} \perp X$.
  \item A presheaf $X \in \psh{\OO_{\leq n}}$ is in the essential image of
    $\iota^{\leq n}_*$ if and only if for all $\omega \in \OO_{< m}$ we have
    $(\mathsf{o}_\omega)_{\leq n} \perp X$, i.e. $X_\omega$ is a singleton.
  \end{enumerate}
\end{proposition}
\begin{proof}
  The first point follows from the fact that $\iota^{\geq m}$ and $\iota^{\leq
    n}$ are fully faithful. The rest are straightforward verifications.
\end{proof}

\begin{notation}
  To ease notations, we sometimes leave truncations implicit, e.g. point (3) of
  the previous proposition can be reworded as: a presheaf $X \in \psh{\OO_{\geq
      m}}$ is in the essential image of $\iota^{\geq m}_*$ if and only if
  $\sfB_{> n} \perp X$.
\end{notation}




\section{Opetopic algebras}
\label{sec:opetopic-algebras}
Let $k \leq n \in \NN$, and recall that $\OO_{n-k, n} \subto \OO$ is the full
subcategory of those opetopes $\omega$ such that $n - k \leq \dim \omega \leq
n$. A \emph{$k$-coloured, $n$-dimensional opetopic algebra}, or $(k,
n)$-opetopic algebra, will be an algebraic structure on a presheaf over
$\OO_{n-k, n}$, whose cells of dimension $n$ are ``operations'' that can be
``composed'' in ways encoded by $(n+1)$-cells\footnote{Recall that an
  $(n+1)$-opetope is precisely a pasting diagram of $n$-opetopes.}. As we will
see in the next \cref{sec:opet-nerve-functors}, the fact that the operations and
relations of a $(k, n)$-opetopic algebra are encoded by opetopes of dimension
$\geq n$ results in the category $\Oalg k n$ of $(k, n)$-opetopic algebras
always having a canonical fully faithful \emph{nerve functor} to the category
$\psh{\OO}$ of opetopic sets (\cref{th:nerve-theorem-O}).

We begin this section by surveying elements of the theory of \emph{parametric
  right adjoint (p.r.a.) monads}. This will be essential to the definition of
the \emph{coloured} $\optPolyFun^n$ monad, which is a generalisation of the
polynomial monad $\optPolyFun^n\colon \Set\slice{\OO_n}\to \Set\slice{\OO_n}$ to
the presheaf category $\psh{\OO_{n-k, n}}$. The category $\oAlg_{k,n}$ algebras
of this new monad is the category of $(k, n)$-opetopic algebras. Then, we
introduce the category $\bbLambda$ of \emph{opetopic shapes}. We 
obtain a reflective adjunction to the category $\oAlg_{k,n} \eqdef
\optPolyFun^n\alg$
\[
  \tau : \psh{\bbLambda} \localisation \oAlg_{k,n} : N,
\]
where the left adjoint is called the \emph{algebraic realisation}, and where the
right adjoint is the \emph{algebraic nerve functor}. Finally, we describe the
previous adjunctions as a Gabriel-Ulmer localisation at \emph{spine} inclusions
(\cref{th:nerve-theorem-lambda}).

Each of the monads $\optPolyFun^n $ on $\psh{\OO_{n-k,n}}$ will be finitary. We
can therefore conclude that each category $\oAlg_{k,n}$ is the category of
algebras of an $\OO_{n-k,n}$\=/sorted theory (\cref{def:C-sorted-theory}). The
cases of interest (and in fact, the \emph{only} real cases of interest) are
$\oAlg_{1,1}$, $\oAlg_{1,2}$ and $\oAlg_{1,3}$, which are the categories $\Cat,
\Opd\pl$ and $\Cmbd_{\mathrm{col}}$ of small categories, planar coloured operads and a coloured
version of Loday's combinads over planar trees \cite{Loday2012a} respectively.

\begin{para}
  [Parametric right adjoint monads]
  \label{para:pra-monads}
  We survey elements of the theory of parametric right adjoint (p.r.a.) monads
  on presheaf categories, which will be essential to the definition and
  description of $(k,n)$-opetopic algebras. A comprehensive treatment of this
  theory can be found in \cite{weber2007familial}.
\end{para}

\begin{definition}
  \label{def:pra}
  If $T : \cC \to \cD$ is a functor, and $\cC$ has a terminal object $1$, then
  $T$ factors as
  \begin{equation}
    \label{eq:pra}
    \begin{tikzcd}
      \cC \simeq \cC \slice 1 \ar[r,"T_1"] &\cD \slice{ T 1} \ar[r] &\cD ,
    \end{tikzcd}
  \end{equation}
  where the second functor is the induced functor between slice categories, and
  the third is the domain functor. We say that $T$ is a \defn{parametric right
    adjoint} (abbreviated p.r.a.) if $T_1$ has a left adjoint $E$.
\end{definition}

\begin{para}
  \label{rem:pra}
  We will immediately restrict \cref{def:pra} to the case where $\cC = \cD =
  \psh A$ for a small category $ A$. If $ A / T 1$ is the category of elements
  of $T 1 \in \psh A$, and using the equivalence $\psh{ A / T 1} \simeq \psh A
  \slice{ T 1}$, the factorisation of \eqref{eq:pra} becomes
  \begin{equation*}
    \begin{tikzcd}
      \psh A
      \ar[r,"T_1"] &\psh{ A / T1}
      \ar[r] &\psh A .
    \end{tikzcd}
  \end{equation*}
  Let $E$ be the left adjoint of $T_1$. Then $T_1$ is the nerve of the
  restriction $E : A / T1 \to \psh A$ of $E$ to the representable presheaves,
  and the usual nerve formula gives
  \[
    (T_1 X)_x = \psh A (E x, X) ,
  \]
  where $X \in \psh A$ and $x \in A / T1$. Therefore, for $a \in A$, we have
  \begin{equation*}
    (T X)_a = \sum_{x \in (T 1)_a} \psh A (E x, X).
  \end{equation*}
  Whence it is clear that the data of the object $T1 \in \psh A$ and of the
  functor $E : A / T1 \to \psh A$ determine the functor $T$ up to isomorphism.
  Let $\Theta_0$ \index{$\Theta_0$|see {$T$-cardinal}} (leaving $T$ implicit) be
  the full subcategory of $\psh A$ that is the image of $E : A / T1 \to \psh A$.
  Objects of $\Theta_0$ are called \defn{$T$-cardinals}.
\end{para}

\begin{definition}
  \label{def:pra-monad}
  A \defn{p.r.a. monad} \index{p.r.a. monad} is a monad $T$ whose endofunctor is
  a p.r.a. and whose unit $\id \to T$ and multiplication $TT \to T$ are
  cartesian natural transformations.
\end{definition}

\begin{remark}
  Any p.r.a. monad $T$ on a presheaf category is a \emph{monad with $\Theta_0$
    as arities} (that we have seen in \cref{para:cellular-arities}) and so we
  can deduce a lot of information about the free-forgetful adjunction $\psh A
  \rightleftarrows T\alg$ and about the category of algebras $T\alg$.
\end{remark}

\begin{notation}
  \label{not:thetaT}
  With slight abuse, let $T : \psh A \to T\alg$ be
  the free $T$-algebra functor. The (identity-on-objects, fully faithful)
  factorisation of the composite $\Theta_0 \subto \psh A \xto{T} T\alg$ will be
  denoted by \index{$\Theta_T$}
  \begin{equation}
    \label{eq:thetaT}
    \begin{tikzcd}
      \Theta_0
      \ar[r,"j_T"] &\Theta_T
      \ar[r,hook,"i_T"] &T\alg .
    \end{tikzcd}
  \end{equation}
  In other words, $\Theta_T$ is the full subcategory of $T\alg$ spanned by free
  algebras over $T$-cardinals.
\end{notation}

\begin{proposition}
  [{\cite[Prop. 4.20]{weber2007familial}}]
  \label{prop:theta0}
  Let $T : \psh A \to \psh A$ be a p.r.a. monad, with $\Theta_0$ the category of
  $T$\=/cardinals. Then the Yoneda embedding $A\subto \psh A$ factors through
  $\Theta_0\subto \psh A$.
  (In other words, representable presheaves are $T$-cardinals.)
\end{proposition}

\begin{proposition}
  \label{th:pra-nerve-theorem}
  Let $T$ be a p.r.a. monad on $\psh A$. Then $T$ has arities $\Theta_0$. Thus:
  \begin{enumerate}
  \item The functors $i_0 : \Theta_0 \to \psh A$ and $i_T : \Theta_T \to T\alg$
    are dense. Equivalently, their nerve functors $N_{i_0} : \psh A \to
    \psh{\Theta_0}$ and $N_{i_T} : T\alg \to \psh{\Theta_T}$ are fully faithful.
  \item The following diagram is an exact adjoint square\footnote{There exists a
      natural isomorphism $N_0 U_T \cong j_T^* N_T$ whose mate ${j_T}_! N_0 \to
      N_T F_T$ is invertible (satisfies the Beck-Chevalley condition).}.
    \[
      \begin{tikzcd}
        \psh A
        \ar[d, hook, "N_0" left]
        \ar[r, shift left = .4em, "F_T"] &
        T\alg
        \ar[d, hook, "N_T"]
        \ar[l, shift left = .4em, "U_T", "\perp" above] \\
        \psh{\Theta_0}
        \ar[r, shift left = .4em, "{j_T}_!"] &
        \psh{\Theta_T}
        \ar[l, shift left = .4em, "j_T^*", "\perp" above]
      \end{tikzcd}
    \]
    In particular, both squares commute up to natural isomorphism.
  \item (Segal condition) A presheaf $X \in \psh{\Theta_T}$ is in the essential
    image of $N_{i_T}$ if and only if $j_T^* X$ is in the essential image of
    $N_{i_0}$.
  \end{enumerate}
\end{proposition}
\begin{proof}
  This is essentially \cite[Prop. 4.22]{weber2007familial}. Point (2) is
  \cite[Prop. 1.9]{Berger2012}, and the Segal condition is \cite[theorem 4.10
  (2)]{weber2007familial}.
\end{proof}

\begin{corollary}
  \label{coroll:algebra-lifting}
  Let
  \[
    \sfJ_T \eqdef {j_T}_! \sfJ_A = \left\{ {j_T}_! \epsilon_\theta : {j_T}_!
      i_! i^* \theta \to {j_T}_! \theta \mid \theta \in \Theta_0 - \im i
    \right\} ,
  \]
  where $\epsilon$ is the counit of the adjunction $i_! \dashv i^*$. Then a
  presheaf $X \in \psh{\Theta_T}$ is in the essential image of $N_{i_T}$ if and
  only if $\sfJ_T \perp X$. As a consequence, the left adjoint $\psh{\Theta_T}
  \to T\alg$ of $N_{i_T}$ (namely the left Kan extension of $i_T$ along the Yoneda
  embedding) is a Gabriel-Ulmer localisation
  \[
    \psh{\Theta_T}\localisation T\alg \simeq\sfJ_T^{-1} \psh{\Theta_T}
  \]
\end{corollary}


\begin{para}
  \label{rem:uncoloured-zn}
  Recall the definition of the polynomial monad $\optPolyFun^n$ from
  \cref{para:opetopes-definition}. If $X = (X_\psi \mid \psi \in \OO_n)$ is in
  $\Set / \OO_n$, and if $\omega \in \OO_n$, then
  \[
    (\optPolyFun^n X)_\omega = \sum_{\substack{\nu \in \OO_{n+1} \\ \tgt \nu =
        \omega}} \prod_{[p] \in \nu^\nodesymbol} X_{\src_{[p]} \nu} .
  \]
  Under the equivalence $\Set / \OO_n \simeq \psh{\OO_n}$, this formula can be
  rewritten as
  \[
    (\optPolyFun^n X)_\omega = \sum_{\substack{\nu \in \OO_{n+1} \\ \tgt \nu =
        \omega}} \psh{\OO_n} (S [\nu], X) ,
  \]
  where $S [\nu]$ is the truncated spine of $\nu$.

  We will extend the polynomial monad $\optPolyFun^n$
  on $\Set / \OO_n = \psh{\OO_n}$ to a p.r.a. monad
  on $\psh{\OO_{n-k, n}}$, where $k \leq n$. This new setup will encompass
  more known examples than the uncoloured case (see \cref{prop:algebra-table}).
  For instance, recall that the polynomial monad $\optPolyFun^2$ on $\Set /
  \OO_2 \cong \Set / \NN$ is exactly the monad of planar operads. The extension
  of $\optPolyFun^2$ to $\psh{\OO_{1, 2}}$ will retrieve \emph{coloured} planar
  operads as algebras. Similarly, the polynomial monad $\optPolyFun^1$ on $\Set$
  is the free-monoid monad, which we would like to vary to obtain ``coloured
  monoids'', namely small categories.

  The first step of this construction is to define $\optPolyFun^n$ as a p.r.a.
  functor, namely an endofunctor $\optPolyFun^n$ on $\psh{\OO_{n-k, n}}$ such that
  in the sequence below, $\optPolyFun^n_1$ is a right adjoint:
  \[
    \begin{tikzcd}
      \psh{\OO_{n-k, n}}
      \ar[r,"\optPolyFun^n_1"] &\psh{\OO_{n-k, n} / \optPolyFun^n 1}
      \ar[r] &\psh{\OO_{n-k, n}} .
    \end{tikzcd}
  \]
  Following \cref{rem:pra}, it suffices to define its value $\optPolyFun^n 1$ on
  the terminal presheaf, and to specify a functor $E : \OO_{n-k, n} /
  \optPolyFun^n 1 \to \psh{\OO_{n-k, n}}$.
\end{para}
\begin{definition}
    \label{def:coloured-zn}
    Define $\optPolyFun^n 1 \in \psh{\OO_{n-k, n}}$ as
    \[
        (\optPolyFun^n1)_\psi
        \eqdef \{*\},
        \qquad\qquad
        (\optPolyFun^n1)_\omega
        \eqdef
        \left\{ \nu \in \OO_{n+1} \mid \tgt \nu = \omega \right\} ,
    \]
    where $\psi \in \OO_{n-k, n-1}$ and $\omega \in \OO_n$, along with the
    obvious restriction maps. We now define a
    functor $E : \OO_{n-k, n} / \optPolyFun^n 1 \to
    \psh{\OO_{n-k, n}}$. On objects, for $* \in (\optPolyFun^n1)_\psi$ and
    $\nu \in (\optPolyFun^n1)_\omega$, let\footnote{Note that in
    \cref{eq:coloured-zn:E}, the presheaves $O [\psi]$ and $S [\nu]$ are
    considered in $\psh{\OO_{n-k, n}}$, but the truncations are left
    implicit.}
    \begin{equation}
        \label{eq:coloured-zn:E}
        E (*) \eqdef O [\psi],
        \qquad\qquad
        E (\nu) \eqdef S [\nu] .
    \end{equation}
    On morphisms, $E$ takes face maps to the canonical inclusions.

    Thus we can define the
    functor $\optPolyFun^n_1 : \psh{\OO_{n-k, n}} \to
    \psh{\OO_{n-k, n} / \optPolyFun^n 1}$ as the nerve functor of $E$, namely
    $\optPolyFun^n_1 = N_E$. We now
    recover the endofunctor $\optPolyFun^n$ explicitly (as in
    \cref{rem:pra}): for $\psi \in \OO_{n-k, n-1}$ we have
    $(\optPolyFun^n X)_\psi \cong X_\psi$, and for $\omega \in \OO_n$, we
    recover a formula similar to the uncoloured case (\cref{rem:uncoloured-zn})
    \[
        (\optPolyFun^n X)_\omega
        \cong
        \sum_{\substack{\nu \in \OO_{n+1} \\ \tgt \nu = \omega}}
            \psh{\OO_{n-k, n}} (S [\nu], X) .
    \]
  \end{definition}

\begin{example}
  Let us unfold \cref{def:coloured-zn} in the case $n = 1$ and $k = 1$. Here,
  $\psh{\OO_{0, 1}} = \psh{\GG_1}$ is the category of directed graphs, whose
  terminal object $1$ is the graph with one vertex and a loop. The graph
  $\optPolyFun^1 1$ also has one vertex, but this time, it has an many loops as
  there are $2$-opetopes, namely one loop per element in $\NN$. The category of
  elements $\OO_{0, 1} / \optPolyFun^1 1$ looks like this:
    \[
        \begin{tikzcd}
            *
                \ar[d, "{s_0, t_0}" left]
                \ar[dr, "{s_1, t_1}" sloped, bend left]
                \ar[drr, "{s_2, t_2}" sloped, bend left]
                \ar[drrrr, "{s_m, t_m}" sloped, bend left]
            & & & & & \\
            0 & 1 & 2 & \cdots & m & \cdots
        \end{tikzcd}
    \]
    where $*$ corresponds to the vertex of $\optPolyFun^1 1$, the numbers on
    the second row correspond to its vertices, and the morphisms are the
    inclusions of $*$ as the source or target of these vertices. The functor $E
    : \OO_{0, 1} / \optPolyFun^1 1 \to \psh{\OO_{0, 1}}$ maps
    $*$ to the graph with one vertex and no edges, and maps $m$ to the linear
    graph with $m$ consecutive edges:
    \[
        E (*) = \left( \bullet \right),
        \qquad\qquad
        E (m) = \left(
            \bullet \to \bullet \to \bullet
            \to \cdots \to \bullet
        \right) .
    \]
    On morphisms, $E (s_n)$ (resp. $E (t_n)$) is the inclusion of $\bullet$ as
    the first (resp. as the last) vertex of $E (m)$. Then, for $X \in
    \psh{\OO_{0, 1}}$, the graph $\optPolyFun^1 X$ has the same vertices as $X$,
    but its edges are paths in $X$. In other words, $\optPolyFun^1 :
    \psh{\OO_{0, 1}} \to \psh{\OO_{0, 1}}$ is the free category monad.
  \end{example}

\begin{para}
  Recall from \cref{def:pra-monad} that a p.r.a. monad is a monad $T$ whose unit
  $\id \to T$ and multiplication $T T \to T$ are cartesian, and such that its
  underlying functor is a p.r.a. We now endow $\optPolyFun^n$ with the structure
  of a p.r.a. monad over $\psh{\OO_{n-k, n}}$. We first specify the unit and
  multiplication $\eta_1 : 1 \to \optPolyFun^n 1$ and $\mu_1 : \optPolyFun^n
  \optPolyFun^n 1 \to \optPolyFun^n 1$ on the terminal object $1$, and extend
  them to cartesian natural transformations
  (\cref{lemma:coloured-zn-structure}). Next, we check that the required monad
  identities hold for $1$ (\cref{lemma:coloured-zn-law}), which automatically
  gives us the desired monad structure on $\optPolyFun^n$.
\end{para}

\begin{definition}
    \label{def:partition-opetopes}
    Let $\OO_{n+2}^{(2)}$\index{o@$\OO_{n+2}^{(2)}$} be the set of
    $(n+2)$-opetopes \emph{of uniform height $2$}, namely of the form
    \[
        \ytree{\alpha} \biggraft_{[[p]]} \ytree{\beta_{[p]}} ,
    \]
    with $\alpha, \beta_{[p]} \in \OO_{n+1}$ and $[p]$ ranging over
    $\alpha^\nodesymbol$.
\end{definition}

\begin{proposition}
    \label{prop:znzn}
    If $X \in \psh{\OO_{n-k, n}}$, then $(\optPolyFun^n \optPolyFun^n X)_{<
    n} = X_{<n}$, and if $\omega \in \OO_n$, then
    \[
        (\optPolyFun^n \optPolyFun^n X)_\omega
        \cong
        \sum_{\substack{
            \xi \in \OO_{n+2}^{(2)} \\ \tgt \tgt \xi = \omega
        }} \psh{\OO_{n-k, n}} (S [\tgt \xi], X) .
    \]
\end{proposition}
\begin{proof}
    Take $\omega \in \OO_n$, and $x \in (\optPolyFun^n \optPolyFun^n
    X)_\omega$, say $x : S [\nu] \to \optPolyFun^n X$, where $\tgt
    \nu = \omega$. For $[p_i] \in \nu^\nodesymbol$, write $x_i \eqdef x [p_i] :
    S [\nu_i] \to X$, where $\tgt \nu_i = \src_{[p_i]} \nu$.
    Informally, $x$ is a ``pasting diagram of pasting diagrams'' of $X$, namely a
    pasting diagram of the $x_i$'s, which are themselves pasting diagrams in
    $X$. The goal is to assemble the $x_i$'s in a single pasting diagram $\Phi
    (x)$. Let
    \[
        \xi \eqdef \ytree{\nu} \biggraft_{[[p_i]]} \ytree{\nu_i} ,
    \]
    and note that $\tgt \tgt \xi = \tgt \src_{[]} \xi = \tgt \nu = \omega$ by
    \condition{Glob1}. We now define a map $\Phi (x) : S [\tgt \xi]
    \to X$. Note that leaf addresses of $\xi$ are of the form
    $[[p_i][l]]$, where $[l] \in \nu_i^\leafsymbol$, thus node addresses of
    $\tgt \xi$ are of the form $\readdress_\xi [[p_i][l]]$. Let
    \[
        \Phi (x) \left( \readdress_\xi [[p_i][l]] \right)
        \eqdef x_i \left( \readdress_{\nu_i} [l] \right) .
    \]
    The construction of $\Phi (x)$ provides a map
    \[
        \Phi : (\optPolyFun^n \optPolyFun^n X)_\omega
        \to \sum_{\substack{
            \xi \in \OO_{n+2}^{(2)} \\ \tgt \tgt \xi = \omega
        }} \psh{\OO_{n-k, n}} (S [\tgt \xi], X)
    \]
    whose inverse we now construct. Let $\xi \in \OO_{n+2}^{(2)}$, say
    \[
        \xi = \ytree{\alpha} \biggraft_{[[p]]} \ytree{\beta_{[p]}} ,
    \]
    be such that $\tgt \tgt \xi = \omega$, and take $y : S [\tgt \xi]
    \to X$. Write $\nu \eqdef \tgt \xi$. As noted in
    \cref{def:partition-opetopes}, $\xi$ exhibits a partition of $\nu$ into
    subtrees, and let $\Psi (y) : S [\alpha] \to \optPolyFun^n X$
    map $[p]$ to the restriction of $y$ to the subtree $\beta_{[p]}$ of $\nu$.
    It is routine verification to check that $\Phi$ and $\Psi$ are mutually
    inverse.
\end{proof}

\begin{definition}
    \label{def:zn-laws-on-1}
    We now define $\eta_1 : 1 \to \optPolyFun^n 1$ and $\mu_1 :
    \optPolyFun^n \optPolyFun^n 1 \to \optPolyFun^n 1$, the monads
    laws of $\optPolyFun^n$, on the terminal presheaf $1 \in \psh{\OO_{n-k,
    n}}$. In dimension $< n$, they are the identity. Let $\omega \in \OO_n$.
    Recall from \cref{def:coloured-zn} that $(\optPolyFun^n 1)_\omega =
    \left\{\nu \in \OO_{n+1} \mid \tgt \nu = \omega \right\}$, and by
    \cref{prop:znzn},
    \begin{equation}
        \label{eq:znzn1}
        (\optPolyFun^n \optPolyFun^n 1)_\omega
        = \left\{
            \xi \in \OO^{(2)}_{n+2}
            \mid \tgt \tgt \xi = \omega
        \right\} .
    \end{equation}
    Now, let $(\eta_1)_\omega$ map the unique element of $1_\omega$ to
    $\ytree{\omega} \in (\optPolyFun^n 1)_\omega$, and let $(\mu_1)_\omega$ map
    $\xi \in (\optPolyFun^n \optPolyFun^n 1)_\omega$ to $\tgt \xi \in
    (\optPolyFun^n 1)_\omega$.
\end{definition}

\begin{remark}
    Let $X \in \psh{\OO_{n-k, n}}$, and consider the terminal map $! :
    X \to 1$. The map $\optPolyFun^n ! : (\optPolyFun^n
    X)_\omega \to (\optPolyFun^n 1)_\omega$ simply maps a pasting
    diagram $f : S [\nu] \to X$ (where $\tgt \nu = \omega$) to its
    shape $\nu$.
\end{remark}

\begin{lemma}
    \label{lemma:coloured-zn-structure}
    Let $X \in \psh{\OO_{n-k, n}}$, and consider the terminal map $! : X
    \to 1$. To alleviate notations, write $p \eqdef \optPolyFun^n !
    : \optPolyFun^n X \to \optPolyFun^n 1$. There exist maps
    $\eta_X : X \to \optPolyFun^n X$ and $\mu_X : \optPolyFun^n
    \optPolyFun^n X \to \optPolyFun^n$ such that the following
    squares are cartesian:
    \begin{equation}
        \label{eq:zn-laws}
        \squarediagram
            {X}{\optPolyFun^n X}{1}{\optPolyFun^n 1 ,}
            {\eta_X}{!}{p}{\eta_1}
        \qquad\qquad
        \squarediagram
            {\optPolyFun^n \optPolyFun^n X}{\optPolyFun^n X}
                {\optPolyFun^n \optPolyFun^n 1}{\optPolyFun^n 1 .}
            {\mu_X}{\optPolyFun^n p}{p}{\mu_1}
    \end{equation}
    In particular, the maps $\eta_X$ and $\mu_X$ assemble into cartesian
    natural transformations $\eta : \id \to \optPolyFun^n$ and $\mu
    : \optPolyFun^n \optPolyFun^n \to \optPolyFun^n$.
\end{lemma}
\begin{proof}
    All morphisms are identities in dimension $<n$, so it suffices to check
    that both squares are cartesian in dimension $n$.
    \begin{enumerate}
        \item If $P$ is the pullback
        \[
            \pullbackdiagram
                {P}{\optPolyFun^n X}{1}{\optPolyFun^n 1,}
                {}{!}{p}{\eta_1}
        \]
        then for $\omega \in \OO_n$ we have
        \[
            P_\omega
            = \left\{
                x \in \optPolyFun^n X
                \mid p (x) = \ytree{\omega}
            \right\}
            = \psh{\OO_{n-k, n}} (S [\ytree{\omega}], X)
            = X_\omega ,
        \]
        as $S [\ytree{\omega}] = O [\omega]$.

        \item Let $P$ be the bullback
        \[
            \pullbackdiagram
                {P}{\optPolyFun^n X}
                    {\optPolyFun^n \optPolyFun^n 1}{\optPolyFun^n 1,}
                {}{}{p}{\mu_1}
        \]
        and let $\omega \in \OO_n$. By definition, and with \cref{eq:znzn1},
        $P_\omega$ is the set of all pairs $(\xi, x)$, where $\xi \in
        \OO_{n+2}^{(2)}$ is such that $\tgt \tgt \xi = \omega$, $x : S [\nu]
        \to X$ is such that $\tgt \nu = \omega$, and subject to the
        constraint that $\tgt \xi = \nu$. By \cref{prop:znzn}, it is clear that
        $P_\omega \cong (\optPolyFun^n \optPolyFun^n X)_\omega$.
        \qedhere
    \end{enumerate}
\end{proof}

\begin{lemma}
    \label{lemma:coloured-zn-law}
    The following diagrams commute:
    \[
        \begin{tikzcd}
            \optPolyFun^n 1
                \ar[r, "\eta_{\optPolyFun^n 1}"]
                \ar[dr, equal] &
            \optPolyFun^n \optPolyFun^n 1
                \ar[d, "\mu_1"] &
            \optPolyFun^n 1
                \ar[l, "\optPolyFun^n \eta_1"']
                \ar[dl, equal] \\
            &
            \optPolyFun^n 1 , &
        \end{tikzcd}
        \qquad\qquad
        \squarediagram
            {\optPolyFun^n \optPolyFun^n \optPolyFun^n 1}
                {\optPolyFun^n \optPolyFun^n 1}{\optPolyFun^n \optPolyFun^n 1}
                {\optPolyFun^n 1 .}
            {\optPolyFun^n \mu_1}{\mu_{\optPolyFun^n 1}}{\mu_1}{\mu_1}
    \]
\end{lemma}
\begin{proof}
    Recall from \cref{def:coloured-zn} that for $X \in \psh{\OO_{n-k, n}}$,
    $(\optPolyFun^n X)_{<n} = X_{<n}$. Thus all diagrams commute trivially in
    dimension $<n$.
    \begin{enumerate}
        \item Let $\omega \in \OO_n$ and $\nu \in \optPolyFun^n 1_\omega$,
        namely $\nu \in \OO_{n+1}$ such that $\tgt \nu = \omega$. Then
        \begin{align*}
            \mu_1 \eta_{\optPolyFun^n 1} (\nu)
            &= \mu_1 \left(
                \ytree{\ytree{\tgt \nu}} \graft_{[[]]} \ytree{\nu}
            \right)
                & \text{see \cref{def:zn-laws-on-1}} \\
            &= \tgt \left(
                \ytree{\ytree{\tgt \nu}} \graft_{[[]]} \ytree{\nu}
            \right)
                & \text{see \cref{def:zn-laws-on-1}} \\
            &= \ytree{\tgt \nu} \subst_{[]} \nu
                & \text{by \cref{prop:opetope-target}} \\
            &= \nu ,
        \end{align*}
        and similarly, if $\left\{ [p_1], \ldots \right\} = \nu^\nodesymbol$,
        \begin{align*}
            \mu_1 (\optPolyFun^n \eta_1) (\nu)
            &= \mu_1 \left(
                \ytree{\nu} \biggraft_{[[p_i]]}
                \ytree{\ytree{\src_{[p_i]} \nu}}
            \right)
                & \spadesuit \\
            &= \tgt \left(
                \ytree{\nu} \biggraft_{[[p_i]]}
                \ytree{\ytree{\src_{[p_i]} \nu}}
            \right)
                & \spadesuit \\
            &= \left( \nu
                \subst_{[p_1]} \ytree{\src_{[p_1]} \nu} \right)
                \subst_{[p_2]} \ytree{\src_{[p_2]} \nu}
                \: \cdots
                & \text{by \cref{prop:opetope-target}} \\
            &= \nu ,
        \end{align*}
        where $\spadesuit$ follows from \cref{def:zn-laws-on-1}.

        \item Akin to \cref{prop:znzn}, one can show that elements of
        $\optPolyFun^n \optPolyFun^n \optPolyFun^n 1_\omega$ are
        $(n+2)$-opetopes $\xi$ of uniform height $3$ such that $\tgt \tgt \xi =
        \omega$. Let $\xi$ be such an opetope, and write it as
        \[
            \xi
            =
            \ytree{\alpha} \biggraft_{[[p_i]]} \underbrace{\left(
                \ytree{\beta_i} \biggraft_{[[q_{i, j}]]} \ytree{\gamma_{i, j}}
            \right)}_{A_i \eqdef}
            =
            \underbrace{\left(
                \ytree{\alpha} \biggraft_{[[p_i]]} \ytree{\beta_i}
            \right)}_{B \eqdef}
                \biggraft_{[[p_i][q_{i, j}]]} \ytree{\gamma_{i, j}}
        \]
        where $\alpha, \beta_i, \gamma_{i, j} \in \OO_n$, $[p_i]$ ranges over
        $\alpha^\nodesymbol$ and $[q_{i, j}]$ over $\beta_i^\nodesymbol$. Then
        \begin{align*}
            \mu_1 (\optPolyFun^n \mu_1) (\xi)
            &= \mu_1 (\optPolyFun^n \mu_1) \left(
                \ytree{\alpha} \biggraft_{[[p_i]]} A_i
            \right) \\
            &= \mu_1 \left(
                \ytree{\alpha} \biggraft_{[[p_i]]} \ytree{\tgt A_i}
            \right) \\
            &= \tgt \left(
                \ytree{\alpha} \biggraft_{[[p_i]]} \ytree{\tgt A_i}
            \right) \\
            &= \tgt \left(
                \ytree{\alpha} \biggraft_{[[p_i]]} A_i
            \right)
                & \text{by \cref{prop:opetope-target}} \\
            &= \tgt \left(
                B \biggraft_{[[p_i][q_{i, j}]]} \ytree{\gamma_{i, j}}
            \right)
                & \text{by definition} \\
            &= \tgt \left(
                \ytree{\tgt B} \biggraft_{[\readdress_B [[p_i][q_{i, j}]]]} \ytree{\gamma_{i, j}}
            \right)
                & \text{by \cref{prop:opetope-target}} \\
            &= \mu_1 \mu_{\optPolyFun^n 1} (\xi) .
        \end{align*}
        \qedhere
    \end{enumerate}
\end{proof}

\begin{proposition}
    The cartesian natural transformations $\mu$ and $\eta$ (whose components
    are defined in \cref{def:zn-laws-on-1,lemma:coloured-zn-structure}) give
    $\optPolyFun^n$ a structure of p.r.a. monad on $\psh{\OO_{n-k, n}}$.
\end{proposition}
\begin{proof}
    This is a direct consequence of
    \cref{lemma:coloured-zn-law,lemma:coloured-zn-structure}.
\end{proof}

\begin{remark}
    Clearly, when $k = 0$, we recover the polynomial monad $\optPolyFun^n$ on $\Set /
    \OO_n$.
\end{remark}

\begin{remark}
 \label{rem:zn-cardinals-spines} 
  Note from~\cref{def:coloured-zn} that the \defn{$\optPolyFun^n$-cardinals} are
  precisely the representable opetopes in $\OO_{n-k,n}$ and the spines
  $S[\nu]$ for all $\nu\in\OO_{n+1}$.
\end{remark}

We come to the essential definition of this chapter.

\begin{definition}
  \label{def:opetopic-algebra}
  A \defn{$k$-coloured, $n$-dimensional opetopic algebra}
  \index{opetopic!algebra} is an algebra of the p.r.a. monad $\optPolyFun^n$ on
  $\psh{\OO_{n-k, n}}$. We write $\Oalg k n$ \index{alg@$\Oalg k n$|see
    {opetopic algebra}} for the Eilenberg-Moore category $\optPolyFun^n\alg$.
\end{definition}

\begin{proposition}
  \label{prop:oalg-locally-presentable}
  The monad $\optPolyFun^n$ on $\psh{\OO_{n-k,n}}$ is finitary. The category
  $\oAlg_{k,n}$ is therefore locally finitely presentable.
\end{proposition}
\begin{proof}
  Every $\optPolyFun^n$-cardinal is a finite colimit of representables, thus is
  finitely presentable. This, together with~\cref{th:pra-nerve-theorem}, implies
  that $\optPolyFun$ has $\fin_{\OO_{n-k,n}}$ as arities.
\end{proof}

\begin{corollary}
  $\Oalg kn$ is the category of $\Set$\=/models of an
  $\OO_{n-k,n}$\=/sorted theory.
\end{corollary}
\begin{proof}
  From
  \cref{lem:O-lfd-cat},~\cref{thm:classification-dep-alg-theories}
  and~\cref{prop:oalg-locally-presentable}.
\end{proof}

\begin{proposition}
  \label{prop:algebra-table}
  Up to equivalence, and for small values of $k$ and $n$ with $k \leq n$, the
  category $\Oalg k n$ is given by the following table\footnote{If $k > n$, then
    $\Oalg k n = \Oalg n n$.}:
  \begin{center}
    \begin{tabular}{c|cccc}
      $k \backslash n$ & $0$ & $1$ & $2$ & $3$ \\
      \hline
      $0$ & $\Set$ & $\Monoid$ & $\Opd\pl^1$ & $\Cmbd$ \\
      $1$ & & $\Cat$ & $\Opd\pl$ & $\Cmbd_{\mathrm{col}}$ \\
      $2$ & & & $\Oalg 2 2$ & $\Oalg 2 3$ \\
      $3$ & & & & $\Oalg 3 3$
    \end{tabular}
  \end{center}
  where $\Monoid$ is the category of monoids, $\Opd\pl^1$ of uncoloured planar
  ($\Set$\=/)operads, $\Opd\pl$ of coloured planar operads, and $\Cmbd$
  (respectively, $\Cmbd_{\mathrm{col}}$) \index{combpt@$\Cmbd$|see {combinad}}
  of combinads \index{combinad} (respectively, coloured combinads) over the
  combinatorial pattern of planar trees \cite{Loday2012a}.
\end{proposition}
\begin{proof}
  [Proof (sketch)] Let us first treat the case where $k = 0$.
  \begin{enumerate}
  \item If $n = 0$, then $\optPolyFun^0$ is by definition the identity functor
    on $\Set / \OO_0 = \Set$, thus $\optPolyFun^0$-algebras bear no structure,
    and are simply sets.
  \item The polynomial monad $\optPolyFun^1 = (\optPolyFun^0)^+$ is isomorphic
    to
    \[
      \polynomialfunctor {\{ \optOne \}}{\NN_<}{\NN}{\{ \optOne \}}
      {\src}{}{\tgt}
    \]
    where for $m \in \NN$, $\NN_< (m) \eqdef \left\{ 0, 1, \ldots, m-1
    \right\}$. The result follows by \cite[Example 1.9]{Gambino2013}.
  \item The functor $\optPolyFun^2 : \Set / \NN \to \Set / \NN$ maps a signature
    $X = (X_m \mid m \in \NN) \in \Set / \NN$ to the set of trees whose nodes
    are adequately decorated by elements of $X$, namely it is the free planar
    operad monad.
  \item A $\optPolyFun^4$-algebra is a set of ``planar trees'' (namely an
    element of $\Set / \OO_3$) with an suitable notion of substitution, which is
    structure encapsulated in the notion of $PT$-combinad.
  \end{enumerate}
  Let us now consider higher values of $k$.
  \begin{enumerate}
  \item Assume $k = n = 1$. Then $\psh{\OO_{0, 1}}$ is the category of graphs,
    and a $\optPolyFun^1$ maps a graph to its graph of paths. A
    $\optPolyFun^1$-algebra is just a graph with an adequate notion of
    composition of paths, namely a category.
  \item Similarly, in the case $k = 1$ and $n = 2$, the category $\psh{\OO_{1,
        2}}$ is the category of signatures whose inputs and output of functions
    are typed. Extending the reasoning of the case $k = 0$, it is easy to see
    that a $\optPolyFun^2$-algebra is a coloured planar operad. \qedhere
  \end{enumerate}
\end{proof}

\begin{para}
  We now describe how the ``ordinary'' nerve theorem
  (\cref{th:pra-nerve-theorem}) for the p.r.a. monad $\optPolyFun^n$ (for every
  $(k,n)$) implies that $\Oalg kn$ is
  a Gabriel-Ulmer localisation of $\psh{\bbLambda_{k, n}}$ at a set of
  \emph{algebraic} spine inclusions, where $\bbLambda_{k, n} \eqdef
  \Theta_{\optPolyFun^n}$ (see~\cref{not:thetaT}).
\end{para}

\begin{definition}
  \label{def:opetopic-shape}
  By \cref{def:pra-monad,def:coloured-zn}, the category $\Theta_0$ of
  $\optPolyFun^n$-cardinals is the full subcategory of $\psh{\OO_{n-k, n}}$
  spanned by the representables $O [\omega]$, where $\omega \in \OO_{n-k, n}$,
  and the spines $S[\nu]$, where $\nu \in \OO_{n+1}$
  (\cref{rem:zn-cardinals-spines}). The category of \defn{opetopic shapes}
  $\bbLambda_{k, n}$\index{$\bbLambda$|see {opetopic shape}} is the full
  subcategory $\Theta_{\optPolyFun^n}$ of $\Oalg k n$ obtained from the
  identity-on-objects, fully faithful factorisation (see \cref{not:thetaT})
  \[
    \Theta_0\xto{h} \bbLambda_{k,n} \subto \Oalg kn .
  \]
\end{definition}

\begin{notation}
  \label{conv:implicit-k-n}
  Throughout the rest of this section, we will frequently fix parameters $k \leq n \in \NN$
  \emph{implicitly}, and suppress them in notation whenever it is unambiguous.
  For example, we write $\bbLambda$ instead of $\bbLambda_{k, n}$, $\optPolyFun$
  instead of $\optPolyFun^n$, $\oAlg$ instead of $\Oalg k n$, etc.
\end{notation}

\begin{para}
  \label{def:algebraic-realisation-lambda}
  Recall from \cref{prop:oalg-locally-presentable} that $\oAlg$ is cocomplete.
  From the ordinary nerve theorem (\cref{th:pra-nerve-theorem}), the
  inclusion $\bbLambda \subto \oAlg$ is dense. We thus have a reflective adjunction
  \index{$\tau$|see {algebraic realisation}} \index{n@$N$|see {nerve}}
  \[
    \tau : \psh{\bbLambda} \localisation \oAlg : N .
  \]
  The left adjoint is called the \defn{algebraic
    realisation}, and the right adjoint (the nerve
  functor) is called the \defn{(algebraic) nerve}.
\end{para}

\begin{example}
  \label{ex:opetopic-shape}
  \begin{enumerate}
  \item Take $n = k = 1$. By \cref{prop:algebra-table}, $\Oalg 1 1 = \Cat$, and
    $\bbLambda_{1, 1}$ is the full subcategory of $\Cat$ spanned by $[m] =
    \optPolyFun^1 O[\optInt{m}]$, where $m \in \NN$. Therefore, $\bbLambda_{1,
      1} = \bDelta$, the simplex category. The algebraic realisation $\tau_{1, 1}
    : \psh\bDelta \to \Cat$ is just the realisation of a simplicial set into a
    category, and its right adjoint $N_{1, 1}$ is the usual simplicial nerve
    functor.
  \item Likewise, $\bbLambda_{1, 2}$ is the category of coloured operads
    generated by trees, thus it is the planar version $\bOmega\pl$ of Moerdijk
    and Weiss's category of dendrices $\bOmega$. The functor $N_{1, 2}$ is the
    \emph{dendroidal nerve} of \cite[Sec. 4]{Moerdijk2007}, and $\tau_{1, 2}$ is
    its left adjoint (in \emph{op. cit.}, they are written $N_d$ and $\tau_d$
    respectively).
  \end{enumerate}
\end{example}

\begin{para}
  [Algebraic spines] We will export the notion of \emph{spine} from $\psh{\OO}$
  to the category $\psh{\bbLambda}$. As we shall see, we will be able to
  characterise $\oAlg$ as the Gabriel-Ulmer localisation of $\psh\bbLambda$ at
  the set $\bbSigma$ of \emph{algebraic} spine inclusions.
\end{para}

\begin{definition}
  \label{def:spine-lambda}
  For $\nu \in \OO_{n+1}$, write $\lambda = \optPolyFun S [\nu]$, and let $S
  [\lambda]$, the \defn{(algebraic) spine} \index{spine} of the opetopic shape
  $\lambda$, be the colimit
  \[
    S [\lambda] \eqdef h_! S [\nu] = \colim \left( \OO_{n-k, n} / S[\nu] \to
      \OO_{n-k, n} \xto{h} \bbLambda \xhookrightarrow{\yoneda} \psh{\bbLambda}
    \right) .
  \]
  Let $\sfS_{\lambda} : S[\lambda] \subto \lambda$ be the \emph{(algebraic)
    spine inclusion of $\lambda$}\index{spine!inclusion}, and let
  $\bbSigma$ be the set of spine inclusions in $\psh{\bbLambda}$:
  \[
    \bbSigma \eqdef \left\{ \sfS_{\lambda} : S [\lambda] \subto \lambda \mid \nu \in
      \OO_{n+1} \right\} .
  \]
\end{definition}

\begin{example}
  If $k = n = 1$, then $\bbLambda_{1, 1} = \bDelta$, and the $(n+1)$-opetopes are
  the opetopic integers. For $m \in \NN$, the diagram $\OO_{0, 1} / S
  [\optInt{m}] \to \OO_{0, 1}$ is
  \[
    \begin{tikzcd} [column sep = small]
      &
      \optOne
      &
      &
      \optOne
      &
      &
      \optOne
      &
      \cdots
      &
      \optOne
      &
      &
      \optOne
      \\
      \optZero \ar[ur, "\src_*" above left]
      &
      &
      \optZero \ar[ur, "\src_*" above left] \ar[ul, "\tgt" above right]
      &
      &
      \optZero \ar[ur, "\src_*" above left] \ar[ul, "\tgt" above right]
      &
      &
      &
      &
      \optZero \ar[ur, "\src_*" above left] \ar[ul, "\tgt" above right]
      &
      &
      \optZero  \ar[ul, "\tgt" above right]
    \end{tikzcd}
  \]
  where there are $m$ instances of $\optOne$. By definition, $\optPolyFun
  \optZero = \bDelta [0]$ and $\optPolyFun \optOne = \bDelta [1]$. Further,
  $\optPolyFun \src_* = d_1$ and $\optPolyFun \tgt = d_0$. Thus, if $\lambda
  \eqdef \optPolyFun S [\optInt{m}]$, then $S [\lambda]$ is the colimit of the
  following diagram in $\psh\bDelta$:
  \[
    \begin{tikzcd} [column sep = tiny]
      &
      \bDelta [1]
      &
      &
      \bDelta [1]
      &
      \cdots
      &
      \bDelta [1]
      &
      &
      \bDelta [1]
      \\
      \bDelta [0] \ar[ur, "d_1" above left]
      &
      &
      \bDelta [0] \ar[ur, "d_1" above left] \ar[ul, "d_0" above right]
      &
      &
      &
      &
      \bDelta [0] \ar[ur, "d_1" above left] \ar[ul, "d_0" above right]
      &
      &
      \bDelta [0]\ar[ul, "d_0" above right]
    \end{tikzcd}
  \]
  Therefore, $S [\lambda]$ is the simplicial spine $S [m]$.
\end{example}

\begin{theorem}
  [Nerve theorem for $\bbLambda$]
  \label{th:nerve-theorem-lambda}
  \begin{enumerate}
  \item The functor $\tau : \bbLambda \to \oAlg$ is dense, or equivalently, the
    nerve $N : \oAlg \to \psh{\bbLambda}$ is fully faithful.
  \item A presheaf $X \in \psh\bbLambda$ is in the essential image of $N$ if and
    only if $\bbSigma \perp X$.
  \item (Segal condition) The reflective adjunction $\tau : \psh{\bbLambda}
    \rightleftarrows \oAlg : N$ exhibits $\oAlg$ as the Gabriel-Ulmer
    localisation of $\psh{\bbLambda}$ at the spine inclusions: $\oAlg \simeq
    \bbSigma^{-1} \psh\bbLambda$.
  \end{enumerate}
\end{theorem}
\begin{proof}
  \begin{enumerate}
  \item This is \cref{th:pra-nerve-theorem}.

  \item Recall that $\Theta_0$ denotes the category of $\optPolyFun^n$-cardinals
    (\cref{rem:pra,def:opetopic-shape}). Consider the composite
    \[
      \OO_{n-k, n} \xhookrightarrow{i} \Theta_0 \xto{h} \bbLambda .
    \]
    The category $\Theta_0 - \im i$ is composed of exactly the (opetopic) spines
    $S [\nu]$, for $\nu \in \OO_{n+1}$. By~\cref{coroll:algebra-lifting}, any
    $X\in\psh\bbLambda$ is in the essential image of the nerve if and only if
    for every $\nu$, we have $h_!S[\nu]\perp X$. But these are exactly the
    algebraic spines (\cref{def:spine-lambda}).
  \item Follows from (2). \qedhere
  \end{enumerate}
\end{proof}

\begin{remark}
  \cref{th:nerve-theorem-lambda} is a general form of the well-known results
  that $\Cat$ (the case $k = n = 1$) and $\Opd\pl$ (the case $k =
  1$, $n = 2$) have fully faithful nerve functors to $\psh\bDelta$ and
  $\psh\bOmega\pl$ \cite[Example 4.2]{Moerdijk2007} respectively, exhibiting them
  as Gabriel-Ulmer localisations of the respective presheaf categories at a set
  of \emph{spine inclusions}, also called \emph{Grothendieck-Segal colimits}. I
  do not think that the corresponding result for $\Cmbd_{\mathrm{col}}$ (the
  case $k=1, n=3$) exists in the literature.
\end{remark}

\section{Opetopic nerve functors for opetopic algebras}
\label{sec:opet-nerve-functors}
In~\cref{sec:opetopic-algebras}, we saw that $\Oalg kn$ was the category of
models in $\Set$ of an $\OO_{n-k,n}$\=/sorted theory. We will now show that we
can construct all opetopic algebras \emph{directly from opetopic sets}, by means
of an \emph{opetopic realisation} functor $h_{k, n} : \psh{\OO} \to \Oalg k n$.
Every opetopic realisation (for all $(k,n)$) will have a fully faithful
accessible right adjoint called the \emph{opetopic nerve functor}. That is to
say, every $\Oalg kn$ is the category of models of an \emph{idempotent}
$\OO$\=/sorted theory.

This construction is based on the following heuristic: given $X \in
\psh{\OO}$, we shall interpret its cells of dimension $\leq n$ as ``generators''
for some object in $\Oalg kn$, and its $(n+1)$\=/cells as ``relations'', while
its $(n+2)$\=/cells will be composites of relations or ``coherences''. This will
together give an ``opetopic presentation'' of an object of $\Oalg kn$, and every
object in $\Oalg kn$ admits a unique presentation.

The first step to implement this is to extend the free functor $\optPolyFun^n
\colon \OO_{n-k, n} \to \bbLambda$ to a functor from $\OO_{n-k, n+2}$.
Informally, the image of an $(n+1)$-opetope represents an algebra with
essentially one relation, and the image of an $(n+2)$-opetope is an algebra,
also with essentially a single relation, but which is presented with many
smaller composable relations (see \cref{ex:doth} for an illustration). Thus,
realisations of $(n+1)$-opetopes are the relations enforcing ``unique
composites'' in opetopic algebras, while realisations of $(n+2)$-opetopes
enforce the coherences ``associativity of composition''.

Then, in \cref{def:algebraic-realisation-O}, the realisation $h_{k, n}$ for
opetopes is defined as a composite of left adjoints
\[
  \psh{\OO} \stackrel{(-)_{n-k,n+2}}{\localisation} \psh{\OO_{n-k, n+2}}
  \rightleftarrows \psh{\bbLambda_{k, n}} \stackrel{\tau_{k, n}}{\localisation}
  \Oalg k n .
\]
Finally, we will show that each of the previous adjunctions is a Gabriel-Ulmer
localisation, and so is the composite. To avoid clutter, we follow
\cref{conv:implicit-k-n}, omitting $k$ and $n$ when possible, thus $\bbLambda =
\bbLambda_{k, n}$, $\oAlg = \Oalg k n$, $\optPolyFun = \optPolyFun^n$, etc.

\begin{definition}
  \label{def:doth}
  There is a canonical functor $\OO_{n-k, n} \to \bbLambda$, mapping an opetope
  $\omega$ to $\optPolyFun O [\omega]$. We now extend it to a functor $h:
  \OO_{n-k, n+2} \to \bbLambda$\index{h@$h$}. On objects, it is given by
  \begin{align*}
    h : \OO_{n-k, n+2} &\to \bbLambda \\
    \omega &\longmapsto \begin{cases}
      \optPolyFun O [\omega]
      & \text{if $\dim \omega \leq n$,} \\
      \optPolyFun S [\omega]
      & \text{if $\dim \omega = n+1$,} \\
      \optPolyFun S [\tgt \omega]
      & \text{if $\dim \omega = n+2$.}
    \end{cases}
  \end{align*}
  We now specify $h$ on morphisms. Since it extends the natural functor
  $\OO_{n-k, n} \to \bbLambda$, it is enough to consider morphisms in $\OO_{n,
    n+2}$, so take $\nu \in \OO_{n+1}$ and $\xi \in \OO_{n+2}$.
  \begin{enumerate}
  \item For $[p] \in \nu^\nodesymbol$, let $ h \left( \src_{[p]} \nu
      \xto{\src_{[p]}} \nu \right) \eqdef \optPolyFun \left( O [\src_{[p]} \nu]
      \xto{\src_{[p]}} S [\nu] \right) $.

  \item In order to define $ h \left( \tgt \nu \xto{\tgt} \nu \right) = \left(
      \optPolyFun O [\tgt \nu] \xto{h \tgt} \optPolyFun S [\nu] \right) $, it is
    enough to provide a morphism $O [\tgt \nu] \to \optPolyFun S [\nu]$, namely
    a cell in $\optPolyFun S [\nu]_{\tgt \nu}$. Let it be $\left( S [\nu]
      \xto{\id} S [\nu] \right) \in \optPolyFun S [\nu]_{\tgt \nu}$.

  \item Let $ h \left( \tgt \xi \xto{\tgt} \xi \right) = \left( \optPolyFun S
      [\tgt \xi] \xto{h \tgt} \optPolyFun S [\tgt \xi] \right) $ be the identity
    map.

  \item
    \label{def:doth-item-4}
    Let $[p] \in \xi^\nodesymbol$. In order to define a morphism of
    $\optPolyFun$-algebras
    \[
      h \left( \src_{[p]} \xi \xto{\src_{[p]}} \xi \right) = \left( \optPolyFun
        S [\src_{[p]} \xi] \xto{h \src_{[p]}} \optPolyFun S [\tgt \xi] \right),
    \]
    it is enough to provide a morphism $h \src_{[p]} : S[\src_{[p]}\xi] \to
    \optPolyFun S [\tgt \xi]$ in $\psh{\OO_{n-k, n}}$, which we now construct.
    \begin{enumerate}
    \item Using~\cref{eq:node-decomposition}, $\xi$ decomposes as
      \[
        \xi = \alpha \graft_{[p]} \ytree{\src_{[p]} \xi} \biggraft_{[[q_i]]}
        \beta_i ,
      \]
      for some $\alpha, \beta_i \in \OO_{n+2}$, and where $[q_i]$ ranges over
      $(\src_{[p]} \xi)^\nodesymbol$. The leaves of any $\beta_i$ are therefore
      a subset of the leaves of $\xi$. More precisely, a leaf address $[l] \in
      \beta_i^\leafsymbol$ corresponds to the leaf $[p[q_i]l]$ of $\xi$. This
      defines an inclusion $f_i : S [\tgt \beta_i] \to S [\tgt \xi]$ that maps
      the node $\readdress_{\beta_i} [l] \in (\tgt \beta_i)^\nodesymbol$ to
      $\readdress_\xi [p[q_i]l] \in (\tgt \xi)^\nodesymbol$.

    \item Note that by definition, the map $f_i$ is an element of
      \[
        \psh{\OO_{n-k, n}} (S [\tgt \beta_i], S [\tgt \xi]) \:\subseteq\:
        \optPolyFun S [\tgt \xi]_{\tgt \tgt \beta_i} ,
      \]
      and since $\tgt \tgt \beta_i = \tgt \src_{[]} \beta_i = \edg_{[p[q_i]]}
      \xi$ (by \condition{Glob1} and \condition{Inner}), we have $f_i \in
      \optPolyFun S [\tgt \xi]_{\edg_{[p[q_i]]} \xi}$.

    \item Together, the $f_i$ assemble into the required morphism $h \src_{[p]}
      : S [\src_{[p]} \xi] \to \optPolyFun S [\tgt \xi]$, that maps the node
      $[q_i] \in (\src_{[p]} \xi)^\nodesymbol$ to $f_i$. So in conclusion, we
      have
      \begin{align*}
        h \src_{[p]} :
        S [\src_{[p]} \xi]
        &\to
          \optPolyFun S [\tgt \xi] \\
        (h \src_{[p]}) [q_i] :
        S [\tgt \beta_i]
        &\to
          S [\tgt \xi] \\
        \readdress_{\beta_i} [l]
        &\longmapsto
          \readdress_\xi [p [q_i] l] ,
      \end{align*}
      for $[q_i] \in (\src_{[p]} \xi)^\nodesymbol$ and $[l] \in
      \beta_i^\leafsymbol$.
    \end{enumerate}
  \end{enumerate}
  This defines $h$ on object and morphisms, and functoriality is
  straightforward.
\end{definition}

\begin{example}
  \label{ex:doth}
  Consider the case $k = n = 1$, so that $h = h_{1, 1}$ is a functor $\OO_{0, 3}
  \to \bbLambda_{1, 1} \cong \bDelta$. In low dimensions, we have $h \optZero =
  [0]$, $h \optOne = [1]$, and $h \optInt{m} = [m]$ with $m \in \NN$, since $h$
  is $\optPolyFun$ in this case. For instance,
  \[
    h \optInt{3} = h \left( \tikzinput[.8]{opetope-2-graphical}{3.numbers}
    \right) = [3]
  \]
  is the category with $3$ generating morphisms, and the $2$-cell of
  $\optInt{3}$ just witnesses their composition.

  Consider the following $3$-opetope $\xi$:
  \[
    \xi = \ytree{\optInt{3}} \graft_{[[*]]} \ytree{\optInt{2}} \graft_{[[**]]}
    \ytree{\optInt{1}} = \left( \tikzinput[.8]{opetope-3-graphical}{ex1.numbers}
    \right)
  \]
  Then $h \xi = \optPolyFun S [\tgt \xi] = \optPolyFun S [\optInt{4}] = [4]$.
  This result should be understood as the poset of points of $\xi$ (represented
  as dots in the pasting diagram above) ordered by the topmost arrows. The
  $2$-dimensional faces of $\xi$ provide several relations among the generating
  arrows, and the $3$-cell is a witness of the composition of those relations.

  Take the face map $\src_{[]} : \optInt{3} \to \xi$, corresponding to the
  trapezoid at the base of the pasting diagram. Then $h \src_{[]}$ maps points
  $0$, $1$, $2$, $3$ of $h \optInt{3} = [3]$ to points $0$, $1$, $3$, $4$ of $h
  \xi$, respectively. In other words, it ``skips'' point $2$, which is exactly
  what the pasting diagram above depicts: the $[]$-source of $\xi$ does not
  touch point $2$ (the topmost one). Likewise, the map $h \src_{[[**]]} : [1] =
  h \optInt{1} \to h \xi$ maps $0$, $1$ to $0$, $1$, respectively.

  Consider now the target map $\tgt : \optInt{4} \to \xi$. Since the target face
  touches all the points of $\xi$ (this can be checked graphically, but more
  generally follows from \condition{Glob2}), $h \tgt$ should be the identity map
  on $[4]$, which is precisely what the definition gives.
\end{example}

\begin{remark}
  Recall that, as a p.r.a.~monad on a presheaf category, $\optPolyFun$ has an
  associated ``generic-free'' factorisation system on the category $\bbLambda$
  \cite[Example 4.21]{weber2007familial}. We note that in \cref{def:doth}, the
  functor $h : \OO_{n-k,n+2} \to \bbLambda$ takes a source map $\src_{[]} :
  \src_{[]}\xi \to \xi$ (where $\xi \in \OO_{n+2}$) to a \emph{generic} morphism
  $\optPolyFun S[\src_{[]}\xi] \to \optPolyFun S[\tgt\xi]$. This motivates part
  \ref{def:doth-item-4} of \cref{def:doth}. Namely, $h$ sends a morphism
  $\src_{[p]} : \src_{[p]}\xi \to \xi$ to a generic-free composite $\optPolyFun
  f \circ h\src_{[]} : \optPolyFun S[\src_{[p]}\xi] \to \optPolyFun S[\tgt\nu]
  \to \optPolyFun S[\tgt \xi]$ where $\nu = \ytree{\src_{[p]} \xi}
  \biggraft_{[[q_i]]} \beta_i$ is the maximal subtree of $\xi$ that ``begins''
  at the node $[p]$, and where $f:S[\tgt\nu]\to S[\tgt\xi]$ is the inclusion of
  the leaves of the subtree $\nu$ into the leaves of $\xi$.
\end{remark}

\begin{definition}
  \label{def:algebraic-realisation-O}
  With slight abuse of notation, let $h : \psh{\OO} \rightleftarrows \oAlg :
  M$ be the composite adjunction
  \[
    \psh{\OO} \stackrel{(-)_{n-k,n+2}}{\localisation} \psh{\OO_{n-k, n+2}}
    \stackrel{h_!}{\rightleftarrows} \psh{\bbLambda}
    \stackrel{\tau}{\localisation} \oAlg .
  \]
  The left adjoint $h$ is the \defn{opetopic realisation functor}, and the right
  adjoint $M$ is the \defn{opetopic nerve functor}.
\end{definition}

\begin{remark}
  The first adjunction of the composite is just a truncation, and does not carry
  any information; the part between $\psh{\OO_{n-k, n+2}}$ and $\oAlg$ is
  actually what implements the $n$-cells of a presheaf as operations, and
  $(n+1)$-cells as relations. The $(n+2)$-cells represent coherences among
  relations (e.g. associativity of composition in categories) and cannot be
  discarded, namely one cannot obtain a reflective adjunction of the
  form $\psh{\OO_{n-k, n+1}} \rightleftarrows \oAlg$.
\end{remark}

\begin{remark}
  We now have a commutative triangle of adjunctions:
  \begin{equation}
    \label{eq:triangle-adjunction}
    \begin{tikzcd} [column sep = small, row sep = large]
      & \oAlg \arrow[dl, shift left = .4em, "M" below right,
      "\scriptstyle{\bot}" {sloped, above}] \arrow[dr, shift right = .4em, "N"
      below left, "\scriptstyle{\bot}" {sloped, above}] &
      \\
      \psh{\OO} \arrow[ur, shift left = .4em, "h" above left] \arrow[rr, shift
      left = .4em] & & \psh{\bbLambda} , \arrow[ul, shift right = .4em, "\tau"
      above right] \arrow[ll, shift left = .4em, "\scriptstyle{\bot}" {sloped,
        above}]
    \end{tikzcd}
  \end{equation}
  The notation $h$ might seem a bit overloaded, but its meaning is quite simple:
  it always takes an opetopic set and produces an algebra. If that opetopic set
  is a representable opetope in $\OO_{n-k, n+2}$, then it falls within the scope
  of \cref{def:doth}, and the output algebra is in fact an opetopic shape,
  namely in $\bbLambda$.
\end{remark}

\begin{para}
  [Diagrammatic morphisms] We will now prove various (technical) facts about the
  functor $h : \OO_{n-k, n+2} \to \bbLambda$ of \cref{def:doth}, eventually
  leading to \cref{lemma:diagramatic-lemma}, stating that all morphisms in
  $\bbLambda$ admit a good ``geometric decomposition'' (see
  \cref{def:diagrammatic-morphism,ex:diagrammatic-morphism}). This result will
  be crucial for the \emph{nerve theorem for $\OO$} (\cref{th:nerve-theorem-O}).
\end{para}
\begin{definition}
  \label{def:diagrammatic-morphism}
  Let $\nu_1, \nu_2 \in \OO_{n+1}$. A morphism $f : h \nu_1 \to h \nu_2$ in
  $\bbLambda$ is \emph{diagrammatic}\index{diagrammatic morphism} if there
  exists an opetope $\xi \in \OO_{n+2}$ and a node address $[p] \in
  \xi^\nodesymbol$ such that $\src_{[p]} \xi = \nu_1$, $\tgt \xi = \nu_2$, and
  $f = (h \tgt)^{-1} (h \src_{[p]})$. This situation is summarized by the
  following diagram, called a \emph{diagram of $f$}\index{diagram}:
  \[
    \frac{
      \begin{tikzcd} [ampersand replacement = \&]
        \&
        \xi \\
        \nu_1 \ar[ur, sloped, near end, "\src_{[p]}"] \& \nu_2 \ar[u, "\tgt"]
      \end{tikzcd}
    }{
      \begin{tikzcd} [ampersand replacement = \&]
        h \nu_1 \ar[r, "f"] \& h \nu_2 .
      \end{tikzcd}
    }
  \]
\end{definition}

\begin{example}
  \label{ex:diagrammatic-morphism}
  Consider the case $k = n = 1$ again, and recall from
  \cref{ex:opetopic-shape} that in this case, $\bbLambda = \bDelta$.
  Consider the map $f : [2] \to [3]$ in $\bDelta$, where $f (0)
  = 0$, $f (1) = 1$, and $f (2) = 2$. In other words, $f = d^3$ is the
  3\textsuperscript{rd} coface map. Taking $\xi$ as on the left, we obtain a
  diagram of $f$ on the right:
  \[
    \xi
    = \ytree{\optInt{2}} \graft_{[[*]]} \ytree{\optInt{2}}
    = \left( \tikzinput[.8]{opetope-3-graphical}{classic} \right) ,
    \qquad\qquad
    \frac{
      \begin{tikzcd} [ampersand replacement = \&]
        \&
        \xi \\
        \optInt{2} \ar[ur, sloped, near end, "\src_{[[*]]}"] \&
        \optInt{3} \ar[u, "\tgt"]
      \end{tikzcd}
    }{
      \begin{tikzcd} [ampersand replacement = \&]
        {} [2] \ar[r, "f"] \& {} [3]
      \end{tikzcd}
    }
  \]
  Consider now a non injective map $g : [2] \to [1]$ where $g (0) = g (1) = 0$
  and $g (2) = 1$. In other words, $g = s^0$ is the 0\textsuperscript{th}
  codegeneracy map. Taking $\xi'$ as on the left, we obtain a diagram of $g$ on
  the right:
  \[
    \xi'
    = \ytree{\optInt{2}} \graft_{[[*]]} \ytree{\optInt{0}}
    = \left( \tikzinput[.8]{opetope-3-graphical}{degen2} \right) ,
    \qquad\qquad
    \frac{
      \begin{tikzcd} [ampersand replacement = \&]
        \&
        \xi' \\
        \optInt{2} \ar[ur, sloped, near end, "\src_{[]}"] \&
        \optInt{1} \ar[u, "\tgt"]
      \end{tikzcd}
    }{
      \begin{tikzcd} [ampersand replacement = \&]
        {} [2] \ar[r, "g"] \& {} [1]
      \end{tikzcd}
    }
  \]
  \cref{lemma:diagramatic-composite-diagramatic} below states that diagrammatic
  morphisms are stable under composition, and these two examples seem to
  indicate that all simplicial cofaces and codegeneracies are diagrammatic. One
  might thus expect all morphisms of $\bDelta$ to be in the essential image of
  $h_{1, 1} : \OO_{0, 3} \to \bDelta$. This is indeed true, not just for
  $\bDelta$, but also for every $\OO_{n-k,n+2}\to\bbLambda$
  (\cref{prop:h-surjective}).
\end{example}

\begin{lemma}
  \label{lemma:diagramatic-composite-diagramatic}
  If $f_1$ and $f_2$ are diagrammatic as on the left, the diagram on the right
  is well-defined, and is a diagram of $f_2 f_1$.
  \[
    \frac{
      \begin{tikzcd} [ampersand replacement = \&]
        \& \xi_1 \& \xi_2 \\
        \nu_1 \ar[ur, sloped, near end, "\src_{[p_1]}"] \&
        \nu_2
        \ar[u, "\tgt"]
        \ar[ur, sloped, near end, "\src_{[p_2]}"] \&
        \nu_3 \ar[u, "\tgt"]
      \end{tikzcd}
    }{
      \begin{tikzcd} [ampersand replacement = \&]
        h \nu_1 \ar[r, "f_1"] \&
        h \nu_2 \ar[r, "f_2"] \&
        h \nu_3 ,
      \end{tikzcd}
    }
    \qquad\qquad
    \frac{
      \begin{tikzcd} [ampersand replacement = \&]
        \&
        \mathclap{\xi_2 \subst_{[p_2]} \xi_1}\phantom{\xi_3} \\
        \nu_1 \ar[ur, sloped, "\src_{[p_2 p_1]}"]
        \&
        \nu_3 \ar[u, "\tgt"]
      \end{tikzcd}
    }{
      \begin{tikzcd} [ampersand replacement = \&]
        h \nu_1 \ar[r, "f_2 f_1"] \& h \nu_3
      \end{tikzcd}
    }
  \]
\end{lemma}
\begin{proof}
  It is a simple but lengthy matter of unfolding the definition of $h$. First,
  note that
  \begin{align*}
    \tgt (\xi_2 \subst_{[p_2]} \xi_1)
    &= \tgt \tgt (\ytree{\xi_2} \graft_{[[p_2]]} \ytree{\xi_1})
    & \text{by \cref{prop:opetope-target}} \\
    &= \tgt \src_{[]} (\ytree{\xi_2} \graft_{[[p_2]]} \ytree{\xi_1})
    & \text{by \condition{Glob2}} \\
    &= \tgt \xi_2 = \nu_3 .
  \end{align*}
  Using \cref{eq:node-decomposition}, we decompose $\xi_1$ as
  \begin{equation}
    \label{eq:diagramatic-composite-diagramatic:xi1}
    \xi_1
    =
    \alpha_1 \graft_{[p_1]} \ytree{\nu_1} \biggraft_{[[q_i]]} \beta_i,
  \end{equation}
  where $[q_i]$ ranges over $\nu_1^\nodesymbol$. If $\beta_i^\leafsymbol =
  \left\{ [l_{i, j}] \mid j \right\}$, then $\xi_1^\leafsymbol =
  \left\{[p_1[q_i]l_{i, j}] \mid i, j \right\}$, and so we have
  $\nu_2^\nodesymbol = (\tgt \xi_1)^\nodesymbol = \left\{ \readdress_{\xi_1}
    [p_1[q_i]l_{i, j}] \mid i, j \right\}$. Using \cref{eq:node-decomposition}
  again, we decompose $\xi_2$ as
  \begin{equation}
    \label{eq:diagramatic-composite-diagramatic:xi2}
    \xi_2
    =
    \alpha_2 \graft_{[p_2]} \ytree{\nu_2}
    \biggraft_{[\readdress_{\xi_1} [p_1[q_i]l_{i, j}]]} \gamma_{i, j}
  \end{equation}
  and write
  \begin{align*}
    \xi_2 \subst_{[p_2]} \xi_1
    &= \left(
      \alpha_2 \graft_{[p_2]} \ytree{\nu_2}
      \biggraft_{[\readdress_{\xi_1} [p_1[q_i]l_{i, j}]]} \gamma_{i, j}
      \right) \subst_{[p_2]} \xi_1
    & \text{see \eqref{eq:diagramatic-composite-diagramatic:xi1}} \\
    &= \alpha_2 \graft_{[p_2]} \xi_1
      \biggraft_{[p_1[q_i]l_{i, j}]} \gamma_{i, j}
    & \text{see \cref{para:substitution}} \\
    &= \alpha_2 \graft_{[p_2]} \left(
      \alpha_1 \graft_{[p_1]} \ytree{\nu_1} \biggraft_{[[q_i]]} \beta_i
      \right) \biggraft_{[[q_i]l_{i, j}]} \gamma_{i, j}
    & \text{see \eqref{eq:diagramatic-composite-diagramatic:xi2}} \\
    &= \left( \alpha_2 \graft_{[p_2]} \alpha_1 \right) \graft_{[p_2 p_1]}
      \ytree{\nu_1} \biggraft_{[[q_i]]} \underbrace{\left(
      \beta_i \biggraft_{[l_{i, j}]} \gamma_{i, j}
      \right)}_{\delta_i}
    & \text{rearranging terms.}
  \end{align*}
  Applying the definition of $h$ we have, for $[q_i] \in \nu_1^\nodesymbol$,
  $[l_{i, j}] \in \beta_i^\leafsymbol$, and $[r] \in \gamma_{i,
    j}^\leafsymbol$,
  \begin{align}
    h \src_{[p_2 p_1]} :
    S [\nu_1]
    &\to
      \optPolyFun S [\nu_3]
      \nonumber \\
    (h \src_{[p_2 p_1]}) [q_i] :
    S [\tgt \delta_i]
    &\to
      S [\nu_3]
      \nonumber \\
    \readdress_{\delta_i} [l_{i, j} r]
    &\longmapsto
      \readdress_{\zeta} [p_2 p_1 [q_i] l_{i, j} r] ;
      \label{eq:diagramatic-composite-diagramatic:a} \\
    h \src_{[p_1]} :
    S [\nu_1]
    &\to
      \optPolyFun S [\nu_2]
      \nonumber \\
    (h \src_{[p_1]}) [q_i] :
    S [\tgt \beta_i]
    &\to
      S [\nu_2]
      \nonumber \\
    \readdress_{\beta_i} [l_{i, j}]
    &\longmapsto
      \readdress_{\xi_1} [p_1 [q_i] l_{i, j}] ;
      \label{eq:diagramatic-composite-diagramatic:b} \\
    h \src_{[p_2]} :
    S [\nu_2]
    &\to
      \optPolyFun S [\nu_3]
      \nonumber \\
    (h \src_{[p_2]}) (\readdress_{\xi_1} [p_1 [q_i] l_{i, j}]) :
    S [\tgt \gamma_{i, j}]
    &\to
      S [\nu_3]
      \nonumber \\
    \readdress_{\gamma_{i, j}} [r]
    &\longmapsto
      \readdress_{\xi_2} [p_2 \:\: \readdress_{\xi_1} [p_1 [q_i] l_{i, j}] \:\: r] .
      \label{eq:diagramatic-composite-diagramatic:c}
  \end{align}
  Thus,
  \begin{align*}
    &(h \src_{[p_2 p_1]}) ([q_i])
      (\readdress_{\delta_i} [l_{i, j} r]) \\
    &= \readdress_{\zeta} [p_2 p_1 [q_i] l_{i, j} r]
    & \text{by \eqref{eq:diagramatic-composite-diagramatic:a}} \\
    &= \readdress_{\xi_2} [p_2 \:\: \readdress_{\xi_1} [p_1 [q_i] l_{i, j}] \:\: r]
    & \spadesuit \\
    &= (h \src_{[p_2]})
      (\readdress_{\xi_1} [p_1 [q_i] l_{i, j}])
      (\readdress_{\gamma_{i, j}} [r])
    & \text{by \eqref{eq:diagramatic-composite-diagramatic:c}} \\
    &= (h \src_{[p_2]})
      \left(
      (h \src_{[p_1]}) ([q_i]) (\readdress_{\beta_i} [l_{i, j}])
      \right)
      (\readdress_{\gamma_{i, j}} [r])
    & \text{by \eqref{eq:diagramatic-composite-diagramatic:b}} \\
    &= \left( h \src_{[p_2]} \:\: \cdot \:\: h \src_{[p_1]} \right)
      ([q_i]) (\readdress_{\delta_i} [l_{i, j} r]) ,
    & \diamondsuit
  \end{align*}
  where equality $\spadesuit$ comes from the monad structure on $\optPolyFun$,
  and $\diamondsuit$ from the definition of the composition in $\bbLambda$ when
  considered as the Kleisli category of $\optPolyFun$. \qedhere
\end{proof}

\begin{lemma}
  \label{lemma:elementary-embeddings-diagramatic}
  \begin{enumerate}
  \item Let $\nu \in \OO_{n+1}$, $\omega \eqdef \tgt \nu$, and $\xi \eqdef
    \ytree{\ytree{\omega}} \graft_{[[]]} \ytree{\nu}$. Note that $\nu = \tgt
    \xi$. The following is a diagram of $h \tgt : h \omega \to h \nu$:
    \[
      \frac{
        \begin{tikzcd} [ampersand replacement = \&]
          \& \xi \\
          \ytree{\omega} \ar[ur, sloped, near end, "\src_{[]}"] \&
          \nu \ar[u, "\tgt"]
        \end{tikzcd}
      }{
        \begin{tikzcd} [ampersand replacement = \&]
          h \omega \ar[r, "h \tgt"] \& h \nu .
        \end{tikzcd}
      }
    \]
  \item Let $\beta, \nu \in \OO_{n+1} = \tree \optPolyFun^{n-1}$, and $i : S
    [\beta] \to S [\nu]$ a morphism of presheaves. Then $i$ corresponds to an
    inclusion $\beta \subto \nu$ of $\optPolyFun^{n-1}$ trees, mapping the node at
    address $[q]$ to $[pq]$, where $[p] \eqdef i ([]) \in \nu^\nodesymbol$ is
    the address of the image of the root node. Write $\nu = \bar{\beta}
    \subst_{[p]} \beta$, for an adequate $\bar{\beta} \in \OO_{n+1}$, and let
    $\xi \eqdef \ytree{\bar{\beta}} \graft_{[[p]]} \ytree{\beta}$. Note that
    $\nu = \tgt \xi$ by \cref{prop:opetope-target}. The following is a diagram
    of $h i$:
    \[
      \frac{
        \begin{tikzcd} [ampersand replacement = \&]
          \& \xi \\
          \beta \ar[ur, sloped, near end, "\src_{[p]}"] \&
          \nu \ar[u, "\tgt"]
        \end{tikzcd}
      }{
        \begin{tikzcd} [ampersand replacement = \&]
          h \beta \ar[r, "h i"] \& h \nu .
        \end{tikzcd}
      }
    \]
  \end{enumerate}
\end{lemma}
\begin{proof}
  Tedious but straightforward matter of unfolding \cref{def:doth}.
\end{proof}

\begin{lemma}
  [Diagrammatic lemma]
  \label{lemma:diagramatic-lemma}
  Let $\nu, \nu' \in \OO_{n+1}$ with $\nu$ non degenerate, and $f : h \nu \to h
  \nu'$ be a morphism in $\bbLambda$. Then $f$ is diagrammatic.
\end{lemma}
\begin{proof}
  Let us first sketch the proof. The idea is to proceed by induction on $\nu$.
  The case $\nu = \ytree{\psi}$ for some $\psi \in \OO_n$ is fairly simple. In
  the inductive case, we essentially show that $f$ exhibits an inclusion $\nu
  \subto \nu'$ of $\optPolyFun^{n-1}$-trees by constructing an $(n+1)$-opetope
  $\bar{\nu}$ such that $\nu' = \bar{\nu} \subst_{[q]} \nu$. Thus by
  \cref{lemma:elementary-embeddings-diagramatic}, the following is a diagram of
  $h f$:
  \[
    \frac{
      \begin{tikzcd} [ampersand replacement = \&]
        \& \xi \\
        \nu \ar[ur, sloped, near end, "\src_{[[q_1]]}"] \&
        \nu' \ar[u, "\tgt"]
      \end{tikzcd}
    }{
      \begin{tikzcd} [ampersand replacement = \&]
        h \nu \ar[r, "f"] \& h \nu' ,
      \end{tikzcd}
    }
  \]
  where $\xi \eqdef \ytree{\bar{\nu}} \graft_{[[q_1]]} \ytree{\nu}$.
  Let us now dive into the details. As advertised, the proof proceeds by
  induction on $\nu$, which by assumption is not degenerate.
  \begin{enumerate}
  \item Assume $\nu = \ytree{\psi}$ for some $\psi \in \OO_n$. Then
    \[
      \bbLambda (h \ytree{\psi}, h \nu')
      = \bbLambda (\optPolyFun S [ \ytree{\psi} ] ,
      \optPolyFun S [ \nu' ] )  \\
      \cong (\optPolyFun S [ \nu ])_\psi .
    \]
    Thus $f$ corresponds to a unique morphism $\tilde{F} : S [ \nu'' ]
    \to S [ \nu' ]$, for some $\nu'' \in \OO_{n+1}$
    such that $\tgt \nu'' = \psi$, and $f$ is the composite
    \[
      h \ytree{\psi} = h \psi \xto{h \tgt} h \nu''
      \xto{\optPolyFun \tilde{F}} h \nu' .
    \]
    Those two arrows are diagrammatic by
    \cref{lemma:elementary-embeddings-diagramatic}, and by
    \cref{lemma:diagramatic-composite-diagramatic}, so is $f$.

  \item By induction, write $\nu = \nu_1 \graft_{[l]}
    \ytree{\psi_2}$ for some $\nu_1 \in \OO_{n+1}$, $[l] \in
    \nu_1^\leafsymbol$, and $\psi_2 \in \OO_n$. Write $\psi_1 \eqdef \tgt
    \nu_1$, and $\nu_2 \eqdef \ytree{\psi_2}$. Then $f$ restricts as $f_i$,
    $i = 1, 2$, given by the composite $h \nu_i \to h \nu
    \xto{f} h \nu'$.

    Let $[l']$ be the edge address of $\nu'$ (or equivalently, the
    $(n-1)$-cell of $S [\nu'] \subseteq h \nu'$) such that $\edg_{[l']}
    \nu' = f (\edg_{[l]} \nu)$. Then $\nu'$ decomposes as $\nu' = \beta_1
    \graft_{[l']} \beta_2$, for some $\beta_1, \beta_2 \in \OO_{n+1}$ (in
    particular, $\beta_1$ and $\beta_2$ are sub $\optPolyFun^{n-1}$-trees
    of $\nu'$), and $f_1$ and $f_2$ factor as
    \[
      \triangleURdiagram
      {h \nu_i}{h \beta_i}{h \nu',}
      {\bar{F}_i}{f_i}{b_i}
    \]
    where $b_i$ correspond to the subtree inclusion $\beta_i
    \subto \nu'$. By induction, $\bar{F}_i$ is diagrammatic, say
    with the following diagram:
    \[
      \frac{
        \begin{tikzcd} [ampersand replacement = \&]
          \& \xi_i \\
          \nu_i \ar[ur, sloped, near end, "\src_{[p_i]}"] \&
          \beta_i \ar[u, "\tgt"]
        \end{tikzcd}
      }{
        \begin{tikzcd} [ampersand replacement = \&]
          h \nu_i \ar[r, "\bar{F}_i"] \& h \beta_i ,
        \end{tikzcd}
      }
    \]
    and thus $\beta_i$ can be written
    as $\beta_i = \bar{\nu}_i \subst_{[q_i]} \nu_i$, for some $\bar{\nu}_i \in
    \OO_{n+1}$ and $[q_i] \in \bar{\nu}_i^\nodesymbol$. In the case $i = 2$,
    note that $\beta_2 = \bar{\nu}_2 \subst_{[q_2]} \nu_2 = \bar{\nu}_2
    \subst_{[q_2]} \ytree{\psi_2} = \bar{\nu}_2$.

    On the one hand we have
    \begin{align*}
      \edg_{[l']} \nu'
      &= f (\edg_{[l]} \nu)
      & \text{by definition of $[l']$} \\
      &= f_1 (\edg_{[l]} \nu_1)
      & \text{since $\nu = \nu_1 \graft_{[l]} \ytree{\psi_2}$} \\
      &= b_1 \bar{F}_1 (\edg_{[l]} \nu_1)
      & \text{since $f_1 = b_1 \bar{F}_1$} \\
      &= b_1 (\edg_{[q_1 l]} \beta_1)
      & \text{since $\beta_1 = \bar{\nu}_1 \subst_{[q_1]} \nu_1$} \\
      &= \edg_{[q_1 l]} \nu ,
    \end{align*}
    showing $[l'] = [q_1 l]$, and thus that $\bar{\nu}_1$ is of the form
    \begin{equation}
      \label{eq:diagramatic-lemma:barnu1}
      \bar{\nu}_1
      =
      \mu_1 \graft_{[q_1]} \ytree{\psi_1}
      \biggraft_{[[r_{1, j}]]} \mu_{1, j} ,
    \end{equation}
    where $[r_{1, j}]$ ranges over $\psi_1^\nodesymbol -
    \{\readdress_{\nu_1} [l] \}$, and $\mu_1, \mu_{1, j} \in \OO_{n+1}$.
    On the other hand,
    \begin{align*}
      \edg_{[l']} \nu'
      &= f (\edg_{[l]} \nu)
      & \text{by definition of $[l']$} \\
      &= f_2 (\edg_{[]} \nu_2)
      & \text{since $\nu = \nu_1 \graft_{[l]} \ytree{\psi_2}$} \\
      &= b_2 \bar{F}_2 (\edg_{[]} \nu_2)
      & \text{since $f_1 = b_2 \bar{F}_2$} \\
      &= b_2 (\edg_{[q_2]} \beta_2)
      & \text{since $\beta_2 = \bar{\nu}_2 \subst_{[q_2]} \nu_2$} \\
      &= \edg_{[l']} \nu' ,
    \end{align*}
    showing $[q_2] = []$, and so $\src_{[]} \beta_2 = \src_{[]} \bar{\nu}_2 =
    \psi_2$, and we can write $\beta_2$ as
    \begin{equation}
      \label{eq:diagramatic-lemma:beta2}
      \beta_2
      =
      \ytree{\psi_2} \biggraft_{[[r_{2, j}]]} \mu_{2, j} ,
    \end{equation}
    where $[r_{2, j}]$ ranges over $\psi_2^\nodesymbol$, and $\mu_{2, j}
    \in \OO_{n+1}$. Finally, we have
    \begin{align*}
      \nu'
      &= \beta_1 \graft_{[l']} \beta_2
        = (\bar{\nu}_1 \subst_{[q_1]} \nu_1) \graft_{[l']} \beta_2 \\
      &= \left(
        \mu_1 \graft_{[q_1]} \nu_1
        \biggraft_{\readdress^{-1}_{\nu_1} [r_{1, j}]} \mu_{1, j}
        \right) \graft_{[l']} \left(
        \ytree{\psi_2}
        \biggraft_{[[r_{2, j}]]} \mu_{2, j}
        \right)
      & \text{by \eqref{eq:diagramatic-lemma:barnu1} and \eqref{eq:diagramatic-lemma:beta2}} \\
      &= \left( \left(
        \mu_1 \graft_{[q_1]}
        \underbrace{\nu_1 \graft_{[l]} \ytree{\psi_2}}_{= \nu}
        \right)
        \biggraft_{[q_1] \cdot \readdress^{-1}_{\nu_1} [r_{1, j}]}
        \mu_{1, j}
        \right)
        \biggraft_{[l'[r_{2, j}]]} \mu_{2, j}
      & \text{rearranging terms} \\
      &= \bar{\nu} \subst_{[q_1]} \nu ,
    \end{align*}
    for some $\bar{\nu}' \in \OO_{n+1}$. Finally, by
    \cref{lemma:elementary-embeddings-diagramatic}, we have a diagram of $h
    f$, where $\xi \eqdef \ytree{\bar{\nu}} \graft_{[[q_1]]} \ytree{\nu}$:
    \[
      \frac{
        \begin{tikzcd} [ampersand replacement = \&]
          \& \xi \\
          \nu \ar[ur, sloped, near end, "\src_{[[q_1]]}"] \&
          \nu' \ar[u, "\tgt"]
        \end{tikzcd}
      }{
        \begin{tikzcd} [ampersand replacement = \&]
          h \nu \ar[r, "f"] \& h \nu' .
        \end{tikzcd}
      }
    \]
    \qedhere
  \end{enumerate}
\end{proof}

\begin{lemma}
  \label{lemma:doth-isomorphisms}
  \begin{enumerate}
  \item If $\omega \in \OO_{n-1}$, then $h$ maps $\tgt \tgt : \omega
    \to \itree{\omega}$ to an identity.
  \item If $\omega \in \OO_n$, then $h$ maps $\src_{[]} : \omega
    \to \ytree{\omega}$ to an identity.
  \item If $\omega \in \OO_{n+2}$, then $h$ maps $\tgt : \tgt \omega
    \to \omega$ to an identity.
  \end{enumerate}
\end{lemma}
\begin{proof}
  By inspection of \cref{def:doth}.
\end{proof}
We will show that every morphism in $\bbLambda$ is in the
image of $h$.

\begin{proposition}
  \label{prop:h-surjective}
  The functor $h : \OO_{n-k,n+2} \to \bbLambda$ is 
  surjective on morphisms.
\end{proposition}
\begin{proof}
  Let $\omega, \omega' \in \OO_{n-k, n+2}$.
  \begin{enumerate}
  \item If $\dim \omega, \dim \omega' < n-1$, then by \cref{def:doth}, $h
    \omega = \omega$ and $h \omega' = \omega'$ as presheaves over
    $\OO_{n-k, n+2}$, and thus
    \[
      \bbLambda (h \omega, h \omega')
      = \psh{\OO_{n-k, n+2}} (\omega, \omega')
      = \OO (\omega, \omega') .
    \]

  \item Assume that $\dim \omega < n-1$ and $\dim \omega' \geq n-1$. We
    first show that $O [\omega']_{< n-1} = (h \omega')_{< n-1}$ by
    inspection of \cref{def:doth}. If $\dim \omega' \leq n$, then the
    claim trivially holds. If $\dim \omega' = n+1$, then $h \omega' =
    \optPolyFun S [\omega']$, and
    \begin{align*}
      (h \omega')_{< n-1}
      &= (\optPolyFun S [\omega'])_{< n-1} \\
      &= S [\omega']_{< n-1}
      & \text{see \cref{def:coloured-zn}} \\
      &= O [\omega']_{< n-1} .
    \end{align*}
    The case where $\dim \omega' = n+2$ is proved similarly. Thus, $O
    [\omega']_{< n-1} = (h \omega')_{< n-1}$, and in particular, $O
    [\omega']_\omega = (h \omega')_\omega$. Finally,
    \begin{align*}
      \bbLambda (h \omega, h \omega')
      &\cong \psh{\OO_{n-k, n+2}} (\omega, h \omega') \\
      &= \psh{\OO_{n-k, n+2}} (\omega, \omega')
      & \spadesuit \\
      &= \OO (\omega, \omega') ,
    \end{align*}
    where $\spadesuit$ results from the observation above.

  \item If $\dim \omega \geq n-1$ and $\dim \omega' < n-1$, then
    $\bbLambda (h \omega, h \omega') = \emptyset$.

  \item Lastly, assume $\dim \omega, \dim \omega' \geq n-1$. By
    \cref{lemma:doth-isomorphisms}, we may assume that $\dim \omega = \dim
    \omega' = n+1$. If $\omega$ is non degenerate, then by
    \cref{lemma:diagramatic-lemma}, every morphism in $\bbLambda (h \omega,
    h \omega')$ is diagrammatic, thus in the essential image of $h$. Assume
    that $\omega$ is degenerate, say $\omega = \itree{\phi}$ for some $\phi
    \in \OO_{n-1}$. Akin to point (2), by inspection of \cref{def:doth},
    one can prove that $O [\omega']_\phi = (h \omega')_\phi$. Finally,
    \begin{align*}
      \bbLambda (h \omega, h \omega')
      &\cong \bbLambda (h \phi, h \omega')
      & \text{by \cref{coroll:y-i:spine-local}} \\
      &\cong \OO (\phi, \omega')
      & \spadesuit ,
    \end{align*}
    where $\spadesuit$ results from the observation above.
    \qedhere
  \end{enumerate}
\end{proof}

\begin{remark}
  It is worthwhile to note that $h : \OO_{n-k,n+2} \to \bbLambda$ is \emph{not}
  full. Take for example $n = k = 1$, so that $h$ is a functor $\OO_{0, 3} \to
  \bDelta$. Let $a, b \in \NN$, $a \neq b$, and consider the corresponding
  opetopic integers $\optInt{a}, \optInt{b} \in \OO_2$. Since they are different
  but have the same dimension, $\OO (\optInt{a}, \optInt{b}) = \emptyset$, but
  of course, $\bDelta (h \optInt{a}, h \optInt{b}) = \bDelta ([a], [b])$ is not
  empty. The diagrammatic lemma says that if $a \neq 0$, then a morphism in
  $\bDelta ([a], [b])$ can be recovered as the image of a face map of
  $\optInt{a}$ in some $3$-opetope whose target is $\optInt{b}$.
\end{remark}

\begin{para}
  [Opetopic nerve theorem]
  Recall from \cref{coroll:algebra-lifting} that we have a reflective adjunction
  \[
    \tau : \psh \bbLambda \localisation \oAlg : N
  \]
  that exhibits $\oAlg$ as the Gabriel-Ulmer localisation of $\psh{\bbLambda}$
  at the set $\bbSigma$ of spine inclusions. This is the ``ordinary'' {nerve
    theorem for $\bbLambda$}.

  We now have enough tools to prove a similar result for $\psh{\OO}$. The
  strategy is to study the adjunction $h_! : \psh{\OO_{n-k, n+2}}
  \rightleftarrows \psh{\bbLambda} : h^*$, and to show that it preserves 
  orthogonality classes of spine inclusions. It follows that it lifts to an
  adjunction $\sfS_{n+1, n+2}^{-1} \psh{\OO_{n-k, n+2}} \rightleftarrows
  \bbSigma^{-1} \psh{\bbLambda} \simeq \oAlg$, that we then prove is an adjoint
  equivalence.
\end{para}
\begin{lemma}
  \label{lemma:adjunction-to-equivalence}
  Let $F : \cA \rightleftarrows \cB : U$ be an adjunction, with unit $\eta :
  1_\cA \to UF$ and counit $\epsilon : FU \to 1_\cB$. Let $\cA'$ and $\cB'$ be
  full subcategories of $\cA$ and $\cB$ respectively. Suppose that:
  \begin{enumerate}
  \item $F \cA' \subseteq \cB'$ and $U \cB' \subseteq \cA'$,
  \item for all $a \in \cA'$, $\eta_a$ is an isomorphism, and for all $b \in
    \cB'$, $\epsilon_b$ is an isomorphism.
  \end{enumerate}
  Then the adjunction lifts to an adjoint equivalence $F : \cA' \rightleftarrows
  \cB' : U$. In particular, if $\cA'$ (resp. $\cB'$) is a Gabriel-Ulmer
  localisation at a class of morphisms $\sfK$ (resp. $\sfK'$), then condition (1)
  above can be replaced with:
  \begin{enumerate}
  \item [(1')] for all $a \in \cA$, if $\sfK \perp a$, then $\sfK' \perp Fa$,
    and for all $b \in \cB$, if $\sfK' \perp b$, then $\sfK \perp U b$.
  \end{enumerate}
\end{lemma}

\begin{proposition}
  \label{prop:h!-spines}
  The functor $h_! : \psh{\OO_{n-k, n+2}} \to \psh{\bbLambda}$ (see
  \cref{def:algebraic-realisation-O}) takes the set $\sfS_{n+1} $ of spine
  inclusions of $n+1$\=/opetopes in $\psh{\OO_{n-k, n+2}}$ to the set $\bbSigma$ of
  morphisms in $\subseteq \psh{\bbLambda}$ (the algebraic spine inclusions,
  \cref{def:spine-lambda}), and takes morphisms in $\sfS_{n+2}$ to $\bbSigma$-local
  isomorphisms.
\end{proposition}
\begin{proof}
  \begin{enumerate}
  \item Let $\nu \in \OO_{n+1}$, and recall from \cref{def:spine-lambda} that
    $\OO_{n-k, n} / S[\nu]$ is the category of elements of $S[\nu]$. We have
    \begin{align*}
      h_! S [\nu]
      &= h_! \colim_{\psi \in \OO_{n-k, n} / S[\nu]} O [\psi] \\
      &\cong \colim_{\psi \in \OO_{n-k, n} / S[\nu]} h_! O [\psi] \\
      &= \colim_{\psi \in \OO_{n-k, n} / S[\nu]}
        \yoneda_\bbLambda (h \psi) \\
      &= S [h \nu]
      & \text{see \cref{def:spine-lambda}.}
    \end{align*}

  \item For $\xi \in \OO_{n+2}$, the inclusion $S[\tgt\xi] \to S[\xi]$ is a
    relative $\sfS_{n+1}$-cell complex by
    \cref{lemma:opetopes-technical:spines}. Since $h_!$ preserves colimits, and
    since $h_! \sfS_{n+1} = \bbSigma$, we have that $h_!(S[\tgt \xi] \to S[\xi])$ is
    a relative $\bbSigma$-cell complex, and thus an $\bbSigma$-local isomorphism. In the
    square below
    \[
      \squarediagram
      {h_! S[ \tgt \xi ]}{h_! S[\xi]}{h_! O[\tgt \xi]}{h_! O[\xi]}
      {}{h_! S_{\tgt \xi}}{h_! S_\xi}{h_! \tgt }
    \]
    the top arrow is an $\bbSigma$-local isomorphism, the right arrow is in $\bbSigma$
    by the previous point, and the bottom arrow is an isomorphism by definition.
    By 3-for-2, we conclude that $h_! S_\xi$ is an $\bbSigma$-local isomorphism.
    \qedhere
  \end{enumerate}
\end{proof}

\begin{lemma}
  \label{lemma:some-spans}
  Let $X \in \psh{\OO_{n-k, n+2}}$ be such that $\sfS_{n+1, n+2} \perp X$, and
  take $\omega \in \OO_{n-k, n+2}$. The following are spans of isomorphisms:
  \begin{enumerate}
  \item for $\psi \in \OO_{n-1}$,
    \[
      \bbLambda(h \omega, h \psi) \times X_\psi
      \xleftarrow{\id \times \tgt \tgt}
      \bbLambda(h \omega, h \psi) \times X_{\itree{\psi}}
      \xto{\bbLambda(h \omega , h \tgt \tgt) \times \id}
      \bbLambda(h \omega, h \itree{\psi}) \times X_{\itree{\psi}} ;
    \]

  \item for $\psi \in \OO_n$,
    \[
      \bbLambda (h \omega, h \psi) \times X_\psi
      \xleftarrow{\id \times \src_{[]}}
      \bbLambda (h \omega, h \psi) \times X_{\ytree\psi}
      \xto{\bbLambda(h \omega , h \src_{[]}) \times \id}
      \bbLambda (h \omega, h \ytree{\psi}) \times X_{\ytree\psi} ;
    \]

  \item for $\psi \in \OO_{n+2}$,
    \[
      \bbLambda (h \omega, h \tgt \psi) \times X_{\tgt \psi}
      \xleftarrow{\id \times \tgt}
      \bbLambda (h \omega, h \tgt\psi) \times X_\psi
      \xto{\bbLambda (h \omega, h \tgt) \times \id}
      \bbLambda (h \omega, h \psi) \times X_\psi .
    \]
  \end{enumerate}
\end{lemma}
\begin{proof}
  Follows from \cref{lemma:doth-isomorphisms}.
\end{proof}

\begin{lemma}
  \label{lemma:h-low-dimensions}
  Let $\omega \in \OO_{n-k,n+2}$. If $\psi\in\OO_{n-k, n-1}$, then $\bbLambda (h
  \omega, h \psi) \cong \OO_{n-k, n+2} (\omega, \psi) $.
\end{lemma}
\begin{proof}
  Easy verification.
\end{proof}

\begin{remark}
  In the case $n=k=1$, the previous lemma just says that the morphisms
  $\bDelta([0],[m])$ are the same as graph morphisms.
\end{remark}

\begin{proposition}
  \label{lemma:h!-unit}
  Let $X \in \psh{\OO_{n-k,n+2}}$. If $\sfS_{n+1, n+2} \perp X$, then the unit
  $\eta_X : X \to h^* h_! X$ is an isomorphism.
\end{proposition}
\begin{proof}
  It suffices to show that for each $\omega \in \OO_{n-k, n+2}$, the following
  map is a bijection:
  \[
    X_\omega \xto{\eta_X} h^* h_! X_\omega = \int^{\psi \in \OO_{n-k, n+2}}
    \bbLambda (h \omega, h \psi) \times X_\psi .
  \]
  If $\omega \in \OO_{n-k, n-1}$, then $h \omega = O [\omega]$, and $\bbLambda
  (h \omega, h -) \cong \OO_{n-k, n+2} (\omega, -)$. Thus,
  \begin{align*}
    h^* h_! X_\omega
    &= \int^{\psi \in \OO_{n-k, n+2}}
      \bbLambda (h\omega, h\psi) \times X_\psi
    & \text{by definition} \\
    &\cong \int^{\psi \in \OO_{n-k, n+2}}
      \OO_{n-k, n+2} (\omega, \psi) \times X_\psi
    & \text{since $\dim \omega \leq n-1$} \\
    &\cong X_\omega
    & \text{by the density formula.}
  \end{align*}
  Assume how that $\dim \omega \geq n$. We construct an inverse of $\eta_X$ via
  a cowedge $\bbLambda (h \omega, h -) \times X_{-} \stackrel{\cdot \:
    \cdot}{\to} X_\omega$.
  \begin{enumerate}
  \item Assume $\omega \in \OO_n$. By \cref{lemma:some-spans}, it suffices to
    consider the case $\psi \in \OO_{n+1}$. To unclutter notations, write $\cC
    \eqdef \psh{\OO_{n-k, n+2}}$. We have the sequence of morphisms
    \begin{align*}
      \bbLambda (h \omega, h \psi) \times X_\psi
      &\xto{\cong} (
        \sum_{\substack{\nu \in \OO_{n+1} \\ \tgt \nu = \omega}}
      \cC (S [\nu], S[\psi])
      ) \times \cC (S[\psi], X)
      & \spadesuit \\
      &\xto{\mathrm{comp.}}
        \sum_{\substack{\nu \in \OO_{n+1} \\ \tgt \nu = \omega}}
      \cC (S[\nu], X) \\
      &\xto{\cong}
        \sum_{\substack{\nu \in \OO_{n+1} \\ \tgt \nu = \omega}}
      X_\nu
      & \spadesuit \\
      &\xto{\tgt} X_\omega ,
    \end{align*}
    where $\spadesuit$ follow from the assumption that $\sfS_{n+1} \perp X$. It
    is straightforward to verify that this defines a cowedge whose induced map
    is the required inverse.

  \item Assume $\omega \in \OO_{n+1}$. If $\omega$ is degenerate, say $\omega =
    \itree{\phi}$ for some $\phi \in \OO_{n-1}$, then $\bbLambda (h \omega, h -)
    \cong \bbLambda(h \phi, h -)$ and we are in a case we have treated before.
    So let $\omega$ be non-degenerate. By
    \cref{lemma:some-spans,lemma:h-low-dimensions}, we may suppose $\psi \in
    \OO_{n, n+1}$. Recall that for every $f \in \bbLambda(h\omega, h\psi) $, the
    diagrammatic \cref{lemma:diagramatic-lemma} computes a $\xi \in \OO_{n+2}$
    and $[p] \in \xi^\nodesymbol$ such that $\src_{[p]} \xi = \omega$, $\tgt \xi
    = \psi$ and $h \src_{[p]} \cong f$. By
    \cref{lemma:comparison:tgt-spine-local}, the target fully faithful $\tgt :
    \psi \to \xi$ is an $\sfS_{n+1, n+2}$-local isomorphism, and by assumption,
    $\sfS_{n+1, n+2} \perp X$. Therefore, we have an isomorphism $\tgt : X_\xi
    \to X_\psi$, which gives rise to a map
    \begin{align*}
      \bbLambda(h \omega, h \psi) \times X_\psi
      &\to X_\omega \\
      (f, x)
      &\longmapsto \src_{[p]} \tgt^{-1} x .
    \end{align*}
    It is straightforward to verify that this assignment defines a cowedge,
    whose associated map is the required inverse.

  \item Assume $\omega \in \OO_{n+2}$. Then by definition of $h$, $\bbLambda(h
    \omega, h -) \cong \bbLambda(h \tgt \omega, h -)$, and this is the case we
    have just treated. \qedhere
  \end{enumerate}
\end{proof}

\begin{corollary}
  \label{coroll:h!-spines}
  Let $X \in \psh{\OO_{n-k, n+2}}$. If $\sfS_{n+1, n+2} \perp X$, then
  $\bbSigma \perp h_! X$.
\end{corollary}
\begin{proof}
  Recall from \cref{prop:h!-spines} that $\bbSigma = h_!
  \sfS_{n+1}$. Let $\nu \in \OO_{n+1}$. To unclutter notations, write
  $\cC \eqdef \psh{\OO_{n-k, n+2}}$. We have
  \begin{align*}
    \psh{\bbLambda} (h_! \nu, h_! X)
    &\cong \cC (\nu, h^* h_! X)
    & \text{since } h_! \dashv h^* \\
    &\cong \cC (\nu, X)
    & \text{by \cref{lemma:h!-unit}} \\
    &\cong \cC (S [\nu], X)
    & \text{since } S_\nu \perp X \\
    &\cong \cC (S [\nu], h^* h_! X)
    & \text{by \cref{lemma:h!-unit},} \\
    &\cong \psh{\bbLambda} (h_! S [\nu], h_! X)
    & \text{since } h_! \dashv h^*
  \end{align*}
  and by construction, this isomorphism is the precomposition by $h_! S_\nu$.
  Therefore, $h_! S_\nu \perp h_!X$.
\end{proof}

\begin{notation}
  \label{not:coend-tensor-notation}
  Let $C$ be a small category, $X : C\op \to \Set$, and $Y : C \to \Set$. Recall
  the description of the coend $\int^c X c \times Y c$ as a quotient in $\Set$:
  \[
    \int^{c \in C} X c \times Y c = \frac{\sum_{c \in C} X c \times Y c}{\sim}
  \]
  where for $f : c \to d$, $x \in X d$, $y \in Y c$, we have an identification
  \[
    \left( x, Yf (y) \right) \sim \left( Xf (x), y \right) .
  \]
  We will write $u \otimes v$ for the equivalence class of a pair $(u, v) \in X
  c \times Y c$ will be denoted by. With slight abuse of notations, we will
  write the equivalence relation $\sim$ above as an identity $x \otimes f(y) =
  f(x) \otimes y$.
\end{notation}

\cref{lemma:h!-unit} provides one half of the equivalence between $\oAlg$ and
the Gabriel-Ulmer localisation $\sfS_{n+1,n+2}^{-1} \psh{\OO_{n-k,n+2}}$. The following
proposition will provide the other.

\begin{proposition}
  \label{lemma:h!-counit}
  Let $Y \in \psh{\bbLambda}$. If $\bbSigma \perp Y$, then the counit map $\epsilon_Y
  : h_! h^*Y \to Y$ is an isomorphism.
\end{proposition}
\begin{proof}
  We have to prove that for each $\lambda \in \bbLambda$, the map
  \begin{equation}
    \label{eq:coend-1}
    h_! Y_\lambda
    = \int^{\psi \in \OO_{n-k, n+2}}
    \bbLambda (\lambda, h \psi) \times Y_{h \psi}
    \xto{(\epsilon_Y)_\lambda} Y_\lambda
  \end{equation}
  is a bijection. Consider the map
  \[
    s : Y_\lambda
    \to \int^{\psi \in \OO_{n-k, n+2}}
    \bbLambda (\lambda, h\psi) \times Y_{h \psi}
  \]
  mapping $y \in Y_\lambda$ to $\id_\lambda \otimes y$ (see
  \cref{not:coend-tensor-notation}). It is well-defined, as $h$ is surjective on
  objects, and it is easy to verify that $s (y)$ it is independent of the choice
  of an antecedent $h \nu = \lambda$. Note that $s$ is a section of
  $(\epsilon_Y)_\lambda$, and we proceed to prove that $s$ is surjective. In
  other words, we show that that every element $f \otimes y$, with $f \in
  \bbLambda (\lambda, h \psi)$ for some $\psi \in \OO_{n-k, n+2}$ and $y \in
  Y_\lambda$, is equal to an element of the form $\id_\lambda \otimes y'$, for
  some $y' \in Y_\lambda$.
  \begin{enumerate}
  \item Assume $\lambda = h \phi$ for some $\phi \in \OO_{n-k, n-1}$. Then
    $\bbLambda(\lambda, h \psi) = \OO_{n-k, n+2} (\phi, \psi)$, and $f \otimes y
    = \id_\phi \otimes f (y)$ has the required form.

  \item Assume $\lambda = h \nu = h S[\nu]$ for some $\nu \in \OO_{n+1}$. If
    $\nu$ is degenerate, say $\nu = \itree{\phi}$, then by
    \cref{lemma:doth-isomorphisms}, $h \nu = h \phi$, so we fall into the
    previous case. Thus, we may assume that $\nu$ is not degenerate. Further, by
    \cref{lemma:some-spans,lemma:h-low-dimensions}, we may consider only the
    case where $\psi \in \OO_{n+1}$. By \cref{lemma:diagramatic-lemma}, $f$
    admits a diagram, say
    \[
      \frac{
        \begin{tikzcd} [ampersand replacement = \&]
          \&
          \xi \\
          \nu \ar[ur, sloped, near end, "\src_{[p]}"] \&
          \psi \ar[u, "\tgt"]
        \end{tikzcd}
      }{
        \begin{tikzcd} [ampersand replacement = \&]
          h \nu \ar[r, "f"] \& h \psi ,
        \end{tikzcd}
      }
    \]
    namely $f \cong h \src_{[p]}$. We then have $f \otimes y =
    \id_{\ytree{\omega}} \otimes (h \src_{[p]}) (y)$.

  \item Assume $\lambda = h \omega$ for some $\omega \in \OO_n$. By
    \cref{lemma:doth-isomorphisms}, $h \omega = h \ytree{\omega}$, and we fall
    into the previous case. \qedhere
  \end{enumerate}
\end{proof}

\begin{para}
  \label{def:u-nu-adjunction}
  We write the Gabriel-Ulmer localisation of $\psh{\OO_{n-k, n+2}}$ at the set
  of spine inclusions $\sfS_{n+1, n+2}$ as
  \[
    u : \psh{\OO_{n-k, n+2}} \localisation \sfS_{n+1, n+2}^{-1} \psh{\OO_{n-k,
        n+2}} : N_u .
  \]

  Recall from \cref{th:nerve-theorem-lambda} that we have an
  adjunction $\tau \dashv N$ that exhibits $\oAlg$ as the Gabriel-Ulmer localisation
  $\bbSigma^{-1} \psh{\bbLambda}$. We are now well-equipped to prove that $\oAlg$ is
  equivalent to the localized category $\sfS_{n+1, n+2}^{-1} \psh{\OO_{n-k,
      n+2}}$.
\end{para}

\begin{proposition}
  \label{lemma:h!-equivalence}
  The adjunction $h_! : \psh{\OO_{n-k, n+2}} \rightleftarrows \psh{\bbLambda} : h^*$
  restricts to an adjoint equivalence $\tilde{H}_! \dashv \tilde{H}^*$, as shown
  below.
  \[
    \begin{tikzcd}
      \sfS_{n+1, n+2}^{-1} \psh{\OO_{n-k, n+2}}
      \ar[d, "N_u"']
      \ar[r, shift left = .4em, "\tilde{H}_!"] &
      \bbSigma^{-1} \psh{\bbLambda} \simeq \oAlg
      \ar[d, "N"]
      \ar[l, shift left = .4em, "\tilde{H}^*", "\perp" above] \\
      \psh{\OO_{n-k, n+2}}
      \ar[r, shift left = .4em, "h_!"] &
      \psh{\bbLambda} .
      \ar[l, shift left = .4em, "h^*", "\perp" above]
    \end{tikzcd}
  \]
\end{proposition}
\begin{proof}
  We check the conditions of \cref{lemma:adjunction-to-equivalence}.
  \begin{enumerate}
  \item By \cref{prop:h!-spines}, for all $Y \in \oAlg \simeq \bbSigma^{-1} \psh{\bbLambda}$, we have that $h_! \sfS_{n+1,n+2}
    \perp N_uY$, or equivalently, that $\sfS_{n+1,n+2} \perp h^*N_u Y$.
    Thus $h^* N_u$ factors through the Gabriel-Ulmer localisation $\sfS_{n+1, n+2}^{-1}
    \psh{\OO_{n-k, n+2}}$. Next, by \cref{coroll:h!-spines}, $h_! N_u$
    factors through $\oAlg$.
  \item By \cref{lemma:h!-unit}, if $X \in \sfS_{n+1, n+2}^{-1}
    \psh{\OO_{n-k, n+2}}$, then the unit map $\eta_X : X \to
    h^* h_! X$ is an isomorphism, and dually, by \cref{lemma:h!-counit}, if
    $Y \in \bbSigma^{-1} \psh{\bbLambda}$, then the counit map $\epsilon_Y$ is an
    isomorphism.
    \qedhere
  \end{enumerate}
\end{proof}

\begin{definition}
  \label{def:sfaa}
  Recall the definition of $\sfO$ from~\cref{notn:empty-incl-O}, and let $\sfA =
  \sfA_{k, n} \eqdef \sfO_{< n - k} \cup \sfS_{\geq n + 1}$\index{a@$\sfA$}.
\end{definition}

\begin{theorem}
  [Opetopic nerve theorem]
  \label{th:nerve-theorem-O}
  The adjunction $h : \psh{\OO} \rightleftarrows \oAlg : M$ is reflective, and
  exhibits $\oAlg$ as the Gabriel-Ulmer localisation $\sfA^{-1} \psh{\OO}$ of
  the category $\psh{\OO}$ of opetopic sets at the set of morphisms $\sfA$.
\end{theorem}
\begin{proof}
  Recall from \cref{def:algebraic-realisation-O} that $h$ is the composite
  \[
    \psh{\OO}
    \xto{(-)_{n-k, n+2}} \psh{\OO_{n-k, n+2}}
    \xto{h_!} \psh{\bbLambda}
    \xto{\tau} \oAlg ,
  \]
  and by \cref{lemma:h!-equivalence}, it is isomorphic to the composite
  \[
    \psh{\OO}
    \xto{(-)_{n-k, n+2}} \psh{\OO_{n-k, n+2}}
    \xto{u} \sfS_{n+1, n+2}^{-1} \psh{\OO_{n-k, n+2}}
    \xto{\simeq} \oAlg .
  \]
  The truncation $(-)_{n-k, n+2}$ is the Gabriel-Ulmer localisation at $\sfO_{<
    n-k} \cup \sfB_{> n+2}$. By definition, $u$ is the Gabriel-Ulmer
  localisation at $\sfS_{n+1, n+2}$. Therefore, $h$ is the localisation at
  $\sfO_{< n-k} \cup \sfS_{n+1, n+2} \cup \sfB_{> n+2}$, which by
  \cref{lemma:opetopes-technical:h-lift} is the Gabriel-Ulmer localisation at
  $\sfA$.
\end{proof}

\begin{theorem}
 \label{thm:OAlg-models-of-idempotent-O-theory} 
  For every $k\leq n$, the category $\Oalg kn$ of $k$\=/coloured,
  $n$\=/dimensional opetopic algebras is the category of $\Set$\=/models of an
  idempotent $\OO$\=/sorted theory.
\end{theorem}
\begin{proof}
  Since every map in $\sfA$ has finite domain and codomain, the fully faithful
  right adjoint $M\colon \Oalg kn\subto\psh\OO$ preserves filtered colimits.
  Thus the idempotent monad $Mh$ is finitary, and we conclude
  by~\cref{thm:classification-dep-alg-theories}.
\end{proof}

\begin{corollary}
 \label{cor:Cat-Opd-Cmbd-O-theories} 
  The categories $\Cat,\Opd\pl,\Cmbd$ of small categories, planar
  coloured operads, and Loday's $PT$\=/combinads are all categories of
  $\Set$\=/models of idempotent $\OO$\=/sorted theories.   
\end{corollary}

\section{Homotopy-coherent opetopic algebras}
\label{sec:homot-coher-opetopic-algebras}
We have seen that every category $\Oalg kn$ of opetopic algebras (such as the
categories $\Cat$, $\Opd\pl$ and $\Cmbd$ of small categories, coloured planar
operads and $PT$\=/combinads) can each be described as the category of models in
$\Set$ of an {idempotent $\OO$-sorted theory}. Moreover, we have
characterised their fully faithful nerve functors to $\psh\OO$ by means of
projective sketches on the the category $\OO$ of opetopes.

In this section, we will show that the nature of these sketches allows us
to use a general construction due to Horel (\cite{horel2015model}, generalised
in \cite{HorelCaviglia2016rigidification}) to obtain model structures on
$\s\Oalg kn$ (simplicial opetopic algebras) as well as on $\sPsh\bbLambda_{k,n}$
(simplicial presheaves on $\bbLambda_{k,n}$) that model what we will call
\emph{homotopy\=/coherent opetopic algebras}.


We thus show that, at least for the categories $\Oalg kn$, there exist model
structures on their categories of simplicial objects that (ostensibly) model the
same structures up to homotopy. This can be seen as a very partial answer to the
general conjecture in~\cref{sec:htpical-models-C-cxl-cats} (see
\cref{para:rigidification-htpy-D-algs-conjecture}).

We begin by summarising Horel's construction
(following \cite{HorelCaviglia2016rigidification}) in the case of simplicial presheaves
(this is purely for the convenience of the reader, as the technique 
is relatively new). The reader familiar with it is encouraged to
skip ahead to \cref{para:Lambda-spaces}.

\begin{para}
  [Finite connected sketches]
  \label{para:finite-connected-sketches}
  We will call a pair $(A,K)$ of a small category $A$ and a set $K$ of cocones
  in $A$ a \defn{(projective) sketch}.\footnote{We will diverge from usual
    terminology by considering cocones instead of cones, and contravariant
    presheaves instead of covariant presheaves.} A \emph{model} of $(A,K)$ in
  any category $\cC$ is a functor $A\op\to\cC$ taking every cocone in $K$ to a
  limit cone in $\cC$, and a morphism of models in $\cC$ is just a natural
  transformation. We write $\psh A_K$ for the category of models of $(A,K)$ in
  $\Set$. The full inclusion $\psh A_K\subto \psh A$ is a right adjoint.

  Every cocone $k\in K$ is a pair $(f_k\colon D_k\to A,f_k\to \Delta_a)$ of a diagram
  in $A$ and a natural transformation to the constant diagram on an object $a\in
  A$. Every $k\in K$ thus gives a map $j_k\colon\mquote\colim f_k\to \relHom A
  -a$ of presheaves from the formal colimit in $\psh A$ of the diagram $f_k$ to
  the representable on $a$. If we write $S_K\eqdef\{j_k\mid k\in K\}$ for the
  set of these maps, then the Gabriel-Ulmer localisation of $\psh A$ at $S_K$ is
  precisely the reflective adjunction $\psh A\localisation\psh A_K$. The full
  inclusion $\psh A_K\subto \psh A$ is therefore exactly the nerve functor
  associated to the composite $A\subto \psh A \to \psh A_K$ of the left adjoint
  with the Yoneda embedding of $A$.

  We will say that a sketch $(A,K)$ is \defn{finite} (respectively,
  \defn{connected}) if each diagram in $K$ is indexed by a finite (respectively,
  connected) category. If $(A,K)$ is a finite (respectively, connected) sketch,
  then the image of the dense functor $A\to \psh A_K$ consists of $\omega$-small
  (respectively, connected) objects---equivalently, its nerve functor $\psh
  A_K\subto \psh A$ preserves $\omega$-filtered colimits (respectively,
  arbitrary coproducts).
\end{para}
\begin{remark}
  The cocones in $K$ of a sketch $(A,K)$ are already colimit-cocones in $A$ if
  and only if the functor $A\to \psh A_K$ is fully faithful. As in
  \cite{horel2015model,HorelCaviglia2016rigidification}, we will require this
  condition throughout, and so we assume it holds for every sketch we consider.
\end{remark}

We fix a finite connected (``fully faithful'') sketch $(A,K)$.

\begin{para}
  Consider the reflective adjunction $F:\psh A\localisation\s\psh A_K:U$. By
  taking simplicial objects as in~\cref{para:simpl-algebras-sheaves}, we obtain
  a reflective simplicial adjunction $\tilde F:\sPsh A\localisation \s\psh
  A_K:\tilde U$ between simplicial presheaves on $A$ and models of $(A,K)$ in
  $\SSet$ (equivalently, simplicial objects in $\psh A_K$). Since the right
  adjoint $U$ preserves filtered colimits and arbitrary coproducts, and since
  colimits of simplicial objects are calculated levelwise, the simplicially
  enriched right adjoint $\tilde U$ preserves filtered colimits and
  \emph{tensors} with arbitrary simplicial sets.
\end{para}
\begin{lemma}
  [{\cite[Lem. 4.1]{horel2015model}, \cite[Lem.
    5.1]{HorelCaviglia2016rigidification}}]
  \label{lem:Horel-main-lemma}
  Let $X$ be in $\psh A_K\subto \s\psh A_K$ (thus $X$ is a constant simplicial
  object) and let $i\colon K\mono L$ be a monomorphism in $\SSet$. Then $\tilde
  U$ takes every cocartesian square in $\s\psh A_K$ of the form below to a
  cocartesian square in $\sPsh A$.
  \[
    \begin{tikzcd}[sep=scriptsize]
      K\otimes X\ar[d,tail,"i\otimes 1_X"'] \ar[r]
      &X'\ar[d]\\
      L\otimes X\ar[r]
      &Y'\pomark
    \end{tikzcd}
  \]
\end{lemma}
\begin{proof}
  Since colimits in $\s\psh A_K$ are calculated levelwise, for every
  $[n]\in\Delta$, we have a cocartesian square in $\psh A_K$ as on the left
  below. Since the nerve functor $U\colon \psh A_K\subto \psh A$ preserves coproducts, we can
  apply it to obtain a square in $\psh A$ as on the right below. Then, it suffices to show
  that this square is cocartesian. 
  \[
    \begin{tikzcd}[sep=scriptsize]
      \coprod_{K_n} X\ar[d,tail] \ar[r]
      &X'_n\ar[d]\\
      \coprod_{L_n} X\ar[r]
      &Y'_n\pomark
    \end{tikzcd}
    \qquad\qquad
    \begin{tikzcd}[sep=scriptsize]
      \coprod_{K_n} UX\ar[d,tail] \ar[r]
      &UX'_n\ar[d]\\
      \coprod_{L_n} UX\ar[r]
      &UY'_n
    \end{tikzcd}
  \] If we write $J_n\eqdef L_n-K_n$, then in the square on the left above we
  have an isomorphism $Y'_n\cong X'_n\amalg (\coprod_{J_n}X)$. Since $U$
  preserves coproducts, we have $UY'_n \cong UX'_n\amalg (\coprod_{J_n}UX)$. But
  this is just the pushout
  \[
    UX'_n\coprod_{\amalg_{K_n}UX} (\amalg_{L_n}UX) . \qedhere
  \]
\end{proof}

Recall the \emph{projective} global model structure on $\sPsh A$
  from~\cref{para:simpl-presheaves}.
  \begin{proposition}
    [{\cite[Cor. 4.2]{horel2015model}, \cite[Prop.
      5.2]{HorelCaviglia2016rigidification}}]
    \label{prop:Horel-right-adj-pres-cofs}
    Let $\tilde FI\proj$ and $\tilde FJ\proj$ be the images under $\tilde F$ of
    the sets of generating projective cofibrations and trivial cofibrations.
    Then $\tilde U$ sends pushouts of maps in $\tilde FI\proj$ (respectively,
    $\tilde FJ\proj$) to projective cofibrations (respectively, trivial
    cofibrations).
\end{proposition}
\begin{proof}
  This follows from~\cref{lem:Horel-main-lemma}, since both $\tilde U$ and
  $\tilde F$ preserve tensors, since the maps in $I\proj$ and $J\proj$ are
  exactly of the form $K\otimes a\mono L\otimes a$ for $a\in A$, and since the
  sketch is ``fully faithful''.
\end{proof}

\begin{proposition}
  [{\cite[Thm 5.1]{horel2015model}}{\cite[Prop.
    5.3]{HorelCaviglia2016rigidification}}] The right\=/transferred model
  structure along $\tilde U\colon \s\psh A_K\subto \sPsh A\proj$ exists and is
  combinatorial. Moreover, the right adjoint $\tilde U$ preserves cofibrations.
\end{proposition}
\begin{proof}
  The first claim follows from
  \cref{prop:right-induced-model-structure,prop:Horel-right-adj-pres-cofs}. For
  the second, note that $\tilde U$ preserves filtered colimits and sends
  pushouts of generating cofibrations to cofibrations, thus it preserves
  relative $\tilde FI\proj$\=/cell complexes. Since any functor preserves
  retracts, $\tilde U$ preserves all cofibrations.
\end{proof}

\begin{proposition}
  [{\cite[Prop. 5.5]{HorelCaviglia2016rigidification}}]
 \label{prop:Horel-cofib-unit-iso} 
  Let $X\to Y$ be a cofibration in $\sPsh A\proj$. If the unit $X\to\tilde
  U\tilde FX$ is an isomorphism, then so is the unit $Y\to\tilde U\tilde FY$.
\end{proposition}
\begin{proof}
  Since isomorphisms are closed under retracts and colimits (in the arrow
  category), it suffices to prove the case when $X\to Y$ is a pushout of a
  generating projective cofibration, which is of the form below.
  \[
    \begin{tikzcd}[sep=small]
      K\otimes a\ar[r]\ar[d,tail] &X\ar[r,"\eta"',"\cong"] \ar[d]
      &\tilde U\tilde FX\ar[d]\\
      L\otimes a\ar[r] &Y\pomark\ar[r,"\eta"] &\tilde U\tilde FY
    \end{tikzcd}
  \]
  By \cref{lem:Horel-main-lemma}, the outer square is cocartesian, thus so is
  the right square.
\end{proof}
\begin{proposition}
    [{\cite[Prop. 5.7]{HorelCaviglia2016rigidification}}]
    \label{prop:Horel-rigidification}
    Let $(A,K)$ be a finite connected ``fully faithful'' projective sketch. Then
    the adjunction $\tilde F:\s\psh A_K\rightleftarrows\sPsh A\proj:\tilde U$ is
    a Quillen equivalence.
\end{proposition}
\begin{proof}
  By~\cref{prop:quillen-eq-iff-unit-weq}, it is necessary and sufficient that
  for every cofibrant $X$ in $\sPsh A\proj$, the unit $X\to \tilde U\tilde FX$
  be a global weak equivalence. But $\tilde U$ and $\tilde F$ preserve tensors,
  hence by \cref{prop:Horel-cofib-unit-iso}, $X\to \tilde U\tilde FX$ is an
  isomorphism.
\end{proof}

\begin{para}
  [$\bbLambda$\=/spaces]
  \label{para:Lambda-spaces}
  Recall from \cref{th:nerve-theorem-lambda} that $\Oalg kn$ is the
  Gabriel\=/Ulmer localisation of $\psh\bbLambda$ at the set $\bbSigma$ of
  algebraic spine inclusions (\cref{def:spine-lambda}). Each such
  $\sfS_\lambda\colon S[\lambda]\subto \lambda$ is a subrepresentable whose
  domain $S[\lambda]$ is colimit of representables indexed by the category of
  elements $\OO/S[\nu]$ of the spine of an $(n+1)$\=/opetope $\nu$.
\end{para}

\begin{lemma}
 \label{lem:spine-is-connected} 
  Let $\nu\in\OO$ be of dimension $n\geq 1$. Then the category of elements
  $\OO/S[\nu]$ is connected.
\end{lemma}
\begin{proof}
  We treat the cases $n=1$ and $n>1$ separately. When $n=1$, then $\omega =
  \optOne$ and $S[\omega]=O[\optZero]$, which is representable. When $n>1$,
  $\omega$ is a tree of $(n-1)$\=/opetopes. We proceed by induction. In the base
  case, $\omega = \itree{\psi}$ for some $(n-2)$\=/opetope $\psi$, and
  $S[\omega]=O[\psi]$, which is representable. In the induction step, we use
  \cref{lemma:opetopes-technical:spine-pushout} and the fact that a pushout (or
  any connected colimit) of connected categories is a connected category.
\end{proof}

\begin{proposition}
   \label{prop:existence-horel-model-struc} 
  The commutative triangle of simplicial adjunctions (obtained from
  \cref{eq:triangle-adjunction} by taking simplicial objects)
  \[
    \begin{tikzcd} [column sep = small, row sep = large]
      & \s\oAlg \arrow[dl,shift left=.4em,hook',"\tilde M"below right,
      "\scriptstyle{\bot}" {sloped, above}] \arrow[dr,hook, shift right = .4em,
      "\tilde N" below left, "\scriptstyle{\bot}" {sloped, above}] &
      \\
      \sPsh{\OO}\proj \arrow[ur, shift left = .4em, "\tilde h" above left]
      \arrow[rr,"h_!"above, shift left = .4em]
      & & \sPsh{\bbLambda}\proj 
      \arrow[ul, shift right = .4em, "\tilde\tau" above right] \arrow[ll, shift
      left = .4em,"h^*"below, "\scriptstyle{\bot}" {sloped, above}]
    \end{tikzcd}
  \]
  is a commutative triangle of Quillen adjunctions of simplicial model
  categories, where the model structure on $\s\oAlg$ is right-transferred along
  $\tilde N$ as well as along $\tilde M$. The adjunction
  $\tilde\tau\dashv\tilde N$ is a Quillen equivalence.
\end{proposition}
\begin{proof}
  By~\cref{para:Lambda-spaces,lem:spine-is-connected}, $\oAlg$ is the category
  of models of a finite connected ``fully faithful'' projective sketch on
  $\bbLambda$. Applying~\cref{prop:Horel-rigidification}, the
  right\=/transferred model structure on $\s\oAlg$ along $\tilde N$ exists, and
  $\tilde\tau\dashv\tilde N$ is a Quillen equivalence. Since $h_!\dashv h^*$ is
  a Quillen adjunction between projective model structures, the diagram is a
  triangle of Quillen adjunctions. Finally, since $h\colon \OO\to\bbLambda$ is
  surjective on objects (\cref{prop:h-surjective}), $h^*$ preserves and reflects
  projective fibrations and global weak equivalences, so the model structure on
  $\oAlg$ is also right-transferred along $\tilde M$.
\end{proof}

\begin{remark}
  The adjunction $\tilde h\dashv\tilde M$ is \emph{not} a Quillen equivalence.
  Were it so, then by \cref{prop:quillen-eq-iff-unit-weq}, the unit on every
  representable $\OO\subto\sPsh\OO$ would be a global weak equivalence, but it
  is relatively easy to show that this is not the case.
\end{remark}

\begin{definition}
  The model structure for \defn{homotopy-coherent $(k,n)$\=/opetopic algebras}
  is the left Bousfield localisation of the projective model structure
  $\sPsh\bbLambda\proj$ at the set $\Sigma$ of algebraic spine inclusions
  (\cref{def:spine-lambda}). We write it as $\sPsh\bbLambda\proj^l$.
\end{definition}

\begin{definition}
  The \defn{local model structure} on $\s\oAlg$ is the left Bousfield
  localisation of the right-transferred model structure
  of~\cref{prop:existence-horel-model-struc} at the set
  $\LL\tilde\tau(\Sigma)$ that is the image of $\Sigma$ under the left derived
  functor of $\tilde\tau\colon \sPsh\bbLambda\proj\to\s\oAlg$. We write it as
  $\s\oAlg^l$.
\end{definition}

\begin{theorem}
  [Rigidification of homotopy-coherent opetopic algebras]
  \label{thm:rigidification-Oalg}
  The adjunction
  $\tilde\tau :\sPsh\bbLambda\proj^l\localisation\s\oAlg^l:\tilde N$ is a
  Quillen equivalence.
\end{theorem}
\begin{proof}
  This follows from the Quillen equivalence
  of~\cref{prop:existence-horel-model-struc} and a general fact about left
  Bousfield localisations (\cref{para:left-bousfield-localisation}).
\end{proof}

\begin{para}
  [Segal $\bbLambda$\=/spaces]
  \label{para:segal-Lambda-spaces}
  When $\oAlg = \Cat$ ($k=n=1$), recall that $\bbLambda=\bDelta$. Then it is
  easily seen that the model structure $\sPsh\bDelta\proj^l$ is the projective
  version of Rezk's model structure for \emph{Segal spaces}\cite[Thm
  7.1]{rezk2001model}.

  When $\oAlg = \Opd\pl$ ($k=1,n=2$), recall that $\bbLambda=\bOmega\pl$ is the
  planar version of Moerdijk-Weiss's dendroidal category. Then it is easily seen
  that the model structure $\sPsh(\bOmega\pl)\proj^l$ is the projective version
 of Cisinski\=/Moerdijk's model structure for \emph{Segal dendroidal
    spaces} \cite[Prop. 5.5]{cisinski2013dendroidal}. 
\end{para}



\part{Localisations of locally presentable \texorpdfstring{\oo}{infinity}-categories}
\label{part:part2}

\chapter*{Preliminaries}

Throughout \cref{part:part2}, our goal will be to prove theorems in
\oo\=/category ($(\infty,1)$\=/category) theory. Despite the
``model-independent'' theory of \oo\=/categories still being nascent, I have
made the choice to work entirely within it. The payoff is twofold:
\begin{enumerate}
\item The theory of locally presentable \oo\=/categories perfectly subsumes the
  theory of locally presentable 1\=/categories. Thus the reader unfamiliar with
  \oo\=/categories has only to take ``category'' to mean ``1\=/category''
  throughout, and to consider that every space is a set, in order to recover
  precise results\footnote{Even when restricted to the particular case of
    locally presentable $1$-categories, the definition of \emph{pre-modulator}
    and the results of \cref{chap:SOA} seem (to the best of my knowledge) to be
    new.} for locally presentable 1\=/categories, especially those results that
  are (implicitly or explicitly) used in \cref{part:part1}.
\item The motivated reader familiar with quasicategories as models for
  \oo\=/categories should \emph{in principle} be able to reconstruct for
  themselves all the arguments in \cref{part:part2} within quasicategories. I
  have included some remarks that refer to material from \cite{LurieHT,LurieHA} and
  other sources where necessary. On the other hand, the reader not of
  quasicategorical persuasion should still be able to interpret at least the
  statements of the main theorems in the model of their predilection.
\end{enumerate}

\section*{Change of vocabulary}
This shift in point of view permits (demands) a change in vocabulary. We will
say \emph{category} instead of \emph{\oo\=/category} and \emph{space} instead of
\emph{\oo\=/groupoid}. In particular, we will say \emph{topos} for what is
called an \emph{\oo\=/topos} in \cite{LurieHT}. If necessary, we will explicitly
use the terms ``1\=/category'' and ``1\=/topos''. An \emph{$n$\=/topos} is the
full subcategory of $n$\=/truncated objects of a topos.

We write $\cS$ for the category of spaces and $\CAT$ for the category of
categories. We denote the space of morphisms between two objects $x$ and $y$ of
a category $\cC$ by $\fun x y$, or by $\relmap \cC x y$ if the category $\cC$
needs to be recalled. We will say that a map in a category is \emph{invertible}
or an \emph{equivalence} if it admits both a left and a right inverse, and is a
\emph{monomorphism} if it is $(-1)$\=/truncated (its diagonal is invertible).
For any functor $f\colon \cA\to\cB$ and any object $x$ in $\cB$, we write
$\cA\relcomma fx$ and $x\relcomma f\cA$ for the over and under categories. For a
functor $g\colon \cC\to\cB$, we write $f\comma g$ (or sometimes $\cA\comma\cC$)
for the comma category.

We write $\Psh C\eqdef\cS^{C\op}$ for the category of presheaves of spaces on a
category $C$. We say \emph{locally presentable} or \emph{presentable category}
for the notion of \emph{presentable \oo\=/category} of \cite{LurieHT}.

We say (\emph{op})\emph{fibration} for the notion of (\emph{co})\emph{Cartesian fibration} of
\cite{LurieHT}. We say \emph{fibration} (respectively, \emph{opfibration})
\emph{in groupoids} for the notion of \emph{right} (respectively, \emph{left})
\emph{fibration} from \cite{LurieHT}.


\chapter{Factorisation systems}
\label{chap:fact-systems}

\section{Arrow categories}
\label{sec:arrow-categories}
The \emph{arrow} is the 1\=/category $I=\{ 0<1 \}$, and the \defn{arrow category} of
$\cC$ is the functor category $\cC\arr\eqdef \cC^I$. We write the domain (source) and
codomain (target) functors as $\s{},\t{}\colon\cC\arr\to\cC$ respectively.

\begin{para}[Cartesian closure]
  \label{para:arrow-cat-cartesian-closure}

  Let $(\cC,\times,1,\intmap--)$ be a cartesian closed category. For morphisms
  $f\colon A\to B,g\colon X\to Y$ in $\cC$, their product in the arrow category $\cC\arr$ is the
  map $f\times g\colon A\times X \to B\times Y$ below.
  \[
    \begin{tikzcd}
      A\times X \ar[r,"A\times g"] \ar[d,"f\times X"'] \ar[dr,"f\times g"]
      &A\times Y \ar[d,"f\times Y"]\\
      B\times X\ar[r,"B\times g"]
      &B\times Y
    \end{tikzcd}
  \]
  If the cartesian square
  \[
    \begin{tikzcd}[sep=scriptsize]
      \intmap fg \pbmark \ar[r] \ar[d,"\intbracemap fg"']
      &\intmap AX \ar[d]\\
      \intmap BY \ar[r] &\intmap AY
    \end{tikzcd}
  \] exists in $\cC$, then $\intbracemap fg$, along with the map $\intbracemap
  fg \times f\to g$ below, is the internal hom of
  $f$ and $g$ in $\cC\arr$. We have an equivalence of spaces
  $\relmap{\cC\arr}h{\intbracemap fg} \simeq \relmap{\cC\arr}{f\times h}g$ for
  every $f,g,h \in \cC\arr$.
  \[
    \begin{tikzcd}[sep=scriptsize]
      \intmap fg\times A \ar[r] \ar[d] \ar[dd,bend right=20,start anchor=south
      west,end anchor=north west,"\intbracemap fg\times f"'] 
      &\intmap AX\times A \ar[r] \ar[d]
      &X\ar[d] \ar[dd,bend left,"g"]\\
      \intmap BY\times A \ar[r] \ar[d]
      &\intmap AY\times A \ar[r]
      &Y\ar[d,equals]\\
      \intmap BY\times B\ar[rr]
      &&Y
    \end{tikzcd}
  \]
  Hence if $\cC$ has finite limits, then $(\cC\arr,\times,1_1,\intbracemap--)$ is
  cartesian closed, and $\intbracemap--$ is an enrichment
  of $(\cC\arr,\times,1_1)$ over itself. Namely, for every $f,g,h\in\cC\arr$, we
  have an equivalence $\intbracemap{f\times h}g \simeq \intbracemap
  h{\intbracemap fg}$ in $\cC\arr$. In particular, the arrow category
  $\cS\arr$ of the category of spaces is cartesian closed. We write its
  internal hom as $\bracemap--$.

  For any category $\cC$ and maps $f\colon A\to B,g\colon X\to Y$ in $\cC$,
  their \defn{external hom} $\bracemap fg\colon \relmap{\cC\arr} fg \to
  \relmap\cC BY$ is the left vertical map in the canonical cartesian square
  below in $\cS$.
  \[
    \begin{tikzcd}[sep=scriptsize]
      \relmap{\cC\arr} fg \pbmark \ar[r] \ar[d,"\bracemap fg"']
      &\relmap\cC  AX \ar[d]\\
      \relmap\cC  BY \ar[r] &\relmap\cC  AY
    \end{tikzcd}
  \]
\end{para}

\begin{remark}
  Recall that a category $\cC$ is cartesian closed if it has finite 
  products and if $(\cC,\times,1)$ is a \emph{closed} (symmetric) monoidal
  category \cite[4.1.1.15]{LurieHA}. Following \cite[4.1.1.16]{LurieHA},
  the universal property of $\intbracemap fg \times f\to g$ can be checked at
  the level of homotopy categories. Following \cite[4.2.1.32]{LurieHA},
  $(\cC,\times,1)$ is enriched over itself if and only if it is cartesian
  closed.
\end{remark}

\begin{para}[Pushout-product in spaces]
 \label{para:arrow-cat-pushout-product} 
 The arrow $I$ has a symmetric monoidal structure $(I,\otimes,1)$ given by
 $0\otimes0 = 0\otimes 1=1\otimes 0 = 0$ and $1\otimes 1=1$. By Day convolution,
 there is a symmetric monoidal structure on the arrow category $\cS\arr$ of the
 cartesian symmetric monoidal category $(\cS,\times,1)$. The monoidal product
 $-\pp-$ is called the \defn{pushout-product}, and its unit is $\emptyset \to
 1$. By the formula for Day convolution, the pushout-product of two maps
 $f\colon A\to B$ and $g\colon X\to Y$ is the cocartesian gap map of
 the outer (cartesian) square in $\cS$ below.
 \[
    \begin{tikzcd}[sep=scriptsize]
      A\times X \ar[r] \ar[d]
      &A\times Y \ar[d] \ar[ddr,bend left]\\
      B\times X \ar[r] \ar[drr,bend right=15]
      &B\times X \coprod_{A\times X} A\times Y \pomark \ar[dr,dashed,"f\pp g"]\\
      &&B\times Y
    \end{tikzcd}
  \]
  The category $(\cS,\times,1)$ is cartesian closed; since $\cS$ is locally
  presentable, this is equivalent to the product functor $-\times-$ preserving colimits
  (in $\cS$) in each variable. Hence, since $-\pp-$ is a Day convolution, it too
  preserves colimits (in $\cS\arr$) in each variable. Since $\cS\arr$ is also locally
  presentable, this implies that the symmetric monoidal category
  $(\cS\arr,\pp,\emptyset\to 1)$ is closed, and its internal hom $\pbh--$ is
  called the \defn{pullback-hom}. The pullback-hom of $f\colon A\to B$ and
  $g\colon X\to Y$ is the cartesian gap map of the outer square in $\cS$ below.
  \[
    \begin{tikzcd}[sep=scriptsize]
      \fun BX \ar[drr,bend left=20,shift left=1.5,"\fun Bg"] \ar[ddr,bend
      right,"\fun fX"']
      \ar[dr,dashed,"\pbh{f}{g}"]\\
      &\fun fg \pbmark \ar[r,"\bracemap fg"] \ar[d]
      &\fun BY \ar[d,"\fun fY"]\\
      &\fun AX \ar[r,"\fun Ag"'] &\fun AY
    \end{tikzcd}
  \]
  The pullback-hom gives $(\cS\arr,\pp,\emptyset\to 1)$ an enrichment over itself,
  namely, we have an equivalence $\pbh {f\pp g}h\simeq \pbh f {\pbh
    gh}$ for every $f,g,h\in\cS\arr$.
\end{para}

\begin{para}[Internal pushout-product]
  \label{para:internal-pushout-prod}
  When $(\cC,\otimes,\bbone,\intmap--)$ is a locally presentable, symmetric
  monoidal closed category, then the same constructions as
  in~\cref{para:arrow-cat-pushout-product} define a symmetric monoidal closed
  structure $(\cC\arr,\boxtimes,0\to\bbone,\intpbh--)$ on the arrow category of
  $\cC$. We will call this the \defn{internal pushout-product} and
  \defn{internal pullback-hom} in $\cC\arr$.
\end{para}

\begin{para}[Tensor, enrichment, cotensor]
  \label{para:arrow-cat-tensor-enrichment}
  Let $\cC$ be a locally presentable category. Recall that this implies that $\cC\arr$ is
  also locally presentable. As $\cC$ is cocomplete, it is canonically
  tensored over $\cS$, and we write $\times\colon \cS\times \cC\to\cC$ for the
  functor that is cocontinuous in each variable. For any space $A$ and any
  $X\in\cC$, $A\times X = \coprod_A X$ is the colimit in $\cC$ of the constant
  diagram on $X$ indexed by the space $A$. The tensor $\times$ is compatible
  with the enrichment of $\cC$ over $\cS$, namely, we have an equivalence of
  spaces $\relmap \cC {A\times X}Y \simeq \relmap \cS A {\relmap \cC XY}$ for
  every space $A$ and every $X,Y$ in $\cC$. Moreover, for every space $A$, the
  functor $A\times-\colon\cC\to\cC$ preserves colimits and (since $\cC$ is
  locally presentable) thus has a right
  adjoint; hence $\cC$ is also cotensored over $\cS$.

  Next, $\cC\arr$ is locally presentable and is canonically tensored over
  $(\cS\arr,\times,1_1)$ via
  \begin{align*}
    \times\colon\cS\arr\times\cC\arr&\to\cC\arr\\
    (f\colon A\to B,g\colon X\to Y)&\mapsto f\times g\colon A\times X\to B\times Y . 
  \end{align*}
  This functor is also cocontinuous in each variable. Since $\cS\arr$ is locally
  presentable, for every $f\in\cC\arr$, the functor $-\times f\colon
  \cS\arr\to\cC\arr$ has a right adjoint, and it is easily verified that this is
  given by the \defn{external hom} $\bracemap f-$ (see
  \cref{para:arrow-cat-cartesian-closure}). This defines an enrichment of
  $\cC\arr$ over $(\cS\arr,\times,1_1,\bracemap--)$, and we have an equivalence
  $\bracemap {f\times g}h \simeq \bracemap f{\bracemap gh}$ in $\cS\arr$ for
  every $f\in\cS\arr$ and $g,h\in\cC\arr$. Moreover $\cC\arr$ is also cotensored
  over $(\cS\arr,\times,1_1)$, since for every $f\colon A\to B$ in $\cS\arr$,
  the functor $f\times-\colon \cC\arr\to\cC\arr$ has a right adjoint (as
  $\cC\arr$ is also locally presentable).

  Finally, $\cC\arr$ is also canonically tensored over the symmetric monoidal
  category $(\cS\arr,\pp,\emptyset\to 1)$, via the \defn{external
    pushout-product} $\pp\colon \cS\arr\times\cC\arr\to\cC\arr$ that takes
  $f\colon A\to B$ in $\cS$ and $g\colon X \to Y$ in $\cC$ to the cocartesian
  gap map $f\pp g$ in $\cC$ below.
  \[
    \begin{tikzcd}[sep=scriptsize]
      A\times X \ar[r] \ar[d]
      &A\times Y \ar[d] \ar[ddr,bend left]\\
      B\times X \ar[r] \ar[drr,bend right=15]
      &B\times X \coprod_{A\times X} A\times Y \pomark \ar[dr,dashed,"f\pp g"]\\
      &&B\times Y
    \end{tikzcd}
  \]
  Since $\cS\arr$ is locally presentable, $\cC\arr$ is enriched over
  $(\cS,\pp,\emptyset\to 1,\pbh--)$ via the \defn{external pullback-hom}, which
  takes $f\colon A\to B$ and $g\colon X \to Y$ in $\cC$ to the cartesian gap map
  $\pbh f g$ below in $\cS$.
  \[
    \begin{tikzcd}[sep=scriptsize]
      \relmap {\cC} BX \ar[drr,bend left=20,shift left=1.5] \ar[ddr,bend right]
      \ar[dr,dashed,"\pbh{f}{g}"]\\
      &\relmap {\cC\arr} fg \pbmark \ar[r,"\bracemap fg"] \ar[d]
      &\relmap {\cC} BY \ar[d]\\
      &\relmap {\cC} AX \ar[r] &\relmap {\cC} AY
    \end{tikzcd}
  \]
  We thus have an equivalence $\pbh{f\pp g}h \simeq \pbh f {\pbh gh}$ for every
  $f\in\cS\arr$ and $g,h\in\cC\arr$. Moreover, since $\cC\arr$ is locally
  presentable, $\cC\arr$ is also cotensored over $(\cS,\pp,\emptyset\to 1)$.
\end{para}

\begin{lemma}[Absorption properties of invertible maps]
\label{lem:pp-pbh-iso-absorption}
Let $f$ and $g$ be two maps in $\cC$ and let $u$ be a map in $\cS$.
\begin{enumerate}
\item If $u$ or $f$ is invertible then so is $u\pp f$.
\item If $f$ or $g$ is invertible then so is $\pbh g f$.
\end{enumerate}
\end{lemma}
\begin{proof}
  Easy verification using the previous diagrams.
\end{proof}

\begin{remark}
  Day convolution is defined in \cite[2.2.6.17, 4.8.1.13]{LurieHA}.
\end{remark}

\begin{remark}
  Recall from \cite[4.8.1]{LurieHA} that $\cS$ is the unit of the symmetric
  monoidal category $\CAT\pres$ of presentable categories (and
  cocontinuous functors). Hence every presentable category is canonically
  tensored (i.e. is a module) over (the commutative monoid in $\CAT\pres$)
  $(\cS,\times,1)$. Next, $(\cS\arr,\times,1_1)$ and $(\cS\arr,\pp,\emptyset\to
  1)$ are commutative monoids in $\CAT\pres$. The cocontinuous functor
  $\cS\to\cS\arr$ that takes a space $A$ to its identity map $1_A$ (and that is
  uniquely determined by the unit $1_1$ of $(\cS\arr,\times,1_1)$) is a
  symmetric monoidal functor from $(\cS,\times,1)$ to $(\cS\arr,\times,1_1)$.
  Similarly, the cocontinuous functor $\cS\to\cS\arr$ taking $A$ to $\emptyset
  \to A$, and uniquely determined by the unit $\emptyset\to 1$ of
  $(\cS\arr,\pp,\emptyset\to 1)$, is symmetric monoidal. Thus the tensor product
  $\cS\arr\otimes\cC\simeq \cC\arr$ of presentable
  categories has the structure of a free $(\cS\arr,\times,1_1)$\=/module as
  well as a free $(\cS\arr,\pp,\emptyset\to 1)$\=/module \cite[4.2.4.8]{LurieHA},
  which are the tensor structures of~\cref{para:arrow-cat-tensor-enrichment}.
  Cotensors and enrichment follow from \cite[4.2.1.33]{LurieHA}. 
\end{remark}

\begin{remark}
 \label{rem:pp-pbh-cocomplete-cats} 
  Recall that the the category $\CAT\pres$ is a full symmetric monoidal
  subcategory of the symmetric monoidal category $\CAT\cc$ of cocomplete
  categories and cocontinuous functors. We note that when $\cC$ is a
  cocomplete (but not necessarily presentable) category, we still have the
  equivalence $\cS\arr\otimes \cC \simeq \cC^{I\op}\simeq\cC\arr$, and all of
  \ref{para:arrow-cat-tensor-enrichment} holds for $\cC\arr$, \emph{except} for
  the existence of cotensor structures.
\end{remark}

\begin{remark}
  For \emph{any} commutative monoid $\cV$ in $\CAT\pres$ (namely, a closed
  symmetric monoidal presentable category), the category of $\cV$\=/enriched,
  $\cV$\=/cocomplete (\ie closed under all weighted colimits) categories is
  equivalent to the category $\cV\Mod(\CAT\cc)$ of $\cV$\=/modules in $\CAT\cc$.
  Moreover, the category of \emph{presentable} $\cV$\=/enriched categories
  is equivalent to $\cV\Mod(\CAT\pres)$ (a proof of this is
  \cite[A.3.8]{mazelgee2021universal}, and the proof for $\cV\Mod(\CAT\cc)$ is
  the same). Hence we could start with \emph{any} commutative monoid
  $(\cC,\otimes,\bbone,\intmap--)$ in $\CAT\pres$ and describe the
  pushout-product and pullback-hom structures for arrow $\cC$-enriched
  categories in exactly the same way, using the internal pushout-product and
  pullback-hom of $\cC\arr$.
\end{remark}

\begin{para}[Pullback-hom for arbitrary arrow categories]
 For any small category $\cC$, post-composing with the Yoneda embedding
$\cC\subto\Psh\cC$ gives a full inclusion of arrow categories $\cC\arr\subto
\Psh\cC\arr$. As $\Psh\cC$ is locally presentable, it has an external
pullback-hom (\cref{para:arrow-cat-pushout-product}), namely it is enriched over
$(\cS\arr,\pp,\emptyset\to 1,\pbh--)$. Hence so is $\cC\arr$, as it is a full
subcategory (\cite[4.2.1.34]{LurieHA}). Further, it is easily seen that $\pbh fg $ is the cartesian gap map
$\relmap\cC
{B}{X}\to\relmap{\cC\arr}fg$  for every $f\colon A\to B,g\colon X\to Y$ in
$\cC\arr$.
\end{para}

\begin{remark}
  For any (not necessarily small) category $\cC$, we have a full inclusion
  $\cC\subto\tilde\cC$ (that may even be supposed colimit-preserving) of $\cC$
  into a cocomplete category $\tilde\cC$ (see \cite[5.3.6.2]{LurieHT}).
  By~\cref{rem:pp-pbh-cocomplete-cats}, $\tilde\cC\arr$ is enriched over
  $(\cS\arr,\pp,\emptyset\to 1,\pbh--)$ via the external pullback hom
  $\pbh--\colon(\tilde\cC\arr)\op\times\tilde\cC\arr\to\cS\arr$. Hence, by
  restriction, so is $\cC\arr$.
\end{remark}

\begin{remark}
 \label{rem:pbh-preserves-limits} 
 For any category $\cC$, the formula $\pbh fg \colon \relmap\cC{\t f}{\s
   g}\to\relmap{\cC\arr}fg$ and the fact that $\t,\s\colon\cC\arr\to\cC$
 preserve all limits and colimits, together tell us that the enrichment functor
 $\pbh--\colon (\cC\arr)\op\times\cC\arr \to \cS\arr$ preserves limits in each
 variable. In other words, conical limits and colimits in $\cC\arr$ for the
 enrichment $\pbh--$ coincide with ordinary limits and colimits in $\cC\arr$.
\end{remark}

\begin{para}[Diagonal, codiagonal]
  \label{para:diag-codiag}
  The diagonal of a map $f\colon A\to B$ in a category $\cC$ is the map $\Delta
  f\colon A\to A\times_BA$ (should the pullback of $f$ along itself exist). When
  they exist, the \defn{iterated diagonals} of $f$ are defined as
  $\Delta^0f\eqdef f$ and $\Delta^{n+1}f\eqdef \Delta(\Delta^n f)$. Dually, its
  \defn{iterated codiagonals} $\nabla^n f$ are the iterated diagonals in
  $\cC\op$.\footnote{In a $1$\=/category, all iterated diagonals $\Delta^nf$ and
    codiagonals $\nabla^nf$ are isomorphisms for $n\geq 2$.}

  For $n\geq -1$, let $S^n$ denote the object of $\cS$ that is the
  $n$-dimensional sphere, with $S^{-1}=\emptyset$, and let $s^n\colon S^n\to 1$
  denote the canonical map to the point. An easy calculation (as
  unreduced suspensions, see \cite[Example 3.2.1]{ABFJblakersmassey}) shows that
  $s^m\pp s^n = s^{n+m+1}$.

  A further easy calculation then shows that for any map $f\in\cC\arr$, the codiagonal
  $\nabla^n f$ is the
  tensor $s^{n-1}\pp f$, and dually, the diagonal $\Delta^n f$ is the cotensor
  $\pbh{s^{n-1}}f$. We thus have equivalences $\Delta^n \pbh fg \simeq \pbh
  f{\Delta^ng} \simeq \pbh{\nabla^nf}g$ in $\cS\arr$.
\end{para}

\begin{definition}
  Recall that a map $f\colon A\to B$ in $\cS$ is a \defn{cover} (or an \defn{effective
    epimorphism}) if it induces a surjection $\pi_0(A) \to \pi_0(B)$ on
  connected components.
\end{definition}

\begin{lemma}[Whitehead completion]
\label{lem:Whitehead}
A map $w$ in $\cS$ is invertible if and only if all its iterated diagonals $\Delta^n w =\pbh {s^{n-1}} w$ are covers.
\end{lemma}

\begin{proof}
Recall that for $-2\leq k\leq \infty$ a map $w$ in $\cS$ is $k$\=/connected if and only if all maps $\Delta^n w$ are covers for $n\leq k+1$~\cite[6.5.1.18]{LurieHT}.
The result follows from the fact that $\infty$\=/connected maps in $\cS$ are invertible by Whitehead's theorem.
\end{proof}

    
    
    




\section{Orthogonality}

\begin{definition}
  \label{def:orthogonality}
    Two arrows $f:A\to B$ and $g:X\to Y$ of a category $\cC$ are said to be \defn{orthogonal}
  if the map $\pbh f g$ in $\cS$ is invertible. We write this relation as
  $f\perp g$. 
\end{definition}
It says that for any commutative square $\alpha:f\to g$ as below,
  there exists a unique diagonal lift $\t f \to \s g$ making the triangles
  commute.
\[
  \begin{tikzcd}[sep=scriptsize] 
    \s f \ar[r,"\alpha_\s"] \ar[d,"f"']
    &\s g\ar[d,"g"] \\
    \t f \ar[r,"\alpha_\t"'] \ar[ur,dashed,"\exists!"]
    &\t g
  \end{tikzcd}
\]

\begin{para}[Orthogonal systems]
  \label{para:orth-systems}
Given a diagram $W\to \cC\arr$ we define (replete) full subcategories of
$\cC\arr$ as follows: $W^\bot \eqdef \{g\, |\, \forall f
\in \im W, f\perp g\}$ and $^\bot W \eqdef \{f\, |\, \forall g \in \im W,
f\perp g\}$. Note that the generating data $W\to
\cC\arr$ can be an arbitrary diagram, indexed by a space or any category. 

A pair $(\cL,\cR)$ of full subcategories of $\cC\arr$ is an
\defn{orthogonal system on $\cC$} if we have $\cL = {^\bot}\cR$ and $\cR = \cL^\bot$. For any
diagram $W\to \cC\arr$, the pair $(\cL_W,\cR_W) := \left(^\bot (W^\bot),
  W^\bot\right)$ is an orthogonal system, said to be \emph{generated by $W$}. The pair
$(\cL_W,\cR_W)$ depends only on the essential image of the functor $W\to
\cC\arr$. 
\end{para}

\begin{lemma}
  \label{lem:prop-os}
  Let $(\cL,\cR)$ be an orthogonal system.
  \begin{enumerate}
  \item
    \label{lem:prop-os:1}
    $\cL$ is absorbing for $\pp$: for $u\in\cS\arr$ and $f\in\cL$, the tensor
    $u\pp f$ is in $\cL$.
  \item
    \label{lem:prop-os:1-dual}
    $\cR$ is absorbing for $\pbh--$: for $u\in\cS\arr$ and $f\in\cR$, the
    cotensor $\pbh u f$ is in $\cR$.
  \item
    \label{lem:prop-os:2}
    $\cL$ is stable by colimits in $\cC\arr$.
  \item
    \label{lem:prop-os:2-dual}
    $\cR$ is stable by limits in $\cC\arr$.
  \item
    \label{lem:prop-os:3}
    $\cL$ is stable under right-cancellation (if $f,gf\in\cL$, then $g\in\cL$).
  \item
    \label{lem:prop-os:3-dual}
    $\cR$ is stable under left-cancellation (if $g,gf\in\cR$, then $f\in\cR$).
  \end{enumerate}
\end{lemma}

\begin{proof}
By duality, it is enough to prove the statements for the left class $\cL$.

\noindent \ref{lem:prop-os:1} For any $g$, we have $\pbh {u\pp f} g = \pbh u
{\pbh f g}$. By the absorption properties of \cref{lem:pp-pbh-iso-absorption},
this map is invertible when $\pbh fg$ is so. Since $\cL = {^\bot}\cR$, this
proves the first assertion.

\noindent \ref{lem:prop-os:2}
Let $f_i$ be a diagram in $\cL$. We have $\pbh {\colim_i f_i} g = \lim_i \pbh
{f_i} g$ (see \ref{rem:pbh-preserves-limits}).
For any $g\in \cR$, each $\pbh {f_i} g$ is invertible, and so is their limit.
Thus $\colim_i f_i$ is in $\cL$.

\noindent \ref{lem:prop-os:3}
For any composable maps $f\colon A\to B$ and $g\colon B\to C$ in $\cC$, the
obvious left-hand square below in $\cC\arr$ is cocartesian (and cartesian). For
any $h\in\cC\arr$, the
right-hand square is thus cartesian in $\cS\arr$. If $f,gf\in\cL$ and $h\in\cR$,
then $\pbh gh$ is a limit of equivalences, thus is an equivalence.
\[
  \begin{tikzcd}[sep=scriptsize]
    f \ar[r]\ar[d] \pbmark
    &gf \ar[d]\\
    1_B\ar[r]
    &g \pomark
  \end{tikzcd}
  \qquad
  \begin{tikzcd}[sep=scriptsize]
    \pbh g h \ar[r]\ar[d] \pbmark
    &\pbh {gf} h \ar[d]\\
    \pbh {1_B}h\ar[r]
    &\pbh f h
  \end{tikzcd}
\]\qedhere
\end{proof}

\begin{definition}
  \label{def:weak-orthogonality}
  Two arrows $f:A\to B$ and $g:X\to Y$ of a category $\cC$ are said to be
  \defn{weakly orthogonal} if the map $\pbh f g$ in $\cS$ is a cover. We write
  this relation as $f\pitchfork g$.
\end{definition}
It says that for any commutative square $\alpha:f\to g$ as below, there exists a
(non-unique) diagonal lift $\t f \to \s g$ making the triangles commute.
\[
  \begin{tikzcd}[sep=scriptsize] 
    \s f \ar[r,"\alpha_\s"] \ar[d,"f"']
    &\s g\ar[d,"g"] \\
    \t f \ar[r,"\alpha_\t"'] \ar[ur,dashed,"\exists"]
    &\t g
  \end{tikzcd}
\]

\begin{para}[Weak orthogonal systems]
  \label{para:wk-orth-systems}
  As in \cref{para:orth-systems},
  given any diagram $W\to \cC\arr$ we define (replete) full subcategories
  $W^\pitchfork$ and $^\pitchfork W$ of $\cC\arr$ using the relation
  $\pitchfork$. We say that a
  pair $(\cA,\cB)$ of full subcategories of $\cC\arr$ is a \defn{weak orthogonal
    system} if $\cB = {^\pitchfork}\cA$ and $\cA = \cB^\pitchfork$. For any
  diagram $W\to \cC\arr$, the pair $(\cA_W,\cB_W) := \left(^\pitchfork
    (W^\pitchfork), W^\pitchfork\right)$ is a weak orthogonal system, said to be
  \emph{generated by $W$}. The pair $(\cA_W,\cB_W)$ depends only on the
  essential image of the functor $W\to \cC\arr$.
\end{para}

\begin{remark}
  The relation $\pitchfork$ is exactly that of \cite[Def. 1.4.1]{LurieDAGX}.  
\end{remark}

\begin{lemma}
  \label{lem:prop-wos}
  Let $(\cA,\cB)$ be a weak orthogonal system. Then $\cA$ is stable by discrete
  sums (coproducts indexed by a set), cobase change (pushout along an arbitrary
  map) and transfinite composition.
\end{lemma}
\begin{proof}
  Stability under discrete sums and cobase change follows since covers in $\cS$
  are stable under discrete products and base change. Transfinite composition
  follows from an argument similar to that of \cite[Prop.
  1.2.9]{LurieDAGX}.
\end{proof}

\begin{remark}
  \label{rem:wos-defect-sum}
  The previous statement about discrete sums cannot, \emph{a priori}, be
  extended to more general sums, namely colimits indexed by a space. Let
  $(\cA,\cB)$ be a weak orthogonal system, $f$ be a map in $\cB$, and $w: I \to
  \cA$ a family of maps in $\cA$ indexed by a space $I$. Although each map $\pbh {w_i} f$
  is a cover, the product $\prod_{i\in I} \pbh {w_i} f$ is generally not a
  cover, unless $I$ is a set. As a counter-example, consider the weak orthogonal
  system on $\cS$ generated by the map $\emptyset\to 1$. A map is in the right
  class if and only if it is a cover. The sum $\coprod_I(\emptyset\to 1)$ indexed by any
  space $I$ is the map $\emptyset\to I$, however not all covers are weakly right
  orthogonal to $\emptyset\to S^1$.
\end{remark}

The following lemma relates weak orthogonality to orthogonality in a category
$\cC$ in which all maps have diagonals and codiagonals. 
\begin{lemma}
  \label{lem:wos-with-codiag-is-os}
  A pair $(\cL,\cR)$ is an orthogonal system if and only if it is a weak orthogonal system
  whose left class is stable by codiagonal.
\end{lemma}
\begin{proof}
  Let $(\cL,\cR)$ be an orthogonal system. By \cref{para:diag-codiag} and
  \cref{lem:prop-os}, $\cL$ is stable by codiagonal, and dually, $\cR$ is stable
  by diagonal. We need to prove that $(\cL,\cR)$ is a weak orthogonal system,
  namely that $\cL^\pitchfork = \cR$ and $\cL ={^\pitchfork}\cR$. Clearly, $\cR
  \subset \cL^\pitchfork$ and $\cL\subset {^\pitchfork}\cR$. Let $g$ be a map in
  $\cL^\pitchfork$, then for any $f\in\cL$, the map $\pbh fg $ is a cover in
  $\cS$. Since $\cL$ is stable by codiagonal, all maps $\pbh {\nabla^nf} g
  \simeq \Delta^n \pbh fg$ are covers and by \cref{lem:Whitehead} the map $\pbh
  fg$ is invertible, thus $\cL^\pitchfork=\cR$. A dual argument proves
  ${^\pitchfork}R=\cL$.

  Conversely, let $(\cL,\cR)$ be a weak orthogonal system with $\cL$ stable by
  codiagonal. We have that $\Delta^n \pbh fg \simeq \pbh {\nabla^nf} g \simeq
  \pbh f {\Delta^n g}$ is a cover for any $g$ in $\cR$ and $f$ in $\cL$. Hence
  $\cR\subset \cL^\bot$ and obviously $\cL^\bot \subset \cL^\pitchfork = \cR$,
  thus $\cR = \cL^\bot$. A dual argument proves $\cL= {^\bot} \cR$.
\end{proof}

\begin{definition}[Accessibility]
  \label{def:os-wos-accessibility}
  Let $\kappa$ be a regular cardinal. A (weak) orthogonal system is
  \defn{$\kappa$\=/accessible} if it is generated by $W\to\cC\arr$
  where $W$ is a $\kappa$\=/small category and each $w\in\im W$ is a
  $\kappa$\=/small object of $\cC\arr$. A (weak) orthogonal system is
  \defn{accessible} if it is $\kappa$\=/accessible for some $\kappa$.
\end{definition}


\section{Factorisation systems}
\begin{definition}
 \label{def:fact-system} 
  A \defn{factorisation system} on a category $\cC$ is the data of an orthogonal
  system $(\cL,\cR)$ such that, for every map $f$ in $\cC$, there exists a
  factorisation of $f$ as $\rho(f)\lambda(f)$, namely a commutative triangle as
  below, with $\lambda(f)\in \cL$ and $\rho(f)\in \cR$.
  \[
    \begin{tikzcd}[sep=scriptsize]
      A \ar[rr,"f"] \ar[dr,"\lambda(f)"']
      &&B\\
      &C\ar[ur,"\rho(f)"']
    \end{tikzcd}
  \]
  The factorisation is only assumed to exist and is a property of the orthogonal
  system. However, the factorisation can be proven to be unique (precisely, the
  category of such factorisations of any $f$ is a contractible groupoid). This
  defines functors $f\mapsto \lambda(f)$ and $f\mapsto \rho(f)$ that are,
  respectively, right and left adjoints to the inclusions $\cL\subset \cC\arr
  \supset \cR$. In particular, $\cR$ is a reflective subcategory of $\cC\arr$.

  A factorisation system is \defn{accessible} if its underlying orthogonal
  system is accessible.
\end{definition}

\begin{remark}
  \label{rem:HTT-fact-systems}
  Factorisation systems are treated in \cite[5.2.8]{LurieHT}. The condition of
  stability under retracts in \cite[5.2.8.8]{LurieHT} is superfluous; it follows
  from~\cref{lem:prop-os} since a retract is both a limit and a colimit in
  $\cC\arr$.
\end{remark}

\begin{definition}
  \label{def:wfs}
  A \defn{weak factorisation system} on a category $\cC$ is the data of a weak
  orthogonal system $(\cA,\cB)$ such that, for every map $f$ in $\cC$, there
  exists a factorisation of $f$ as $ba$, with $a\in \cA$ and $b\in \cB$.

  The factorisation is once again only assumed to exist, and is a property of
  the weak orthogonal system. However in this case, the factorisation is not
  unique and not functorial in general.

  A weak factorisation system is \defn{accessible} if its underlying weak
  orthogonal system is accessible.  
\end{definition}

\begin{para}[Modalities and lex modalities]
 \label{para:modalities-lex-modalities} 
  A factorisation system $(\cL,\cR)$ is called \defn{stable}, or a
  \defn{modality}, if the factorisation is stable under base
  change (pullback along
  arbitrary maps)~\cite{ABFJblakersmassey,RSSmodalities}. It is easily seen that this is
  equivalent to asking that both classes $\cL$ and $\cR$ be stable under base
  change. Since the right class $\cR$ is always stable under base change, the
  condition is in fact only on $\cL$.

  A modality $(\cL,\cR)$ is \defn{left-exact} (\defn{lex} for short) if
  the factorisation is stable under finite limits in $\cC\arr$. Once again,
  since the right class $\cR$ is always stable under limits, the condition is in
  fact only on the class $\cL$.

  A number of characterizations of lex modalities is given in~\cite[Thm
  3.1]{RSSmodalities}, but not the one we have just used as a definition. We will
  show that our definition is equivalent to~\cite[Thm 3.1(x)]{RSSmodalities}.
  Recall that a morphism $\alpha\colon f\to g$ in $\cC\arr$ is a \emph{cartesian map} if the corresponding
  square in $\cC$ is cartesian.
\end{para}

\begin{proposition}
  \label{prop:stable-lex-FS}

  Let $(\cL,\cR)$ be a factorisation system on $\cC$.
  \begin{enumerate}
  \item \label{prop:stable-lex-FS:enum1} $(\cL,\cR)$ is a modality if and only
    if the reflection $\rho:\cC\arr \to \cR$ preserves cartesian maps.
  \item \label{prop:stable-lex-FS:enum2} $(\cL,\cR)$ is a lex modality if and
    only if the reflection $\rho:\cC\arr \to \cR$ is left-exact. In particular,
    a lex modality is a modality.
  \end{enumerate}
\end{proposition}

\begin{proof}
  \noindent \ref{prop:stable-lex-FS:enum1} Let $f\colon A\to B$ and $f'\colon A'\to B'$ be
  two maps in $\cC$ and let $\alpha\colon f'\to f$ be a cartesian map between them.
  The factorisation of $f'$ and $f$ induces a diagram
  \[
    \begin{tikzcd}[sep=scriptsize]
      A'  \ar[d,"\lambda(f')"'] \ar[rr] \ar[rrd,phantom,"\lambda(\alpha)"]
      && A \ar[d,"\lambda(f)"]
      \\
      C'  \ar[d,"\rho(f')"'] \ar[rr] \ar[rrd,phantom,"\rho(\alpha)"]
      && C \ar[d,"\rho(f)"]
      \\
      B'  \ar[rr] 
      && B .
    \end{tikzcd}
  \]
  Then, the factorisation system is a modality if and only if for every
  such $\alpha\colon f'\to f$, both
  $\lambda(\alpha)$ and $\rho(\alpha)$ are cartesian squares in $\cC$. But by the
  cancellation property of cartesian squares, this is equivalent to
  $\rho(\alpha)$ being cartesian only.

  \smallskip
  \noindent \ref{prop:stable-lex-FS:enum2} Let $f_i\colon A_i\to B_i$ be a finite
  diagram in $\cC\arr$ and $f= \lim f_i\colon A\to B$. We have the canonical
  diagram below.
  \[
    \begin{tikzcd}[sep=scriptsize]
      A  \ar[d,"\lambda(f)"'] \ar[rr, equal]
      && \lim A_i \ar[d,"\lim \lambda(f_i)"]
      \\
      C  \ar[d,"\rho(f)"'] \ar[rr]
      && \lim C_i \ar[d,"\lim \rho(f_i)"]
      \\
      B  \ar[rr, equal] 
      && \lim B_i
    \end{tikzcd}
  \]
  The factorisation system is left-exact if and only the canonical map $C\to \lim C_i$ is
  an equivalence. But this is the same as $\rho (\lim f_i ) \simeq \lim \rho(f_i)$.

  For the last statement, let $f\colon X\to Y$ and $f'\colon X'\to Y'$ be two maps in $\cC$
  and $\alpha\colon f'\to f$ be a morphism in $\cC\arr$. Remark that $\alpha$ is
  cartesian if and only if the following square is cartesian in $\cC\arr$.
  \[
    \begin{tikzcd}[sep=scriptsize]
      f' \ar[r,"\alpha"] \ar[d]
      & f\ar[d] \\
      1_{Y'} \ar[r]
      & 1_{Y}
    \end{tikzcd}
  \]
  Since the reflection $\rho:\cC\arr \to \cR$ is left-exact and codomain-preserving, the image
  under $\rho$ of this square is cartesian in $\cC\arr$ and of the same type. Thus a lex
  modality is a modality.
\end{proof}

\begin{para}[Enriched factorisation systems]
  \label{para:enriched-FS}

  When $(\cC,\otimes,\bbone,\intmap--)$ is a locally presentable, symmetric
  monoidal closed category, the internal pushout-product and
  pullback-hom define a symmetric closed monoidal structure
  $(\cC\arr,\boxtimes,0 \to \bbone,\intpbh--)$
  (see~\cref{para:internal-pushout-prod}). We define the \defn{enriched
    orthogonality} relation $f\iperp g$ by the condition that the internal
  pullback-hom $\intpbh f g$ be an invertible map in $\cC$. The enriched and
  external orthogonality are related by
  \[
    f \iperp g \quad\Leftrightarrow\quad \forall X\in\cC,\ (X\otimes f) \perp g
    \quad\Leftrightarrow\quad \forall X\in\cC,\ f \perp \intmap X g.
  \]
  
  An \defn{enriched orthogonal system} is a pair $(\cL,\cR)$ of full
  subcategories of $\cC\arr$ such that $\cR=\cL^\intperp$ and
  $\cL={^\intperp}\cR$. An (external) orthogonal system $(\cL,\cR)$ on $\cC$ is
  enriched if and only if $\cL$ is stable under $X\otimes-$ for any $X$ in
  $\cC$, if and only if $\cR$ is stable by $\intmap X-$ for any $X$ in $\cC$.
  Every enriched orthogonal system is an orthogonal system. An enriched
  orthogonal system is \defn{accessible} if it is of the form $\left(^\intperp
    (W^\intperp), W^\intperp\right)$ for $W\to \cC\arr$ a small diagram of
  internally small objects in $\cC\arr$.

  An \defn{enriched factorisation system} is an enriched orthogonal system
  $(\cL,\cR)$ such that, for every map $f$ in $\cC$, there exists a
  factorisation $f= \rho(f)\lambda(f)$ where $\lambda(f)\in \cL$ and $\rho(f)\in
  \cR$. Every enriched factorisation system is a factorisation system. An
  enriched factorisation system is \defn{accessible} if its underlying enriched
  orthogonal system is accessible.

  When $(\cC,\times,1,\intmap--)$ is a locally presentable \emph{cartesian}
  closed category, the enriched orthogonality relation also satisfies
  \[
    f \iperp g \quad\Ra\quad \forall X\in\cC,\ f\perp X\times g,
  \]
  since $X\times g$ is the base change of $g\colon A\to B$ along  $X\times
  B\to B$. Thus for any enriched factorisation system $(\cL,\cR)$ on $\cC$,
  the factorisation of $f:A\to B$ is always stable by base change
  along projections $X\times B\to B$ (weaker than modalities, see
  \cref{table:strengthFS}).
\end{para}

\begin{center}	
\begin{table}[htbp]
\renewcommand{\arraystretch}{1.5}
\begin{tabularx}{.9\textwidth}{
|>{\hsize=.8\hsize\linewidth=\hsize\centering\arraybackslash}X
|>{\hsize=1.2\hsize\linewidth=\hsize\centering\arraybackslash}X
|>{\hsize=1\hsize\linewidth=\hsize\centering\arraybackslash}X
|}
\hline
$\cC$ is any category 
    & {\sl External orthogonality}
    $w\perp f$ 
    & Factorisation systems
\\
\hline
$\cC$ is a cartesian closed category 
    & {\sl Enriched orthogonality}
    $w \iperp f \ \iff\ \forall X, (X\times w) \perp f$ 
    & Enriched FS 
    
    (stable by product)
\\
\hline
$\cC$ is a locally cartesian closed category with 1
    & {\sl Fiberwise orthogonality} 
    $w\fwperp f \ \iff\ w'\perp f$\qquad\qquad
    $\forall$ base change $w'\to w$ 
    & Modality 
    
    (FS stable by base~change)
\\
\hline
$\cC$ is a topos
    & {\sl Lex orthogonality}\qquad\qquad
    for $W\subset \cC\arr$,\qquad\qquad 
    $W\lexperp f\ \iff\ \lim_iw_i\perp f$\qquad\qquad 
    $\forall$ finite diagram $I\to W$
    & Lex~modality/ lex~localisation\qquad 
    (FS stable by finite limits)
\\
\hline
\end{tabularx}
\medskip
\caption{Strength of orthogonality relations and of factorisation systems}
\label{table:strengthFS}
\end{table}
\end{center}


\section{Factorisation systems and localisations}
\label{sec:localisations}

\begin{para}[Localisations]
  \label{para:localisation}
  Let $\cC$ have a terminal object $1$. If $(\cL,\cR)$ is an orthogonal system
  on $\cC$, we write $\cR_1$ for the full subcategory of $\cC$ on
  objects $X$ such that $X\to 1$ is in $\cR$.

  When $(\cL,\cR)$ is a \emph{factorisation} system, then for any $X$ in $\cC$,
  the factorisation of the map $X\to 1$ produces a reflection of $\cC$ into
  $\cR_1$. Since $\cR$ is stable under left-cancellation, we have full
  inclusions $\cR_1\arr\subset \cR\subset \cC\arr$. We
  define $\cL^{3/2} \eqdef {^\bot}(\cR_1\arr)$, and $\cR^{2/3}\eqdef(\cL^{3/2})^\bot$.
  Thus $(\cL^{3/2},\cR^{2/3})$ is an orthogonal system on $\cC$. In particular,
  $\cL^{3/2}$ is closed under colimits in $\cC\arr$.

  For every $f\in\cC\arr,g\in\cR_1\arr$, the reflection $P\colon
  \cC\localisation \cR_1\, \colon\! \iota$ provides an equivalence
  $\pbh{Pf}{g}\simeq\pbh{f}{g}$ in $\cC\arr$. Hence $\cL^{3/2}\subset\cC\arr$ is
  exactly the class of maps that are inverted by $P\colon \cC\to\cR_1$, and as
  such, satisfies the ``2~out~of~3'' condition. Thus the reflection $P\colon
  \cC\to \cR_1$ is the \defn{localisation} of $\cC$ at $\cL^{3/2}$. Namely, for
  any category $\cD$, pre-composition with $P$ is a full inclusion of functor
  categories $\fun{\cR_1}\cD\subset\fun\cC\cD$ with image those functors $\cC\to\cD$
  that invert every map in $\cL^{3/2}$. 
  
  Moreover, since $\cL\subset \cL^{3/2}$ and since for any map $f\colon A\to B$
  in $\cC$,
  the top and bottom maps in the following square
  in $\cC$ are in $\cL$, we conclude that $\cL^{3/2}$ is the \defn{saturation}
  of $\cL$ (the closure of $\cL$ under the ``2~out~of~3'' condition).
  \[
    \begin{tikzcd}[sep=scriptsize]
      A \ar[r,"\in\cL"] \ar[d,"f"']
      & PA\ar[d,"Pf"] \\
      B \ar[r,"\in\cL"]
      & PB
    \end{tikzcd}
  \]
\end{para}


\begin{proposition}
 \label{prop:mod-lccc} 
  Let $\cC$ be locally cartesian closed, and let $(\cL,\cR)$ be
  a modality on $\cC$. Then for every $g\colon X\to Y$ in $\cR$ and every
  $h\colon Y\to Y'$, the map $\Pi_hg \in \cC\comma {Y'}$ is in $\cR$ (where
  $\Pi_h$ is right adjoint to base change $h^*\colon \cC\comma Y'\to\cC\comma Y$). 
\end{proposition}
\begin{proof}
  For any $f\in\cL$, the map $\pbh f {\Pi_hg}$ in $\cS$ is invertible if and
  only if each of its fibres over $1\in\cS$ is invertible (or \emph{contractible}).
  But every such fibre is also a fibre of $\pbh {f'} {g}$ where $f'$ is a base
  change of $f$. We conclude since $\cL$ is stable under base change.\qedhere
\end{proof}
\begin{corollary}
 \label{cor:mod-lccc} 
  Under the hypotheses of~\cref{prop:mod-lccc}, if in addition $\cC$ has a
  terminal object, then the reflective subcategory $\cR_1\subset \cC$ is locally cartesian
  closed. 
\end{corollary}

\begin{para}[Reflective factorisation systems]
  \label{para:reflective-FS}
  If $\cC$ has a terminal object $1$, we shall say that a factorisation system
  $(\cL,\cR)$ is \defn{reflective} (following \cite{cassidyhebertkelly1985}) if
  the class $\cL$ satisfies the ``2~out~of~3'' condition (namely,
  $\cL=\cL^{3/2}$). Since left classes of orthogonal systems are always stable
  by composition and right cancellation (\cref{lem:prop-os}), this condition is
  equivalent to $\cL$ having the left cancellation condition.

  When $\cC$ has finite limits and $(\cL,\cR)$ is a \emph{modality}, then the
  orthogonal system $(\cL^{3/2},\cR^{2/3})$ from \cref{para:localisation} is a
  reflective factorisation system. In effect, the
  $(\cL^{3/2},\cR^{2/3})$-factorisation $X\to Z\to Y$ of a map $X\to Y$ is given
  by the diagram
  \[
    \begin{tikzcd}
      X \ar[r]\ar[rd]\ar[rr,bend left=30,"\lambda(X\to 1)"]  & Z \ar[r] \ar[d] \pbmark & PX \ar[d] \\
      & Y \ar[r,"\lambda(Y\to 1)"'] & PY.
    \end{tikzcd}
  \]
  The maps $X\to PX$ and $Y\to PY$ are in $\cL$ by definition of $P$. The map
  $Z\to PX$ is in $\cL$ because $\cL$ is stable by base change. Hence, the map
  $X\to Z$ is in $\cL^{3/2}$. The map $Z\to Y$ is in $\cR^{2/3}$ since it is a
  base change of the map $PX\to PY$ in $\cR_1$. The existence of this
  factorisation also implies that for any $f$, we have $f\in\cR^{2/3}$ if and
  only if $f$ 
  is a base change of some $g\in\cR_1\arr$.

  Recall that a reflective localisation $P:\cC\localisation \cD$ is
  \emph{semi-left-exact} in the sense of \cite[Thm 4.3]{cassidyhebertkelly1985}
  or \emph{locally cartesian} in the sense of \cite{gepnerkock2017univalence} if
  $P$ preserves base change along maps in $\cD$ (viewed as maps in $\cC$). The
  following result also appears as \cite[Prop. 5.1]{cherubiniRijke2020} in the
  internal homotopy type theory of a topos, however it holds more generally.
\end{para}

\begin{proposition}
  \label{prop:mod-slex}
  Let $\cC$ have a terminal object $1$ and let $(\cL,\cR)$ be a modality on
  $\cC$. Then the associated localisation $P:\cC\localisation \cR_1$
  is semi-left-exact.
\end{proposition}
\begin{proof}
  Any map in $\cR_1$ is of the type $PX\to PY$ for some map $X\to Y$ in $\cC$.
  We consider a pullback square
  \begin{equation}
    \tag{\ensuremath{\star}}
    \label{eqn:mod-slex}
    \begin{tikzcd}[sep=scriptsize]
      X' \ar[r,"f"]\ar[d] \pbmark & PX \ar[d] \\
      Y' \ar[r,"g"'] & PY.
    \end{tikzcd}
  \end{equation}
  and its factorisation for the modality $(\cL,\cR)$
  \[
    \begin{tikzcd}[sep=scriptsize]
      X' \ar[r,"\lambda(f)"]\ar[d] \pbmark
      & X'' \ar[r,"\rho(f)"]\ar[d] \ar[rd,near end,"P(\ref{eqn:mod-slex})" description,
      phantom] \pbmark
      & PX \ar[d] \\
      Y' \ar[r,"\lambda(g)"']
      & Y'' \ar[r,"\rho(g)"']
      & PY.
    \end{tikzcd}
  \]
  Both squares are cartesian by
  \cref{prop:stable-lex-FS}\ref{prop:stable-lex-FS:enum1}. Since the maps $PX\to
  1$ and $PY\to 1$ are in $\cR$ by definition of $P$, the objects $X''$ and
  $Y''$ are equivalent to $PX'$ and $PY'$ respectively. Since $P$ is
  idempotent, the square $P(\star)$ is the image of the square
  (\ref{eqn:mod-slex}) under $P$. Hence $P$ preserves such cartesian
  squares.
\end{proof}


\begin{corollary}[cf. {\cite[Prop.
    1.4]{gepnerkock2017univalence}}, {\cite[6.1.3.15]{LurieHT}}]
  \label{cor:mod-univ-colims}
  
  Under the hypotheses of \cref{prop:mod-slex}, if in addition $\cC$ is
  cocomplete, finitely complete, and has universal colimits, then the same is
  true of $\cR_1$.
\end{corollary}
\begin{proof}
  Colimits in $\cR_1$ are calculated by applying $P$ to the colimit calculated
  in $\cC$ (of the same diagram). Finite limits in $\cR_1$
  are simply those calculated in $\cC$.

  For
  universality of colimits, let $X_{\bullet}\colon I\to \cR_1$ be a diagram, let
  $X$ be its colimit in $\cC$, and let $Y\to PX$ be a map in $\cR_1$. By
  pullback, we obtain a diagram $Y_{\bullet}\colon I\to\cR_1$ and a cartesian
  transformation $\alpha\colon Y_{\bullet}\to X_{\bullet}$, as below.
  \[
    \begin{tikzcd}[sep=scriptsize]
      Y_i \ar[r]\ar[d,"\alpha_i"'] \pbmark
      & Y'\pbmark \ar[r]\ar[d]
      & Y \ar[d] \\
      X_i \ar[r]
      & X \ar[r,"\lambda(X)"']
      & PX.
    \end{tikzcd}
  \] Since colimits are universal in $\cC$, we have $Y'=\colim_iY_i$ in $\cC$.
  Finally we have $PY'\simeq Y$ by~\cref{prop:mod-slex}.
\end{proof}

\begin{proposition}[cf. {\cite[Thm 3.1]{RSSmodalities}}]
\label{prop:lex-loc=lex-mod}
 Let $\cC$ have finite limits and let $(\cL,\cR)$ be a modality on
  $\cC$. Then $(\cL,\cR)$ is a reflective factorisation system if and only if it
  is a lex modality.  
\end{proposition}
\begin{proof}
  For all $f\colon A\to B,g\colon B\to C$, the following square is cartesian in $\cC\arr$. Hence if $\cL$
  is stable under finite limits, it is stable by left-cancellation. 
  \[
    \begin{tikzcd}[sep=scriptsize]
      f \ar[r]\ar[d] \pbmark
      &gf \ar[d]\\
      1_B\ar[r]
      &g \pomark
    \end{tikzcd}
  \]
  Conversely, if $\cL = \cL^{3/2}$, then for any cartesian square in $\cC$ as on
  the left below, in the diagram on the right, the back face is cartesian (since
  $(\cL,\cR)$ is a modality) and side faces are cartesian (since $(\cL,\cR)$ is
  reflective). Thus the front face is cartesian, so $P$ preserves finite limits.
  \[
    \begin{tikzcd}
      A \ar[r]\ar[d,"f"'] \pbmark
      &X \ar[d,"g"]\\
      B\ar[r]
      &Y
    \end{tikzcd}
   \qquad\qquad 
    \begin{tikzcd}[sep=small]
      A'\ar[rr]\ar[rd] \ar[dd,"\rho f"']
      && X' \ar[rd] \ar[dd,"\rho g"' near start]\\
      &PA
      && PX \ar[from=ll, crossing over]\ar[dd, "Pg"]\\
      B \ar[rd] \ar[rr]
      && Y \ar[rd]\\
      &PB \ar[from = uu, crossing over,"Pf"' near start] \ar[rr]
      && PY
    \end{tikzcd}
  \]
  Hence if $h_i\colon C_i\to D_i$ is a finite diagram in $\cC\arr$ and $h= \lim
  h_i\colon C\to D$, then it is easy to see that $\rho (h) =\lim_i \rho (h_i)$ since it is
  the base change of $P(h) = \lim_i P(h_i)$. We conclude
  by~\cref{prop:stable-lex-FS}\ref{prop:stable-lex-FS:enum2}.  
\end{proof}

\chapter{Kelly's small object argument and the plus-construction for sheaves}
\label{chap:SOA}

The original \emph{small object argument}
(\cite[\textsection{3}.L3]{quillen1967}) is a construction that, under suitable
``smallness'' conditions, allows an accessible weak orthogonal system on a
$1$\=/category to be enhanced into a weak factorisation system. Variants of this
construction, as in Gabriel-Ulmer \cite[Satz 8.5]{gabrielulmer1971} and Kelly
\cite[Thms 10.1 and 11.3]{kelly1980} allow an accessible orthogonal system on a
locally presentable $1$\=/category to be enhanced into a factorisation system.
From now on, we will use the sobriquet \defn{small object argument}
(abbreviated \defn{SOA}) to refer to the generalisation of these
constructions to \oo\=/categories.

For a topological space, and more generally for a Grothendieck topology $J$ on a
$1$\=/category $C$, the \emph{sheafification} of a presheaf $X\in\psh C$ is a
particular case of the previous situation---since (following
\cref{sec:localisations}) it is exactly the factorisation of $X\to 1$ in the
left-exact modality $({}^\bot(J^\bot),J^\bot)$ on $\psh C$. In this case, there
is a construction on $X$ (\cite[Exp. II, 3.0.5]{SGA4-1}) that we will call the
``plus-construction'', that can be iterated to factor $X\to 1$.

In this chapter, we will show that Kelly's SOA admits a generalisation to
locally presentable \oo-categories. Moreover, we will show that under mild
conditions, that can be satisfied without loss of generality, Kelly's SOA
simplifies to a construction that is \emph{exactly} the plus\=/construction in
the particular case of a lex modality generated by a Grothendieck topology. We
use this simplification (the theory of \emph{pre-modulators}) to define
\emph{modulators} and \emph{lex modulators}, which are a good
generalisation of Grothendieck topologies to \oo-categories (since they capture
all accessible semi-left-exact and left-exact localisations).

\begin{notation}
 \label{notn:plus-construction-W-notation} 
  For a map $f$ in $\cC$, we recall the notation $\s f$ and $\t f$ to refer to
  the source and the target of $f$, that is, we write $f\colon\s f\to \t f$. We
  will use the terms ``source'' and ``domain'' (respectively, ``target'' and
  ``codomain'') interchangeably.

  For a diagram $W\to \cC\arr$, we will say (with slight abuse) that a map $w$
  in $\cC$ is ``in $W$'' if it is in the image of $W\to\cC\arr$. For any
  $f\in\cC\arr$, we write $W\comma f$ and $f\comma W$ for the over and under
  comma categories. For any $X\in\cC$, we write $W\comma X$ and $X\comma W$ for
  $W\comma 1_X$ and $1_X\comma W$. Remark that these are just $W\relcomma\t X$
  and $X\relcomma\s W$ for the source and target functors $\s,\t\colon W\to\cC$.
\end{notation}

\begin{remark}
  The results of \cref{sec:modulators-modalities,sec:lex-modulators} are largely
  due to M. Anel. They appear in \cite[Sec. 3]{anelLS2020small} in more detail
  than here. I have slightly changed and abridged their presentation from
  \emph{op. cit.} and have corrected a few errors. I have included them here to
  show the usefulness of pre-modulators (\cref{defn:pre-modulator}) and the
  general plus-construction (\cref{thm:plus-construction}).
\end{remark}

\section{Orthogonal factorisation after Kelly}
\label{sec:SOA-Kelly}
We use the pushout-product and pullback\=/hom to describe a version of the SOA
that constructs the accessible factorisation system $(\cL_W,\cR_W)$ generated by
a set of maps $W$ of a locally presentable category. This argument is a slight
modification of the one in~\cite[\textsection{10,11}]{kelly1980}.

\begin{para}
  It is useful to study the case of the orthogonal system
  $({}^\bot(w^\bot),w^\bot)$ generated by a single map $w$ of a category
  $\cC$. Let $f$ be any map in $\cC$. In the original argument of Quillen, when
  constructing an approximation $\Quillenmap f$ to a \emph{weak} factorisation,
  lifts against $w$ are added freely to $f$, independently of any existing lifts.
  Moreover, if several lifts exist, they are not identified. These observations
  suggest that in order to obtain unique lifts, we replace the lifting problem
  \[
    \begin{tikzcd}[sep=small]
      \fun{w}{f}\times \s w \ar[r]\ar[d] & \s f\ar[d, "f"]\\
      \fun{w}{f}\times \t w \ar[r] & \t f
    \end{tikzcd}
  \]
  with the lifting problem
  \[
    \begin{tikzcd}[sep=small]
      \fun{\t w}{\s f}\times \t w\coprod_{\fun{\t w}{\s f}\times \s w}\fun{w}{f}\times
      \s w \ar[r]\ar[d]
      & \s f\ar[d, "f"]\\
      \fun{w}{f}\times \t w \ar[r]
      & \t f.
    \end{tikzcd}
  \]
  The space $\fun{w}{f}$ is that of all squares between $w$ and $f$, with or
  without lifts, while $\fun{\t w}{\s f}$ parametrises all squares with lifts.
  Altogether, this new lifting problem contains the old one (the right component
  of the pushout), but enforces existing lifts (left component) to be identified
  with the corresponding free lift (middle component). The primary observation is
  that this new square is nothing but the counit square
  \[
    \pbh w f \pp w \to f .
  \]
  Then, by adjointness, we see that the previous square has a (not necessarily
  unique) lift if and only if the identity square below has a lift, which is
  exactly $w\perp f$.
  \[
    \begin{tikzcd}
      \pbh wf \ar[r,equals] &\pbh wf 
    \end{tikzcd}
  \]
\end{para}

\begin{remark}
  In~\cite[\textsection{10,11}]{kelly1980}, Kelly uses the sum $\fun{\t w}{\s
    f}\times \t w\coprod \fun{w}{f}\times \s w$ instead of the pushout $\fun{\t
    w}{\s f}\times \t w\coprod_{\fun{\t w}{\s f}\times \s w}\fun{w}{f}\times \s
  w$. The two constructions turn out to be equivalent in the context of
  1-categories, but the correct formula for $\infty$\=/categories does need the
  pushout. The difference concerns the uniqueness of the lifts.
\end{remark}

\begin{remark}
  Quillen's SOA uses a generating {\em set} of maps, but the variant for
  orthogonal factorisation allows for any small \emph{diagram} of maps $W\to
  \cC\arr$. For such a diagram, we will consider the single map that is the
  coend \footnote{See \cite[\textsection{3}]{anel2018exponentiable},
    \cite{Haugseng2021coends} for a treatment of coends.}
  \[
    W_K (f) \eqdef \int^{w\in W} \pbh w f \pp w .
  \]
  This is nothing but the domain of the density comonad of the functor
  $W\to\cC\arr$ for the enrichment $\pbh--$.
\end{remark}


\begin{lemma}
  [Functorial lifting]
  \label{lem:caracortho}
  Let $W\to \cC\arr$ be a diagram of arrows in $\cC$. The following conditions
  are equivalent:
  \begin{enumerate}[label={(\alph*)}, leftmargin=*]
  \item \label{lem:caracortho-1} For all $w$ in $W$, the map $\pbh{w}{f}$ is
    invertible.
  \item \label{lem:caracortho-2} For all $w$ in $W$, the identity map of
    $\pbh{w}{f}$ has a diagonal lift.
    \[
      \begin{tikzcd}[sep=small]
        \fun{\t w} {\s f}	\ar[r, equal] \ar[d] 		& \fun{\t w}{\s f}	\ar[d]\\
        \fun w f \ar[r, equal] \ar[ru, dashed]& \fun w f.
      \end{tikzcd}
    \]

  \item \label{lem:caracortho-3} For all $w$ in $W$, the map $\pbh{w}{f} \pp w
    \to f$ has a diagonal lift.
  \item \label{lem:caracortho-4} The single map $\epsilon:\int^{w\in W}
    \pbh{w}{f}\pp w\to f$ has a diagonal lift.
  \end{enumerate}
\end{lemma}
\begin{proof}
  The equivalence \ref{lem:caracortho-1} $\iff$ \ref{lem:caracortho-2} is
  straightforward. The equivalence \ref{lem:caracortho-2} $\iff$
  \ref{lem:caracortho-3} is direct by adjunction. Notice that the lift in
  \ref{lem:caracortho-2} is not asked to be unique (but by uniqueness of
  inverses, it will be automatically) and that we use this fact in the
  equivalence \ref{lem:caracortho-2} $\iff$ \ref{lem:caracortho-3}.

  If all maps $\pbh{w}{f}$ are invertible, then
  by~\cref{lem:pp-pbh-iso-absorption} so are the maps $\pbh{w}{f}\pp w$ and the
  coend $\int^w \pbh{w}{f}\pp w$. This proves \ref{lem:caracortho-1} $\Ra$
  \ref{lem:caracortho-4}. Reciprocally, we simply use the canonical maps $\pbh w
  f \pp w \to \int^w \pbh w f \pp w$ to show that \ref{lem:caracortho-4} $\Ra$
  \ref{lem:caracortho-3}.
\end{proof}

\begin{para}
  [Kelly's construction]
 \label{para:kelly-construction} 
  The commutative square $\epsilon$ of 
  \cref{lem:caracortho}\ref{lem:caracortho-4} is
  \begin{equation}
    \tag{$!\boxslash$}
    \label{eqn:liftcoend!}
    \begin{tikzcd}[sep=scriptsize]
      \int^w\left( \fun{\t w}{\s f}\times \t w\coprod_{\fun{\t w}{\s f}\times \s w}
        \fun{w}{f}\times \s w \right) \ar[r]\ar[d,"W_K (f)"']
      & \s f\ar[d, "f"]\\
      \int^w \fun{w}{f}\times \t w \ar[r] & \t f.
    \end{tikzcd}
  \end{equation}
  The construction then proceeds as in \cite{kelly1980}. The map $\Kellymap
  f:\Kelly f\to \t f$ is defined as the cogap map of the square
  (\ref{eqn:liftcoend!})
  \begin{equation}
    \tag{$k$\=/constr}
    \label{eqn:kconstruction}
    \begin{tikzcd}[sep=scriptsize]
      \int^w\left( \fun{\t w}{\s f}\times \t w\coprod_{\fun{\t w}{\s f}\times \s w}
        \fun{w}{f}\times \s w \right) \ar[r]\ar[d,"W_K (f)"']
      & \s f\ar[d, "u(f)"'] \ar[ddr, bend left, "f"]&\\
      \int^w \fun{w}{f}\times \t w \ar[r, "\ell"] \ar[rrd, bend right=20]
      & \Kelly f \ar[rd, "\Kellymap f"'] \pomark &\\
      && \t f.
    \end{tikzcd}
  \end{equation}
  Let $\t_W f$ be the identity map of the codomain of $W_K (f)$, that is, of
  $\int^w \fun{w}{f}\times \t w$. In the presentation of~\cite[Prop.
  9.2]{kelly1980}, the functor $f\mapsto \Kellymap f$ is defined via the pushout
  in $\cC\arr$
  \[
    \begin{tikzcd}[sep=small]
      W_K (f) \ar[r] \ar[d] & f \ar[d]
      \\
      \t_W f \ar[r] & \Kellymap f \pomark .
    \end{tikzcd}
  \]
  It is convenient to look at this diagram in terms of the
  lift
  \begin{equation}
    \tag{$k$\=/lift}\label{eqn:k-lift}
    \begin{tikzcd}[sep=scriptsize]
      \int^w\left( \fun{\t w}{\s f}\times \t w\coprod_{\fun{\t w}{\s f}\times \s w}
        \fun{w}{f}\times \s w \right) \ar[r]\ar[d]
      & \s f\ar[d, "f"] \ar[r,"u(f)"]
      & \Kelly f\ar[d, "\Kellymap f"]\\
      \int^w \fun{w}{f}\times \t w  \ar[r] \ar[rru, dashed, "\ell"']
      & \t f \ar[r, equal] & \t f.
    \end{tikzcd}
  \end{equation}
  The map $\Kellymap f$ can then be understood as the best approximation of $f$ on the
  right that admits a lift from the map $\int^w\pbh{w}{f} \pp w$.

  We have a factorisation of $f$ as $\Kellymap f u(f):\s f\to \Kelly f\to \t f$.
  Because $\cL$ is absorbing for $\square$ and stable by colimits in $\cC\arr$
  (\cref{lem:prop-os}), the map $\int^w\fun{w}{f} \pp w$ is in $\cL$ and so is
  its cobase change $u(f):X\to \Kelly f$. The map $\Kellymap f:\Kelly f\to \t f$
  need not be in $\cR$ but we now show that it will be after a transfinite
  iteration.
\end{para}

\begin{theorem}
  [{Kelly's SOA \cite[Thm 11.5]{kelly1980}}]
  \label{thm:!SOA2}
  Let $\cC$ be a cocomplete category and $W\to \cC\arr$ be a small diagram of
  arrows with small domains and codomains, and let $f$ be a map in $\cC$. Then,
  the transfinite iteration of~(\ref{eqn:kconstruction}) on $f$ converges to a
  factorisation of $f$ for the orthogonal system $(\cL_W, \cR_W)$ generated by
  $W$. Thus $(\cL_W,\cR_W)$ is an orthogonal factorisation system.
\end{theorem}

\begin{proof}
  The natural transformation $u:f\to \Kellymap f$ defines a transfinite sequence
  \[
    f\xto{u(f)} \Kellymap f\xto{u(\Kellymap f)} \Kellymapn2 f\xto{u(\Kellymapn2
      f)} \dots
  \]
  We will show that this sequence converges. By hypothesis, all maps $w\in W$
  are small (namely, their sources and targets are small) and $W$ is a small
  category. Thus, we can fix a regular cardinal $\kappa$ majoring the size of
  all $w$. The map $u^\kappa(f)$ is in $\cL_W$ because $\cL_W$ is stable by
  transfinite compositions. If the map $\Kellymapn \kappa f$ is
  in $\cR_W=W^\bot$, then $\Kellymapn \kappa f \simeq \Kellymapn
  {\kappa+1} f$. We will show that it is.

  We define $\Kellyn \kappa f\eqdef\s \Kellymapn \kappa f$. For any $w$ in $W$, $\pbh{w}{\Kellymapn \kappa
    f}$ is invertible if and only if the identity map of $\pbh{w}{\Kellymapn
    \kappa f}$ has a lift
  \[
    \begin{tikzcd}[sep=scriptsize]
      \fun{\t w}{\Kellyn \kappa f}\ar[r, equal]\ar[d]
      & \fun{\t w}{\Kellyn \kappa f} \ar[d]\\
      \fun{w}{\Kellymapn \kappa f}\ar[r, equal]\ar[ru, dashed]
      & \fun{w}{\Kellymapn \kappa f}.
    \end{tikzcd}
  \]
  Since $w$ is $\kappa$\=/small and $\kappa$ is $\kappa$\=/filtered, we have
  \[
    \pbh{w}{\Kellymapn \kappa f} = \colim_{\lambda<\kappa} \pbh{w}{\Kellymapn \kappa f}.
  \]
  Also, for a fixed $w$, the map $\ell$ of the square (\ref{eqn:k-lift})
  provides a diagonal lift for the square $\pbh w f \pp w \to \Kellymap f$. By
  adjunction, it also provides a diagonal lift $\ell(f)$ for the square $\pbh w
  f \to \pbh w {\Kellymap f}$.
  \[
    \begin{tikzcd}[column sep=scriptsize]
      \fun {\t w}{\s f}	\ar[r]\ar[d]					& \fun{\t w}{\Kelly f} \ar[d]\\
      \fun w f		\ar[r]\ar[ru, dashed, "\ell(f)"]	& \fun{w}{\Kellymap f}.
    \end{tikzcd}
  \]
  By iteration, we have a diagram
  \[
    \begin{tikzcd}[sep=scriptsize]
      \fun{\t w}{\Kellyn \lambda f}\ar[r]\ar[d]
      & \fun{\t w}{\Kellyn {\lambda+1} f}\ar[r]\ar[d]
      &\dots \ar[r]
      &  \fun{\t w}{\Kellyn \kappa f}\ar[r, equal]\ar[d]
      & \fun{\t w}{\Kellyn \kappa f} \ar[d]\\
      \fun{w}{\Kellymapn \lambda f}\ar[r]\ar[ru, "\ell^\lambda"description]
      & \fun{w}{\Kellymapn {\lambda+1} f}\ar[r] \ar[ru, "\ell^{\lambda+1}"description]
      & \dots \ar[r]
      &  \fun{w}{\Kellymapn \kappa f}\ar[r, equal]\ar[ru, dashed, "\ell"]
      & \fun{w}{\Kellymapn \kappa f}
    \end{tikzcd}
  \]
  where the maps $\ell^\lambda \eqdef \ell(\Kellymapn \lambda f)$ provide the
  desired lift $\ell$ at the limit.
\end{proof}

\begin{remark}
  \label{rem:Kelly-comparison}
  Because of its focus on reflective subcategories,
  \cite[\textsection{10}]{kelly1980} only considers the adjunction
  \[
    \begin{tikzcd}[sep=scriptsize]
      \s (-\pp w):\cS\arr \ar[r, shift left]
      \ar[from=r, shift left]
      & \cC : \pbh w-
    \end{tikzcd}
  \]
  where $\cC$ is identified with the full subcategory of $\cC\arr$ on maps $X\to
  1$. However, the construction in \cref{para:kelly-construction} relies on the
  adjunction
  \[
    \begin{tikzcd}[sep=scriptsize]
      -\pp w :\cS\arr \ar[r, shift left]
      \ar[from=r, shift left]
      & \cC\arr : \pbh w-.
    \end{tikzcd}
  \]
  More precisely, the relevant adjunction for a diagram $W\to \cC\arr$ of arrows
  is the enriched realisation-nerve adjunction
  \[
    \begin{tikzcd}
      &W \ar[rd]\ar[ld]\\
      \fun {W\op} {\cS\arr} \ar[rr, shift left=1.6,"{\int^w(-\pp w)}"] \ar[from=rr,
      shift left=1.6,"{w\mapsto \pbh w-}","\bot"']
      && \cC\arr \ar[loop right, "W_K"]
    \end{tikzcd}
  \]
  built from the diagram $W\to \cC\arr$ using the enrichment $\pbh--$ of
  $\cC\arr$ over $(\cS\arr, \pp)$. The comonad $W_K \colon f\mapsto \int^w \pbh
  w f \pp w$ is then the enriched density comonad of the diagram $W\to \cC\arr$.
\end{remark}


\begin{corollary}
  \label{cor:+-fix=R}
  Under the hypotheses of \cref{thm:!SOA2}, a map $f$ is in $W^\bot$ if and only
  if the canonical map $f\to \Kellymap f$ is invertible.
\end{corollary}
\begin{proof}
  If $f$ is in $W^\bot$, the maps $\pbh w f$ are invertible and so is the map
  $W_K (f)=\int^w \pbh w f \pp w$. Thus $f\simeq \Kellymap f$. Reciprocally, if
  $f\simeq \Kellymap f$, then by \cref{lem:caracortho}\ref{lem:caracortho-4}
  and \ref{eqn:k-lift}, $f$ is in $W^\bot$.
\end{proof}

\section{Pre-modulators and the plus construction}
\label{sec:plus-construction}

The purpose of this section is to connect the $k$\=/construction of
\cref{thm:!SOA2} to the plus-construction involved in sheafification~\cite[Exp.
II, 3.0.5]{SGA4-1}. Recall the enrichment of $\cC\arr$ over the cartesian closed
category $(\cS\arr,\times)$ (see \cref{para:arrow-cat-cartesian-closure}).
\[
  \bracemap w f : \fun w f \to \fun {\t w}{\t f} 
\]
\begin{definition}[The plus-construction]
  \label{defn:plus-construction}
  Let $f$ be a map in $\cC$ and $W\to \cC\arr$ a diagram of maps. We define the
  \defn{plus-construction} of the map $f$ as the map
  \begin{align*}
    \Plusmap f &\eqdef\int^w \bracemap w f \times \t w 
    \\
               &= \int^w \fun w f \times \t w \to \int^w \fun {\t w}{\t f} \times \t w
    \\
               &= \colim_{W\commaindex f} \t w \to \colim_{W \relcomma\t \t f} \t w
  \end{align*}
  where $W \comma{\t f}$ is the comma category of $\t\colon W \to \cC$ over
  the object $\t f$. We write $\Plus f$ for the domain of $\Plusmap f$. The
  second iteration of the plus-construction will be written $\Plusplusmap f$ and
  $\Plusplus f$.
\end{definition}

\begin{remark}
  Remark that $\Plus -$ is the left Kan extension of $\t\colon
  W\to\cC\arr\to\cC$ along $W\to \cC\arr$, and that $\t \Plusmap {(-)}$ is the
  left Kan extension of $\t\colon W\to\cC\arr\to\cC$ along itself. The map
  $\Plusmap {(-)}$ is the transformation associated to the left Kan
  extension of $\Plus-$ along $\t\colon\cC\arr\to\cC$.
\end{remark}

\begin{para}
  [Grothendieck topologies]
\label{para:Groth-top-plus-constr}
  We recall the usual ``plus-construction'' in the
  context of sheafification of presheaves. Let $C$ be a small category, let $\cC
  = \Psh C$, and let $W\subset \cC\arr$ be the full subcategory of covering
  sieves $R\to x$ of some Grothendieck topology on $C$. For $x$ in $C$, we let
  $W(x)$ be the category of covering sieves of $x$ (\ie the fiber of the target
  functor $\t\colon W\to C$ at $x$). Recall that for a presheaf $F$, the presheaf
  $\Plusoriginal F$ is defined by
  \begin{align*}
    C\op &\tto \cS \\
    x &\mto \Plusoriginal F(x) = \colim_{R\to x\in W(x)\op}\fun R F.
  \end{align*}
  Let us show that this definition coincides with~\cref{defn:plus-construction}
  in the particular case of a Grothendieck topology.
\end{para}

\begin{lemma}
  \label{lem:plus-is-plus-for-topology}
  Let $C,\cC,W,F$ be as in~\cref{para:Groth-top-plus-constr}.
  Then the presheaf $\Plusoriginal F$ is equivalent to the domain $\Plus {F\to 1}$ of
  the plus-construction $\Plusmap {(F\to 1)}$ for the diagram $W\subset\cC\arr$.
\end{lemma}
\begin{proof}
  The domain of $\Plusmap {(F\to 1)}$ is the presheaf
  \[
    H_F :x \mapsto \fun x {\colim_{W\relcomma\s F} \t w} = \colim_{W\relcomma\s
      F} \fun x {\t w},
  \]
  where $W\relcomma\s F$ is the comma category of $\s\colon W\to\cC$ over $F$.
  We need to construct an equivalence $\Plusoriginal F \simeq H_F$.
  We consider the category $W(x,F)$ that is the category of elements of the
  functor
  \begin{align*}
    W\relcomma\s F &\tto \cS   \\
    (sw\to F) &\mto \fun x {\t w}.
  \end{align*}
  Its objects are diagrams
  \[
    \begin{tikzcd}[sep=small]
      & R \ar[r]\ar[d,"w"']  & F\ar[d] \\
      x \ar[r]    & y \ar[r]        & 1
    \end{tikzcd}
  \]
  where $w\colon R\to y$ is a covering sieve. Its morphisms are the natural
  transformations that are the identity on $x$ and $F$. Remark that
  $W(x,F)\simeq (x\relcomma\t W)\relcomma\s F$. Recall that the \emph{external
  groupoid} of a category $C$ (obtained by localising at all maps in $C$) is the
  colimit $\abs C\eqdef\colim_C1$ in $\cS$. Then, we have that
  $H_F(x)=\colim_{W\relcomma\s F} \fun x {\t w} = \abs{W(x,F)}$.

  Let $W(x)/ F$ be the category of elements of the functor
  \begin{align*}
    W(x)\op &\tto \cS   \\
    R\to x &\mto \fun R F.
  \end{align*}
  We have that $\Plusoriginal F(x) = \abs{W(x)/ F}$. Remark that $W(x)/F\simeq
  W(x)\relcomma\s F$.
  There is an adjunction $G:W(x)/ F \rightleftarrows W(x,F):D$ such
  that
  \[
    G\left(
      {\begin{tikzcd}[sep=small]
          R \ar[r]\ar[d]  & F\ar[d] \\
          x \ar[r] & 1
        \end{tikzcd}}
    \right)
    \quad = \quad
    {\begin{tikzcd}[sep=small]
        & R \ar[r]\ar[d]  & F\ar[d] \\
        x \ar[r,equal]  & x \ar[r]
        & 1
      \end{tikzcd}}
  \]
  \[
    D\left(
      {\begin{tikzcd}[sep=small]
          & R \ar[r]\ar[d]  & F\ar[d] \\
          x \ar[r]    & y \ar[r]           & 1
        \end{tikzcd}}
    \right)
    \quad = \quad
    {\begin{tikzcd}[sep=small]
        R\times_yx \ar[r]\ar[d]  & F\ar[d] \\
        x \ar[r]           & 1
      \end{tikzcd}}
  \]
  which is induced by the colocalisation $W(x)\rightleftarrows (x\relcomma\t W)$
  (which is due to $\t\colon W\to C$ being a fibration). So we obtain an
  equivalence of external groupoids $\abs{W(x)/ F} \simeq \abs{W(x,F)}$ which is
  natural in $x$ and in $F$.
\end{proof}

\begin{remark}
  \label{rem:plus-is-plus}
  The proof of~\cref{lem:plus-is-plus-for-topology} relies implicitly on the
  stability under base change of covering sieves (in order to construct the
  right adjoint $D$) and on the fact that $W$ contains the identity of
  generators, that is, representable functors (to have the natural map $F\to
  \Plusoriginal F$). The facts that the covering sieves are monomorphisms and
  local is irrelevant.~\cref{lem:plus-is-plus-for-topology} holds for any small
  full subcategory $W\to {\Psh C}\arr$ such that the codomain functor $\t\colon
  W\to {\Psh C}\arr \to \Psh C$ has values in $C$, $W$ contains all the identity
  maps of $C$, and $W$ is stable by base change along maps of $C$ (namely,
  $\t\colon W\to C$ is a fibration). This motivates the subsequent definitions of a
  pre-modulator (\cref{defn:pre-modulator}) and a modulator
  (\cref{defn:modulator}).
\end{remark}

\begin{remark}
  \label{rem:values+construction-in-presheaves}
  Given any diagram $W\to {\Psh C}\arr$ whose codomains are in $C$ and a map $f$
  in $\Psh C$, let $C\comma \Plus f$ be the category of elements of the presheaf
  $\Plus f=\s (\Plusmap f)$. The associated fibration in groupoids $C\comma
  {\Plus f}\to C$ is the right part of the (cofinal, fibration in groupoids)
  factorisation of the functor $\t\colon W\comma f\to C$ sending a map $w\to f$
  to $\t w$.
  \[
    \begin{tikzcd}[sep=small]
      &C\comma \Plus f\ar[rd, "\text{fib. in gpds}"]\\
      W\comma f \ar[rr,"\t"] \ar[ru,"\text{cofinal}"] && C
    \end{tikzcd}
  \]
  The value of the presheaf $\Plus f$ at some $x$ in $C$ is the external
  groupoid of the category $W(x,f)\simeq (x\relcomma\t W)\comma f$ whose objects
  are diagrams
  \[
    \begin{tikzcd}[sep=small]
      & \s w \ar[r]\ar[d,"w"'] & \s f\ar[d,"f"] \\
      x \ar[r]    & \t w \ar[r]            & \t f.
    \end{tikzcd}
  \]
  When $W$ is stable by base change, exactly the same proof as in
  \cref{lem:plus-is-plus-for-topology} gives
  \[
    \Plus f(x) \simeq \colim_{W(x)\op}\fun w f .
  \]
\end{remark}


\begin{para}
  We require some more notation. We set $\Minusmap f\eqdef\int^{w\in W}\!\!
  \pbh w f \pp w$ to be the cocartesian gap map of the square below, and we
  define $a(f),b(f),c(f)$ as below.
  \[
    \begin{tikzcd}[sep=scriptsize]
      \int^w \fun{\t w}{\s f}\times \s w  
      \ar[dd]
      \ar[rrr]
      \ar[rrrdd, "c(f)\eqdef\int^w \pbh w f \times w" description]
      &&& \int^w \fun w f \times \s w \ar[dd, "b(f)\eqdef\int^w \fun w f \times w"]
      \\
      \\
      \int^w \fun{\t w}{\s f}\times \t w
      \ar[rrr,"a(f)\eqdef\int^w \pbh w f \times \t w"'] 
      &&& \int^w \fun w f \times \t w
    \end{tikzcd}
  \]
  This square can also be written using colimits (recall $W \comma X$ from
 \cref{notn:plus-construction-W-notation}). 
  \[
    \begin{tikzcd}[sep=scriptsize]
      \colim_{W \commaindex {\s f}} \s w
      \ar[dd]
      \ar[rr]
      \ar[rrdd, "c(f)" description]
      && \colim_{W\commaindex f} \s w
      \ar[dd, "b(f)"]
      \\
      \\
      \colim_{W \commaindex {\s f}} \t w
      \ar[rr,"a(f)"] 
      && \colim_{W\commaindex f}\t w
    \end{tikzcd}
  \]
  Finally, \cref{fig:main-diag-plus-constr} is going to be central to \cref{thm:plus-is-k}.
\end{para}
\begin{figure}
  \[
    \mathclap{
      \begin{tikzcd}[column sep=scriptsize,ampersand replacement=\&]
        {\displaystyle\colim_{W \commaindex {\s f}} \s w} 
        \ar[rd, "\alpha",start anchor={[xshift=-1.5ex,yshift=1ex]}]
        \ar[r, "\beta"]
        \ar[dddrr, "c(f)=\int^w \pbh w f \times w" description, near start, bend right=60] 
        \&{\displaystyle\colim_{W \commaindex f} \s w} 
        \ar[dddr, "b(f)=\int^w \fun w f \times w" description, bend
        right=15]
        \ar[rd, "\alpha'"]
        \ar[rrd, "\gamma", bend left=10]
        \\
        \&{\displaystyle\colim_{W \commaindex {\s f}}\t w} 
        \ar[ddr, "a(f)=\int^w \pbh w f \times \t w" description, bend right=40] 
        \ar[r, "\beta'", crossing over] 
        \& {\displaystyle\colim_{W \commaindex {\s f}} \t w} 
        \underset{\colim_{W \commaindex {\s f}} \s w}{\coprod}
        {\displaystyle\colim_{W \commaindex f} \s w}
        \ar[dd,"\Minusmap f=\int^w \pbh w f \pp w" description]
        \ar[r, "\delta"']
        \& \s f 
        \ar[dd,"u(f)"  description]
        \ar[dddd,"f",bend left=50]
        \\
        \&\&{}
        \\
        \&\& \Plus f={\displaystyle\colim_{W \commaindex f}\t w} 
        \ar[dd, "\Plusmap f = \int^w \bracemap w f \times \t w" description]
        \ar[r, "\varepsilon"]
        \& \Kelly f\pomark
        \ar[dd,"\Kellymap f" description]
        \\
        {}
        \\
        \&\&{\displaystyle\colim_{W \commaindex{\t f}} \t w} 
        \ar[r, "\zeta"]
        \&\t f 
      \end{tikzcd}
    }
  \]
  \caption{Decomposition of $\Plusmap{(-)}$ and $\Kellymap-$.}
  \label[figure]{fig:main-diag-plus-constr}
\end{figure}

\begin{para}
  [A coincidence condition for $\Plusmap {(-)}$ and $\Kellymap -$] Inspired by
  the particular case of \cref{lem:plus-is-plus-for-topology,rem:plus-is-plus},
  we will give conditions for the construction $\Kellymap f$ of \cref{thm:!SOA2}
  to coincide with $\Plusmap f$. We shall do so by proving that, under suitable
  conditions on the diagram $W\to\cC$, all the Greek-lettered maps
  of~\cref{fig:main-diag-plus-constr} are invertible.

  So far the only condition on $\cC$ has been that it is a cocomplete category.
  We have also assumed the maps in $W$ to be small. The proof of the coincidence
  of the $k$\=/ and plus-constructions will require the stronger assumption that
  $\cC$ is locally presentable. Since any object is small in a locally
  presentable category, this will remove the smallness assumption on the maps in
  $W\to \cC\arr$ when applying \cref{thm:!SOA2}, as long as $W$ is still a small
  indexing category.
\end{para}

\begin{definition}[Pre-modulator]
  \label{defn:pre-modulator}
  Let $\cC$ be a locally presentable category with a fixed generating small
  subcategory $C\subset \cC$. A \defn{pre-modulator} is a diagram $W\to \cC\arr$
  such that:
  \begin{enumerate}
  \item \label[axiom]{enum:modulator:1} $W$ is a small category,
  \item \label[axiom]{enum:modulator:2} $W\to \cC\arr$ is fully faithful,
  \item \label[axiom]{enum:modulator:3} the codomains of the maps in $W$ are all in
    $C$,
  \item \label[axiom]{enum:modulator:4} the inclusion $C\subto \cC\arr$ sending a
    generator to its identity map factors through $W\subto \cC\arr$.
  \end{enumerate}
  We often leave implicit the choice of the generating category $C$ when
  considering pre-modulators, since \cref{enum:modulator:3,enum:modulator:4}
  imply that $C$ is just the image of $\t\colon W\to \cC$.
\end{definition}

\begin{remark}
  \label{rem:def-admissibility}
  Consider the comma category $\cC \comma C$. The functor $\t\colon
  \cC\comma C\to C$ admits a fully faithful right adjoint $1\colon C\subto \cC\comma C$
  sending a generator to its identity map. The adjunction $\t\colon \cC\comma
  C\localisation C:1$ is then a reflective localisation. A pre-modulator can
  equivalently be defined as a small full subcategory $W\subset \cC\comma C$
  such that the previous localisation restricts to a localisation $\t\colon
  W\localisation C:id$.
\end{remark}

\begin{theorem}[The plus-construction]
  \label{thm:plus-is-k}
  \label{thm:plus-construction}
  Let $\cC$ be a locally presentable category, and let $W\to \cC\arr$ be a
  pre-modulator. Then in \cref{fig:main-diag-plus-constr}, we have the
  simplifications
  \[
    u(f) = a(f) = b(f) = c(f) = \Minusmap f \eqnand \Kellymap f = \Plusmap f.
  \]
  In particular, the factorisation of $f$ for the orthogonal system
  $(\cL_W,\cR_W)$ can be obtained by a transfinite iteration of the
  plus-construction. Moreover, a map $f$ is in $\cR$ if and only if
  $f\simeq\Plusmap f$.
\end{theorem}

\begin{proof}
  We prove that all the Greek\=/lettered maps in
  \cref{fig:main-diag-plus-constr} are invertible. First, consider the
  generating small category $C\subset \cC$ relative to which $W$ is a
  pre-modulator. Let $F$ be any object in $\cC$. By hypotheses
  \ref{enum:modulator:2} and \ref{enum:modulator:4} we get a fully faithful
  functor
  \begin{align*}
    C \comma F &\tto W \commaindex F\\
    c\to F & \mto (1_c, c\to F).
  \end{align*}
  Then, using hypothesis \ref{enum:modulator:3} and
  \cref{rem:def-admissibility}, this functor has a left adjoint sending $(w,\t
  w\to F)$ to $\t w\to F$. Recall that a right adjoint functor is always
  cofinal. Since $C$ is a generating subcategory, we have
  \[
    \colim_{W \commaindex {F}} w \ =\ \colim_{C\commaindex {F}} 1_c \ =\ 1_{F}.
  \]
  Applied to $F\eqdef \s f$, we obtain $\alpha = 1_{\s f}$.
  The same argument with $F\eqdef \t f$ gives
  \[
    \colim_{W \commaindex {\t f}} \t w \ =\ \t f,
  \]
  proving that $\zeta$ is invertible, and with $F\eqdef\s f$, gives
  $\gamma\beta=1_{\s f}$.
 
  Next, we prove that $\gamma$ is also the identity map of the object $\s f$.
  Consider the category $C \comma W\comma f$ whose objects are triples
  $(c,w,1_c\to w\to f)$ where $c$ is in $C$, $w$ in $W$, $1_c\to w$ in $W$ and
  $w\to f$ in $\cC\arr$, and whose morphisms $(c,w,1_c\to w\to f) \to
  (c',w',1_{c'}\to w'\to f)$ are pairs $(c\to c',w\to w')$ such that the obvious
  diagram commutes. Remark that $C\comma f = C\comma sf$ and that the functor
  $C\to C\comma W$ sending $c$ to $1_{1_c}\colon 1_c\to 1_c$ induces a functor
  $C\comma sf\to C\comma W\comma f$. This has a right adjoint $h\colon C \comma
  W\comma f \to C \comma {\s f}$ sending $1_c\to w \to f$ to $c\to \s f$ (using
  hypotheses \ref{enum:modulator:2} and \ref{enum:modulator:4}). In particular,
  $h$ is cofinal. Then we have
  \begin{align*}
    \colim_{W\commaindex f}\ \s w &= \colim_{W\commaindex f}\colim_{C\commaindex {\s w}} c
    \\
                                  &= \colim_{C\commaindex W\commaindex f} c
    \\
                                  &= \colim_{C\commaindex {\s f}} c & \text{by cofinality}
    \\
                                  &= \s f,
  \end{align*}
  proving that $\gamma$ is invertible. Hence, since $\gamma\beta=1_{\s f}$, thus
  $\beta$ is invertible. Then, by pushout, so are $\alpha'$ and $\beta'$. Then
  $\delta$ is invertible since $\gamma$ and $\alpha'$ are. By pushout, so is
  $\varepsilon$. Thus all the Greek-lettered maps are invertible, giving the
  identities of the theorem. The last assertion follows from $\Kellymap f =
  \Plusmap f$ and \cref{cor:+-fix=R}.
\end{proof}

\begin{remark}
  There is a concrete payoff to the $k$\=/construction simplifying to the
  plus-construction in the case of a pre\=/modulator. The colimit formula
  defining the plus-construction is better suited to check exactness conditions
  such as stability by base change and left-exactness (see \cref{thm:+modality}
  and \cref{thm:lex+construction}). Moreover, any diagram $W\to\cC\arr$ may be
  replaced with a pre-modulator in an essentially harmless way
  (\cref{prop:modulator-completion}).
\end{remark}

\begin{remark}
  \label{rem:weakening-admissibility}
  A careful reading of the proof of \cref{thm:plus-construction} shows that it
  depends on less than the hypothesis that $W\to \cC\arr$ be fully faithful. The
  proof that $C \comma F \to W \commaindex F$ is cofinal needs only that, for
  any $w$ in $W$, the canonical map $w\to 1_{\t w}$ is in $W$. The proof
  that $h:C \comma W\comma f \to C \comma {\s f}$ is cofinal needs only that
  $\relmap W {1_c} w = \relmap \cC c {\s w}$.
\end{remark}

\begin{proposition}
  [Pre-modulator envelope]
  \label{prop:modulator-completion}
  Any diagram $W\to \cC\arr$ can be completed into a pre-modulator generating
  the same factorisation.
\end{proposition}
\begin{proof}
  Let $W\to \cC\arr$ be a small diagram. Because $W$ is small, we can find a
  small category $C$ of generators of $\cC$ containing the codomains of all maps
  in $W$. We consider the full subcategory $W'\subset\cC\arr$ generated by the
  image of $W$ and the identity maps of $C$. Then $W'$ is a pre-modulator.
  We need to show that $W^\bot = (W')^\bot$. The definition of $(W')^\bot$
  depends only on the objects in the image of $W'$ and not the morphisms between
  them. The objects of $W'$ are those in the image of $W$ and the identity maps
  of $C$. We deduce that $(W')^\bot = W^\bot \cap C^\bot$, and we conclude since
  $C^\bot$ is the whole of $\cC\arr$.
\end{proof}

\begin{remark}
  The $k$\=/constructions associated to $W$ and $W'$ are \emph{a priori} different,
  but they converge to the same endofunctor.
\end{remark}

\section{Modulators and modalities}
\label{sec:modulators-modalities}

In this section, we use \cref{thm:plus-construction} to give sufficient
conditions on $W$ for the factorisation system $(\cL_W,\cR_W)$ to be a modality.
Throughout this section, we fix a locally presentable category $\cC$ with
universal colimits (equivalently, locally presentable and locally cartesian
closed).

\begin{notation}
  Recall that for any $f,g\in\cC\arr$, a map $f\to g$ is \emph{cartesian} if it
  is a cartesian square in $\cC$ (it is a cartesian arrow of the fibration
  $\t\colon\cC\arr\to\cC$).
\end{notation}

\begin{definition}[Modulator]
  \label{defn:modulator}
  Let $C$ be generators for $\cC$. We say that a pre-modulator $W\to \cC\arr$ is
  a \defn{modulator} if
  \begin{enumerate}
    \setcounter{enumi}{4}
  \item \label{enum:modulator:5} the codomain functor $\t:W \to C$ is a
    fibration.
  \end{enumerate}
  Namely, for any $w$ in $W$, any $x$ in $C$ and any map $x\to \t w$, the base
  change $w':\s w\times_{\t w}x\to x$ is in the (essential) image of $W$.
\end{definition}



\begin{lemma}
  \label{lem:BC-cart}
  Given a cartesian square in $\cC\arr$
  \[
    \begin{tikzcd}[sep=small]
      g'\ar[r]\ar[d] \pbmark & f'\ar[d]\\
      g\ar[r] & f
    \end{tikzcd}
  \]
  such that the map $f'\to f$ is cartesian, then the map $g'\to g$ is also
  cartesian.
\end{lemma}
\begin{proof}
  By the cancellation property of cartesian squares.
\end{proof}

\begin{lemma}
  \label{lem:Kan-opfib}
  Let $C$ and $D$ be small categories, $F\colon C\to D$ an opfibration,
  and $G\colon C\to \cC$ a functor with values in a cocomplete category. The
  left Kan extension of $G$ along $F$ is given pointwise by the colimits
  \[
    \lan_FG(d) = \colim_{c\in F(d)} G(c).
  \]
  where $F(d)$ is the fibre of $F$ over $d$.
\end{lemma}
\begin{proof}
  For any $d\in D$, recall that the fibre $F(d)$ of $F$ over $d$ and
  the over category $C\relcomma F d$ are defined by the pullbacks
  \[
    \begin{tikzcd}[sep=small]
      F(d) \ar[r, "\iota_d"]\ar[d]\pbmark  
      & C\relcomma F d \ar[r]\ar[d]\pbmark 
      & C \ar[d,"F"]
      \\
      \{d\} \ar[r]& D\comma d \ar[r]& D
    \end{tikzcd}
  \]
  Since $F$ is an op-fibration,
  the functor $\iota_d:F(d)\to C\relcomma F d$ has a left adjoint. Hence
  $\iota_d$ is cofinal, and the result follows.
\end{proof}

\begin{lemma}
  \label{lem:+cartesian}
  Let $W\to \cC\arr$ be a modulator and $\alpha:g\to f$ a cartesian map in
  $\cC\arr$. Then, for any ordinal $\kappa$, the two following squares are
  cartesian.
  \[
    \begin{tikzcd}[sep=scriptsize]
      \s g\ar[d] \ar[r] \ar[rd, phantom,"(1)"]
      & \s g^{+\kappa}\ar[d] \ar[r,"g^{+\kappa}"] \ar[rd, phantom,"(2)"]
      & \t g \ar[d] \\
      \s f \ar[r]         & \s f^{+\kappa} \ar[r,"f^{+\kappa}"'] & \t f
    \end{tikzcd}
  \]
\end{lemma}
\begin{proof}
  By cancellation for cartesian squares, it is enough to prove that the square
  (2) is cartesian. We proceed by induction on $\kappa$.

  \noindent\textbf{(Base.)} When $\kappa = 1$, we
  need to show that the canonical map $\s g^+ \to \s \Plusmap f \times_{\t f} \t
  g$ is invertible, where $\s g^+ = \colim_{W\commaindex g} \t w$ and where
  \begin{align*}
    \s \Plusmap f \times_{\t f} \t g 
    &= \left(\colim_{W\commaindex f} \t w\right) \times_{\t f} \t g \\
    &= \colim_{W\commaindex f} \left(\t w \times_{\t f} \t g \right) \\
    &= \colim_{W\commaindex f} \t (w \times_fg).
  \end{align*}
  Given $v\to g$ in $W\comma g$, we have a square
  \[
    \begin{tikzcd}[sep=small]
      v \ar[d, equal]\ar[r] & g \ar[d,"\alpha"] \\    
      v \ar[r] & f
    \end{tikzcd}
  \]
  and a map $v\to v \times_fg$ in $W\comma g$ (by \cref{lem:BC-cart} and since
  $W$ is a modulator). The induced transformation
  \[
    \begin{tikzcd}[sep=small]
      W\comma g \ar[rd] \ar[rr]   & \ar[d, "\Ra" description, phantom]
      & W\comma f \ar[ld,"{-\times_fg}"]\\
      & \cC\arr
    \end{tikzcd}
  \]
  is just the unit of the adjunction $\alpha_!: \cC\arr\comma g\rightleftarrows
  \cC\arr\comma f:\alpha^*$, restricted to $W\comma g$ and $W\comma f$. The
  corresponding map between colimits is
  \[
    \zeta:\colim_{W\commaindex g} w \tto \colim_{W\commaindex f} w \times_fg,
  \]
  and the domain of $\zeta$ is the map $\s g^+ \to \s \Plusmap f \times_{\t f}
  \t g$. The result follows if we show that $\zeta$ is invertible in $\cC\arr$.
  We do so by introducing an auxiliary diagram $V \to \cC\arr$ and two morphisms
  of diagrams
  \[
    \begin{tikzcd}[sep=small]
      W\comma g \ar[rrd]
      &\ar[rd, "=" description, phantom]
      & V \ar[d] \ar[rr] \ar[ll] 
      &\ar[ld, "\Ra" description, phantom]
      & W\comma f \ar[lld,"{-\times_fg}"]
      \\
      && \cC\arr
    \end{tikzcd}
  \]
  that will induce equivalences between the colimits. The category $V$ is
  defined as the full subcategory of $(\cC\arr)\arr \comma \alpha$ consisting of
  squares
  \[
    \begin{tikzcd}[sep=small]
      v \ar[d, "\beta"']\ar[r] 
      & g \ar[d,"\alpha"] \\    
      w \ar[r] & f
    \end{tikzcd}
  \]
  where $v$ and $w$ are in $W$ and $\beta$ is cartesian. The diagram $V\to
  \cC\arr$ is given by forgetting to $v$. The obvious functor
  $V\to W\comma g$ has a left adjoint given by
  \[
    v\to g \mto
    \begin{tikzcd}[sep=small]
      v \ar[d, equal]\ar[r] & g \ar[d,"\alpha"] \\    
      v \ar[r] & f
    \end{tikzcd}.
  \]
  Hence it is cofinal and $\colim_V v \simeq \colim_{W\commaindex g} v$.

  There is also an obvious functor $B:V\to W\comma f$ and a morphism of diagrams
  given by the maps $v \to w\times_fg$ (this is the unit of
  $(\cC\arr)\arr\comma\alpha \rightleftarrows \cC\arr\comma f$). The composition
  of the induced map $\colim_V v \to \colim_{W\commaindex f} w\times_fg$ with
  the previous identification $\colim_{W\commaindex g} v \simeq \colim_V v$
  gives back the map $\zeta$.
  
  We show that $B\colon V\to W\comma f$ is an opfibration. Let $w\to w'\to f$ be a map in $W\comma f$ and
  $\beta\colon v\to w$ in the fibre $B(w\to f)$ of $B$ at $w\to f$. Since $\cC$
  has finite limits, the composite map $v\to w\to w'$ in $\cC\arr$ factors
  uniquely as $v\to v'\to w'$ where $\t v=\t v'$ and $\beta'\colon v'\to w'$ is
  cartesian.
  \[
    \begin{tikzcd}[sep=small]
      v \ar[d, "\beta"']\ar[r, dashed] & v' \ar[d,"\beta'", dashed] \\    
      w \ar[r] & w'
    \end{tikzcd}
  \]
  We claim that $\beta'$ is in the fibre $B(w'\to f)$. Since $\t v'=\t v$, $v'$
  is a base change of $w'$, and since $W$ is a modulator, $v'$ is in $W$. Since
  $\t v'=\t v$ and $\alpha\colon g\to f$ is cartesian, the maps $v\to g$ and
  $v'\to w'\to f$ together give a square below.
  \[
    \begin{tikzcd}[sep=small]
      v' \ar[d, "\beta'"']\ar[r] 
      & g \ar[d,"\alpha"] \\    
      w' \ar[r] & f
    \end{tikzcd}
  \]
  The claim follows since $\beta'$ is cartesian, and implies that $B$ is an
  opfibration.

  Now, the fibre $B(w\to f)$ is just the full subcategory of
  $W\comma (w \times_fg)$ consisting of cartesian maps $v\to w \times_fg$.
  %
  Since $B$ is an opfibration, we use $V\to W\comma f \to 1$, decomposition of
  left Kan extensions, and \cref{lem:Kan-opfib} to get
  \[
    \colim_V v = \colim_{W\commaindex f} \colim_{v\underset{\mathrm{cart}}{\to}
      w\times_fg} v.
  \]
  We finally show that $\colim_{v\to w\times_fg} v = w\times_fg$. The object
  $X=\t (w\times_fg)$ is the colimit of the canonical diagram $C\comma X\to
  \cC$. For any $\xi\colon y\to X$ in $C\comma X$, let $w_\xi\colon Z\to
  y$ be the left-hand pullback
  \[
    \begin{tikzcd}[sep=scriptsize]
      Z\ar[d,"w_\xi"']\ar[r]\pbmark
      &\s (w\times_fg) \ar[d,"w\times_fg"']\ar[r]\pbmark
      &\s w \ar[d,"w"] \\
      y \ar[r,"\xi"]
      &\t (w\times_fg) \ar[r]
      &\t w.
    \end{tikzcd}
  \]
  Since $W$ is a modulator, by pasting of pullbacks, $w_\xi$ is in $W$. Thus the
  diagram of the $w_\xi$ is exactly the diagram of cartesian maps $v\to
  w\times_fg$ where $v$ is in $W$. But by universality of colimits, the map
  $w\times_fg$ is the colimit of the maps $w_\xi$, thus $\colim_{v\to
    w\times_fg} v = w\times_fg$. We conclude that
  \[
    \zeta:\colim_{W\commaindex g} v \simeq \colim_V v \to \colim_{W\commaindex
      f} w\times_fg.
  \]
  is invertible.
  
  \noindent\textbf{(Induction.)} If $\kappa$ is not a limit
  ordinal, then the result is a consequence of the previous computation applied
  to $g^{+\kappa-1} \to f^{+\kappa-1}$. If $\kappa$ is a limit ordinal, we have
  $\s g^{+\kappa} = \colim_{\lambda<\kappa}\s g^{+\lambda}$. By induction
  hypothesis, we have $\s g^{+\lambda} = \s f^{+\lambda}\times_{\t f}\t g$ for
  all $\lambda<\kappa$. We conclude by universality of colimits.
\end{proof}


\begin{theorem}[Stable plus-construction]
  \label{thm:+modality}
  Let $\cC$ be a locally presentable locally cartesian closed category with a
  fixed small category of generators $C\subset \cC$. If $W\to \cC\arr$ is a
  modulator, then
  \begin{enumerate}
  \item \label{thm:+modality:enum1} the plus-construction $f\mapsto \Plusmap
    f$ preserves cartesian maps in $\cC\arr$,
  \item \label{thm:+modality:enum2} the factorisation system generated by $W$
    is a modality, and
  \item \label{thm:+modality:enum3} a map $f$ is modal (\ie in $\cR$) if and
    only if $f=\Plusmap f$.
  \end{enumerate}
\end{theorem}
\begin{proof}
  \noindent \ref{thm:+modality:enum1} follows from \cref{lem:+cartesian}.

  \smallskip
  \noindent \ref{thm:+modality:enum2} By \cref{thm:plus-construction}, the
  reflection into the right class is given by $\rho(f) = f^{+\kappa}$ for some
  $\kappa$. The result follows from \cref{lem:+cartesian} and
  \cref{prop:stable-lex-FS}.

  \smallskip
  \noindent \ref{thm:+modality:enum3} This is the final claim of
  \cref{thm:plus-construction}.
\end{proof}

\begin{remark}
  The fact that the plus-construction respects cartesian maps says, intuitively,
  that it is a fibrewise process on the map $f$. In particular, it can be
  defined as an endomorphism $+:\UU\to \UU$ of the universe (the fibration
  $\t\colon\cC\arr\to\cC$) of the category $\cC$.
\end{remark}


\begin{definition}
  \label{defi:base-change-envelope}
  Let $C\subset \cC$ be a generating small category. For any small diagram $W\to
  \cC\arr$, we define the \defn{modulator envelope} $W^{\mathrm{mod}}$ of $W$
  (relative to $C$) as the full subcategory of $\cC\arr$ consisting of the
  identity maps on objects of $C$ and all pullbacks of the maps in $W$ over the
  objects of $C$. The diagram $W^{\mathrm{mod}}\subto \cC\arr$ is always a
  modulator.
\end{definition}

\begin{corollary}[Generation of modalities]
  \label{cor:gen-modality}
  If $\cC$ is a locally presentable locally cartesian closed category with a
  fixed small category of generators $C\subset \cC$, the modality generated by
  an arbitrary family $W\to \cC\arr$ can be constructed by applying
  \cref{thm:!SOA2} to the modulator envelope $W^{\mathrm{mod}}$.
\end{corollary}
\begin{proof}
  By \cref{thm:+modality}, the factorisation system $(\cL,\cR)$ generated by
  $W^{\mathrm{mod}}$ is a modality. We need to prove that it is the smallest
  modality (with respect to the order given by inclusion of left classes) such
  that $W\subset \cL$. We first prove that $W\subset \cL$. Let $w$ be a map in
  $W$. For any $X$ in $C$ we define the diagram
  \begin{align*}
    w(-):C\comma \t w &\tto \cC\arr \\
    (\xi:x\to \t w) &\mto w(\xi):x\times_{\t w} \s w\to x.
  \end{align*}
  By universality of colimits, we have $\colim_\xi w(\xi) = w$. By construction,
  all maps $w(\xi)$ are in $W^{\mathrm{mod}}\subset\cL$. Hence $w\in\cL$ by
  stability under colimits (\cref{lem:prop-os}\ref{lem:prop-os:2}). Finally, let
  $(\cL',\cR')$ be another modality such that $W\subset \cL'$. Then
  $W^{\mathrm{mod}}\subset \cL'$, so $\cR'\subset (W^{\mathrm{mod}})^\bot$ and since
  $\cR=(W^{\mathrm{mod}})^\bot$, thus $\cL\subset \cL'$.
\end{proof}

\section{Lex modulators and lex modalities}
\label{sec:lex-modulators}
\label{sec:sheafification}

We now turn to the problem of generating left-exact localisations, namely lex
modalities (see~\cref{prop:lex-loc=lex-mod}). We stay in the context of a
locally presentable category $\cC$ with universal colimits.

\begin{definition}[Lex modulator]
  \label{defi:lex-modulator}
  Let $C\subset \cC$ be a generating small category. We say that a
  modulator $W\to \cC\arr$ is \defn{left-exact}, or \defn{lex}, if
  \begin{enumerate}
    \setcounter{enumi}{5}
  \item \label{enum:modulator:6} for all $x$ in $C$, the fibre $W(x)$ of
    $\t:W\to C$ is a co-filtered category.
  \end{enumerate}
  We say that a pre-modulator $W\to \cC\arr$ is \defn{stable under finite
    limits} if
  \begin{enumerate}
    [label=(\arabic*$^+$), leftmargin=*]
    \setcounter{enumi}{5}
  \item \label{enum:modulator:6+} the image of $W\to \cC\arr$ is stable under
    finite limits in $\cC\arr$.
  \end{enumerate}
  Condition~\ref{enum:modulator:6+} implies that the generating category $C$
  must have finite limits. A pre-modulator stable under finite limits is a lex
  modulator. That is, together with \ref{enum:modulator:3} and
  \ref{enum:modulator:4}, \ref{enum:modulator:6+} implies both conditions
  \ref{enum:modulator:5} and \ref{enum:modulator:6}.
\end{definition}

\begin{theorem}
  [Lex plus-construction]
  \label{thm:lex+construction}
  Let $\cE$ be an $n$\=/topos ($n\leq \infty$) and let $W\to \cE\arr$ be a lex
  modulator. Then:
  \begin{enumerate}
  \item \label{thm:lex+construction:1} the functor $f \mapsto \Plusmap f$ is
    left-exact as an endofunctor of $\cE\arr$,
  \item \label{thm:lex+construction:2} the factorisation system generated by $W$
    is a lex modality, and
  \item \label{thm:lex+construction:3} A map $f$ is a \emph{relative sheaf}
    (namely, in $\cR$) if and only if $f\simeq\Plusmap f$.
  \end{enumerate}
\end{theorem}

\begin{proof}
  We first prove the theorem for topoi (\oo\=/topoi).

  \smallskip
  \noindent \ref{thm:lex+construction:1} We first consider the case $\cE = \P C$
  of a presheaf topos. Let $f$ be a map in $\cE$, then for each object $x$
  of $C$, we have a map $f(x)$ in $\cS$. The plus-construction is left-exact if
  and only if
  all functors $f\mapsto \Plusmap f(x)$ are left-exact. Since $\t \Plusmap f =
  \t f$ and the functor $\t:\cE\arr \to \cE$ is always left exact, the condition
  reduces to proving that all functors $f\mapsto \Plus f(x)$ are left-exact. Let
  $W(x)$ be the fibre at $x$ of the target functor $\t\colon W\to C$. By
  \cref{lem:plus-is-plus-for-topology} and \cref{rem:plus-is-plus}, we have
  \[
    \Plus f(x) \quad=\quad \colim_{W(x)\op}\fun w f.
  \]
  By assumption, $W(x)\op$ is filtered, thus the functor $f\mapsto \Plus f(x)$
  is left-exact.

  Any topos $\cE$ is a left-exact localisation $P:\P C \localisation \cE:\iota$,
  for some small subcategory $C\subset\cE$ of generators. Both functors $P$ and
  $\iota$ are left-exact. For a map $f$ in $\cE$, the plus-construction is
  defined by a colimit indexed by $W\comma f$. The fully faithful functor
  $\iota\colon \cE\subto \P C$ induces a fully faithful functor
  $\iota\colon\cE\arr \subto {\P C}\arr$ and the category $W\comma f$ computed
  in $\cE\arr$ or ${\P C}\arr$ is the same. Hence, the plus-construction of
  $\cE$ is related to that of $\P C$ by the formula
  \[
    \Plusmap f = P\big((\iota f)^+\big).
  \]
  Since the three functors $P$, $\iota$ and the plus-construction of $\P C$ are
  left-exact (since $W$ is also a lex modulator in $\P C$), so is the
  plus-construction in $\cE$.

  \smallskip
  \noindent \ref{thm:lex+construction:2} Since filtered colimits of left-exact
  functors are still left-exact, the transfinite iterations $(-)^{+\kappa}$ are
  also left-exact. We deduce that the reflector $\rho$ is left-exact. Then the
  result follows from \cref{prop:stable-lex-FS}.

  \smallskip
  \noindent \ref{thm:lex+construction:3} This is the final claim of
  \cref{thm:plus-construction}.

  \smallskip
  This finishes the proof for \oo\=/topoi. The strategy is similar for
  $n$\=/topoi. The presheaf category $\P C = \fun {C\op} \cS$ needs only to be
  replaced by the category $\cP_n(C) = \fun {C\op} {\cS\truncated {n-1}}$, where
  $\cS\truncated {n-1}$ is the category of $(n-1)$\=/groupoids.
\end{proof}

\begin{lemma}
  \label{lem:BC-cart2}
  Given a cartesian map $f'\to f$ in $\cC\arr$, the obvious square
  \[
    \begin{tikzcd}[sep=small]
      f'\ar[r]\ar[d] \pbmark & f\ar[d]\\
      1_{\t f'}\ar[r] & 1_{\t f}
    \end{tikzcd}
  \]is cartesian.
\end{lemma}
\begin{proof}
  A straightforward computation.
\end{proof}

\begin{proposition}
  \label{prop:lex-mod-4-loc-lex}
  Let $f^*:\cE\to \cF$ be an accessible left-exact localisation of a topos
  $\cE$. Then there exists a lex modulator presenting $f^*$.
\end{proposition}
\begin{proof}
  Let $W\subset \cE\arr$ be the subcategory of arrows inverted by $f^*$. Since
  $f^*$ is a left-exact localisation, $W$ is stable by finite limits in
  $\cE\arr$. By assumption, there exists a regular cardinal $\kappa$ such that
  $W\subto\cE\arr$ is an accessible subcategory. Let $W(\kappa)\subset W$ be
  subcategory of $\kappa$\=/compact objects, choosing $\kappa$ such that $W$ is
  $\kappa$\=/accessible, $W(\kappa)$ is stable by finite limits and such that
  any object of $C$ is $\kappa$\=/compact. For such a $\kappa$, we set $W'
  \eqdef W(\kappa)\cap \P C \comma C$. Let us show that it is a lex modulator
  with respect to the generators $C$. First, $W'$ contains all equivalences
  between objects of $C$. Thus the stability by finite limits of $W$ and
  \cref{lem:BC-cart2} implies that $W'$ is a modulator. Then, for any $x$ in
  $C$, we need to prove that $W'(x)$ is co-filtered. Let $D\colon I\to W'(x)$ be
  a finite diagram. Its limit in $W'$ has codomain $\lim_Ix = x^{\abs I}$ (where
  $\abs I = \colim_I1$), which need not be an object of $C$, since we have not
  assumed $C$ to be stable by finite limits. But the base change of this limit
  along the diagonal $x\to x^{\abs I}$ is an element in $W'(x)$, giving a cone
  over $D$.
\end{proof}

\begin{para}
  [Lex localisation of truncated objects]
  \label{para:sheaf-truncated}
  We give a condition under which the plus-construction associated to a lex
  modulator converges in $(n+2)$ steps on $n$\=/truncated objects.
\end{para}

\begin{definition}[Mono-saturation]
  \label{lem:mono-saturated}
  Let $\cE$ be an $n$\=/topos ($n\leq \infty$) and $C\subset \cE$ a generating
  category. Let $W\subto \cE\arr$ be a lex modulator and let $(\cL,\cR)$ be the
  corresponding factorisation system. We say that $W$ is \defn{mono-saturated} if
  \begin{enumerate}
    \setcounter{enumi}{6}
  \item \label{enum:modulator:7} any monomorphism in $\cL$ with codomain in $C$
    is in $W$.
  \end{enumerate}
  For example, any Grothendieck topology on $\P C$ is mono-saturated. Any
  $n$\=/topos is well-powered, so a lex modulator can always be completed to a
  mono-saturated lex modulator with the same lex modality.
\end{definition}

The following lemma is the primary reason for \cref{prop:truncated-sheaf}.

\begin{lemma}
  \label{lem:truncated-sheaf}
  Let $\cE$ be an $n$\=/topos ($n\leq \infty$) and $W\subto \cE\arr$ a
  mono-saturated lex modulator. Let $(\cL,\cR)$ be the factorisation system
  generated by $W$. Then, for any monomorphism $m$ in $\cL$, $m^+$ is
  invertible. More generally, if $w$ is in $\cL$ and is $n$\=/truncated, then
  $w^+$ is $(n-1)$\=/truncated.
\end{lemma}
\begin{proof}
  Recall that a map $f\colon A\to B$ in $\cC$ is a monomorphism if and only if $f
  \simeq f\times_{1_B} f$ in $\cC\arr$. Because the plus-construction preserves
  codomains and is left-exact, it preserves monomorphisms. Hence, in order to
  prove that $m^+$ is invertible it is enough to prove that it has a section.

  Let $C$ be the generating category for $\cC$ and let $x$ be in $C$. Since $W$ is
  mono-saturated, any base change $m_\xi$ of $m$ along $\xi\colon x\to \t m$ is in
  $W$. Because colimits are universal in $\cC$, $m$ is the colimit of the
  $m_\xi$. By construction of $m^+$, we have lifts
  \[
    \begin{tikzcd}
      \s m_\xi \ar[r] \ar[d] & \s m \ar[r]\ar[d,"m" near end] &\s m^+ \ar[d,"m^+"]\\
      \t m_\xi \ar[r] \ar[rru,dashed, bend left=05] & \t m \ar[r,equal] &\t m
    \end{tikzcd}
  \]
  These lifts are unique since $m^+$ is a mono. Passing to the colimit, they
  define a lift
  \[
    \begin{tikzcd}
      \s m \ar[r]\ar[d,"m"'] &\s m^+ \ar[d,"m^+"]\\
      \t m\ar[r,equal] \ar[ru,dashed]&\t m.
    \end{tikzcd}
  \]

  For the second statement, recall that $f$ is $n$\=/truncated if and only if
  $\Delta^{n+2} f$ invertible, if and only if $\Delta^{n+1} f$ is a
  monomorphism. Left-exactness of $\Plusmap{(-)}$ and the previous statement
  imply that $\Delta^{n+1} (\Plusmap f)=\Plusmap{(\Delta^{n+1}f)}$ is
  invertible, hence $\Plusmap f$ is $(n-1)$\=/truncated.
\end{proof}

\begin{para}
  We need another lemma. Any two maps $f\colon A\to B$ and $g:B\to C$ in a
  category $\cC$ define a commutative cube
  \[
    \begin{tikzcd}[sep=small]
      A\ar[rr,equal]\ar[rd,"f" near end] \ar[dd,"f"'] && A \ar[rd,"f"] \ar[dd]\\
      &B && B \ar[from=ll, crossing over, equal]\ar[dd, "g"]\\
      B \ar[rd, equal] \ar[rr, "g" near start] && C \ar[rd, equal]\\
      &B \ar[from = uu, crossing over, equal] \ar[rr, "g"] && C
    \end{tikzcd}
  \]
  where both top and bottom faces are cartesian and cocartesian. This gives the
  following result.
\end{para}

\begin{lemma}
  \label{lem:facto-bicartesian}
  Viewed from above, the previous cube defines a square in $\cC\arr$ which is
  both cartesian and cocartesian.
  \[
    \begin{tikzcd}[sep=scriptsize]
      f \ar[r]\ar[d] \pbmark &gf \ar[d]\\
      1_B\ar[r]&g \pomark
    \end{tikzcd}
  \]
\end{lemma}

\begin{lemma}
  \label{lem:+facto}
  Let $\cE$ be an $n$\=/topos ($n\leq \infty$) and $W\subto \cE\arr$ a
  mono-saturated lex modulator. Let $f\colon A\to B$ be a map in $\cE$ and
  $\rho(f)\lambda(f):A\to M\to C$ its factorisation for the left-exact modality
  generated by $W$. Then, we have $\rho(\Plusmap f) = \rho(f)$ and
  $\lambda(\Plusmap f) = \lambda(f)^+$.
\end{lemma}
\begin{proof}
  Using \cref{lem:facto-bicartesian} for $f=\rho(f)\lambda(f)$ and the
  left-exactness of the plus-construction, we get a cartesian square
  \[
    \begin{tikzcd}[sep=scriptsize]
      \lambda(f)^+ \ar[r]\ar[d] \pbmark & \Plusmap f\ar[d]\\
      (1_M)^+\ar[r]&\rho(f)^+.
    \end{tikzcd}
  \]
  We have $\rho(f)^+ = \rho(f)$ since $\Plusmap{(-)}$ fixes the class $\cR$.
  This implies that $M$ is also the middle object of the factorisation of
  $\Plusmap f$. Using $(1_M)^+ = 1_M$ and \cref{lem:facto-bicartesian} for
  $\Plusmap f=\rho(\Plusmap f)\lambda(\Plusmap f)$, we get $\lambda(\Plusmap
  f)=\lambda(f)^+$.
\end{proof}

\begin{proposition}[Sheafification of truncated objects]
  \label{prop:truncated-sheaf}
  Let $\cE$ be an $n$\=/topos ($n\leq \infty$) and $W\subto \cE\arr$ a
  mono-saturated lex modulator. Then, if $f$ is an $m$\=/truncated map ($m\leq
  n$), we have $f^{+m+2} = \rho(f)$.
\end{proposition}
\begin{proof}
  By \cref{lem:+facto}, $f^{+m+2} = \rho(f^{+m+2})\lambda(f^{+m+2}) =
  \rho(f)\lambda(f)^{+m+2}$. Then, by \cref{lem:truncated-sheaf},
  $\lambda(f)^{+m+2}$ is an equivalence and $f^{+m+2} = \rho(f)$.
\end{proof}

\begin{remark}
  \label{rem:truncated-sheaf}
  Notice that we do not need $W\to \cC\arr$ to consist only of monomorphisms
  (although it needs to be mono-saturated). \Cref{prop:truncated-sheaf} works
  for all accessible left-exact localisations, topological or not. However, the
  fact that the modulator is lex is crucial.
\end{remark}


\begin{table}[htbp]
  \begin{center}	
    \caption{Summary of the conditions for the plus-construction}
    \label{table:cdt+}
    \smallskip
    \renewcommand{\arraystretch}{2}
    \begin{tabularx}{\textwidth}{
      |>{\hsize=.7\hsize\linewidth=\hsize\centering\arraybackslash}X
      |>{\hsize=.4\hsize\linewidth=\hsize\centering\arraybackslash}X
      >{\hsize=1.3\hsize\linewidth=\hsize\raggedright\arraybackslash}X
      |>{\hsize=1.6\hsize\linewidth=\hsize\centering\arraybackslash}X
      |}
      \hline
      $W\to \cC\comma C$
      & \multicolumn{2}{c|}{{\em Condition}}
      & {\em Property}
      \\
      \hline
      Pre-modulator
      & \ref{enum:modulator:1}-\ref{enum:modulator:4} 
      & $\t:W\rightleftarrows C:id$ is a reflective~localisation
      & The plus- and $k$\=/constructions coincide (\cref{thm:plus-construction}).
      \\
      \hline
      Modulator
      & \ref{enum:modulator:1}-\ref{enum:modulator:5} 
      &
        {\begin{tikzcd}[ampersand replacement=\&,sep=small,cramped]
            W\ar[rd,"\t"']\ar[r,hook] \& \cC\comma C \ar[d,"\ \t"]\\
            \& C    
          \end{tikzcd}}
      is a sub-fibration
      & The plus-construction generates modalities (\cref{thm:+modality}).
      \\
      \hline
      Lex modulator
      & \ref{enum:modulator:1}-\ref{enum:modulator:6} 
      & $W\subto \cC\comma C$ is a sub-fibration and fibres of $\t\colon W\to C$ are co-filtered
      & The plus-construction is left-exact (\cref{thm:lex+construction}).
      \\
      \hline
      Mono-saturated lex modulator
      & \ref{enum:modulator:1}-\ref{enum:modulator:7} 
      & $W$ contains all the monos in $\cL_W$ whose codomain is in $C$
      & The plus-construction converges in $(n+2)$ steps on $n$\=/truncated objects (\cref{prop:truncated-sheaf}).
      \\
      \hline 
    \end{tabularx}
  \end{center}
\end{table}



\bibliographystyle{halpha}
\bibliography{./biblio}
\end{document}